\documentclass[a4paper, 12pt]{book}
\usepackage[utf8]{inputenc} 
\usepackage{enumerate, theorem}
\usepackage{amsmath, amsfonts, amssymb}
\usepackage[height=22.5cm, width=15cm]{geometry}
\usepackage[all]{xy}
\usepackage{mathrsfs, yfonts}
\usepackage{graphicx}
\DeclareMathOperator{\M}{\frak{m}}
\DeclareMathOperator{\C}{\mathbb{C}}
\DeclareMathOperator{\R}{\mathbb{R}}
\DeclareMathOperator{\N}{\mathbb{N}}
\newcommand{\parag}[1]{\paragraph{\sc{#1.}} }

\DeclareMathOperator{\id}{id}

\DeclareMathOperator{\Sym}{Sym}

\DeclareMathOperator{\vol}{vol}

\newtheorem{thm}{Theorem}[subsection]
\newtheorem{defn}[thm]{Definition}
\newtheorem{cor}[thm]{Corollary}
\newtheorem{prop}[thm]{Proposition}
\newtheorem{lemma}[thm]{Lemma}

\begin{document}

\title{ Cycles of finite type} 
\author{Daniel Barlet\footnote{Institut Elie Cartan, G\'eom\`{e}trie,\newline
Universit\'e de Lorraine, CNRS UMR 7502   and  Institut Universitaire de France.} and J\'on Magn\'usson\footnote{Department of Mathematics, School of Engineering and Physical Sciences, \newline University of Iceland.} }

\maketitle

\tableofcontents

\parag{AMS Classification} 
32-02, 32-C15- 32 C25, 32 H 02, 32 H 04, \\ 32 K 05,  32 K 12.

\parag{Key words} Finite Type Cycle \ -- Quasi-Proper Map \ -- Holomorphic Fiber Map \\ 
-- Geometric $f$-Flattening \ -- Meromorphic Quotients.

\parag{Abstract}
The aim of this book is to show that the use of $f$-analytic families of finite type cycles (cycles having finitely many irreducible components, but not compact in general) in a given complex space  may be useful in complex geometry, despite the fact that the corresponding functor is not, in general, representable, in contrast to the compact case. This study leads to the notion of strongly quasi-proper map which is characterized by the existence of a geometric $f$-flattening which is a generalization of the Geometric Flattening Theorem for proper holomorphic maps. As applications we prove an existence theorem for meromorphic quotients of reduced complex spaces and a generalization of the classic Stein factorization.

\newpage
\section*{\centerline{FOREWORD}}

{\bf About different notions of flatness for a surjective holomorphic map}

\bigskip

\hfill{\it Quelques platitudes en hors-d'oeuvre}

\bigskip

Classically three notions of "flatness" have been considered in complex geometry  for  a holomorphic surjective map $\pi : M \to N$ between irreducible complex spaces 
\begin{enumerate}
\item In {\em algebraic} sense: the map $\pi$ is flat if  the sheaf $\mathcal{O}_M$ is flat as a module over $\pi^*(\mathcal{O}_N)$.
\item In {\em geometric} sense: the map $\pi$  is flat if there exists an analytic family of $n$-cycles (where $n := \dim M - \dim N$)  $(X_y)_{y \in N}$ parametrized by $N$ such that for each $y$ in $N$ we have
$\pi^{-1}(y) = \vert X_y\vert$, and also that $X_y$ is reduced for generic $y$ in $N$.
\item In {\em topological} sense:  the map $\pi$ is flat if it is an open map.
\end{enumerate}
Note that $1)$ implies $2)$\footnote{this implication is not trivial; see \cite{[BM.2]} chapter X.} which implies $3)$.\\
Of course, as such flat maps are considered in the respective points of view as an "ideal situation" an important question is to know if a given holomorphic surjective map can be "slightly transform" to such  a flat map. This leads to the question of existence of  a flattening.\\
To find a flattening in the sense $1)$, $2$ or $3)$  for $\pi$ means for us  to find a (proper) modification $\tau : \tilde{N} \to N$ such that the strict transform $\tilde{\pi} : \tilde{M} \to \tilde{N}$ of $\pi$ by $\tau$ (so $\tilde{M}$ is the union of irreducible components of $M\times_N\tilde{N}$ which dominate $M$) is flat in the sense of $1)$, $2)$ or $3)$.

\parag{The proper case} In his  famous paper \cite{[H.75]} H. Hironaka solves the problem of finding an "algebraic" flattening for any proper surjective holomorphic map  with a nice (succession of blow-ups with smooth centers) modification.\\
Few years after, D. Barlet in \cite{[B.78]} proved the existence of a canonical "geometric" flattening  for any proper surjective morphism (but with few information on the corresponding modification). \\
In the proper case, the "topological flattening" is not very interesting as it is already obtained by the previous results.
\parag{Non proper case} The first attempt without the  proper hypothesis (but of course with some  rather restrictive assumptions)  was made  in the sense of $3)$ by H.Grauert (see \cite{[G.86]}) in the middle of the  eighties and continued by his student B. Siebert (see \cite{[Si.93]} and \cite{[Si.94]}) in the  early nineties.  But it seems that their results are not so easy to use.

\parag{The quasi-proper case} At the end of the nineties  D. Mathieu consider the setting  of {\em geometric $f$-flat} maps which is defined by adding to  the  geometrically flat condition the assumption that $\pi$ is {\bf quasi-proper}. Then, due to the fact that the strict transform of a quasi-proper map by a modification is not always quasi-proper, he introduces an extra condition on a holomorphic  quasi-proper surjective map, which is called now-days "{\bf strongly quasi-proper map}", and then proves an existence theorem for a $f$-geometric flattening  (meaning that the flattening produced is a quasi-proper geometrically flat map)  for this class of holomorphic maps (see \cite{[Mt.00]}).\\

\parag{Finite cycles's space} Then these notions of quasi-proper (and strongly quasi-proper) maps put in light the fact that to consider the space $\mathcal{C}_n^{loc}(M)$ of all (closed) $n$-cycles in a given complex space $M$ was not the good point of view for the problem of existence of a geometric flattening in the quasi-proper case. This leads to introduce the space $\mathcal{C}_n^f(M)$ of finite type $n$-cycles (so cycles with finitely many irreducible components) with a topology which is stronger than the topology induced by $\mathcal{C}_n^{loc}(M)$, and to define the corresponding  notion of $f$-analytic family of finite type $n$-cycles in $M$ parametrized by a reduced complex space. In fact the only new condition on an analytic family $(X_s)_{s \in S}$  of $n$-cycles in $M$ to be an $f$-analytic is that its set-theoretic graph 
$$\vert G\vert := \{ (s, x) \in S \times G \  / \ x \in \vert X_s\vert \} $$
is quasi-proper over $S$ (of course this implies that each cycle $X_s$ is of finite type, but it asks more).
Then these tools allow to give a very simple reformulation of D. Mathieu result.

\begin{thm} 
Any strongly quasi-proper  surjective holomorphic map $\pi : M \to N$ between irreducible complex spaces admits a  global canonical $f$-flattening.
\end{thm}

In conclusion, this point of view makes appear a new kind of "flat" morphisms, the geometrically $f$-flat maps defined as follows: a holomorphic surjective map is $f$-geometrically flat if it is quasi-proper and equidimensional and if there exists a $f$-analytic family $(X)_{y \in N}$ of finite type $n$-cycles in $M$ such that $\vert X_y\vert = \pi^{-1}(y), \forall y \in N$ with $X_y $ reduced for generic $y \in N$. Then  one of the  characterization of strongly quasi-proper maps  which are  given in the present monograph is that a holomorphic surjective map $\pi : M \to N$ admits a  geometric $f$-flattening if and only if it is strongly quasi-proper. But we also show that the class of strongly quasi-proper holomorphic maps has several interesting stability properties and may be useful, for instance, to produce meromorphic quotients.

\chapter*{Introduction}

The space of compact analytic cycles has been used successfully for the last decades to study proper holomorphic maps  in complex geometry, see \cite{[B.99]}. It appears that one important idea is to consider, for a  holomorphic proper surjective map $\pi : M \to N$ between two irreducible complex spaces, its fiber map which is a meromorphic family of compact $n$-cycles ($n := \dim M - \dim N$) in $M$ parametrized by $N$. This fiber map is obtained by sending the generic point $y$ in $N$ to the reduced $n$-cycle given by the set-theoretic fiber $\pi^{-1}(y)$ of $\pi$ at $y$. Of course, the fact that this  family of compact $n$-cycles in $M$ is given by a meromorphic map (in the usual sense:  between reduced complex spaces) from $N$ to $\mathcal{C}_{n}(M)$ (which is, in a natural way, a reduced complex space locally of finite dimension) is an important tool for many applications (see  \cite{[BM.1]} ch.IV or \cite{[BM.2]}  for instance).\\
We shall explain, in this article, that, despite the fact that for non compact cycles the functor which associates to a reduced complex space $S$ the set of analytic families of $n$-cycles in $M$ parametrized by $S$ is, in general, neither representable in the category of reduced complex spaces (see nevertheless  \cite{[B.17]}) nor in the category of Banach analytic sets\footnote{This functor was introduced in  \cite{[B.75]} (see also \cite{[BM.2]}) but we shall consider a different one here which co\"incides with the classical one for compact $n$-cycles. It associates to a reduced complex space $S$ the set of $f$-analytic families of (finite type) $n$-cycles in $M$ parametrized by $S$.}, it is interesting and useful to consider  fiber maps (holomorphic and meromorphic) for a rather large class of surjective holomorphic maps $\pi : M \to N$ between irreducible complex spaces.\\
Let us explain now what are the problems coming from the non compactness of the cycles we are looking at.\\
The  "local charts" used  to classify $n$-cycles nearby a given $n$-cycle $X_0$ in a complex space $M$ are given by  $n$-scales on $M$ adapted to $X_0$. They allow to obtain Banach analytic classifying sets for local  analytic families of $n$-cycles in $M$. But when we consider non compact cycles, a finite collection of such $n$-scales cannot control what happens globally in $M$ even near $X_0$. And considering  countably many scales is not compatible with the Banach analytic setting. So we have to face two problems:
\begin{enumerate}
\item 
How can we  ensure, for a given cycle $X_0$, that near infinity in $M$ there are no "extra"  irreducible components belonging to  cycles arbitrarily close to $X_{0}$ ?\\
Note that such far away "extra" irreducible components  cannot be detected  with a finite number of scales, which, by definition, stay in a compact subset of $M$.
\item 
Assuming that no ``extra irreducible component''  appears near infinity for nearby cycles in a neighborhood of a given finite type cycle $X_{0}$ (see the definition below), how to control all these  nearby cycles  with only a finite number of  scales adapted to $X_{0}$ ?
\end{enumerate}

The first problem is solved by restricting ourself to the subset
 $$i : \mathcal{C}_{n}^{f}(M) \hookrightarrow  \mathcal{C}_{n}^{loc}(M)$$
 of {\bf finite type cycles}, meaning that we consider only closed $n$-dimensional cycles in $M$ having only {\bf finitely many irreducible components}, and also by choosing a {\bf new topology} on this subset which is stronger than the topology induced by the inclusion $i$, adding to the basis of open sets of the  topology induced by the inclusion $i$ the finite intersections of subsets of the following type
  $$\Omega(W) := \{X \in \mathcal{C}_{n}^{f}(M) \ / \ {\rm any \ irreducible \ component \ of \  X \  meets \ W} \}$$
   where $W$ is a relatively compact open set in $M$.\\
    Let $S$ be a Hausdorff topological space. We shall say that a family of finite type $n$-cycles $(X_{s})_{s \in S}$ in $M$  is   {\bf $f$-continuous} when the corresponding classifying map $\varphi : S \to \mathcal{C}^{f}_{n}(M)$ is continuous. We shall discuss below how to describe $f$-continuous families of finite type $n$-cycles with this new topology.\\
   
   The second problem will be solved by the following rather strong analytic extension result.   
   \begin{thm}\label{cycles}
Let $M$ be a complex space and fix a non negative integer  $n$. Then consider an $f$-continuous family  $(X_{s})_{s\in S}$ of finite type  $n$-cycles in  $M$ parametrized by a reduced complex space  $S$. Fix a point $s_{0}$ in $S$ and assume that there exists an open set $M'$ in $M$ meeting each irreducible component of $\vert X_{s_{0}}\vert$  and such that the family $(X_{s}\cap M')_{s \in S}$ is analytic at  $s_{0}$\footnote{In the usual sense, so using $\mathcal{C}_n^{loc}(M')$.}. Then  the family  $(X_{s})_{s\in S_{0}}$ is  analytic at $s_{0}$\footnote{ In fact in an open neighborhood of $s_0$, see the Analytic Extension Theorem \ref{cycles}.}.
\end{thm}

Note that this result is false in general when the parameter space is a Banach analytic set which is singular and not locally finite dimensional. See a counter-example in \cite{[BM.2]}   Chapter V section  2.4.\\
But clearly,  the previous theorem will be enough in order to solve our second problem at least when we are dealing with a $f$-continuous family of of finite type  $n$-cycles in  $M$ parametrized by a reduced complex space.\\

Let us now return to the signification of the new topology introduced on $\mathcal{C}^{f}_{n}(M)$.\\

The first remark is that, with this topology, the tautological family of finite type $n$-cycles parametrized by $\mathcal{C}_{n}^{f}(M) $ has a set-theoretic graph $\vert G\vert \subset \mathcal{C}_{n}^{f}(M) \times M$ which is {\bf quasi-proper}\footnote{The definition of a quasi-proper map is given below in such a context (not locally compact target space).} on $\mathcal{C}_{n}^{f}(M) $ by the natural projection. Moreover, it is easy to see that a continuous family of $n$-cycles in $M$ parametrized by a Hausdorff topological space $S$ such that all of them are of finite type, (by definition, this is simply a continuous map $\varphi : S \to \mathcal{C}_{n}^{loc}(M)$ taking values in $\mathcal{C}_{n}^{f}(M)$) factorizes as $i\circ \psi$ where $\psi : S \to \mathcal{C}_{n}^{f}(M) $ is continuous for the topology defined above, if and only if the graph of $\varphi$ is quasi-proper on $S$. And then, we shall say that an analytic family\footnote{In the usual sense; see \cite{[BM.1]} Chapter IV section 3.}  of finite type $n$-cycles, parametrized by a reduced complex space $S$, is an {\bf  $f$-analytic family of $n$-cycles} (of finite type) in $M$ when its set-theoretic graph $\vert G \vert \subset S \times M$ is quasi-proper over $S$. \\

This is the first step in defining an {\em analytic structure} (in a weak sense) on the topological space $\mathcal{C}^{f}_{n}(M)$: we say that a  map $\varphi : S \to \mathcal{C}^{f}_{n}(M)$, defined on a Banach analytic set $S$ is {\bf holomorphic} when the corresponding family is a $f$-analytic family of $n$-cycles in $M$. \\

To compare with the case of compact $n$-cycles, recall that an analytic family of $n$-cycles $(X_{s})_{s \in S}$ in $M$, parametrized by a reduced complex space $S$, is a ``proper analytic family of compact $n$-cycles in $M$'' if and only if its  set-theoretic graph $\vert G\vert \subset S \times M$ is proper over $S$. So we see that  we simply replace the properness condition of the graph in the compact cycles case  by the quasi-properness condition  of the graph in the case of finite type $n$-cycles.\\ 

This leads, following the case of compact cycles (see the definition IV.9.1.2 in \cite{[BM.1]} of  geometrically flat proper map), to consider  the morphisms of irreducible complex spaces $\pi : M \to N$ which are surjective, with fibers of pure dimension $n$ such that there exists a holomorphic  $f$-fiber map $\varphi : N \to \mathcal{C}^{f}_{n}(M)$. That is to say that we ask for a $f$-analytic family $(F_{y})_{y\in N}$,  of (finite type) $n$-cycles parametrized by $N$,  with a quasi-proper set-theoretic graph, such that, for {\em any} \, $y \in N$, we have the set-theoretic equality  $\pi^{-1}(y) = \vert F_{y}\vert$ and such that for $y$ generic in $N$ the cycle $F_y$ is reduced. These morphisms we be called {\bf geometrically $f$-flat maps}.\\
Note that, in fact, the condition of quasi-properness of the set-theoretic  graph (which is  isomorphic to $M$ in such a case) is equivalent to the quasi-properness of $\pi$, so  a geometrically $f$-flat map $\pi$ is simply a quasi-proper  map which possesses an $f$-analytic generically reduced  fiber map.\\

Before explaining a major ingredient of this work, let us recall that in \cite{[Ku]} Kuhlmann proved that the image of a {\it semi-proper}\footnote{See Definition \ref{semi-proper 0} below.}  holomorphic map between reduced complex spaces is an analytic subset of the target space\footnote{This is a generalization of the classical Remmert's Direct Image Theorem.}. In \cite{[Mt.99]} D. Mathieu generalized Kuhlmann's theorem to the case where the target space is an open subset of a complex Banach space. Using this theorem and the strong analytic extension given in section IV 3, we prove that the image of a semi-proper holomorphic map $\varphi : S \to \mathcal{U}$ from a reduced complex space $S$ to an open subset $\mathcal{U}$ of $\mathcal{C}_n^f(M)$ is a reduced complex space of $\mathcal{U}$ (see Definition IV 2.1.3 (vii)). This result has its roots in Mathieu's thesis (see \cite{[Mt.99]} and \cite{[Mt.00]}) and a weaker version of our result  was proved and applied in \cite{[B.08]},\cite{[B.13]} and \cite{[B.18]}.

But, in general, a quasi-proper surjective morphism $\pi : M \to N$ between two irreducible complex spaces is only generically equidimensional, so $f$-geometrically flat  on the complement of a closed analytic subset $\Sigma \subset N$ with no interior points in $N$: take for $\Sigma$  the union of the locus of ``big'' fibers (which is a closed analytic subset with no interior points in $N$, thanks to the quasi-properness of $\pi$) and of the subset of non normal points of $N$. Then we have on the dense  Zariski open set $N \setminus \Sigma$   a holomorphic reduced fiber map
$$ \varphi : N \setminus \Sigma \to \mathcal{C}^{f}_{n}(M)$$
and the question is now: 
\begin{itemize}
\item Is this map $\varphi$ meromorphic along $\Sigma$, as it is the case when  $\pi$  is a proper map, replacing  $\mathcal{C}^{f}_{n}(M)$ by $\mathcal{C}_{n}(M)$ ?
\end{itemize}
The answer if definitely NO !\\
 We shall give an example below (see example 1 in  Chapter V section 4). \\
 
  The reason for this is that, unlike $\mathcal{C}_n(M)$, the space $\mathcal{C}_n^f(M)$ is not  a complex space.\\
 Let us precise what we mean by a  meromorphic fiber map in this context. In the situation above of  a quasi-proper surjective morphism $\pi : M \to N$ between two irreducible complex spaces, we want to find a (proper) modification $\tau : \tilde{N} \to N$ with center contained in $\Sigma$  and a holomorphic map \ $\tilde{\varphi} : \tilde{N} \to \mathcal{C}^{f}_{n}(M)$ such that the restriction of  \ $\tilde{\varphi} $ \ to $\tilde{N} \setminus \tau^{-1}(\Sigma) \simeq N \setminus \Sigma$ gives the $f$-fiber map of $\pi$ restricted to $N \setminus \Sigma$. \\
This implies that, defining $\bar \Gamma$ as the closure in $N \times \mathcal{C}^{f}_{n}(M)$ of the graph $\Gamma$ of the  fiber map  $\varphi : N \setminus \Sigma \to \mathcal{C}^{f}_{n}(M)$, the projection $pr :\bar  \Gamma \to N$ is proper\footnote{Let us prove this point, because it lies at  the heart of the problem which comes up when we consider surjective quasi-proper maps with "big" fibers:\\
 To prove the properness of  $pr : \bar \Gamma \to N$, if such a modification $\tau$ exists, consider a compact set $K$ in $N$. Then $\tau^{-1}(K)$ is a compact set in $\tilde{N}$ and then $K\times \tilde{\varphi}(K)$ is compact in $N \times \mathcal{C}^{f}_{n}(M)$. The point is that this compact set contains $pr^{-1}(K)$ because if $(y, X)$ is in $pr^{-1}(K)$ with $y \in \Sigma$ there exists a sequence $(y_{\nu})_{\nu \geq 0}$ in $N \setminus \Sigma$ converging to $y$ such that the sequence  $(X_{\nu}:= \tilde{\varphi}(y_{\nu}))_{\nu \geq 0}$ converges to $X$ in $\mathcal{C}^{f}_{n}(M)$. By taking a sub-sequence if necessary, we may assume that $(y_{\nu})_{\nu}$ converges to $\tilde{y}$ in $\tilde{N}$. Then $X = \tilde{\varphi}(\tilde{y})$ and $\tau(\tilde{y}) = y$ proving our claim.}.\\

In fact, what happens is the fact that the quasi-properness of $\pi$ is not a strong enough condition, in general, to ensure the properness of the projection of $\bar \Gamma$ on $N$ when $N$ has big fibers along $\Sigma$. This means, assuming that $\pi : M \to N$ is a quasi-proper surjection between irreducible complex spaces, that for some $y \in \Sigma$ such that  $\pi^{-1}(y)$ has dimension strictly bigger than $n$, it may happen that, for a sequence  $(y_\nu)_{\nu \in \mathbb{N}}$ in $N \setminus \Sigma$ converging to $y$ the subset 
$\{\varphi(y_\nu) \}$ of $\mathcal{C}_n^f(M)$ is not relatively compact.\\
Another way to explain this phenomenon is to observe  that the strict transform of a quasi-proper surjective map $\pi : M \to N$ between irreducible complex spaces by  a  (proper) modification of $N$ may  no longer be  a quasi-proper map.\\
This leads to the notion of a {\bf strongly quasi-proper} (or SQP map for short) which is a quasi-proper map $\pi : M \to N$ having the property that the closure of of the graph of the reduced  holomorphic fiber map $N \setminus \Sigma \to \mathcal{C}_n^f(M)$ is proper over $N$. Every {\em equidimensional} quasi-proper  map has this property, but this is no longer true, in general, when the map admits big fibers.

 We shall show that, in general,  this notion corresponds exactly to the quasi-proper surjective morphisms admitting a meromorphic $f$-fiber map in the sense described above (but remark that it is not evident that $\bar \Gamma$ is a finite dimensional analytic subset in $N\times \mathcal{C}^{f}_{n}(M)$). \\
  We prove that SQP maps have good functorial properties. Moreover we show that the strict transform of an SQP map by a modification of the target space is an SQP map. But quasi-proper maps do not, in general, have this stability property.\\
 
In fact we shall prove that a holomorphic surjective quasi-proper map $\pi : M \to N$ between irreducible complex spaces is strongly quasi-proper if and only if it admits a geometric $f$-flattening. This means that there exists a (proper) modification of $N$, 
 $\tau : \tilde{N} \to N$, such that the strict transform by $\tau$ of $\pi$, $  \tilde{\pi} : \tilde{M} \to \tilde{N}$, is a geometrically $f$-flat map.\\
Note that in such a case we may compose the $f$-fiber map $\tilde{\varphi}\colon \tilde{N} \rightarrow \mathcal{C}_n^f(\tilde{M})$ with the direct image by the proper projection $p : \tilde{M} \to M$ for finite type $n$-cycles 
$$ p_* : \mathcal{C}_n^f(\tilde{M}) \to \mathcal{C}_n^f(M) $$
in order to obtain a holomorphic map $ \varphi : \tilde{N} \to \mathcal{C}_n^f(M)$ corresponding to the meromorphic $f$-fiber map of $\pi$.\\

We give two applications of this study of SQP maps. The first one consists of proving that a meromorphic equivalence relation on a reduced complex space $M$  admits a meromorphic quotient if it is strongly quasi-proper (in other words if the natural projection of its graph onto  $M$ is an SQP map). The second application is a (optimal) generalization of the classical Stein factorization for a proper holomorphic map to the strongly quasi-proper maps.\\

We give now a brief description of the content of this monograph.\\

In chapter I, after proving the basic results on semi-proper maps we prove the generalization of Kuhlmann's  Direct Image Theorem with values in a Banach open set. This result was originally proved in D. Mathieu's thesis \cite{[Mt.99]}, but the proof we give here is considerably simpler.\\
In Chapter II we generalize the notion of {\bf quasi-proper} map and prove some fundamental results on such maps. To illustrate this notion, we give a proof of the direct image theorem with values in a Banach open set for a quasi-proper map defined on an irreducible complex space which  is a great deal shorter than the proof in the semi-proper case.
This uses an easy generalization of the classical Remmert-Stein  Theorem in a Banach open set. \\
In chapter III we discuss the topology on the space $\mathcal{C}_n^f(M)$ where $M$ is a given complex space. We compare relatively compact sets in $\mathcal{C}_n^{\rm loc}(M)$ and $\mathcal{C}_n^f(M)$  and we give a characterization of relatively compact subsets in these two topological spaces.\\
In chapter IV we introduce the notion of $f$-analytic families of finite type $n$-cycles in a given complex space  $M$ and we define a "weak analytic structure" on the space $\mathcal{C}_n^f(M)$ for any given complex space $M$ and any integer $n \geq 0$. This makes possible to define the notion of a  { \bf meromorphic map}  from a reduced complex space to $\mathcal{C}_n^f(M)$, notion which plays an important role in the following chapters.\\
The chapter V is devoted to $f$-geometrically flat maps ($f$-GF maps for short) and to strongly quasi-proper maps (SQP maps for short). The latter being introduced because the strict transform of a quasi-proper map by a modification of the target space is not quasi-proper in general.\\
 The main results in Chapter VI are the existence of meromorphic quotients  for strongly quasi-proper meromorphic  equivalence relations and a generalization of the classical Stein factorization.

\newpage

\parag{Terminology and Notations}
\begin{itemize}
\item
A {\em complex space} is assumed to be countable at infinity and consequently  second-countable and metrizable.
\item
An {\em analytic subset} of a complex space is assumed to be closed.
\item 
By convention, an  {\em irreducible}  complex space is reduced. An analytic cycle is irreducible when it has exactly one irreducible component which has multiplicity $1$. An irreducible component is non empty by definition. As the empty $n$-cycle $\emptyset[n]$  has no irreducible component, it is not an irreducible cycle.
\item
An  $n$-{\em scale} on a complex space $M$ is a triple $E = (U,B,j)$, where $j$ is a closed holomorphic embedding of an open subset $V$ of $M_{\rm red}$ into an open subset $W$ of a numerical space $\C^m$, $U$ and $B$ are relatively compact open polydiscs in $\C^n$ and $\C^{m-n}$ respectively such that $\bar{U}\times\bar{B} \subseteq W$.

We call $V$ the {\em domain} of the scale and we call $c(E) := j^{-1}(U\times B)$ the {\em center} of the scale.

Such a scale is said to be {\em adapted to} an analytic subset $X$ of $M$ if 
$$
j(X)\cap\left(\bar{U}\times\partial{B}\right) = \emptyset.\footnote{This is a generalization of the usual notion of an {\em adapted $n$-scale} since we skip the condition that $X$ is of pure dimension $n$.}
$$
Note that in this situation, either $X \cap j^{-1}(U \times B)$ is empty, or it has dimension at most equal to $n$. Moreover the map $\pi: j(X)\cap (U\times B) \to U$ is proper and has finite fibers.
\item
We say that $X = \sum\limits_{i\in I}k_iX_i$ is the canonical expression of an $n$-cycle $X$ if $k_i$ are positive integers and  $(X_i)_{i\in I}$ is a locally finite family of irreducible analytic subsets such that $X_i\neq X_j$ for $i\neq j$. 
\item
An open subset of a Banach space will be called a {\em Banach open set}.
\item
For a holomorphic map $\pi\colon M\rightarrow N$ between complex spaces and for  a natural number $k$ we put
$$
\Sigma_k(\pi) := \{x\in M\ /\ \dim_x\pi^{-1}(\pi(x)) \geq k\}
$$
\end{itemize}

\begin{itemize}
\item A {\bf Zariski open}  subset in a reduced complex space $M$  is, by definition, the complement of a (closed) analytic subset in $M$.
\item  
Recall that a {\bf modification} between two reduced complex spaces will always be a proper holomorphic map which induces an isomorphism between two dense Zariski open sets.
\item 
We say that a holomorphic map  $\pi : M \to N$ between two irreducible complex spaces is {\bf dominant} (or that $M$ {\bf dominates} $N$) if the image of $\pi$ contains a non empty open subset of $N$.\\
When $M$ is not irreducible, we say that $\pi$ is {\bf dominant} if every irreducible component of $M$ dominates $N$.
\item  We say that in a  reduced complex space $M$ a subset $T$ is $M$ in {\bf very general} if its complement in $M$  is a countable union of locally closed analytic subsets  with no interior point in $M$. So a countable intersection of very general subset is a very general subset.\\
 Remark that for any open set $M'$ in $M$ the intersection $M'\cap T$ is very general in $M'$ when $T$ is very general in $M$.\\
Also, if  $Z$ is a closed analytic subset with no interior point  in $M$ and assume that $T \subset M \setminus Z$ is very general in $M \setminus Z$. Then $T$ is very general in $M$.\\
\item Conversely, if   for a subset $T$ in $M$ and  for every $x \in M $ there exists an open neighborhood $U$ of $x$ in $M$ such that $U \cap T$ is very general in $U$, then $T$ is very general in $M$ because $M$ has a countable basis of open set. \\
\end{itemize}

\chapter{Semi-proper maps}

A very useful notion in topology is the notion of  a proper map, which is the relative notion of compactness. For instance, in a continuous family of compact cycles $(X_s)_{s \in S}$ in a given complex space $M$, parametrized by a Hausdorff topological space $S$, the projection of the set theoretic graph,
$$\vert G\vert := \{(s,x) \in S \times M \ / \ x \in \vert X_s\vert \}$$
of such a family is assumed to be proper on $S$.\\

As we want to generalize this notion to the case of a family of finite type $n$-cycles in $M$, we shall demand that the projection on $S$  of the set theoretic graph is {\bf quasi-proper}. The classical notion of quasi-proper map, which is discussed in  Chapter II  below, is not purely topological as it takes into account the fact that the fibers of the map are analytic subsets of $M$ in order to consider the irreducible components of the fibers.\\

It is remarkable that there exists a purely topological notion, the semi-properness, such that quasi-proper maps are always semi-proper and which gives a sufficient (topological) condition in order to obtain a Direct Image Theorem in the category of reduced complex spaces. This result was proved by N. Kuhlmann in the early sixties, generalizing Remmert's Direct Image Theorem which corresponds to the proper case. We shall give in section 4 of chapter IV  a generalization of 's result in the case of a semi-proper holomorphic map $f : N \to \mathcal{C}_n^f(M)$ where $M$ and $N$ are reduced complex spaces and where the space $\mathcal{C}_n^f(M)$ is equipped with its "weak Banach analytic set structure" which will be defined in section 2 of chapter IV. This result will be the main tool in several applications, for instance, in the theorem of existence of meromorphic quotients (see section 1 of chapter VI).\\

In the present chapter, after collecting some basic results on semi-proper maps, we prove in section 2  the generalization of 's result when the target space is a Banach open set, which is a crucial step for the case where the target is $\mathcal{C}_n^f(M)$.\\
 The case when the target is a Banach open set  is originally due to D. Mathieu in his thesis, see \cite{[Mt.99]}.\\

\section{Definition and basic properties}

Let us begin by recalling the standard definition of a semi-proper map.

\begin{defn}\label{semi-proper 0}
Let $S$ be a {\em locally compact} Hausdorff space, $T$ be a topological space and $f\colon S \rightarrow T$ be a continuous map. We shall say that $f$ is {\bf\em semi-proper at a point $t_0 \in T$} when there exists a neighborhood $T_0$ of $t_0$ in $T$ and a compact subset $L$ in $S$ such that we have the equality
$$ 
f(S) \cap T_0 = f(L)\cap T_0.
$$
We say that $f$ is {\bf\em semi-proper} when it is  semi-proper at every point in $T$.
\end{defn}

\begin{prop} \label{semi-proper 1} 
 Let $f\colon S \rightarrow T$ be a continuous map from a locally compact Hausdorff space $S$  to a {\em Hausdorff} space $T$. Then the following properties hold true: 
 \begin{enumerate}[(i)]
\item
If $f$ is semi-proper, then the induced map $f^{-1}(X)\rightarrow X$ is semi-proper for every closed subset  $X$  of $T$.
 \item 
The map $f$ is semi-proper if and only if $f(S)$ is closed in $T$ and the induced map $S \rightarrow f(S)$ is semi-proper.
\item
If the induced map $S \rightarrow f(S)$ is semi-proper, $f(S)$ is locally compact.
\item 
If $f$ is semi-proper, then the induced map $f^{-1}(T') \rightarrow T'$ is semi-proper, for every {\em locally closed} subset $T'$ of $T$. 
\end{enumerate}
\end{prop}

 {\sc Proof of $(i)$}
Let $X$ be a closed subset of $T$ and  $t_0$ be a point in \ $X$. Then we take a neighborhood $T_0$ of $t_0$ in $T$ and a compact subset $L$ of $S$ which satisfy $T_0\cap f(L) = T_0\cap f(S)$. Then, as $f^{-1}(X)$ is a closed subset of $S$, the subset $f^{-1}(X)\cap L$ of $f^{-1}(X)$ is compact and consequently we get
$$T_0\cap f(L\cap f^{-1}(X)) = T_0\cap f(L)\cap X = T_0\cap f(S)\cap X = T_0\cap f(f^{-1}(X)).$$
\medskip
{\sc Proof of $(ii)$} 
Suppose first that $f$ is semi-proper. Then, due to (i), it is enough to prove that $f(S)$ is closed in $T$. To do so take a point $t_0$ in $T\setminus f(S)$, an open neighborhood $T_0$ of $t_0$ in $T$ and a compact subset $L$ of $S$ which satisfy  the equality $T_0\cap f(L) = T_0\cap f(S)$. Then $t_0$ is not in the compact subset $f(L)$ of $T$ and, as $T$ is a Hausdorff space, $f(L)$ is a closed subset of $T$. It follows that $T_0\setminus f(L)$ \ is an open neighborhood of $t_0$ which does not intersect $f(S)$. So $f(S)$ is closed in $T$.\\
Conversely, assume that $f(S)$ is closed in $T$ and that the induced map $S \rightarrow f(S)$ is semi-proper. For any $t_0 \in T \setminus f(S)$ the open set $T_0 := T \setminus f(S)$ satisfies the equality $f(S) \cap T_0 = f(\emptyset) \cap T_0 = \emptyset$ and so $f$ is semi-proper at $t_0$.\\
For any  $t_0 \in f(S)$ there exists an open neighborhood $\Theta_0$ of $t_0$ in $f(S)$ and a compact set $L$ in $S$ with $f(S) \cap \Theta_0 = f(L) \cap \Theta_0$. Now choose an open set $T_0$ in $T$ such that $T_0 \cap f(S) = \Theta_0$. It satisfies $T_0 \cap f(S) = f(L) \cap T_0$ and $f$ is semi-proper at $t_0$.

\parag{Proof of $(iii)$}For any $t_0 \in f(S)$ there is an open neighborhood $T_0$ of $t_0$ in $T$ and a compact set $L$ in $S$ with $f(S) \cap T_0 = f(L) \cap T_0$. Then $f(L)$ is a compact neighborhood of $t_0$ in $f(S)$.

\parag{Proof of $(iv)$} 
Let $T'$ be a locally closed subset of $T$. Then there exist a closed subset  $X$ and an open subset $V$ of $T$ such that $T' = X\cap V$. As the subset  $f^{-1}(T') = f^{-1}(X)\cap f^{-1}(V)$ is locally compact it is enough, due to (i), to show that the induced map $f^{-1}(V) \rightarrow V$ is semi-proper. To this end we take an arbitrary point $t_0$ in $V\cap f(S)$ and a compact subset $L$ of $S$ such that $f(L)$ is a neighborhood of $t_0$ in $f(S)$. Since $f(S)$ is a locally compact Hausdorff space there exists a compact neighborhood $C$ of $t_0$ in $V\cap f(S)$. It follows that $f^{-1}(C)\cap L$ is a compact subset of $f^{-1}(V)$ and $f(f^{-1}(C)\cap L) = C\cap f(L)$ is a compact neighborhood  of $t_0$ in $V\cap f(S)$. This shows that the induced map $f^{-1}(V) \rightarrow V\cap f(S)$ is semi-proper so, by (ii), the induced map $f^{-1}(V) \rightarrow V$ is semi-proper since $V\cap f(S)$ is closed in $V$.$\hfill \blacksquare$\\

\begin{prop} \label{semi-proper 1bis} 
 Let $f\colon S \rightarrow T$ be a continuous map from a locally compact Hausdorff space $S$  to a {\em Hausdorff} space $T$. Then the following properties hold true: 
\begin{enumerate}[(a)]
\item
Let $X$  be a closed subset of $T$  which contains  $f(S)$. Then the induced map $S\rightarrow X$ is semi-proper if and only if $f$ is semi-proper.  

\item
Suppose that $f$ is an open map. Then the induced map $S\rightarrow f(S)$ is semi-proper. Moreover, assuming that $f$ is an open map, then $f\colon S\rightarrow T$ is semi-proper if and only if it $f(S)$ is a union of connected components of $T$.
\item
If $f$ is proper then it is semi-proper.
\item
Suppose that $S$ is countable at infinity, $T$ is first countable and that $f$ is a closed map. Then $f$ is semi-proper.
\end{enumerate}
\end{prop}

 Notice that $(d)$  is no longer true if the hypothesis that $S$ is {\em locally compact} is skipped. For instance, if $E$ is an infinite dimensional Banach space, the identity map of $E$ is proper (i.e. a closed map with compact fibers) but does not satisfy the condition of  Definition \ref{semi-proper 0}.

\bigskip

{\sc Proof of $(a)$} 
If $f$ is semi-proper, then the induced map $S\rightarrow X$ is semi-proper thanks to (i) in Proposition \ref{semi-proper 1}. \\
Conversely, if the induced map $S\rightarrow X$ is semi-proper, then (ii) in  Proposition \ref{semi-proper 1} implies that the induced map $S\rightarrow f(S)$ is semi-proper and $f(S)$ is a closed subset of $X$. It follows that $f(S)$ is a closed subset of $T$ and consequently $f$ is semi-proper thanks  again to (ii)  in  Proposition \ref{semi-proper 1}.

\parag{Proof of $(b)$} 
To prove the first assertion we consider a point $t_0$ in $f(S)$. Then we pick a point $s_0$ in $f^{-1}(t_0)$ and a compact neighborhood $L$ of $s_0$ in $S$. Then $f(L)$ is a compact neighborhood of $t_0$ in $T$ since $f$ is an open map. So $f: S \to T$ is semi-proper at $t_0$.\\ 
Let us  prove the second assertion. If  $f: S \to T$ is semi-proper, then $f(S)$ is both open and closed in $T$ and consequently a union of connected components of $T$.\\
Conversely, suppose that $f$ is an open map and  the subset $f(S)$ of $T$ is a union of connected components of $T$. Then $f: S \to f(S)$ is semi-proper by the direct part and $f(S)$  is closed in $T$, so, by (ii) in Proposition \ref{semi-proper 1}, $f$ is semi-proper.

\parag{Proof of $(c)$} 
By definition a proper map is a continuous map which is closed and such that each fibers compact. So $f(S)$ is closed. It is then  enough to prove that $f$ is semi-proper at each point in $f(S)$ thanks to $(ii)$ in Proposition \ref{semi-proper 1}. Choose any $t_0 \in f(S)$ and let $U$ be a relatively compact open set containing the compact set $f^{-1}(t_0)$. Then $F := S \setminus U$ is a closed set in $S$ and so $f(F)$ is closed in $T$. Define $T_0 := T \setminus f(F)$. This open set contains $t_0$ and let $L := \bar U$. Then we have $T_0 \cap f(S) = f(L) \cap T_0$ because if $t$ is in $T_0 \cap f(S)$ there exists a point $s \not\in F$, so $s \in U$, with $f(s) = t$. As $\bar U$ is a compact subset of $S$ the map $f$ is semi-proper at $t_0$.

\parag{Proof of $(d)$} 
Take any $t_0 \in T$ and let $(W_n)_{n \in \mathbb{N}}$  a basis of open neighborhoods of $t_0$ in $T$. Let $(L_m)_{m \in \mathbb{N}}$ be an increasing  exhausting sequence of compact sets in $S$ with $t_0 \in f(L_0)$. Assume that for each $n \geq 0$ there exists a point $t_n = f(s_n)$ with $t_n \in W_n \setminus f(L_n)$. So $s_n$ is not in $ L_n$ and $f(s_n)$ is in $W_n$ and not equal to $t_0$. Then the  set $F := \{s_n, n\geq 0 \}$ is closed but $f(F)$ is not closed because it  closure contains $t_0$ which is not in $f(F)$. Contradiction. So there exists an integer $n_0$ such that $W_{n_0} \cap f(S) = W_{n_0} \cap f(L_{n_0})$ concluding the proof.$\hfill \blacksquare$\\

We shall need two more topological  lemmas.

\begin{lemma}\label{semi-propre 2}
Let S be a locally compact Hausdorff space and T be a Hausdorff space. Let $f\colon S \rightarrow T$ be  a continuous map. Then the semi-properness of $f$ at the point $t_0 \in T$ is equivalent to the following:
\begin{itemize}
\item
There exists an open neighborhood $T_1$ of $t_0$ in $T$ and a subset $L$ of $f^{-1}(T_1)$ such that the map $g\colon L\rightarrow T_1$ induced by $f$ is proper and 
\begin{equation*}
f(S) \cap T_1 \ = \ f(L). \tag{$*$} 
\end{equation*}
\end{itemize}
\end{lemma}

\parag{Proof} Suppose first that $f$ is semi-proper. Then there exists an open neighborhood  $T_0$ of $t_0$ and a compact subset $K$ in $S$ such that 
$$ 
f(S) \cap T_0 = f(K)\cap T_0.
$$
Then $T_1 := T_0$ and $L := K \cap f^{-1}(T_0)$ satisfy condition $(*)$.

\smallskip
Conversely, assume that $f$ satisfies condition $(*)$. As $S$ is locally compact the compact subset $f^{-1}(t_0)\cap L = g^{-1}(t_0)$ admits a compact neighborhood  $C$ in $f^{-1}(T_1)$. Then there exists an open neighborhood $V$ of $t_0$ in $T_1$ such that $f^{-1}(V)\cap L = g^{-1}(V)$ is contained in $C$ since $g$ is a closed map. It follows that $f^{-1}(V)\cap L\subseteq f^{-1}(V)\cap C \subseteq f^{-1}(V)$ and consequently  $V\cap f(L) = V\cap f(C) = V\cap f(S)$.
\hfill $\blacksquare$

\begin{lemma}\label{semi-propre 3} 
Let $S$ and $T$ be first countable Hausdorff spaces with $S$ locally compact  and $f\colon S\rightarrow T$ be a semi-proper map. Then, for every open subset $V$ of $T$, the restriction $ \overline{f^{-1}(V)} \to T$ of $f$ to the subset $\overline{f^{-1}(V)}$ is a semi-proper map.
\end{lemma}

\parag{Proof}By continuity $f(\overline{f^{-1}(V)})\subseteq \overline{V}$ so thanks to (iii) of Proposition I.1.0.2 it is enough to prove that the induced map $\overline{f^{-1}(V)}\rightarrow \overline{V}$ is semi-proper.
Now let $y$ be a point in $\overline{V}$ and take a compact subset $K$ of $S$ such that $f(K)$ is a neighborhood of $y$ in $T$. In order to prove that the restriction of $f$ to $\overline{f^{-1}(V)}$ is semi-proper at $y \in f(S)$ it is enough to show that $f\left(\overline{f^{-1}(V)}\cap K\right)$ contains $\overline{V}\cap{\rm int}(f(K))$. So let $y_0$ be a point in $\overline{V}\cap{\rm int}(f(K))$. Then there exists a sequence $(y_n)_{n\geq 1}$ in $V\cap\,{\rm int}(f(K))$ which converges to $y_0$ and a sequence $(z_n)_{n\geq 1}$ in $f^{-1}(V)\cap K$ such that $f(z_n) = y_n$ for all $n$. By taking a subsequence we may suppose that $(z_n)_{n\geq 1}$ converges to a point $z$ in  $\overline{f^{-1}(V)}\cap K$ since $\overline{f^{-1}(V)}\cap K$ is compact. It follows that $f(z) = y_0$ and as $y_0$ is an arbitrary point in $\overline{V}$ this shows that $\overline{V}\cap{\rm int}(f(K))\subseteq f\left(\overline{f^{-1}(V)}\cap K\right)$. Hence the restriction of $f$ to $\overline{f^{-1}(V)}$ is a semi-proper map
\hfill{$\blacksquare$}

\section{Generalization of 's Direct Image Theorem}

The goal of this section is to prove the following theorem, which is a generalization of Kuhlmann's theorem for semi-proper holomorphic maps between reduced complex spaces. This result was proved by D. Mathieu in his doctoral thesis (see University H. Poincar\'e, Nancy 1999). The main part of the thesis is contained in \cite{[Mt.00]}, but the proof of this result is not given in this article.\\
We present here a proof which is simpler  than  D.Mathieu's proof and our proof  is also self-contained in the sense that it  does not use Kuhlmann's theorem.

\begin{thm}\label{Kuhl-Banach}
Let $M$ be a reduced complex space, $\mathcal{U}$ an open subset of a Banach space $E$ and $\pi \colon M\rightarrow \mathcal{U}$  a semi-proper holomorphic map. Then $\pi(M)$ is a reduced complex subspace\footnote{This means that, endowed with the sheaf of holomorphic functions obtained from $\mathcal{U}$, $\pi(M)$ is a reduced complex space. Moreover, $\pi(M)$ is locally contained in a finite dimensional  sub-manifold of $\mathcal{U}$, thanks to the {\em Enclosability Theorem} (see Theorem III.7.4.1 in \cite{[BM.1]} or \cite{[e]}).} of 
$\mathcal{U}$.
\end{thm}

To prove the theorem we need some preliminary results.

\parag{Notation} For a holomorphic map $\pi$  from a reduced complex space $M$ into a Banach open set, we put
$$
\Sigma_m(\pi) := \{x\in M\ |\ \dim_x\pi^{-1}(\pi(x)) \geq m\}
$$
for every integer $m\geq 0$.

\begin{prop}\label{Kuhl-banach.0}
Let \ $\pi\colon M \rightarrow \mathcal{U}$\ be a holomorphic map from a reduced complex space into a Banach open set. Let \ $x$ \ be a point in \ $M$ \ and put \ $q := \dim_x\pi^{-1}(\pi(x))$. Then there exists an open neighborhood \ $W$ \ of \ $x$ \  in \ $M$ \ and an open neighborhood $\mathcal{V}$ of $\pi(x)$  in $\mathcal{U}$ having the following properties:
\begin{enumerate}[(i)]
\item
\ $\dim_z\pi^{-1}(\pi(z))\leq q$ \ for all \ $z$ \ in \ $W$.
\item
\ $\pi(W\cap\Sigma_q(\pi))$ \ is a reduced complex subspace of $\mathcal{V}$.
\end{enumerate}
\end{prop} 

\parag{Proof} Take a $q-$scale $(U,B,j)$ adapted to $\pi^{-1}(\pi(x))$ near $x$. Let $pr \colon U\times B\rightarrow U$ be the natural projection and consider the holomorphic map
$$
g\colon j^{-1}(U\times B)\longrightarrow U\times \mathcal{U},\qquad z\mapsto (pr(j(z)),\pi(z)).
$$
Then $g^{-1}(g(x))$ is a finite subset of $j^{-1}(U\times B)$ so  there exists an open neighborhood $W_x$ of $x$ in $j^{-1}(U\times B)$, an open connected neighborhood $U_x$ of $pr(j(x))$ in $U$ and an open neighborhood $\mathcal{U}_x$ of $\pi(x)$ in $\mathcal{U}$ having the property that $g$ induces a proper map $h\colon W_x\rightarrow U_x\times\mathcal{U}_x$ with finite fibers (see  Proposition 2.1.6 in \cite{[e]}). It follows that, for each $z$ in $W_x$, the induced map $W_x\cap \pi^{-1}(\pi(z))\rightarrow U_x\times\{\pi(z)\}$ is proper with  finite fibers. Hence the analytic subset $W_x\cap \pi^{-1}(\pi(z))$ of $W_x$ is of dimension at most $q$ for each $z \in W_x$ proving the point $(i)$.\\
 Moreover the dimension of  $W_x\cap \pi^{-1}(\pi(z))$ is  equal to $q$ if and only if the restriction of $h$ to $W_x\cap \pi^{-1}(\pi(z))$  is surjective onto $U_x\times\{\pi(z)\}$. \\

Due to Remmert's Direct Image Theorem, generalized to the case where the target space is a Banach open set  (see \cite{[BM.1]}, ch. III, sect. 7), $h(W_x)$ is a reduced complex subspace of $U_x\times\mathcal{U}_x$. Hence, by shrinking $U_x$ around $\pi(j(x))$ and $\mathcal{U}_x$ around $\pi(x)$ and replacing $W_x$ by the inverse image of $U_x\times\mathcal{U}_x$ by $h$, we may assume that $h(W_x)$ is the zero set of a holomorphic map $\Phi\colon U_x\times\mathcal{U}_x\rightarrow F$ where $F$ is a Banach space. Now take a relatively compact open neighborhood $V$ of $\pi(j(x))$ in $U_x$ and let 
$$\tilde\Phi\colon \mathcal{U}_x\rightarrow {\rm Hol}(\bar{V},F)$$
 be the holomorphic map defined by $\tilde\Phi(\zeta) := \Phi(-,\zeta)$. The vanishing of this map at a point $z \in \mathcal{U}_x$ means that $V\times \{z\}$ is in $h(W_x)$, so that $z = \pi(x) $ with $\pi^{-1}(\pi(x))$ of dimension $q$, which implies  $x \in \Sigma_q(\pi)$. \\
 Then, since $U_x$ is connected, $\pi(W_x\cap\Sigma_q(\pi))$ coincides with $ \tilde\Phi^{-1}(0)$ and consequently 
 $\pi(W_x\cap\Sigma_q(\pi))$ is an analytic subset of $\mathcal{U}_x$. \hfill{$\blacksquare$}

\parag{Remark} The above proof implies in particular that \ $\Sigma_q(\pi)$ \ is a closed subset of \ $M$ \ for every \ $q$, because we proved that if $\dim_{x_0}(\pi^{-1}(\pi(x)) = q$ then there exists an open neighborhood $W$ of $x_0$ on which $\dim_x(\pi^{-1}(\pi(x))$ is at most equal to $q$ for each $x \in W$. So the complement of $\Sigma_{q+1}(\pi)$ in $M$ is open for each $q\geq 0$ (and $\Sigma_0(\pi) = M$).

\begin{prop}\label{Kuhl-banach.1} 
Let $\pi\colon M \rightarrow \mathcal{U}$ be a holomorphic map from a reduced complex space into a Banach open set and let $y$ be a point in $\pi(M)$. Suppose we have a compact subset $K$ of $M$ such that $\pi(K)$ is a neighborhood of $y$ in $\pi(M)$\footnote{Note that the existence of such a $K$ is equivalent to $\pi$ being semi-proper at $y$.} and put 
$$
p := \sup\limits_{x\in K\cap \pi^{-1}(y)} \{\dim_x\pi^{-1}(y)\}.
$$ 
Then there exists an open neighborhood $U$ of $\pi^{-1}(y)\cap K$ in $M$  and an open neighborhood $\mathcal{V}$ of $y$ in $\mathcal{U}$ having the following properties:
\begin{enumerate}[(a)]
\item
\ $\pi^{-1}(\mathcal{V})\cap K\ \subseteq \ U\ \subseteq \ \pi^{-1}(\mathcal{V})$
\item
\ $\pi(U\cap\Sigma_p(\pi))$ is a reduced complex subspace of $\mathcal{V}$.
\end{enumerate}
\end{prop}
\parag{Proof}
Take an open neighborhood $\mathcal{W}$ of $y$ in $\mathcal{U}$ such that $\pi(K)\cap\mathcal{W} = \pi(M)\cap\mathcal{W}$. 

Then by Proposition \ref{Kuhl-banach.0} there exists, for each $x$ in $\pi^{-1}(y)\cap K$, an open neighborhood $W_x$ of $x$ in $M$ such that $W_x\cap\Sigma_p(\pi) = \emptyset$ if $\dim_x\pi^{-1}(y) < p$, and such that $\pi(W_x\cap\Sigma_p(\pi))$ is a reduced complex space of an open neighborhood of $y$ in $\mathcal{U}$ when $\dim_x(\pi^{-1}(y) = p $. As $\pi^{-1}(y)\cap K$ is compact there exist $x_1,\ldots,x_l$ in $\pi^{-1}(y)\cap K$ such that the open set  $W := W_{x_1}\cup\cdots\cup W_{x_l}$ contains $\pi^{-1}(y)\cap K$. For each $j$ let $\mathcal{U}_j$ be an open neighborhood of $y$ in $\mathcal{W}$ such that $A_j := \pi(W_{x_j}\cap\Sigma_p(\pi))$ is a reduced complex subspace of  $\mathcal{U}_j$.  Then there exists an open neighborhood $\mathcal{V}$ of $y$ in  $\mathcal{U}_1\cap\cdots\cap\mathcal{U}_l$ satisfying $\pi^{-1}(\mathcal{V})\cap K\subseteq W_{x_1}\cup\cdots\cup W_{x_l}$. Put $U := \pi^{-1}(\mathcal{V})\cap(W_{x_1}\cup\cdots\cup W_{x_l})$. Then $U$ clearly satisfies condition $(a)$ and 
$$
\pi(U\cap\Sigma_p(\pi)) = \mathcal{V}\cap(A_1\cup\cdots\cup A_l)
$$
is a reduced complex subspace of $\mathcal{V}$; so $U$ satisfies condition $(b)$ also.
\hfill{$\blacksquare$}\\

In the situation of Theorem \ref{Kuhl-Banach}, let $C$ be an irreducible component of $M$. Recall that the maximal rank of $\pi_{|C}$ on the smooth part of $C$ is called the {\bf generic rank of $\pi$ on $C$}.  We observe that the generic rank of $\pi$ on $C$ is equal to
$$
\max\limits_{x\in C}(\dim C - \dim_x\pi^{-1}(\pi(x))). 
$$
The map $\pi$ is said to be of {\bf constant generic rank} if it has the same generic rank on every irreducible component of $M$.\\

It should be noted that if $\pi$ is of constant generic rank $n$ and if  $\pi(M)$ is a reduced complex space,  then $\pi(M)$ is of pure dimension $n$.\\

We shall now prove Theorem \ref{Kuhl-Banach} in a special case.

\begin{thm}\label{Kuhl-Banach.special case}
Under the assumptions of Theorem \ref{Kuhl-Banach}, suppose moreover that $\dim E < \infty$ and that $\pi$ is of constant generic rank. Then $\pi(M)$ is a reduced complex subspace of $\mathcal{U}$. 
\end{thm}

\parag{Proof} As $\pi(M)$ is closed in $\mathcal{U}$ and $\dim E < \infty$ it is enough to show that every point in $\pi(M)$ admits an open neighborhood $W$ in $\mathcal{U}$ such that $\pi(M)\cap W$ is an analytic subset of $W$. 

Let $y$ be a point in $\pi(M)$ and choose a compact subset $K$ of $M$ and an open neighborhood $V$ of $y$ in $\mathcal{U}$ having the property  that $\pi(K)\cap V = \pi(M)\cap V$. Let $M_1$ denotes the union of those irreducible components of $M$ which intersect $K$. Since $M_1$ is a reduced complex space having only finitely many irreducible components the fiber dimension of $\pi_{|M_1}$ is bounded. Moreover, thanks to point $(iii)$  in Proposition \ref{semi-proper 1}, the induced map $\pi^{-1}(V)\rightarrow V$ is semi-proper. But the pair $(V, K)$ gives the semi-properness of $\pi$ at each point in $V$, and as $K$ is contained in $M_1$ we conclude that the induced  map $\pi_1 : M_1 \cap \pi^{-1}(V)\rightarrow V$ is semi-proper. Then to prove that $\pi(M) \cap V$ is an analytic subset in $V$ it is enough to prove this result for the map $\pi_1$.

\smallskip
{\em Hence we may, without loss of generality, assume that the fiber dimension is bounded and we shall prove the theorem by induction on the maximal fiber dimension of the map $\pi$. }

\smallskip
Suppose that all fibers of $\pi$ are $0-$dimensional. Then $\Sigma_0(\pi) = M$ and in the setting of Proposition \ref{Kuhl-banach.1} we have that $\pi(U) = \mathcal{V}\cap \pi(K) =  \mathcal{V}\cap f(M)$ is a pure dimensional analytic subset of $\mathcal{V}$ whose dimension is equal to $\dim M$. Hence in this case $\pi(M)$ is a reduced complex subspace of $\mathcal{U}$.

\smallskip
{\em Now suppose that $f$ is of maximal fiber dimension $p \geq 1$ and suppose the theorem proven for all semi-proper holomorphic maps whose maximal fiber dimension is strictly less than $p$. }

Let $T$ be the  subset of all $y$ in $\pi(M)$ such that there exists a compact subset $K_y$ of $M$ which satisfies the following properties:
\begin{enumerate}
\item[($\alpha$)]
\ $\dim_x\pi^{-1}(y) < p$ for all $x$ in $\pi^{-1}(y)\cap K_y$
\item[($\beta$)]
\ $\pi(K_y)$ is a neighborhood of $y$ in $\pi(M)$.
\end{enumerate}
Observe that ($\alpha$) is equivalent to $\pi^{-1}(y)\cap K_y\cap\Sigma _p(\pi) = \emptyset$. This means that the fiber $\pi^{-1}(y)$ does not meet the compact set $ K_y\cap\Sigma _p(\pi) $. So for an open neighborhood $U$ of $y$  in $\mathcal{U}$, contained in the interior of $\pi(K_y)$  and small enough, we have $\pi(K_y)\cap U = \pi(M)\cap U$ and $\Sigma_p(\pi)\cap \overline{K_y\cap \pi^{-1}(U)} = \emptyset$. So $T$ is an open subset in $\pi(M)$. \\
Moreover the induced map $\pi^{-1}(U)\setminus\Sigma_p(\pi)\rightarrow U$ is semi-proper,  since $K_y\cap \pi^{-1}(U)$ does not meet $\Sigma_p(\pi)$ and is proper over $U$. So we may apply Lemma \ref{semi-propre 2}.\\
But its maximal fiber dimension is less than or equal to $p-1$. By the induction hypothesis it then follows that $\pi(M)\cap U$ is an analytic subset of $U$. Hence every point in $T$ admits an open neighborhood in $\mathcal{U}$ in which $\pi(M)$ is an analytic subset.\\

Now let $y$ be a point in $\pi(M)\setminus T$ and let $K$ be a compact subset of $M$ such that $\pi(K)$ is a neighborhood of $y$ in $\pi(M)$. Then for at least one $x$ in $\pi^{-1}(y)\cap K$ we have $\dim_x\pi^{-1}(y) = p$ and, by Proposition \ref{Kuhl-banach.1}, there exist  an open neighborhood $U$ of $\pi^{-1}(y)\cap K$ in $M$ and an open neighborhood $V$ of $y$ in $\mathcal{U}$ satisfying the conditions 
\begin{itemize}
\item
\ $\pi^{-1}(V)\cap K\ \subseteq \ U\ \subseteq \ \pi^{-1}(V)$
\item
\ $A := \pi(U\cap\Sigma_p(\pi))$ is a an analytic subset of $V$. 
\end{itemize}

Since $ \pi(K\cap\Sigma_p(\pi))\cap V\subset A$ we have $(\pi(M)\setminus T)\cap V\subseteq A$ and consequently 
\begin{equation*}\tag{$*$}
\pi(M)\cap (V\setminus A) = T\cap (V\setminus A).
\end{equation*} 

Let $n$ denote the (constant) generic rank of $\pi$. Then $T$ is either empty, in which case $A = \pi(M)\cap V$ and the proof is completed, or $\pi(M)\cap (V\setminus A)$ is an analytic subset of pure dimension $n$ of $V\setminus A$ and obviously  $\dim A \leq n$. Consider the decomposition  $A =  A_1\cup A_2$ where which  $A_1$ is the union of all $n-$dimensional irreducible components of $A$ and  $A_2$ is the union of the others. Then $(\pi(M)\cap (V\setminus A_1))\setminus A_2$ is an analytic subset of pure dimension $n$ in the open set $V \setminus A_1$. Since $\dim A_2 < n$ it follows, due to the Remmert-Stein Theorem\footnote{See section 3 or  Theorem 2.4.75 in \cite{[e]}}, that  the closure $X$ of $\pi(M)\cap (V\setminus A_1) \setminus A_2$ in $V \setminus A_1$  is an analytic subset of pure dimension $n$ of $V \setminus A_1$ contained in $\pi(M)$. Let $Y$ be the closure of $\pi^{-1}(V \setminus A_1) = \pi^{-1}(X)$ in $\overline{\pi^{-1}(V)}$. Then $Y$ is the union of those irreducible components of $\pi^{-1}(V)$ whose image by $\pi$ are not contained in $A_1$. In particular no irreducible component of $Y$ is send into an irreducible component of $A_1$. Now by Lemma \ref{semi-propre 3}  the restriction $\pi_{|Y}\colon Y\rightarrow V$ is semi-proper and consequently $\pi(Y)$ is a closed subset of $V$. It follows that $\pi(Y)$ is the closure $\bar{X}$ of $X$ in $V $ . As $(V \setminus\bar X)$ is an open set which intersects every irreducible component of $A_1$ the Remmert-Stein Theorem tells us that $\bar X$ is an analytic subset of $V$. Hence $\pi(M)\cap V =  \bar X\cup A_1$ is an analytic subset of $V$.
\hfill{$\blacksquare$}\\

For the proof of Theorem \ref{Kuhl-Banach} we will need some more or less known technical results.

\begin{lemma}\label{Kuhl-banach.3} 
Let $M$ be a reduced complex space, $U$ be an open subset of $\C^n$ and $\pi\colon M\rightarrow U$ be a holomorphic map whose generic rank on every irreducible component of $M$ is at most $n-1$. Then, for every compact subset $K$ of $M$, $\pi(K)$ is a $b$-negligible subset of $U$.
\end{lemma}

\parag{Proof} As only finitely many irreducible components of $M$ meet a given  compact subset of $M$ we may assume that $M$ is irreducible. Let $p$ denote the maximal fiber dimension of $\pi$ and consider the stratification $\Sigma_p(\pi)\subseteq \Sigma_{p-1}(\pi)\subseteq\cdots\subseteq \Sigma_q(\pi) = M$, where $q$ is the  minimal generic fiber dimension of $\pi$ on an irreducible component of $M$. \\
Let $K$ be a compact subset of $M$ and let $y$ be a point in $\pi(K)$.  For each $x$ in $\pi^{-1}(y)\cap K$ there exists, thanks to Lemma \ref{Kuhl-banach.0}, an open neighborhood $W_x$ of $x$ in $M$ and an open neighborhood $U_x$ of $y$ in $U$ such that $A_x := f(W_{x}\cap\Sigma_p(\pi))$ is a nowhere dense analytic subset of  $U_x$. As $\pi^{-1}(y)\cap K$ is compact there exist $x_1,\ldots,x_l$ in $\pi^{-1}(y)\cap K$ such that $W := W_{x_1}\cup\cdots\cup W_{x_l}$ contains $\pi^{-1}(y)\cap K$. Then there exists an open neighborhood $V$ of $y$ in  $U_{x_1}\cap\cdots\cap U_{x_l}$ satisfying $\pi^{-1}(V)\cap K\subseteq W$. Hence we get
$$
V\cap \pi(K\cap\Sigma_p(\pi)) \ = \ \pi(K\cap\Sigma_p(\pi)\cap \pi^{-1}(V))\ \subseteq \ V\cap(A_{x_1}\cup\cdots\cup A_{x_l})
$$
and $V\cap \pi(K\cap\Sigma_p(\pi))$ is a $b$-negligible subset of $V$ since $V\cap(A_{x_1}\cup\cdots\cup A_{x_l})$ is a nowhere dense analytic subset of $V$. As $y$ is an arbitrary point in $\pi(K)$ it follows that $\pi(K\cap\Sigma_p(\pi))$ is $b$-negligible in $U$. \\
Now, let $y$ be a point in $U\setminus \pi(K\cap\Sigma_p(\pi))$. Then there exists an open neighborhood $V$ of $y$ in $U$ such that $\pi^{-1}(V)\cap K\cap\Sigma_p(\pi) = \emptyset$, so $\pi^{-1}(V)\cap K$ has an open neighborhood $M'$ in $M$ where the maximal fiber dimension of $\pi$ is at most $p-1$. With the same reasoning as above we then see that $\pi(K\cap(\Sigma_{p-1}(\pi))\cap V$ is a $b$-negligible subset of $V$. It follows that $\pi(K\cap(\Sigma_{p-1}(\pi))\setminus \pi(K\cap\Sigma_p(\pi))$ is a $b$-negligible subset of $U\setminus \pi(K\cap\Sigma_p(\pi))$ and consequently $\pi(K\cap(\Sigma_{p-1}(\pi))$ a $b$-negligible subset of $U$. Continuing in this way we finally get that $\pi(K\cap(\Sigma_{q}(\pi)) = \pi(K)$ is a $b$-negligible subset of $U$.
\hfill{$\blacksquare$}

\begin{cor}\label{Kuhl-banach.4}
Let $M$ be a reduced complex space, $U$ an open subset of $\C^n$ and $\pi\colon M\rightarrow U$  a holomorphic map. Let $A$ be a nowhere dense analytic subset of $M$ and $K$ be a compact subset of $M$. Then the closure of the subset of all $y$ in $U$ such that one of the irreducible components of $\pi^{-1}(y)$ meets $K$ and is contained in $A$ is $b$-negligible in $U$.
\end{cor}

\parag{Proof} Since only a finite number of irreducible components of $M$ and $A$ intersect $K$ we may suppose $M$ and $A$ both irreducible. If the restriction $\pi_{|A}$ is of rank strictly less than $n$ the result is obvious thanks to Lemma \ref{Kuhl-banach.3}, so suppose it is equal to $n$. Let $q$ denote the generic fiber dimension of $\pi$. Then the generic fiber dimension of $\pi_{|A}$ is at most $q-1$ and consequently the restriction of $\pi$ to  $\Sigma_q(\pi_{|A})$ is of rank at most $n-1$. Thus the subset $\pi(\Sigma_q(\pi_{|A})\cap K)$ of $U$, which obviously contains the subset in question, is $b$-negligible  due to Lemma \ref{Kuhl-banach.3}.
\hfill{$\blacksquare$}\\

Let $E$ be a Banach space as before and let $E = G\oplus F$ be a topological decomposition where $G$ is a finite dimensional subspace and let $U$ be an open connected subset of $G$ and $B$ an open ball in $F$. A subset $X$ in $U\times B$ is a {\bf reduced multigraph} if it satisfies the following conditions.
\begin{itemize}
\item $X$ is closed in $U \times B$.
\item
The projection onto $G$ along $F$ induces a proper\footnote{Recall that proper means that the map is closed with compact fibers. Here $X$ is not assumed to be locally compact. Remark that here the compact fibers are necessarily finite because $X$ is then a locally finite dimensional analytic subset of $U\times B$.} surjection
$$
pr \colon X\longrightarrow U.
$$
\item
There exists a closed $b$-negligible subset $R$ of $U$ such that $X$ is the closure of $X\setminus pr^{-1}(R)$ in $U\times B$ and such  that every point $x\in U\setminus R$ admits an open neighborhood $V$ in $U\setminus R$ having the property that $pr^{-1}(V)$ is a finite union of mutually disjoint graphs of holomorphic maps from $V$ into $B$.
\end{itemize}

A subset $X$ of $E$ is said to be {\bf locally a reduced multigraph} if for every point $x$ in $X$ there exists a topological decomposition $E = G\oplus F$ with $G$ finite dimensional, an open connected subset $U$ of $G$ and an open ball $B$ in $F$ such that $x\in U\times B$ and $X\cap(U\times B)$ is a reduced multigraph of $U\times B$. Such a subset $X$ has a unique structure of a reduced complex space which has the property that the canonical injection $X\hookrightarrow E$ induces a closed holomorphic embedding  in an open set in  $E$ (see Corollary 3.7.21 in \cite{[e]}). 

Moreover, if $\mathcal{U}$ is an open subset of $E$ and $X$ is a closed subset of $\mathcal{U}$, then $X$ is a reduced complex subspace of $\mathcal{U}$ if and only if it is locally a finite union of reduced multigraphs (see {\em loc. cit.}).

\begin{prop}\label{Kuhl-banach.5}
Let $E = G\oplus F$ be a topological decomposition where $G$ is a finite dimensional  subspace, $U$ is an open connected subset of $G$, $B$ an open ball in $F$  and  let $\pi \colon M\rightarrow U\times B$ be a semi-proper holomorphic map of constant generic rank from a reduced complex space $M$ to $U \times B$. Suppose also that the projection onto $U$ induces a proper surjective map with finite fibers
$$
pr \colon \pi(M)\longrightarrow U.
$$
Then every point in $U$ admits an open neighborhood $V$ such that $pr^{-1}(V)$ is a reduced multigraph in $V\times B$\footnote{It can be shown that $\pi(M)$ is in fact globally a reduced multigraph in $U\times B$, but it is not needed here.}.
\end{prop}

\parag{Proof} Denote $n$ the dimension of $G$. Then  $\pi$ and $pr \circ \pi$ are both of  constant generic rank $n$. Let $S(M)$ denote the singular locus of $M$, let $S'$ denote the analytic subset of $M\setminus S(M)$ consisting of those points where $\pi$ is of rank at most $n-1$ and put $S := S'\cup S(M)$. 

\smallskip
Take a point $y$ in $U$.  Then, as $\pi$ is semi-proper, there exists a compact subset $K$ of $M$ such that $\pi(K)$ is  a neighborhood of $pr^{-1}(y)$ in $\pi(M)$. As $pr$ is proper $y$ has an open neighborhood $V$ in $U$ which satisfies $pr^{-1}(V)\subseteq f(K)$. 
Let $R$ denote the closure of the set of all points $v$ in $V$ such that one of the irreducible components of $(pr \circ \pi)^{-1}(v)$ meets $K$ and is contained in $S$. Then, thanks to Corollary \ref{Kuhl-banach.4}, $R$ is $b$-negligible in $V$. Now, take a point $v$ in $U\setminus R$ and let $z_1,\ldots,z_k$ be the mutually distinct points of $pr^{-1}(v)$. Then, for each $z_j$, there exist a point $x_j$ in $\pi^{-1}(z_j)$ such that $x_j$ is a smooth point of $M$ and $pr \circ \pi$ is a submersion at $x_j$. Thus there exist an open neighborhood $W$ of $v$ in $V\setminus R$ and $k$ holomorphic sections $\sigma_1,\ldots,\sigma_k$ of the map $(pr \circ \pi)^{-1}(W)\rightarrow W$. It follows that, shrinking $W$ if necessary, $\pi \circ\sigma_1,\ldots,\pi \circ\sigma_k$ are $k$ holomorphic sections of $pr^{-1}(W)\rightarrow W$ whose images are mutually disjoint. As $pr \circ \pi$ is a dominant map the set $(pr \circ \pi)^{-1}(R)$ is of empty interior in $M$ and consequently the interior of  $pr^{-1}(R)$ in $\pi(M)$ is empty. Thus we can conclude that $pr^{-1}(V)$ is a reduced multigraph in $V\times B$.
\hfill{$\blacksquare$}

\begin{lemma}\label{Kuhl-banach.6}
Under the assumptions of Theorem \ref{Kuhl-Banach} suppose moreover that $\pi$ has maximal generic rank $n$ and let $M_1$ denote the union of all irreducible components of $M$ where $\pi$ is of generic rank $n$. Then the restriction $\pi_{|M_1}\colon M_1\rightarrow \mathcal{U}$ is a semi-proper map.

\end{lemma}
\parag{Proof} 
Pick an arbitrary  point $y$ in $\pi(M_1)$ and choose a compact subset $K$ of $M$ such that $\pi(K)$ contains an open neighborhood $Z$ of $y$ in $f(M)$. We are going to prove that $\pi(K\cap M_1)$ contains $Z\cap \pi(M_1)$. To this end let $M_2$ denote the (finite) union of those irreducible components of $M$ which intersect $K$ and where $\pi$ is of rank at most  $n-1$. Now, if on the contrary $\pi(K\cap M_1)$ does not contain $Z\cap f(M_1)$, then $Z\cap \pi(M_1)\setminus \pi(K\cap M_1)$ is a non-empty open subset  of $\pi(M_1)$. Then, by the constant rank theorem,  there exists an open subset $T$ of the smooth part of $M_1$ such that $\pi(T)$ is a locally closed complex sub-manifold of dimension $n$ in $Z\cap \pi(M_1)\setminus \pi(K\cap M_1)$. Hence the $n$-dimensional manifold $\pi(T)$ is contained in $\pi(M_2)$, which is absurd since the rank of $\pi$ on every irreducible component of $M_2$ is strictly smaller than $n$.
\hfill{$\blacksquare$}\\

Now we are ready to prove Theorem \ref{Kuhl-Banach}.

\parag{Proof of Theorem \ref{Kuhl-Banach}}
Since for every point $y$ in $\pi(M)$ there exists a compact subset $K$ of $M$ such that $\pi(K)$ is a neighborhood of $y$ we may assume that the generic rank of $\pi$ is bounded.

\smallskip
{\em We shall now prove the theorem by induction on the maximal generic rank of $\pi$.} 

\smallskip
Suppose that the map $\pi$ is of generic rank zero. Then $\pi$ is constant on every irreducible component of $M$, and each point in $\pi(M)$ has a neighborhood which only contains a finite number of points in $\pi(M)$ since $\pi$ is semi-proper. As $\pi(M)$ is closed in $\mathcal{U}$ it follows that  $\pi(M)$ is  a $0-$dimensional reduced complex subspace of $\mathcal{U}$.

\smallskip
Now suppose that, for a given integer $n > 0$, the image of every semi-proper holomorphic map from a reduced complex space $M$ into a Banach open set  is a reduced complex subspace of this open set, if the generic rank of the map is at most $n-1$ on every irreducible component of $M$.

Assume that the map $\pi \colon M\rightarrow \mathcal{U}$ has maximal generic rank $n$ and let $M_1$ be the union  of all irreducible components of $M$ where $\pi$ is of generic rank $n$. Then by Lemma \ref{Kuhl-banach.6} the restriction $\pi_{|M_1}\colon M_1\rightarrow \mathcal{U}$ is a semi-proper map. Suppose for the moment that we have shown that $\pi(M_1)$ is a reduced complex subspace of $ \mathcal{U}$ and put $Z := \mathcal{U} \setminus \pi(M_1)$. Then $M_2 := \overline{\pi^{-1}(Z)}$ is the union of those irreducible components of $M$ which are not mapped by $\pi$ into $\pi(M_1)$. Moreover the generic rank of $\pi$ is bounded by $n-1$ on $M_2$ and the restriction $\pi_{|M_2}\colon M_2\rightarrow \mathcal{U}$ is semi-proper due to Lemma \ref{semi-propre 3}. Hence  $\pi(M_2)$ is a reduced complex subspace of $\mathcal{U}$ by  our induction hypothesis and it follows that  $\pi(M) = \pi(M_1)\cup \pi(M_2)$ is a reduced complex subspace of $\mathcal{U}$. This shows that we may, without loss of generality, assume that $\pi$ is of constant generic rank $n$.

\smallskip
Let $E^*$ denote the topological dual space of $E$. Take an arbitrary point  $x$ in $M$ and fix a compact subset $K$ of $M$ such that $\pi(K)$ is a neighborhood of $\pi(x)$ in $f(M)$. As $E^*$ is a Hausdorff space we have 
$$
\{y\in \pi(M)\ |\ l(y) = l(\pi(x)), \ \forall\,l\in E^*\} \ = \ \{\pi(x)\}
$$
and consequently $\pi^{-1}(\pi(x)) = \bigcup\limits_{l\in E^*}(l\circ \pi)^{-1}(\pi(x))$. Thus every $x'$ in $\pi^{-1}(\pi(x))$ has an open neighborhood $W$ such that $W\cap \pi^{-1}(\pi(x))$ is an intersection of finitely many sets of the form $W\cap(l\circ \pi)^{-1}(\pi(x)))$ where $l\in E^*$.  Since $\pi^{-1}(\pi(x))\cap K$ can be covered by finitely many such neighborhoods, there exist $l_1,\ldots, l_N$ in $E^*$ and an open neighborhood $V$ of $\pi^{-1}(\pi(x))\cap K$  in $M$ such that $\pi^{-1}(\pi(x))\cap V$ is the set of common zeros of the holomorphic functions
$l_1\circ \pi - l_1(\pi(x)),\ldots, l_N\circ \pi - l_N(\pi(x))$ in $V$. It follows that the continuous linear map 
$$
L\colon E\longrightarrow \C^N,\qquad y\mapsto (l_1(y),\ldots, l_N(y))
$$
has the property that the fiber over $L(\pi(x))$ of the restriction  $L_{|\pi(M)}$ is the singleton $\{\pi(x)\}$. Hence there exists an open neighborhood $\mathcal{V}$ of $\pi(x)$ in $\mathcal{U}$ and an open neighborhood $W$ of $L(\pi(x))$ in $\C^N$ such that  the induced map $\pi(M)\cap\mathcal{V}\rightarrow W$ is proper. Now $\pi$ is of constant generic rank $n > 0$ so no irreducible component of $M$ is contained in $\pi^{-1}(f(x))$. It follows that, for all $z$ near enough to $L(\pi(x))$, none of the irreducible components of $M$ which meet $K$ is mapped into $L^{-1}(z)$. 
So, by shrinking $W$ around $L(\pi(x))$ and replacing $\mathcal{V}$ by  $\mathcal{V}\cap L^{-1}(W)$,  we obtain that $\pi$ induces a semi-proper holomorphic map $\pi^{-1}(L^{-1}(z))\rightarrow L^{-1}(z)$ of maximal generic rank at most $n-1$ for all $z$ in $W$. By our induction hypothesis $\pi(M)\cap L^{-1}(z)$ is then a reduced complex subspace of $L^{-1}(z)$ of dimension at most $n-1$. Replacing $M$ by $\pi^{-1}(\mathcal{V})$ we have the following situation
$$
\xymatrix{M \ar[rd]_{h\circ \pi} \ar[r]^{\pi\qquad} & \pi(M)\cap\mathcal{V} \ar[d]^{h} \\
{} & W } 
$$
where $h$ denotes the restriction of $L$ to $\pi(M)\cap\mathcal{V}$. Now, $h$ is proper and consequently $h^{-1}(y) = \pi(M)\cap L^{-1}(y)\cap\mathcal{V}$ is a compact  analytic subset of $\mathcal{V}$, for every $y$ in $W$. Hence $h$ has finite fibers. It follows that $h\circ \pi$ is a semi-proper map of  constant generic rank $n$ and $h(\pi(M))$ is a reduced complex subspace of $W$, thanks to Theorem \ref{Kuhl-Banach.special case}.\\
Shrinking $W$ around $h(\pi(x))$ if necessary,  we can find a projection of $\C^N$ onto an $n-$dimensional subspace of $\C^N$ which induces a proper map with finite fibers of $W$ onto an open subset $U$ of the subspace. This enables us to reduce our situation to the situation of Proposition \ref{Kuhl-banach.5} and consequently the proof is complete.
\hfill{$\blacksquare$}

\parag{Comment} As the reader may see, the generalization of the delicate theorem of  to the case where the target is a Banach open set uses the same tools than the generalization of Remmert's Direct Image Theorem (completed by the enclosability theorem) due to \cite{[B-Mz]} (see also \cite{[e]} chapter III sections 7.3 and 7.4).\\
This result will be generalized again to the case where the target space is the space $\mathcal{C}_m^f(P)$ of finite type $m$-cycles in a reduced complex space $P$ endowed with its "weak Banach analytic structure" (see sections 2 and 4 in chapter IV).\\
But the result above is the fundamental key of the use of finite type cycles in the present book.\\

\chapter{Quasi-proper maps}

\section{Definitions and examples}

The classical notion of a quasi-proper map is as follows: Let $\pi : M \to N$ be an holomorphic map between reduced complex spaces. The map $\pi$ is {\bf quasi-proper } at a point $y_0 \in N$ when there exists an open neighborhood $W$ of $y_0$ in $N$ and a compact set $K$ in $M$ such that for all $y \in W$ and every irreducible component $C$ of $\pi^{-1}(y)$ we have $K \cap C \not= \emptyset$. \\
But in our study of families of analytic cycles in a reduced complex space parametrized by a Hausdorff topological space we need to generalize the notion of a quasi-proper map. The main example we will be interested in and where this generalization is needed is the following:\\
Let $g : N \to \mathcal{C}_n^{loc}(M)$ be a continuous map of a Hausdorff topological space $N$ to the (topological) space\footnote{See chapter III section 1 here below for a brief reminder on the topology of $\mathcal{C}_n^{loc}(M)$ or \cite{[e]} ch.IV section 2  for a systematic treatement  of this topology.} of closed $n$-cycles in the complex space $M$. Let $ \vert G\vert \subset N \times M$ be the set theoretic graph of the corresponding continuous family of cycles (so, by definition,  $\vert G\vert := \{(y, x) \in N \times M \ / \ x \in \vert g(y)\vert \}$) and let $pr : \vert G\vert \to N$ be the natural projection. Then $pr$ is continuous and each fiber of $pr$ is, in an obvious way, identified with a closed analytic subset of $M$. So it makes sense to say that $pr$ is quasi-proper, even though neither $G$ nor $N$ is  a reduced complex space.\\
Remark also that, in this case, the Hausdorff  topological space $\mathcal{C}_n^{loc}(M)$ is not locally compact in general. This leads to  the following definition.

\begin{defn}\label{q-p-0}
\begin{enumerate}[(i)]
\item 
A {\bf\em correspondence} is a triple $(M, N, G)$ where  $M$ is a reduced complex space, $N$ a Hausdorff space and $G$ a closed subset of $N \times M$ having the following property: 
\begin{enumerate}
\item[$(*)$]  
For all \ $y$ in $N$,\ $G\cap (\{y\} \times M) = \{y\}\times G_y $,  where $G_y$ is an analytic subset of $M$.  
\end{enumerate}
\item
We shall say that a triple $(M, N, G)$ is {\bf\em quasi-proper at the point $y_0 \in N$}  when  $(M, N, G)$  is a correspondence and when the following condition holds:
\begin{enumerate}
\item[$(**)$]
There exists an open neighborhood $N_0$ of $y_0$ in $N$ and a compact subset $K$ in $M$  such that, for all $y$ in $N_0$, every irreducible component of $G_y$  meets $K$\footnote{Recall that, by definition, the empty set has no irreducible component.}.
\end{enumerate}
We shall say that the triple $(M, N, G)$  is {\bf \em quasi-proper} when it is quasi-proper at every point in $N$.
\end{enumerate}
\end{defn}

\medskip
Let  $(M, N, G)$ be a correspondence, $\pi\colon G\rightarrow N$ be the canonical projection and $C$ be an irreducible component of $G_y$ for some $y$ in $N$. Then we say that $\{y\}\times C$ is an irreducible component of the fiber $\pi^{-1}(y)$.

\medskip
Often we shall  consider  a continuous map $\pi: M \to N$ where $M$ is a reduced complex space and $N$ a Hausdorff space such that any fiber of $f$ is an analytic subset in $M$. Then we shall associate to such a map $\pi$ the triple $(M, N, G_\pi )$, where $G_\pi \subset  N \times M$ is the subset 
$$G_\pi := \{ (y, x) \in N\times M \ / \  y= \pi(x) \} ,$$ 
and we say that $\pi$ is {\bf quasi-proper at a point} $y_0 \in N$ when the triple $(M, N, G_\pi)$ is quasi-proper at $y_0$. In this situation the condition $(**)$ is equivalent  to the existence of  an open neighborhood $N_0$ of $y_0$ in $N$ and a compact set $K$ in $M$ such that every irreducible component of every fiber of $\pi$ over a point in $N_0$ intersects $K$. The map $\pi$ will be called {\bf quasi-proper} if it is quasi-proper at every point in $N$.

\smallskip
When we consider a holomorphic map $\pi:  M \to N$ between two reduced complex spaces, the subset $G_\pi$ defined as above, always satisfies condition $(*)$  and the quasi-properness of $\pi$ at $y_0 \in N$ reduces to the condition $(**)$.
\medskip
Obviously, a   holomorphic map between reduced complex spaces which is proper is quasi-proper. Moreover, a quasi-proper map between reduced complex spaces is semi-proper. 
\parag{Examples} 
\begin{enumerate}
\item Let $M$ and $N$ be reduced complex spaces with $M$ irreducible. Then the natural  projection $pr: M\times N \to N$ is quasi-proper.
\item Let $M$ and $N$ be complex manifolds and let $\pi : M \to N$ be a holomorphic submersion with connected fibers. Then $\pi$ is quasi-proper.
\item Let $M$ and $N$ be reduced complex spaces and let $\pi : M \to N$ be a holomorphic map with finite fibers. Then $\pi$ is quasi-proper if and only if it is proper
\item 
Assume that the triple $(M, N, G)$ is quasi-proper and that $N$ is a singleton $\{y_0\}$. Then condition $(*)$ of Definition \ref{q-p-0} says that $G = \{y_0\}\times X$ where $X$ is a closed analytic subset in $M$ and condition $(**)$ is equivalent to the fact that $X$ has only finitely many irreducible components. So we encounter in this special case the notion of a {\it cycle of finite type} which we shall use intensively in the sequel.
\item 
Let $M$ be a complex space and $n$ be a natural number. Note $\mathcal{C}_n^f(M)$ the set of finite type cycles of dimension $n$. We shall see in section 4.1 that if  
$|G| \subset \mathcal{C}_n^f(M)\times M$ is the set theoretic graph of the tautological family of finite type $n-$cycles in $M$, 
$$ 
|G\vert := \{( X, x) \in \mathcal{C}_n^f(M)\times M \ / \  x \in \vert X\vert \} ,$$
then the topology we shall define in section 4.1  on $\mathcal{C}^f_n(M)$ is precisely chosen in order to make  the triple $(M,  \mathcal{C}_n^f(M), |G|)$ quasi-proper.\\
\end{enumerate}

\begin{prop}\label{quasi-propre 1}
Let $(M, N, G)$ be a correspondence  and $\pi\colon G\rightarrow N$ be the canonical projection. The following properties hold true:
\begin{enumerate}[(i)]
\item 
Assume that $(M, N, G)$ is quasi-proper  at the point $y_0 \in N$. Then $(M, N, G)$ is quasi-proper at each point of an open neighborhood of $y_0$ in $N$, so quasi-properness for a given correspondence  $(M, N, G)$ is a local property on $N$. More precisely, the correspondence  $(M, N, G)$ is quasi-proper if and only if, for every open subset $N'$ of $N$ the induced correspondence $(M, N', G')$, where  $ G' := (N'\times M)\cap G$, is quasi-proper. It follows that, if we have  an open covering $(N_j)_{j\in J}$ of $N$ such that all the correspondences $(M, N_j, G_j)$, where $G_j :=  (N_j \times M)\cap G$, are quasi-proper, then  the correspondence $(M, N, G)$ is quasi-proper. 
\item 
Assume that $(M, N, G)$ is quasi-proper. Then $\pi(G)$ is a closed subset of $N$.
\item 
If $\pi(G)$ is closed in $N$, then $(M, N, G)$  is quasi-proper if and only if the induced  triple $(M, \pi(G), G) $  is quasi-proper. But $\pi(G)$ is not locally compact in general.
\item 
Assume that  $(M, N, G)$ is quasi-proper. Then  for any closed subset $F$ of $G$ which is a union of irreducible components of fibers of $\pi$, the induced correspondence  $(M, N, F)$ is quasi-proper. 
\item
If  $(M, N, G)$ is quasi-proper, then from (i) above we see that  the following condition is satisfied.
\begin{itemize}
\item[{$(@)$}] 
For each compact set $L$ in $N$ there exists a compact set $K(L)$ in $M$  which intersects every irreducible component  of $G_y$ for all $y$ in $L$.
\end{itemize}
In the case where $\pi(G)$ is a closed subset of $N$ and also locally compact (in particular if $N$ is locally compact), then the correspondence $(M, N, G)$ is quasi-proper if and only if it  satisfies the condition $(@)$.
\end{enumerate}
\end{prop}

\parag{Proof} The points $(i)$ and $(iv)$ are left to the reader as easy exercises. \\
Point $(v)$ is obtained by a standard compactness argument.\\
 So we shall only give more details for points $(ii)$ and $(iii)$. \\
Let us first prove (ii). To see that $\pi(G)$ is closed, take a point $y \in N \setminus \pi(G)$ and let $N_y$ and $K_y$ be respectively an open neighborhood of $y \in N$ and a compact set in $M$ such that for any $y' \in N_y$ any irreducible component of $\pi^{-1}(y')$ meets $K_y$. Then remark that the restriction of $\pi$ to $(N\times K_y) \cap G$ is a proper map. So its image is closed, and as $y$ is not in this closed set, for each $y'$ in the open neighborhood $N'_y := N_y \cap (N \setminus \pi(N\times K_y) \cap G)$ of $y$ in $N$ we have $y' \not\in \pi(N \times K_y) \cap G)$ and any irreducible component of $\pi^{-1}(y')$ must meet $\{y'\} \times K_y$. This implies that $\pi^{-1}(y') = \emptyset$ and it follows that the complement of the subset $\pi(G)$  in $N$ is open.\\
 To prove point $(iii)$ we assume $\pi(G)$  closed and $(M, \pi(G), G)$ quasi-proper. As the quasi-properness of the triple $(M, N, G)$ at  $y \in \pi(G)$ is clear, let $y $ be a point in $N \setminus \pi(G)$. Then, as $\pi(G)$ is closed, the choices  $N_y := N \setminus \pi(G)$ and $K_y := \emptyset$  give the quasi-properness.\\
 The  example below shows that $\pi(G)$ is not locally compact, in general.$\hfill \blacksquare$\\
 
 \parag{Example} Put $M := \C\times\C$, $N := {\rm Hol}(\C,\C)$, i.e. the space of all holomorphic maps from $\C$ to $\C$ endowed with the topology of compact convergence, and  let $G:=\{(f,(x,y))\in N\times M\ | \ f(x) = y\}$. Then the triple $(M,N,G)$ is quasi-proper since for any non empty compact set $K$ in $\C$ and any $\varepsilon > 0$, for each $g $ in the neighborhood  $\mathcal{V}(K, \varepsilon) := \{g \in Hol(\C,\C) \ / \ \vert\vert g - f\vert\vert_K < \varepsilon \}$ of $f$  the graph of $g$ meets the compact set 
 $$ \mathcal{K}(f, K, \varepsilon) := \{(z, t) \in \C \times \C\ / \ K \times \{t \in \C \ / \ d(t,f(K)) \leq \varepsilon \}\} .$$
   But ${\rm Hol}(\C,\C)$ is not locally compact. $\hfill \square$ \\

The proof of the following  very easy lemma is left to the reader.

\begin{lemma}\label{recouvre}
Let $(M, N, G)$ be a correspondence and $(M_j)_{j \in J}$ be a {\em finite} family of open sets  in $M$. For each $j$ in $J$  put $G_j := (N \times M_j)\cap G$  and suppose that  the correspondence  $(M_j, N, G_j)$ is quasi-proper. Then  the triple $\left(\cup_{j \in J} M_j, N, \cup_{j \in J} G_j\right)$ is quasi-proper.$\hfill \blacksquare$
\end{lemma}

Now we shall concentrate on the construction of quasi-proper equidimensional holomorphic map, from a given equidimensional holomorphic map $\pi : M \to N$.\\
Our first step will be a local result.

\begin{lemma}\label{restr.2} Let $\pi : M \to N$ be a holomorphic map between irreducible complex spaces and note $n := \dim M - \dim N$. Let $y_0$ be a point in $N$ such that $\pi^{-1}(y_0)$ has pure dimension $n$ and let $K$ be a compact set in  an open set $M_0$ of $M$ such that $\pi(K)$ is a neighborhood of $y_0$. Then there exists a relatively compact open set $M(y_0) \subset M_0$, a relatively compact open set $N(y_0)$ in $N$ containing $y_0$ such that the following property holds:
\begin{enumerate}
\item The restriction of $\pi$ to $M(y_0)$ defines a surjective equidimensional holomorphic map $\pi(y_0) : M(y_0) \to N(y_0)$.
\item There exists a closed $\pi(y_0)$-proper subset $\mathcal{L}(y_0)$ in $M(y_0)$, containing $K \cap M(y_0)$, such that for each $z \in N(y_0)$ each irreducible component of $\pi(y_0)^{-1}(z)$ meets $\mathcal{L}(y_0)$.
\end{enumerate}
\end{lemma}

\parag{Proof} Remark first that under our hypothesis the fibers of $\pi$ has dimension at least equal to $n$. For each point $z \in K$  choose  a pair of $n$-scales on $M_0$,  $E_{z}:= (j_{z}, U_{z}, B_{z})$ and $E'_z = (j_z, U_z, B'_z)$, where $j_z$ is a closed holomorphic embedding of an open subset $W(z)$ of  $M_0$ into an open  neighborhood of $\bar U_z\times \bar B_z$   satisfying the following conditions:
\begin{itemize}
\item 
$ \bar B'_z \subset B_z$
\item 
$z \in j_z^{-1}(U_z\times B'_z)$ 
\item $\vert \pi^{-1}(y_0) \vert \cap \left( j^{-1}(\bar U_z \times(\bar B_z \setminus B_z'))\right) = \emptyset$, $E_z$ is adapted to $\pi^{-1}(y_0)$.
 \end{itemize}
 So each scale $E_z$ is adapted to $\pi^{-1}(y_0)$. \\
 Now construct for each $z \in K$ an open irreducible  neighborhood $N_z$ as follows:
  $$ N_z := N \setminus \pi\big(j_z^{-1}(\bar U_z\times (\bar B_z \setminus B'_z)) \big).$$
  When $z$ describes $K$ the open sets $c(E_z)$ cover $K$ so we can choose  finitely many points $z_1, \dots, z_m$ in $K$ such that  the relatively compact open set $M'_0 := \cup_{h= 1}^m  c(E_{z_h})$ of $M_0$ contains $K$. \\
 Now choose for each $h \in [1, m]$ a point $t_h \in U_z$  and put  $L_{z_h} := j_{z_h}^{-1}(\{t_{z_h}\}\times \bar B '_{z_h})$.\\
 Let $N(y_0)$ be the open neighborhood of $y_0$ in $N$ defined as $N(y_0) := \cap_{h=1}^m N_{z_h}$ and define $M(y_0) := \pi^{-1}(N(y_0))\cap M'_0.$\\
Then we define   $L(y_0) := M(y_0) \cap \cup_{h=1}^m L_{z_h}$.\\
  Now we shall prove that for each $y \in N(y_0)$, each irreducible component $\Gamma$ of  the cycle $\pi^{-1}(y) \cap M(y_0)$ meets the $\pi(y_0)$-proper  set $L (y_0)$ in $M(y_0)$.\\
 By definition such a $\Gamma$ is not empty, so there exists at least one $h \in [1, m]$ such that $\Gamma$ meets $c(E_h)$. As  the scale $E_{z_h}$ is adapted to $\pi^{-1}(y)$ it is adapted to 
 $\Gamma$.  So $\Gamma$ is of  dimension $n$. Moreover, the degree of $\Gamma$ in this scale is positive as we know that $\Gamma$ meets $c(E_{z_h})$. Then $\Gamma$ is a non empty union of irreducible components of the multigraph defined by $\pi^{-1}(y)$ in the adapted scale $E_{z_h}$ and then has to meet the set $L_h$.  So $\Gamma$ meets $L(y_0)$. \\
 The condition we proved  implies that the map $\pi(y_0)$ is quasi-proper and equidimensional, so it is an open map. Then  it is enough  replace $N(y_0)$ by the image by $\pi$ of $M(y_0)$ which is a relatively compact open set in $N$ and to define  $\mathcal{L}(y_0)$  by $ \mathcal{L}(y_0):= L(y_0)  \cup (K \cap M(y_0))$ to complete the proof.$\hfill \blacksquare$\\

To globalize the previous lemma we need the following topological result.

\begin{lemma}\label{topology}
Let $\pi : M \to N$ be a continuous surjective map between two locally compact and paracompact first countable Hausdorff spaces. Assume that $N = N_1 \cup N_2$ where $N_i$ are open subsets of $N$ and let $M_i := \pi^{-1}(N_i), i = 1,2$. Assume that $L_1$ and $L_2$ are respectively closed $\pi$-proper subsets in $M_1$ and $M_2$ which satisfy $\pi(L_1) = N_1$  and $\pi(L_2) = N_2$. Then there exists a closed $\pi$-proper subset $\Lambda$ in $M$ with the following properties (with the notations $L_i(y) := \pi^{-1}(y) \cap L_i, i = 1, 2$ and $\Lambda(y) := \pi^{-1}(y) \cap \Lambda$):
\begin{enumerate}
\item $\Lambda \subset L_1 \cup L_2$;
\item For each  $y \in N_1\setminus N_2$ then $\Lambda(y)$ contains $L_1(y)$.
\item For each  $y \in N_2 \setminus N_1$ then $\Lambda(y)$ contains $L_2(y)$.
\item For each $y \in N_1\cap N_2$ then $\Lambda(y)$ contains $L_1(y)$ or $L_2(y)$.
\end{enumerate}
\end{lemma}

Note that condition on the set $\Lambda$ implies that  $\pi(\Lambda) = N$.

\parag{Proof} First cover $N$ by a locally finite  countable family $(K^\nu), \nu \in \mathbb{N}$, of compact subsets  such that each $K^\nu$ is contained in $N_1$ or in $N_2$ (and may be  in both for some $\nu$). Then define the map $\varphi : \mathbb{N} \to \{1, 2\}$ and $\psi : \mathbb{N} \to \{1, 2\}$ as follows:\\
--$\varphi(\nu) = 1$ if $K^\nu$ is  contained in $N_1$ and $\varphi(\nu) = 2$ if $K^\nu$ is not contained in $N_1$.\\
--$\psi(\nu) = 2$ if $K^\nu$ is contained in $N_2$ and $\psi(\nu) = 1$ if $K^\nu$ is not contained in $N_2$.\\
Then define:\\
-- $\Lambda_1 := \big(\cup_{\varphi(\nu) = 1} \pi^{-1}(K^\nu) \cap L_1\big) \bigcup  \big(\cup_{\varphi(\nu) = 2} \pi^{-1}(K^\nu) \cap L_2\big)$ and also \\
-- $\Lambda_2 := \big(\cup_{\psi(\nu) =2} \pi^{-1}(K^\nu)\cap L_2\big)  \bigcup \big(\cup_{\psi(\nu) = 1}\pi^{-1}(K^\nu)\cap L_1\big)$.\\
As each $\pi^{-1}(K^\nu)\cap L_i$ is a compact set in $M$ and as this family is locally finite, $\Lambda_1$ and $\Lambda_2$ are closed and $\pi$-proper and then $\Lambda := \Lambda_1 \cup \Lambda_2$ is closed and $\pi$-proper and  $\pi(\Lambda) = N$. Moreover, it is clear that $\Lambda$ satisfies the desired properties 1 to 4.$\hfill \blacksquare$\\

\begin{thm}\label{quasi-propre ex.}
Let $\pi : M \to N$ be a surjective and $n$-equidimensional holomorphic map between irreducible complex spaces. Fix  a closed  $\pi$-proper subset $\mathcal{K}$ which is surjective on $N$ and  an open set $M_0$ in $M$ which contains $\mathcal{K}$. Assume that for each $y \in N$   there exists an open neighborhood $N(y)$ of $y$ in $N$, a relatively compact open set $M(y)$ in 
$M_0 \cap \pi^{-1}(N(y))$ containing the subset $\pi^{-1}(N(y)) \cap \mathcal{K}$ such the restriction $\pi(y) : M(y) \to N(y)$ of $\pi$ satisfies the following property
\begin{itemize}
\item There exists a $\pi(y)$-proper subset $L(y) \subset \pi^{-1}(N(y))$ such that  for each $z \in N(y)$, any irreducible component of $\pi^{-1}(z) \cap M(y)$ meets $L(y)$.
\end{itemize}
Then there exists a closed  $\pi$-proper subset $\mathcal{L}$ in $M_0$, an open set $M'$ in $M$ such that $\mathcal{L} \subset M' \subset M_0$ and such that, for each $y \in N$, each irreducible component of $\pi^{-1}(y) \cap M'$ meets $\mathcal{L}$.\\
So the restriction of $\pi' : M' \to N$ is quasi-proper.
\end{thm}

\parag{Proof} Thanks to the previous lemma, it is enough to consider a locally finite  covering  $\overline{N(y_\mu)}, \mu \in \mathbb{N}$  of $N$, extracted from the covering given  by the relatively compact open sets $N(y), y \in N$,  deduced from the hypothesis, and to argue by an induction on $\mu \in \mathbb{N}$ to construct step by step the open sets 
$$ M_p := \cup_{\mu =0}^p  M(y_\mu) $$
and the closed $\pi$-propers subsets $L_p $ in $M_p$ by use of the previous lemma.\\
 We conclude by letting 
 $$M' := \cup_{\mu \geq 0} M(y_\mu) = \cup_{p \geq 0} M_p$$
  and $\mathcal{L} = \lim_{p \to \infty} L_p $ which is closed and $\pi$-proper in $\mathcal{M}$ thanks to the locally finiteness of our covering of $\mathcal{M}$ by the $M(y_\mu)$:
\begin{itemize}
\item for $p$ given, we have to modify $L_{p+j}$ over $N_p := \cup_{\mu = 0}^p \bar N(y_\mu)$ by applying the previous lemma  only for finitely many values of $j \geq 1$ since there are only finitely many $\mu \geq p+1$ for which $\bar N(y_\mu) $ intersect $N_p$. $\hfill \blacksquare$\\
\end{itemize}

Note that, since an equidimensional holomorphic map is open and a quasi-proper map has a closed image, the surjectivity assumption in the previous result is not restrictive.

\section{ Stability properties of quasi-proper maps}

In general the composition of two quasi-proper maps  is not quasi-proper as can be seen from the example below. Nevertheless  quasi-proper maps have some important stability properties with respect to composition.

\parag{Example}
Let $f\colon M\rightarrow\C^2$ be be blow-up of $\{0\}\times\mathbb Z$ in $\C^2$ and $g\colon\C^2\rightarrow\C$ be the natural projection $(z,w)\mapsto w$. Then $f$ is proper and $g$ is quasi-proper, but the composition $g\circ f$ is not quasi-proper.

\begin{lemma}\label{composition}
Let $\pi : M\rightarrow N$ and  $\rho : N\rightarrow S$   be holomorphic maps between reduced complex spaces.
\begin{enumerate}[(i)] 
\item
If $\pi$ is quasi-propre and $\rho$ is proper, then \ $\rho\circ \pi$ \ is quasi-proper.
\item
If $\rho\circ \pi$ is quasi-proper and  $\pi$ is surjective, then $\rho$  is quasi-proper. 
\end{enumerate}
\end{lemma}

\parag{Proof} (i) Every  point $s$ in $S$ has a relatively compact open neighborhood  $V$ in $S$ and  $\rho^{-1}(\bar{V})$ is a compact subset of $N$. As $\pi$ is quasi-proper there exists a compact subset $K$ of $M$ which meets every irreducible component of every fiber of $\pi$ over $\rho^{-1}(\bar{V})$. Hence $K$ intersects every irreducible component of every fiber of $\rho\circ \pi$ over $V$, so  $\rho\circ \pi$ \ is quasi-proper at $s$.

\smallskip
(ii) Let $s$ be a point in $S$. Then there exists an open neighborhood  $V$ of $s$ in $S$ and a  compact subset $K$ of $M$ which meets every irreducible component of every  fiber
 of $\rho\circ \pi$ over $V$. Let us show that $\pi(K)$  intersects every irreducible component of every fiber of \ $\rho$  over $V$. So let $C$ be an irreducible component of $\rho^{-1}(t)$ for some $t$ in $V$. As $\pi$ is surjective there exists an irreducible component $\hat{C}$ of $\pi^{-1}(C)$  such that $\pi(\hat{C}) = C$. As $C$ has non-empty interior in $\rho^{-1}(t)$ the interior of $\hat{C}$ in $\pi^{-1}(\rho^{-1}(t))$ is also non-empty so $\hat{C}$ is an irreducible component of $\pi^{-1}(\rho^{-1}(t))$. It follows that $K\cap\hat{C}\neq\emptyset$ and consequently \ $\pi(K)\cap C\neq\emptyset$. \hfill{$\blacksquare$}\\

The following result shows that quasi-proper triples are stable by {\em base change}.

\begin{prop}\label{stab.0}
Let $(M, N, G)$ be a quasi-proper triple. Then for any continuous map $g\colon S \rightarrow N$ of a Hausdorff space $S$ to $N$, the triple $(M, S, (g \times \id_M)^{-1}(G))$ is quasi-proper.
\end{prop}

\parag{Proof} Let $g\colon S \rightarrow N$ be a continuous map and put $\Gamma := (g\times \id_M)^{-1}(G)$. Then we have $\Gamma_s = G_{g(s)}$ for every $s$ in $S$.
Now take an arbitrary point  $s_0$ in $S$ and choose an open neighborhood $U$ of $g(s_0)$ in $N$ and a compact set $K$ in $M$ such that, for all $y \in U$, every  irreducible component of $G_{y}$   meets $K$. Then, for all $s$ in the open neighborhood $g^{-1}(U)$ of $s_0$ in $S$, every  irreducible component of $\Gamma_s = G_{g(s)}$  meets $K$.
$\hfill \blacksquare$\\

We end this paragraph by proving two more stability results for quasi-proper maps between reduced complex spaces.

\begin{lemma}\label{modif.quasiproper.0}
Let $\pi\colon M\rightarrow N$ be a holomorphic map between reduced complex spaces and $\tau\colon\tilde{M}\rightarrow M$ be a modification whose center is $\pi$-proper. Then $\pi\circ\tau$ is quasi-proper if and only if $\pi$ is quasi-proper.
\end{lemma}
\parag{Proof} 
Suppose that $\pi$ is quasi-proper. Let $C$ denote the center of $\tau$ and let $K$ be a compact subset of $N$. Then there exists a compact subset $L$ of $M$ which intersects every irreducible component of $\pi^{-1}(y)$ for all $y$ in $K$. As $\tau$ is proper and $C$ is $\pi$-proper the subset $\tau^{-1}(L)\cup\tau^{-1}(C\cap\pi^{-1}(K))$ of $\tilde{M}$ is compact. Now take an arbitrary point $y$ in $K$ and an irreducible component $Z$ of $\tau^{-1}(\pi^{-1}(y))$ and let us show that $Z$ intersects the compact subset $\tau^{-1}(L)\cup\tau^{-1}(C\cap\pi^{-1}(K))$. Indeed, if $Z$ does not meet $\tau^{-1}(C\cap\pi^{-1}(K))$, then $Z\subseteq \tilde{M}\setminus \tau^{-1}(C)$ so $\tau(Z)$ is an irreducible component of $\pi^{-1}(y)$. Hence $\tau(Z)\cap L\neq\emptyset$ and consequently $Z\cap\tau^{-1}(L)\neq\emptyset$.

\smallskip
Conversely, if $\pi\circ\tau$ is quasi-proper then $\pi$ is quasi-proper by point  (ii) of Lemma \ref{composition}, since a modification is, by definition, surjective.\hfill{$\blacksquare$}

\begin{lemma}\label{open.quasipropre}
Let $\pi_1 : M_1\rightarrow N$ and $\pi_2 : M_2\rightarrow M_1$ \ be quasi-proper maps between reduced complex spaces and suppose moreover that  $\pi_2$ is an open map. Then the composed map $\pi_1\circ \pi_2$ is quasi-proper.
\end{lemma}
\parag{Proof} Let $K$ be a compact subset of $N$. As $\pi_1$ is quasi-proper  there exists a compact subset $L_1$ of $M_1$ such that ${\rm int}(L_1)$ intersects every irreducible component of $\pi_1^{-1}(y)$ for all $y$ in $K$. As $\pi_2$ is quasi-proper  there exists a compact subset $L_2$ of $M_2$ which intersects every irreducible component of $\pi_2^{-1}(x)$ for all $x$ in $L_1$. \\
We are going to show that every irreducible component of  $\pi_2^{-1}(\pi_1^{-1}(y))$ meets $L_2$ for all $y$ in $K$.\\
Now take an arbitrary point $y$ in $K$. Since the map $\pi_2$ is both open and quasi-proper the same is true for the induced map $\pi_2^{-1}(\pi_1^{-1}(y))\rightarrow\pi_1^{-1}(y)$ and consequently every irreducible component of $\pi_2^{-1}(\pi_1^{-1}(y))$ is mapped surjectively onto an irreducible component of $\pi_1^{-1}(y)$. Hence for any  irreducible component $Z$ of $\pi_2^{-1}(\pi_1^{-1}(y))$ we have $\pi_2(Z)\cap{\rm int}(L_1)\neq\emptyset$ so there exists $x$ in  $\pi_1^{-1}(y)$ such that $Z$ contains an irreducible component of $\pi_2^{-1}(x)$. It follows that $Z$ intersects $L_2$.
\hfill{$\blacksquare$}

\section{Direct Image Theorem for quasi-proper maps: a simple proof}

In this section we give a proof of Theorem \ref{Kuhl-Banach} in the case where the map $\pi$ is {\em quasi-proper} and the space $M$ is irreducible. This  special case of the theorem illustrates the difference between semi-proper and quasi-proper maps. As the reader will see, the proof is much simpler in the quasi-proper case  than in the semi-proper case. For the proof we use a generalization  to an ambient Banach open set of the \lq\lq simple case\rq\rq\ of the classical Remmert-Stein theorem, which is also proved in this section.

\begin{thm}\label{Quasi-proper.Banach}
Let $M$ be a irreducible complex space, $\mathcal{U}$ an open subset of a Banach space $E$ and $\pi\colon M\rightarrow \mathcal{U}$  a quasi-proper holomorphic map. Then $\pi(M)$ is a reduced complex subspace of $\mathcal{U}$.
\end{thm}

\parag{Proof}
As we assume that $M$ is irreducible we shall prove the theorem by induction on the maximal fiber dimension of $\pi$. 

\smallskip
If the maximal fiber dimension of $\pi$ is $0$ the map $\pi$ is quasi-proper with finite fibers and consequently it is a proper map. Then $\pi(M)$ is a reduced complex subspace of $\mathcal{U}$ due to Remmert's Direct Image Theorem\footnote{Here we need only the simple case for proper holomorphic maps with finite fibers.} (see Theorem 3.7.3  in \cite{[e]}) .\\

Now let $q$ be a strictly positive natural number such that the theorem is true if the maximal fiber dimension of $\pi$ is at most $q-1$. Suppose then that $\pi\colon M\rightarrow \mathcal{U}$ satisfies the hypotheses of the theorem and has maximal fiber dimension $q$.  \\

We begin by proving that $\pi(\Sigma_q(\pi))$ is a reduced complex subspace of $\mathcal{U}$. To do so we recall that, by Proposition \ref{Kuhl-banach.0}, $\Sigma_q(\pi)$ is a closed subset of $M$ and, as $\Sigma_q(\pi)$ is a union of irreducible components of fibers of $\pi$, it follows that the induced map $\Sigma_q(\pi)\rightarrow\pi(\Sigma_q(\pi))$ is semi-proper. Consequently   $\pi(\Sigma_q(\pi))$ is a closed subset of $\mathcal{U}$ so it is enough to show that every point $y$ in $\pi(\Sigma_q(\pi))$ has an open neighborhood $\mathcal{V}$ in $\mathcal{U}$ such that $\pi(\Sigma_q(\pi))\cap\mathcal{V}$ is a reduced complex subspace of $\mathcal{V}$. This is consequence of Proposition \ref{Kuhl-banach.1}.

Now, $S := \pi^{-1}(\pi(\Sigma_q(\pi)))$ is an analytic subset of $M$ and to complete the prove we have to consider two cases.

\smallskip
If $S = M$, then $\pi(M) = \pi(\Sigma_q(\pi))$ is a reduced complex subspace of $\mathcal{U}$.

\smallskip
If $S \neq M$, then $M\setminus S$ is an irreducible complex space and the map 
$$
M\setminus S \ \longrightarrow \ \mathcal{U}\setminus\pi(\Sigma_q(\pi)),
$$ 
induced by $\pi$, is quasi-proper. As the maximal fiber dimension of this map is at most $q-1$  its image, $\pi(M\setminus S) = \pi(M)\setminus \pi(S))$, is a reduced complex subspace of $\mathcal{U}\setminus\pi(S)$ by the induction hypothesis. Now $\dim S < \dim (M\setminus S)$ and all the fibers of the induced map $S\rightarrow\pi(S)$ are of dimension $q$. As the maximal fiber dimension of the induced map $M\setminus S\rightarrow  \pi(M)\setminus \pi(S))$ is at most $q-1$ it follows that $\dim\pi(S) < \dim\pi(M\setminus S)$. \\
 It then follows from Theorem \ref{Remmert-Stein.Banach} below that $\pi(M)$ is a reduced complex subspace of $\mathcal{U}$ since $\pi(M)$ is the closure of $\pi(M\setminus S)$ in $\mathcal{U}\setminus\pi(S)$.
\hfill{$\blacksquare$}\\

The following theorem is a generalization of the classical Remmert-Stein theorem in the "easy" case.

\begin{thm}\label{Remmert-Stein.Banach}
Let $\mathcal{U}$ be an open subset of a Banach space $E$, $A$ be a reduced complex subspace of $\mathcal{U}$ and $X$ be an irreducible complex subspace of $\mathcal{U}\setminus A$. Suppose moreover that $\dim X > \dim A$ and that the closure $\bar{X}$ of $X$ in $\mathcal{U}$ is locally compact. Then $\bar{X}$ is a reduced complex subspace of $\mathcal{U}$. 
\end{thm}

For the proof of the theorem we use the following lemma.

\begin{lemma}\label{Remmert.Baire}
Let $E$ be a Banach space of dimension at least $1$ and let $S$ be a countable subset of $E\setminus\{0\}$. Then there exists a closed hyperplane in $E$ which does not contain any point of $S$. 
\end{lemma}

\parag{Proof}
Let $E^*$ be the topological dual of $E$. Then, for every $s \in E$, the subset  $V_s := \{l\in E^*\ /\ l(s) = 0\}$ is  closed and with  empty interior in $E^*$. It follows that the union $\bigcup_{s\in S}V_s$ has non-empty interior in $E^*$ since $E^*$ is a Baire space. Hence, for each $l$ in the dense subset $E^*\setminus\bigcup_{s\in S}V_s$, we have $l(s) \neq 0$ for all $s$ in $S$.
\hfill{$\blacksquare$}

\parag{Proof of Theorem \ref{Remmert-Stein.Banach}}
The case $A = \emptyset$ being trivial, we suppose $A\neq\emptyset$ and put $n:= \dim A$ and $d := \dim X - \dim A$.

\smallskip
Obviously it is enough to show that every point $a$ in $A\cap\bar{X}$  has an open neighborhood $\mathcal{V}$ such that $\bar{X}\cap\mathcal{V}$ is a reduced complex subspace of $\mathcal{V}$.  To this end we fix a point $a$ in $A\cap\bar{X}$ and by translating if necessary we may assume $a = 0$.
Then we split the proof into two steps.

\medskip
{\sc First step.} \ Let us show that there exists a closed vector subspace $H$ of codimension $n+d+1$ in $E$ such that $H\cap(A\cup X) = \{0\}$. To do so we let $S_1$ be a set which contains exactly one point in $X$ and one point in each irreducible component of $A$. As $S_1$ is countable there exists a closed hyperplane $H_1$ in  $E$ such that $H_1\cap(S_1\setminus\{0\}) = \emptyset$ thanks to Lemma \ref{Remmert.Baire}. Hence $\dim A\cap H_1 = n -1$ and $\dim (X\cap H_1) = n+d-1$.\footnote{A priori we can not exclude the possibility that this set is empty, but the proof will bring into light that this is impossible.} Now, let $S_2$ be a set which contains exactly one point in each irreducible component of $X\cap H_1$ and each irreducible component of $A\cap H_1$. Then, by Lemma \ref{Remmert.Baire}, there exists a closed hyperplane $H_2$ in $H_1$ such that  $H_2\cap(S_2\setminus\{0\}) = \emptyset$. Continuing in this way we end up with a closed vector subspace $H_{n+d}$ in $E$ such that $\dim(X\cap H_{n+d}) = 0$ and $A\cap H_{n+d} = \{0\}$. But the set $X\cap H_{n+d}$ is countable so, again due to Lemma \ref{Remmert.Baire}, there exists a closed hyperplane $H$ in $H_{n+d}$ such that $H\cap(A\cup X) = \{0\}$.

\medskip
{\sc Second step.} \ We conclude the proof by showing that $\bar{X}$ is a reduced multigraph in a neighborhood of $0$. To this end we let $H$ be as above, $L_1$ be  a complementary vector subspace to $H$ in $E$ and  $\pi_1\colon E\rightarrow L_1$ denote the projection along $H$ onto $ L_1$. As $\pi^{-1}(0) = \{0\}$ and \ $X\cup A$ is locally compact, there exist connected open neighborhoods, $V_1$ of the origin in $E$ and $U_1$ of the origin in $L_1$, such that $\pi_1$ induces a proper map $V_1\rightarrow U_1$. It follows that $A_1 := \pi_1(A\cap V_1)$ is a reduced complex subspace in $U_1$, thanks to Remmert's Direct Image Theorem, and $\dim  A_1 =n$ since the restriction of $\pi_1$ to $A_1$ has finite fibers. Now $X_1 := X\cap (V_1\setminus \pi_1^{-1}(A_1))$ is a reduced complex subspace of $V_1\setminus \pi_1^{-1}(A_1)$ and the induced map  $V_1\setminus \pi_1^{-1}(A_1)\rightarrow U_1\setminus A_1$ is proper. As $X_1$ is a complex subspace of pure dimension $n+d$ it follows that $\pi_1(X_1)$ is a reduced complex subspace of dimension $n+d$ in $U_1\setminus A_1$. The classical Remmert-Stein theorem then implies that $\pi_1(V_1\cap(X\cup A))$ is a reduced complex subspace of $U_1$.\footnote{In fact a hypersurface.}
A conveniently chosen  linear projection of $L_1$ onto a hyperplane $L$ in $L_1$ then makes $\pi_1(V_1\cap(X\cup A))$ into a branched covering over an open neighborhood of the origin in $L$.  Let $\pi\colon E\rightarrow L$ be the linear projection obtained by composing $\pi_1$ with the projection of $L_1$ onto $L$. Then there exists an open neighborhood $V$ of the origin in $E$ and a connected open neighborhood $U$ of the origin in $L$ such that  $\pi$ induces a surjective proper map $V\cap(X\cup A)\rightarrow U$, whose fibers are all  finite. 

\smallskip
Now $A_0 := \pi(V\cap A)$ is, by Remmert's Direct Image Theorem, an analytic subset of dimension $n$ in $U$ and consequently a $b$-negligible subset of $U$ since $\dim U = n+d$. As  $\pi^{-1}(U\setminus A_0)\cap V\cap\bar{X}$ is clearly equal to  $\bar{X}\cap V$ it follows that $V\cap\bar{X}$ is a reduced multigraph\footnote{See dicussion following Corollary \ref{Kuhl-banach.4} } in $U\times  F$, where $F := \pi^{-1}(0)$.
\hfill{$\blacksquare$}

\parag{Examples} 
Let $H$ be a complex  separable Hilbert space with orthonormal basis $e_1, e_2, \ldots, e_n, \dots $ 
\begin{enumerate}
\item
Let $C$ be the union on the lines through the origin generated by $e_1, e_2, e_3, \dots $. Then $C\setminus \{0\}$ is a one dimensional smooth complex submanifold of $H \setminus \{0\}$ and, as $C$ is not locally compact near $0$, $C$ is not a finite dimensional complex subset in $H$. \\
This example already shows that without the local compactness hypothesis for $X \cup A$ (or  equivalently for $\bar X$) the Remmert-Stein theorem does not hold.\\

\item
The interest of our second example is to show that, even in the case where $X$ is smooth and connected, without the local compactness assumption $\bar X$ may not be a reduced complex subspace of $H$.

Let $\gamma : D \to H$ be the holomorphic map given by   $\gamma(z) := \sum_{\nu=1}^\infty  z^\nu e_\nu $ where $D$ is the open unit disc in $\mathbb{C}$. 

\begin{lemma}
The map $\gamma$ is injective closed and of constant rank $1$. So its image $\Gamma$ is a closed $1-$dimensional connected smooth complex sub-manifold in $H$.
\end{lemma}

\parag{proof} The rank is $1$ because the component of  $\gamma'(z)$ on $e_1$ is equal to $1$. Also $\gamma(z) = \gamma(z')$ implies the equality of the components on $e_1$ so $z = z'$.\\
To see that $\gamma$ is a closed map, let $\gamma(z_n)$ be a sequence in $\Gamma$ converging to a point $y \in H$. Then the component on $e_1$ of $\gamma(z_n)$ converges to $y_1$. So
the sequence $(z_n)$ converges to $y_1 \in \mathbb{C}$. Then for each $\nu \geq 1$ the sequence $(z_n^\nu)$ converges to $y_\nu = y_1^\nu$. As $y $ is in $H$ the series  $(\vert y_1^\nu\vert^2)$ is summable and this implies that $\vert y_1\vert < 1$. Then there exists $ a \in [0, 1[$ such that $\vert z_n\vert \leq a$ for $n$ large enough. So we may extract a subsequence of the sequence $(z_n)$ converging to a point $x \in D$ and then $\gamma(x) = y$.$\hfill \blacksquare$\\

Consider now the cone $C$ over $\Gamma$ and let $\bar C$ the closure of $C$ in $H$.

\begin{lemma}
The subset $\bar C  \setminus \{0\}$ is a closed smooth connected complex sub-manifold of dimension $2$ in $H \setminus \{0\}$. But $\bar C$ is not an analytic subset of finite dimension of $H$.
\end{lemma}

\parag{proof} Let $\varphi : \mathbb{C}\times D \to H$ be the holomorphic map defined by $\varphi(\lambda, z) = \lambda\gamma(z)$. Then it sends the open set $\mathbb{C}^*\times D^*$ into $C \setminus \{0\}$.\\
-- This restriction of the  map $\varphi$ is bijective : \\
the surjectivity is clear. Moreover if $\lambda\gamma(z) = \lambda'\gamma(z')$ with $\lambda,z,\lambda',z'$ non zero, we obtain that for each $\nu$ the equality $(z'/z)^\nu = \lambda/\lambda'$ which implies $z= z'$ and then $\lambda = \lambda'$.\\
-- This restriction has rank $2$ at each point: we have $\partial_\lambda \varphi (\lambda, z)= \gamma(z) \not= 0$ and $\partial_z\varphi(\lambda, z) = \lambda\gamma'(z)$.\\
 But if, for some $\alpha \in \mathbb{C}^*$, we have
$ \lambda\gamma'(z) = \alpha\gamma(z)$ for some $\lambda \in \mathbb{C}^*$  it implies $\lambda \nu z^{\nu-1} = \alpha z^{\nu}$ for each $\nu \geq 1$ and then $\alpha z= \lambda \nu $ for each $\nu \geq 1$. This is impossible as $\lambda$ is not zero.\\
-- $C\setminus \{0\}$ is closed in $H \setminus \mathbb{C}e_1$.\\
 Assume that a sequence $\lambda_n\gamma(z_n)$ converges to $y \in H \setminus  \mathbb{C}e_1$.\\
 If the sequence $(\lambda_n)$ is bounded, we may pass to a sub-sequence which converges to $\lambda_0$. Then $\lambda_n z_n$ converges to $y_1$ and if $\lambda_0 \not= 0$ the sequence  $(\gamma(z_n))$ converges to $y/\lambda_0$ which is in $\Gamma$ because the map $\gamma$ is closed (see  above). So $y $ is in $C$.\\
If $\lambda_0 = 0$ then  $\lambda_n z_n$ converges to  $y_1 = 0$, because $ \vert z_n\vert < 1$ for each $n$, and so  we have, for each $\nu \geq 2$, the sequence $(\lambda_nz_n^\nu = (\lambda_nz_n)z_n^{\nu-1})$  which converges to $y_\nu = 0$.\\
 So $y = 0$ in this case. Contradiction. \\
Now assume that $\vert \lambda_n\vert$ goes to $+\infty$ when $n \to + \infty$. Then  the sequence $(z_n)$ converges  to $0$ and if we assume that $y_1$ is not $0$,  we have again, for each $\nu \geq 2$, the sequence $(\lambda_nz_n^\nu = (\lambda_nz_n)z_n^{\nu-1})$  which converges to $y_\nu = 0$. So $y$ is in $\mathbb{C}e_1$. Contradiction.\\

Note that $\bar C$ is the union $C \cup \mathbb{C}e_1$, because $\mathbb{C}e_1$ is the tangent to $\Gamma$ at the origin. So we have for each $\lambda \in \mathbb{C}$ the equality $\lambda e_1 = \lim_{n \to \infty} \lambda n \gamma(1/n) $.\\

We shall prove now that $\bar C$ is smooth around a point $x_1^0 e_1$ when $x_1^0 \not = 0$.\\
Let $\pi : H \to \mathbb{C}^2$ be the projection on the coordinates of $e_1$ and $e_2$. Take a point $(x_1, x_2)$ in $\mathbb{C}^2$ near to $(x_1^0, 0)$. Then if $x_2 = 0$ the only point in $\pi^{-1}(x_1, 0) \cap \bar C$ is the point $x_1 e_1$. If $x_2 \not= 0$ a point in $\pi^{-1}(x_1, x_2) \cap \bar C$ is in $C$ so  is of the form $\lambda\gamma(z)$ with $\lambda z\not= 0$;  and so we must have $x_1 = \lambda z$ and $x_2 = \lambda z^2$. This implies $z = x_2/x_1$ and $\lambda = x_1^2/x_2$. To prove that $\bar C$ is smooth near $x_1^0e_1$ it is enough to show, as $\bar C$ is the graph of the map $\theta$ defined by:
$$ \theta(x_1, x_2) = (x_1^2/x_2)\gamma(x_2/x_1) \quad {\rm for} \ x_2 \not= 0 \quad {\rm and} \quad  \theta(x_1,0) = x_1 e_1 $$
that $\theta$ is holomorphic near $(x_1^0, 0)$ with $x_1^0 \not= 0$.\\
We have for $z \in D$ the equality  $\gamma(z) = z\delta(z)$ where $\delta(z) := \sum_{\nu= 0}^\infty  z^\nu e_{\nu+1}$ is a holomorphic function on $D$. This allows to write
$ \theta(x_1, x_2) = x_1\delta(x_2/x_1)$ which gives the holomorphy of $\theta $ near $(x_1^0, 0)$ when $x_1^0 \not= 0$. So $\bar C \setminus \{0\} = (C \cup \mathbb{C}e_1) \setminus \{0\}$ is a closed connected $2-$dimensional sub-manifold of $H \setminus \{0\}$.

\smallskip
We have to prove now that $\bar C$ is not a finite dimensional analytic subset in $H$. If this is not true, the Enclosability Theorem (see \cite{[e]} ch.III section 7)  gives us a locally closed finite dimensional complex sub-manifold of $H$ in an open neighborhood $\mathcal{U}$ of $0$ in $H$ which contains $\bar C\cap \mathcal{U}$. But if a cone is contained near $0$ in a complex sub-manifold $W$, then it is contained in the tangent space $V$ of $W$ at $0$. Then $\bar C$ would be contained in the finite dimensional complex vector  space $V$. Then $V$ has to contain each $e_n, n \in \mathbb{N}^*$, because it has to contain $\Gamma$ and then all the Taylor coefficients of the map $\gamma$ at the origin.\\
 Contradiction.$\hfill \blacksquare$\\
\end{enumerate}

\chapter{The space $\mathcal{C}_n^f(M)$}

{\em In this chapter $M$ will always be a reduced complex space.}\\

We begin this section by giving a brief account of the topological space $\mathcal{C}_n^{loc}(M)$ with emphasis on the characterization of compact subsets of $\mathcal{C}_n^{loc}(M)$\footnote{For a detailed discussion see chapter 4 of \cite{[e]} of or \cite {[BM.1]}.}.

\section{Compactness in $\mathcal{C}_n^{loc}(M)$: The bounded local volume property}

The set of all (closed)  $n-$cycles in $M$ is denoted $\mathcal{C}_n^{loc}(M)$ and is endowed with the topology generated by all subsets  $\Omega_k(E)$ defined by
$$ 
\Omega_k(E) := \{ X \in \mathcal{C}_n^{loc}(M) \ / \  E \  {\rm is \ adapted \ to } \ X \ {\rm and} \ \deg_E(X) = k \}.
$$
where $E$ is an $n$-scale (see Terminology in  section I.3) on $M$ and $k$ a natural number. With this topology $\mathcal{C}_n^{loc}(M)$ is a {\em second countable Hausdorff space} (see \cite{[e]} Theorem 4.2.28).

For the study of the fibers of a holomorphic map $\pi$  the subset of $\pi$-relative cycles will be helpful.

\begin{prop}\label{relatif.loc}
Let $\pi\colon M \rightarrow N$ be a holomorphic map between two reduced complex spaces, let $\mathcal{C}_{n}^{\rm loc}(\pi)$ denote the subset of $\mathcal{C}_{n}^{\rm loc}(M)$ consisting of those $n$-cycles which are contained in a fiber of $\pi$\footnote{In other words  $\mathcal{C}_{n}^{\rm loc}(\pi)$ is the set of {\em $\pi$-relative $n-$cycles  in $M$}.} and let $\mathcal{C}_{n}^{\rm loc}(\pi)^*$ denote the subset of all non empty cycles in  $\mathcal{C}_{n}^{\rm loc}(\pi)$, i.e. 
 $\mathcal{C}_{n}^{\rm loc}(\pi)^* := \mathcal{C}_{n}^{\rm loc}(\pi)\setminus\{\emptyset[n]\}$. 
\begin{enumerate}[(i)]
\item
The subset $\mathcal{C}_{n}^{\rm loc}(\pi)$ is a closed subset of $\mathcal{C}_{n}^{\rm loc}(M)$. 
\item
The obvious map $\lambda\colon \mathcal{C}_{n}^{\rm loc}(\pi)^*\rightarrow N$, which associates to each (non empty) $\pi$-relative $n$-cycle the unique point in $N$ whose fiber contains the cycle, is continuous.
\end{enumerate}
\end{prop}
\parag{Proof}
To prove (i)  we take a cycle $X_{0}$ in $\mathcal{C}_{n}^{\rm loc}(M)\setminus\mathcal{C}_{n}^{\rm loc}(\pi)$. Then $|X_0|$ contains two points $x$ and $y$ such that $\pi(x) \neq \pi(y)$ so there exist disjoint open neighborhoods, $U$ of $\pi(x)$ and $V$ of $\pi(y)$, in $N$. Let  $E$ and $E'$ be two $n$-scales adapted to $X_{0}$ such that $x\in c(E)\subseteq \pi^{-1}(U)$ and $y\in c(E')\subseteq \pi^{-1}(V)$. It follows that the degrees $k := \deg_E(X_0)$ and $l := \deg_{E'}(X_0)$ are positive and consequently every $X\in\Omega_k(E)\cap\Omega_l(E')$ intersects more than one fiber of $\pi$. Hence $\Omega_k(E)\cap\Omega_l(E')$ is an open neighborhood of $X_0$ in $\mathcal{C}_{n}^{\rm loc}(M)\setminus\mathcal{C}_{n}^{\rm loc}(\pi)$.

\smallskip
For the proof of (ii) we take $X_0$ in $\mathcal{C}_{n}^{\rm loc}(\pi)^*$, put $y := \lambda(X)$ and fix an open neighborhood  $V$ of $y$ in $N$. Next we take a point $x$ in $|X_0|$ and an $n$-scale $E$ adapted to $X_0$ such that $x\in c(E)\subseteq \pi^{-1}(V)$ and put $k := \deg_E(X_0)$.   Then $k $ is positive so every $X$ in $\Omega_k(E)\cap \mathcal{C}_{n}^{\rm loc}(\pi)$ intersects $\pi^{-1}(V)$  and consequently $\lambda(X)\in V$.
\hfill{$\blacksquare$}\\

We shall now discuss  compactness in $\mathcal{C}_n^{loc}(M)$ and this discussion consists  more or less of  summarizing material from \cite{[e]} or \cite {[BM.1]}.\\ 
We begin by observing that a subset of  $\mathcal{C}_{n}^{loc}(M)$ is compact if and only it is sequentially compact since $\mathcal{C}_{n}^{loc}(M)$ is second-countable. 

\smallskip
The main tool to study compactness in $\mathcal{C}_n^{loc}(M)$ is Bishop's Theorem which gives a very simple characterization of relatively compact subsets in $\mathcal{C}_n^{loc}(M)$ in terms of the {\it local bounded volume property} (see Theorem \ref{BLV 2} below).

\medskip
A {\em continuous hermitian metric} on $M$ is a positive definite continuous  differential $(1,1)$-form on $M$ which is locally induced in local embeddings of $M$ into  open subsets of some affine space by positive definite continuous differential $(1,1)$-forms on these open sets. (See Chapter 3 of \cite{[e]} or of \cite {[BM.1]})

\smallskip
If $X$ is an $n$-cycle in $M$, $W$ a relatively compact open subset of $M$  and $h$ a continuous  hermitian metric on $M$, then
$$ 
{\rm vol}_h(X\cap W) := \int_{X \cap W} h^{\wedge n}
$$
is called the {\em volume of $X$ in $W$} (with respect to $h$).

\begin{defn}\label{BLV 1}
Let $\mathcal{A}$ be a subset of $\mathcal{C}_n^{\rm loc}(M)$. We shall say that $\mathcal{A}$ has the {\bf \em bounded local volume property} (or {\bf\em BLV} property for short) when the following holds:
\begin{itemize}
\item
There exists a continuous hermitian metric $h$ on $M$  such that, for every relatively compact open subset $W$ of $M$, there exists a constant $C(W)$ satisfying
\begin{equation*}\tag{BLV}
\int_{X\cap W} h^{\wedge n} \leq C(W),\qquad\text{for all} \ \ X \in \mathcal{A}.
\end{equation*}
\end{itemize}
\end{defn}

We say that a family $(X_{\lambda})_{\lambda\in\Lambda}$ of $n$-cycles has the BLV property if the subset $\{X_{\lambda}\ /\ \lambda\in\Lambda\}$ of  $\mathcal{C}_n^{\rm loc}(M)$ has the BLV property.

\parag{Remarks}
 \begin{enumerate}[(i)]
\item 
The BLV property is independent of the choice of a continuous  hermitian metric $h$ because if $k$ is another continuous  hermitian metric then, for every relatively compact open subset $W$ of $M$, there exist  two positive  constants
$\gamma(h,k,W)$ and $\Gamma(h,k,W)$ such that 
$$
\gamma(h,k,W) h^{\wedge n} \ \leq\  k^{\wedge n} \ \leq \ \Gamma(h,k,W) h^{\wedge n}
$$ 
on $\bar W$ where the inequality is taken in the sense of P. Lelong. Then for any cycle $X$ it follows that
$$ 
\int_{X\cap W} k^{\wedge n} \ \leq \ \Gamma(h,k,W)\int_{X \cap W} h^{\wedge n}
$$
\item 
A subset $\mathcal{A}$ of  $\mathcal{C}_n^{loc}(M)$ has the BLV property if and only if it satisfies the following condition, which does not involve a choice of a hermitian metric.
\begin{itemize}
\item
For every continuous  Lelong-positive differential $(n,n)$-form $\varphi$ with compact support in $M$ there exists a constant $C(\varphi) > 0$ such that for every $X$ in $\mathcal{A}$ we have
\begin{equation*}
\int_{X} \varphi \  \ \leq \ C(\varphi) \tag{$@$} 
\end{equation*}
\end{itemize}
It is clear that $\mathcal{A}$ satisfies condition ($@$) if it has the BLV property since, for every continuous  Lelong-positive $(n,n)$-form $\varphi$ with compact support in $M$ and every hermitian metric $h$, there exists a constant $C(\varphi, h) > 0$ such that the estimate  $\varphi \ \leq \ C(\varphi, h) h^{\wedge n}$ holds in the sense of P. Lelong. 
Conversely, suppose that $\mathcal{A}$ satisfies condition ($@$) and let $W$ be a relatively compact open subset of $M$. Then, for any continuous function $\sigma$ with compact support in $M$ and values in $[0,1]$  such that $\sigma \equiv 1$ on $W$, we have 
$$ 
\int_{X \cap W} h^{\wedge n} \leq  \int_{X} \sigma h^{\wedge n} 
$$
for every $n-$cycle $X$ in $M$. Hence $\mathcal{A}$ has the BLV property since $\sigma h^{\wedge n}$ is a continuous Lelong-positive $(n,n)$-form with compact support in $M$.\\
Note that any continuous $(n, n)$-form $\varphi$ with compact support may be written as  $\varphi = \psi_1 + i\psi_2$ where $\psi_j$ are real for $j = 1, 2$ continuous and  with compact supports contained in the support of $\varphi$. Moreover any real continuous $(n, n)$-form $\psi$ with compact support  may be written $\psi = \psi_+ - \psi_-$ \ where $\psi_+$ and $\psi_-$ are continuous, positive in the sense of Lelong and with compact supports contained in the support of $\psi$.
\item 
Note that, by Proposition 4.2.17 in \cite{[e]}, the function $X \mapsto \int_{X} \varphi$ is  continuous on $\mathcal{C}_n^{loc}(M)$ for every continuous $(n,n)$-form $\varphi$ with compact support in $M$. It follows that a subset of $\mathcal{C}_n^{loc}(M)$ has the BLV property if and only if its closure  in $\mathcal{C}_n^{\rm loc}(M)$ has the BLV property.
\end{enumerate}

\parag{Example} Let us give an example of the decomposition $\psi = •psi_+ - \psi_-$ mentioned in Remark $(ii)$
\begin{align*} 
&i( fdx_1\wedge d\bar x_2 - \bar f d\bar x_1\wedge dx_2) = i(fdx_1 + dx_2)\wedge \overline{(fdx_1 + dx_2)} -  i(f\bar f dx_1\wedge d\bar x_1 + dx_2\wedge d\bar x_2) \\
& \qquad =  i(dx_1 + \bar f dx_2)\wedge \overline{(dx_1+ \bar f dx_2)} - i(dx_1\wedge d\bar x_1 + f\bar f dx_2\wedge d\bar x_2) 
\end{align*} 

We shall also use the following definition.

\begin{defn} \label{inequal}  For  two $n$-cycles $X$ and $Y$ in $M$, we write $ Y \leq X$ when  every irreducible component $\Gamma$ of $Y$ is an irreducible component of $X$ and the multiplicity of $\Gamma$ in $Y$ is at most equal to the multiplicity of $\Gamma$ in $X$. 
\end{defn}

\smallskip
We leave the proof of the following lemma as an exercise for the reader.

\begin{lemma}\label{order.volume}
Let $X$ and $Y$ be two $n$-cycles in $M$. Then $Y\leq X$ if and only if for every continuous hermitian metric $h$ on $M$ and every relatively compact open subset $W$ of $M$ we have
$$ 
\int_{Y \cap W} h^{\wedge n} \leq  \int_{X \cap W} h^{\wedge n}.
$$
\hfill{$\blacksquare$}
\end{lemma}

A direct consequence of the lemma is that, for  two $n$-cycles $X$ and $Y$ in $M$, we have $Y \leq X$ if and only if 
$$ 
\int_{Y} \varphi \ \leq \ \int_{X} \varphi
$$
for every continuous Lelong-positive $(n,n)-$form $\varphi$ with compact support in $M$.\\

Since, for every continuous $(n,n)$-form $\varphi$ with compact support in $M$, the function  $X \mapsto \int_{X} \varphi$ is  continuous on $\mathcal{C}_n^{loc}(M)$  the following result is an immediate consequence of Lemma \ref{order.volume}.

\begin{cor}\label{order.volume.1}
Let  $(Y_\nu)$ and $(X_\nu)$ be sequences in $ \mathcal{C}_n^{loc}(M)$ which converge respectively to $Y$ and $X$ and satisfy $Y_\nu \leq X_\nu$ for all $\nu$. Then we have $Y \leq X$.
\hfill{$\blacksquare$}
\end{cor}

Let us recall here Theorem 4.2.69 in \cite{[e]} which  is a rather simple  consequence of Bishop's Theorem.

\begin{thm}\label{BLV 2}
A subset $\mathcal{A}$ in $\mathcal{C}_n^{loc}(M)$ is relatively compact if and only if it has the bounded local volume property.
\hfill{$\blacksquare$}
\end{thm}

\begin{cor}\label{utile}
Let $(X_\nu)_{\nu \in \mathbb{N}}$ be a  sequence of non empty cycles in $\mathcal{C}_n^{loc}(M)$ converging to a cycle $X$. Choose for each $\nu $ an irreducible component $\Gamma_\nu$ of $X_\nu$. Then  $\{\Gamma_\nu\ /\  \nu \in \mathbb{N}\}$  is a relatively compact subset of $\mathcal{C}_n^{loc}(M)$ and every cycle $Y$ which is a limit of a subsequence of the sequence $(\Gamma_\nu)$ satisfies $Y \leq X$.
\end{cor}

\parag{Proof} By Lemma \ref{order.volume} the set $\{\Gamma_\nu\ /\  \nu \in \mathbb{N}\}$ has the bounded local volume property so it is relatively compact in $\mathcal{C}_n^{loc}(M)$ due to Theorem \ref{BLV 2}, and every cycle $Y$ which is a limit of a subsequence of the sequence $(\Gamma_\nu)$ satisfies $Y \leq X$ thanks to Corollary \ref{order.volume.1}.
\hfill{$\blacksquare$}\\

We give now some examples.

\parag{First example} For each integer $n \geq 1$  the $0$-cycle $X_n:= \left\{\frac 1n, \frac 1{(n-1)}, \dots, 1\right\}$ in $\C$ is compact and contained in the relatively compact open subset $\{z \in \mathbb{C} \ / \ |z| < 2\}$ of $\C$. The volume of $X_n$ is $n$ so it goes to infinity when $n$ goes to infinity and consequently $\{X_n\ /\ n\geq 1\}$ is a discrete closed subset of $\mathcal{C}_0^{loc}(\C)$.
\hfill{$\square$}

\parag{Second example} For each integer $n \geq 1$ consider the $0$-cycle $X_n:= n.\{0\}$ in $\C$. These cycles have  $\{0\}$ as support, but  $\{X_n\ /\ n\geq 1\}$ is  a discrete closed subset of $\mathcal{C}_0^{loc}(\C)$ since the volume of $X_n$ tends to infinity as $n$ tends to infinity.
\hfill{$\square$}\\

Even though these two simple examples are not very interesting they show what is going on when a sequence of cycles does not have bounded local volume property; namely that \lq\lq local branches\rq\rq\ of the cycles are piling up somewhere when counted with their multiplicities.\\

Our next example, where we give a sequence of irreducible cycles (in fact smooth and connected) whose local volume is not bounded, is much more interesting.

\parag{Third example} For every $n \in \mathbb{N}^*$ let $X_n$ be the irreducible $1$-cycle in $M := \C^2$ defined by $z_2 = z_1^n$.
 Choose a real number $\delta > 0$ and put
$$ W := \{ 1 < \vert z_1\vert < 100\} \times \{ 1+\delta < \vert z_2 \vert < 1+2\delta \} .$$
Then $W$ is a relatively compact open set in $M$. Let $h := \frac{i}{2}(dz_1\wedge d\bar z_1 + dz_2 + d\bar z_2)$ be the standard K\"ahler form on $\mathbb{C}^2$. Then the following lemma shows that the volume of $X_n$ in $W$ with respect to  $h$  goes to infinity when  $n$ goes to infinity.\\
Note that on the open set $V := \{\vert z_1 \vert < 1\}\times \mathbb{C}$ the sequence $X_n\cap V$ converges in $\mathcal{C}_1^{\rm loc}(\C^2)$ to the $1-$cycle $\{z_2 = 0\}$.\\
 It is easy to deduce from the lemma below that  the local volume property is not satisfied near any point of the unit circle.
\begin{lemma}\label{13/11}
For $n$ large enough we have
$$
{\rm vol}_h(X\cap W) = \int_{W\cap X_n} \frac{i}{2}(dz_1\wedge d\bar z_1 + dz_2 + d\bar z_2) = \pi\delta(3\delta+2)n + O(1) 
$$
when $n$ goes to $+\infty$.
\end{lemma}

\parag{proof} For $n \geq 1$ we have $(1+2\delta)^{1/n}  \leq  1 + 2\delta$  so defining $z_1 := \rho e^{i\theta}$ we obtain
\begin{align*}
& {\rm vol}_h(W\cap X_n) = 2\pi\int_{(1+\delta)^{1/n}}^{(1+2\delta)^{1/n}}  (\rho + n^2\rho^{2n-1}) d\rho =  2\pi\big(\frac{\rho^2}{2} + \frac{n\rho^{2n}}{2}\big)_{(1+\delta)^{1/n}}^{(1+2\delta)^{1/n}} \\
& \qquad = \pi\Big((1+2\delta)^{2/n} - (1+\delta)^{2/n}\Big) + \pi n\Big((1+2\delta)^2 - (1+\delta)^2\Big) \\
& \qquad  = \pi(3\delta^2 +2\delta)n + O(1) 
\end{align*} 
and  the proof is completed.
$\hfill \blacksquare$

\parag{Comment} This third example is not so easy to understand because in the real world (we mean in $\mathbb{R}^2$), the length of the  corresponding curves in $W \cap \mathbb{R}^2$ are bounded. The reason for this comes from the fact that in the real world the equation $x^n = y$ has at most two roots for each given $y$; so the curve $y = x^n$ near the rectangle $]1, 100[\times]1+\delta,1+2\delta[$ \ for $\delta > 0$ has a nice limit given  by the sub-manifold $\{x = 1\}$ in this rectangle. But in the complex world the trace of $X_n$ on $W$ has $n$ branches which goes to the real hyper-surface $\{\vert z_1\vert = 1\}$ in $W$.

\parag{Fourth example} Consider the Hironaka's example which is described in \cite{[BM.1]} page 433  or \cite{[e]} p.444 (see also \cite{[B.15]}  which explains the construction of  this example with more details). Then, in this compact  complex connected $3$-dimensional manifold $Z$ (not projective but birational to $\mathbb{P}_3(\mathbb{C}$)), we have an analytic family of compact $1$-cycles parametrized by a smooth compact connected complex curve $T$ such that for a value $t_0 \in T$ we have the cycle $A$ and for an other value $t_1$ we have the value $A+B$ where $A$ and $B$ are two distinct smooth non empty $1$-cycles meeting at one point. Then it is possible to find a continuous map $ \gamma\colon ]0, 1] \to T$ such that, in the continuous family of compact cycle in $Z$ parameterized by $]0, 1]$ via $\gamma$, the cycle $\gamma(1/n)$ is equal to $A + nB$. Then, near the point in $A\cap B$ the local volume of such a continuous family of (compact) cycles in $Z$ goes to infinity.

\parag{Comment} There are two important examples of families of $n$-cycles which have the bounded local volume property.
\begin{itemize}
\item 
The connected  components of the space $\mathcal{C}_n(M)$ of compact $n$-cycles in a K\"ahler space are compact (see \cite{[e]} Corollary 2.7.26 in chapter IV or  \cite{[BM.2]} chapter XII for more details) and consequently every subset of an irreducible component of $\mathcal{C}_n(M)$ has the BLV property.
\item 
The theorem \ref{BLV 3} below which says that if $\pi\colon M \to N$ is a dominant holomorphic map  between irreducible complex spaces the family of general fibers of $f$ satisfies the local bounded volume property.
\end{itemize}
Note that in a projective complex space (i.e. a compact complex sub-space of some $\mathbb{P}_N(\mathbb{C})$) the degree of the cycles is locally constant in a continuous family, so the local bounded volume property holds when the parameter is connected. This is a special case of the first example above as a projective complex space is always a K\"ahler space, the degree being the volume for a suitable K\"ahler form defined by the Fubini-Study K\"ahler metric induced from an embedding in some $\mathbb{P}_N(\mathbb{C})$.
\hfill{$\square$}\\

The following theorem is an easy  consequence of the main result  in \cite{[B.78]} (see also Theorem 3.6.6 in  \cite{[e]}).

\begin{thm}\label{BLV 3}
Let $\pi\colon M \rightarrow N$ be a holomorphic map between irreducible complex spaces. Assume that $\pi$ is of generic rank equal to the dimension of $N$ and define  $n := \dim M - \dim N$. Let $N'$ be the subset of those points in  $N$ where the fiber of $\pi$ has dimension $n$. Then the subset $\{\pi^{-1}(y)\ /\  y \in N' \}$ of $\mathcal{C}_n^{loc}(M)$ has the bounded local volume property.
$\hfill \blacksquare$
\end{thm}

As a consequence of the previous two theorems, we obtain the following result, which will be of great importance in the sequel.

\begin{cor}\label{BLV 4}
Let $\pi\colon M \rightarrow N$ be a surjective holomorphic map between irreducible complex spaces and define  $n := \dim M - \dim N$. Let $N'$ be the subset of those points in  $N$ where the fiber of $\pi$ has dimension $n$\footnote{Observe that under these hypotheses $N'$ is dense in $N$ (see Proposition 2.4.60 in \cite{[e]})}. Let $\varphi\colon N' \rightarrow \mathcal{C}_n^{loc}(M)$ be the map given by $\varphi(y) := \pi^{-1}(y)$ (here $\pi^{-1}(y)$ is a reduced cycle),  let $\Gamma'$ be the graph of $\varphi$ and let $\Gamma$ be the closure  of   $\Gamma'$ in $N \times \mathcal{C}_n^{loc}(M)$. Then the natural projection $\tau\colon \Gamma \rightarrow N$ is proper (and consequently surjective).
\end{cor}

\parag{proof} The set $\varphi(N')$ is relatively compact in $\mathcal{C}_n^{\rm loc}(M)$, thanks to Theorem \ref{BLV 3} and Theorem \ref{BLV 2}, so $\Gamma$ is a closed subset of $N\times \overline{\varphi(N')}$ which is proper over $N$. It follows that $\tau\colon \Gamma \rightarrow N$ is proper.
$\hfill \blacksquare$\\

In contrast with the definition of a {\em strongly quasi-proper map} (see Chapter V  below)  where for a quasi-proper dominant holomorphic map $\pi$ we require that the closure in $N \times \mathcal{C}_n^{f}(M)$ of the  graph of the maximal  reduced fiber map is proper over $N$, condition which gives a non trivial restriction on such a map $\pi$, we see that  when we take  the closure of the graph  inside  $N \times \mathcal{C}_n^{loc}(M)$ the properness over $N$  is automatic !\\

 The previous corollary can be formulated in terms of  $\pi$-relative cycles.

\begin{cor}\label{fiber.map.loc}
Let $\pi\colon M \rightarrow N$ be a holomorphic map between two irreducible complex spaces and put $n := \dim M -\dim N$. Let $N'$ be a subset of $N$ on which there exists a continuous map   $\varphi\colon N'\rightarrow\mathcal{C}_{n}^{\rm loc}(\pi)^*$  such that $\varphi(y)$ is  the reduced cycle  equal to  $\pi^{-1}(y)$ for all $y$ in $N'$. Then the closure 
$\overline{\varphi(N')}$ in $\mathcal{C}_{n}^{\rm loc}(\pi)$  is a compact subset in $\mathcal{C}_{n}^{\rm loc}(\pi)$ \end{cor}

\parag{Proof} Since $\varphi$ is a continuous section on $N'$ of the natural  continuous map
$\lambda\colon \mathcal{C}_{n}^{\rm loc}(\pi)^*\rightarrow N$, defined in Proposition \ref{relatif.loc}, $\overline{\varphi(N')}$ is  compact thanks to Theorem \ref{BLV 3}.
\hfill{$\blacksquare$}\\

 If $\pi$ is equidimensional and if $N$ is normal, there always  exists  (see \cite{[e]} Corollary 4.3.13 ) a continuous map $\psi : N \to \mathcal{C}_n^{\rm loc}(\pi)$ which satisfies $\vert \psi(y)\vert = \pi^{-1}(y)$ for each $y \in N$ and a dense  set $N'$ of points such  $\psi(y)$ is reduced. In this case the closure of $\varphi(N')$ in $\mathcal{C}_n^{\rm loc}(\pi)$ is equal to $\psi(N) \cup \{\emptyset[n]\}$ when $N$ is not compact. So it is the Alexandroff compactification of $N$, as the natural map $\mathcal{C}_n^{\rm loc}(\pi)^* \to N$ gives a continuous inverse of $\psi$.\\

\section{Topology of  $\mathcal{C}_n^f(M) $}

In the sequel we shall consider $n$-cycles {\bf of finite type}, in other words the (closed) $n$-cycles in $M$ which have only finitely many irreducible components\footnote{This corresponds to  maps ${\rm Irr_n}(M) \rightarrow \mathbb{N}$ having finite support, where ${\rm Irr_n}(M)$ denotes the set of all non empty analytic subsets of dimension $n$ in $M$ (See \cite{[e]} Ch.4, Def. 4.1.1). In other words the $n-$cycles of finite type in $M$ form the free abelian monoid with basis  ${\rm Irr_n}(M)$}. They form a subset $\mathcal{C}_n^f(M)$ of $\mathcal{C}_n^{loc}(M)$ and we denote
$$ 
i : \mathcal{C}_n^f(M) \longrightarrow \mathcal{C}_n^{loc}(M)
$$
the natural injection.

\smallskip
For every relatively compact open set $W$ in $M$ we put 
$$ 
\Omega(W) := \{X \in \mathcal{C}_n^f(M) \ / \  {\rm each \ irreducible \ component \ of } \ X \ {\rm meets} \ W \} 
$$
and we endow $\mathcal{C}_n^f(M)$ with the coarsest topology which contains all such sets and makes the the injection $i$ continuous. So this is the topology generated by the sets $\Omega(W)$ and the sets 
$$
\Omega_k^f(E) := \Omega_k(E)\cap\mathcal{C}_n^f(M),
$$
where $W$ ranges over all relatively compact open subsets of $M$, $E$ ranges over all $n$-scales on $M$ and $k$ ranges over all natural integers.

\smallskip
For relatively compact open subsets $W_1, \dots, W_k$ in $M$ we put  
$$
\Omega(W_1, \dots, W_k) := \bigcap_{j \in [1, k]}\Omega(W_j).
$$

\begin{lemma}\label{Denom.}
For any complex space $M$ and any integer $n$ the topology of the space $\mathcal{C}_{n}^{f}(M)$ has a countable basis.
\end{lemma}

\parag{Proof} This an easy consequence of the analogous result for the topology of $\mathcal{C}_{n}^{loc}(M)$ which is proved in section 4.2.4 of \cite{[e]}.

\smallskip
As $M$ is second countable and locally compact its topology has a countable basis of relatively compact open subsets.
Let us fix such a basis and let $(U_n)_{n\in \mathbb{N}}$ be the countable family of all finite unions of sets belonging to this basis. Then, for any $n$-cycle $X$ in $\mathcal{C}_{n}^{f}(M)$ and any relatively compact open subset $V$ in $M$ satisfying $X\in\Omega(V)$, there exists $n$ in $\mathbb{N}$ such that $U_n$ is contained in $V$ and intersects every irreducible component of $X$, i.e. $X\in\Omega(U_n)\subset\Omega(V)$. It follows that, if $(\mathcal{U}_{\nu})_{\nu \in \mathbb{N}}$ is a countable basis for the topology of $\mathcal{C}^{loc}_{n}(M)$, then the  family $(\Omega(U_n)\cap \mathcal{U}_{\nu})_{(j,\nu)\in\mathbb{N}^2}$ is a countable subbasis for the topology of $\mathcal{C}_{n}^{f}(M)$.
$\hfill \blacksquare$

\begin{lemma}\label{1}
Let $M$ be a complex space and $n$ a natural number. Let $W$ be a relatively compact open set  in $M$ and  $V$ be an open subset of $W$. Then, $\Omega(V)$ is an open subset of $\Omega(W)$ in the topology induced by $\mathcal{C}_{n}^{\rm loc}(M)$.
\end{lemma}

\parag{Proof} It is enough to show that, for every $X$ in $\Omega(V)$, there exists an open neighborhood  $\mathcal{U}$ of $X$ in $\mathcal{C}_{n}^{loc}(M)$ such that 
$$
\mathcal{U} \cap \Omega(V) = \mathcal{U}\cap \Omega(W).
$$
Let us prove this by contradiction and assume that the result is not true. Then there exists a sequence $(X_\nu)_{\nu \geq 0}$ in $\Omega(W)\setminus\Omega(V)$ which converges to $X$ in $\mathcal{C}_n^{loc}(M)$. No $X_\nu$ is the empty $n$-cycle, because $\emptyset[n] \in\Omega(V)$, so for each $\nu$ at least one  irreducible component $\Gamma_\nu$ of $X_\nu$  does not meet $V$. Then, by Corollary \ref{utile}, there exists  a subsequence of the sequence $(\Gamma_\nu)$ which converges in $\mathcal{C}_n^{loc}(M)$ to a non empty cycle $Y$ which satisfies $Y \leq X$. Note that $Y$ is not empty because each $\Gamma_\nu$ meets $W$ which is relatively compact and so $Y$ contains at least one point in $\bar W$.  Moreover  each  irreducible component $\Gamma$ of $Y$  does not meet $V$ because each $\Gamma_\nu$ lies in the closed set $M \setminus V$. This  contradicts the fact that $Y \leq X$ since every irreducible component of $X$ intersects $V$.
$\hfill \blacksquare$\\

The following corollary is an obvious consequence of  Lemma \ref{1}.

\begin{cor}\label{2}
Let $M$ be a complex space and let $(W_{m})_{m \in \mathbb{N}}$ be an exhaustive sequence of relatively compact open subsets in $M$. Then any open set in the topology of $\mathcal{C}_{n}^{f}(M)$ is a union of some open sets  \  $ \mathcal{U} \cap \Omega(W_{m})$ where $\mathcal{U}$ is an open set in $\mathcal{C}_{n}^{loc}(M)$ and $m \in \mathbb{N}$.
\hfill{$\blacksquare$}
\end{cor}

Observe that Lemma \ref{Denom.} is a direct consequence of Corollary \ref{2} since the latter gives a simple way to construct a countable basis for the topology of $\mathcal{C}_{n}^{f}(M)$ from a countable basis for the topology of $\mathcal{C}_{n}^{loc}(M)$ and an exhaustion of $M$ by compact subsets. Note that the countable basis of the topology of $\mathcal{C}_n^f(M)$ obtained in this corollary is {\em a priori} \lq\lq smaller\rq\rq\  than the one given in Lemma \ref{Denom.}.

\begin{defn}
We say that a family of $n$-cycles in a complex space $M$, parameterized by a topological  Hausdorff space $S$, is {\bf\em $f$-continuous} if its classifying map induces a continuous map from $S$ into $\mathcal{C}_n^f(M)$.
\end{defn}

We end this paragraph by establishing a necessary and sufficient condition for a family of $n$-cycles in a complex space $M$, parametrized by a first countable Hausdorff space $S$,  to be $f$-continuous.

\begin{lemma}\label{$f-$continuous}
Let $(X_s)_{s\in S}$ be a family of $n$-cycles in a complex space $M$, parameterized by a first-countable Hausdorff space $S$, and put 
$$
G := \{(s, x)\in S\times M\ /\ x\in |X_s|\}.
$$ 
Then $(X_s)_{s\in S}$ is $f$-continuous if and only if:
\begin{enumerate}[(i)]
\item  
The family $(X_s)_{s\in S}$ is continuous in $\mathcal{C}_n^{loc}(M)$. 
\item  
The correspondence $(M,S,G)$ is quasi-proper\footnote{ In the case where $S$ is locally compact this condition is equivalent to the quasi-properness of the natural projection $G \rightarrow S$.}. 
\end{enumerate}
\end{lemma}

\parag{Proof} Since the case $S = \emptyset$ is trivial we may assume $S\neq\emptyset$.

\smallskip
Remark also that,  if there exists $s_0 \in S$ such that $X_{s_0}$ has infinitely many irreducible components, then  the correspondence $(M,S,G)$ cannot be quasi-proper at $s_0$, 
so we may assume that the  classifying map for the family $(X_s)_{s\in S}$ factorizes through a map  $\varphi\colon S \rightarrow \mathcal{C}_n^f(M)$. \\
Hence it is enough to prove that \ $\varphi$ \ is continuous if and only if $(M,S,G)$ is quasi-proper, assuming that the composition of $\varphi$ and the canonical injection  of $ \mathcal{C}_n^f(M)$ into $\mathcal{C}_n^{\rm loc}(M)$ is continuous.\\ 
Suppose first that $\varphi$ is continuous and let $s_0$ be a point in $S$. To show that $(M,S,G)$ is quasi-proper at $s_0$ let us take a relatively compact open subset $U$ of $M$ which intersects every irreducible component of the cycle $\varphi(s_0)$. Then $\varphi^{-1}(\Omega(U))$ is an open neighborhood of $s_0$ in $S$ and every irreducible component of $\varphi(s)$ intersects the compact subset $\bar{U}$ of $M$ for all $s\in \varphi^{-1}(\Omega(U))$. Thus $(M,S,G)$ is quasi-proper at $s_0$.\\
Now suppose that  $(M,S,G)$ is quasi-proper and 
let $(W_m)$ be an exhaustive sequence of relatively compact open subsets of $M$. Then due to Corollary \ref{2} the map $\varphi$ is continuous if and only if $\varphi^{-1}(\Omega(W_m))$ is an open subset of $S$ for all $m$. 

Let us fix $m$ and take a point $s_0$ in $\varphi^{-1}(\Omega(W_m))$. Then, as $(M,S,G)$ is quasi-proper, there exists a neighborhood \ $S_0$ of $s_0$ in $S$ and a compact subset $K$ of $M$ such that every irreducible component of $\varphi(s)$ intersects $K$ for all $s\in S_0$. But $(W_m)$ being an exhaustive sequence of relatively compact open subsets of $M$ there exists an integer $m_1 > m$ such that $W_{m_1}$ contains $K\cup W_m$. It follows that $\varphi^{-1}(\Omega(W_{m_1}))$ is a neigborhood of $s_0$ in $S$ and consequently $\varphi^{-1}(\Omega(W_m))$ is also a neighborhood of $s_0$ in $S$ since $\varphi^{-1}(\Omega(W_m))$ is an open subset of $\varphi^{-1}(\Omega(W_{m_1}))$ thanks to Lemma \ref{1}. It follows that $\varphi^{-1}(\Omega(W_m))$ is an open subset of $S$.$\hfill \blacksquare$

\newpage

\section{Compactness in $\mathcal{C}_n^f(M)$}

\subsection*{Introduction}

The image in $\mathcal{C}_n^{\rm loc}(M)$ of the closure of a  relatively compact subset  $\mathcal{A}$ of $\mathcal{C}_n^f(M)$ is  compact  so it has the BLV property. Hence $\mathcal{A}$ has also the BLV property. But an important difference between the topological spaces $\mathcal{C}_n^{loc}(M)$ and $\mathcal{C}_n^f(M)$ is the following:\\
 A sequence of non-empty $n$-cycles can converge to the empty $n$-cycle $\emptyset[n]$ in $\mathcal{C}_n^{loc}(M)$, and this can lead to many annoying problems. Fortunately this phenomenon  can not occur in the topological space $\mathcal{C}_n^f(M)$ since $\emptyset[n]$ is an isolated point in $\mathcal{C}_n^f(M)$. In fact the singleton $\{\emptyset[n]\}$ is a closed subset of $\mathcal{C}_n^f(M)$ since the space is Hausdorff, but to see that it is also open we recall that the empty $n$-cycle does not have any irreducible component (an irreducible component being non empty by definition) and consequently $\Omega(\emptyset) = \{\emptyset[n] \}$. In fact $\emptyset[n]$ is an isolated point in the open set  $\Omega(V)$ of $\mathcal{C}_n^f(M)$ for every relatively compact open subset $V$ of $M$.

\medskip
To get a better understanding of  the relative compactness of a subset in $\mathcal{C}_n^f(M)$ we introduce two more notions. The first one, called the {\it no escape to infinity} property, comes from the following result (see for instance Corolloray III 3.3.6 below):

\begin{itemize} 
\item Let $(X_\nu)_{\nu \in \mathbb{N}}$ be a sequence in $\mathcal{C}_n^f(M)$ of non empty cycles  which converges in $\mathcal{C}_n^{f}(M)$  to a  (finite type) cycle $X$. Let $(\Gamma_\nu)_{\nu \in \mathbb{N}}$ be a sequence obtained by choosing for each $\nu \in \mathbb{N}$ an irreducible component of the cycle $X_\nu$. Then $\{\Gamma_\nu\ /\ \nu \in \mathbb{N}\}$ is a relatively compact subset of $\mathcal{C}_n^f(M)$. Hence there exists a subsequence of the sequence  $(\Gamma_{\nu})_{\nu \in \mathbb{N}}$ which converges   in $\mathcal{C}_n^{f}(M)$ to a {\bf non empty} cycle $Y$ which, thanks to Corollary \ref{order.volume.1}, satisfies $Y \leq X$.
\end{itemize} 
Then, if a subset $\mathcal{A}$ of $\mathcal{C}_n^f(M)$ is relatively compact  in $\mathcal{C}_n^f(M)$, it satisfies the following property:
\begin{itemize}
\item 
No sequence $(\Gamma_\nu)_{\nu \in \mathbb{N}}$ of irreducible components of cycles in $\mathcal{A}$ can {\em escape to infinity}, in other words no such sequence satisfies:
\begin{equation*} 
\forall K \quad {\rm compact \ in} \ M \quad \exists \nu_K \ {\rm such \ that} \ \forall \nu \geq \nu_K \quad \Gamma_\nu \cap K = \emptyset .\tag{EI}
\end{equation*}
\end{itemize}
We observe that such a sequence $(\Gamma_\nu)_{\nu \in \mathbb{N}}$ satisfies condition (EI) if and only if it goes to infinity in $M$ as a sequence of closed subsets. It is also easy to see that the sequence $(\Gamma_\nu)_{\nu \in \mathbb{N}}$ escapes to infinity if and only if it converges to the empty $n$-cycle in the topology of $\mathcal{C}_n^{\rm loc}(M)$.\\
We  shall define the $NEI$ property as the negation of the condition  $EI$ above (see Definition III 3.2 below).\\

The second notion we are going to introduce is simpler. Consider a subset $\mathcal{A}$ of  $\mathcal{C}_n^f(M)$  and let  $(X_\nu)_{\nu \in \mathbb{N}}$ be a sequence in $\mathcal{A}$ which converges in $\mathcal{C}_n^{loc}(M)$ to a cycle $X$ which has infinitely many irreducible components. Then $\mathcal{A}$ cannot be relatively compact in $\mathcal{C}_n^f(M)$.
So a relatively compact subset $\mathcal{A}$ of $\mathcal{C}_n^f(M)$ satisfies the following property, called the {\it finite type limit} property:
\begin{itemize}
\item
The closure of $\mathcal{A}$ in $\mathcal{C}_n^{loc}(M)$ is contained in $\mathcal{C}_n^f(M)$.\hfill(FTL)
\end{itemize}
The main goal of this section is to prove Theorem \ref{BLV+NEI + FTL} below which gives a  characterization of relatively compact subsets in $\mathcal{C}_n^f(M)$ in terms of the three properties, BLV,  NEI and FTL.

\subsection{No Escape to Infinity property}

\begin{defn}\label{NEI 1}
We say that a subset $\mathcal{A}$ of $\mathcal{C}_n^{\rm loc}(M)$  has the {\bf\em no escape to infinity} property (or {\bf\em NEI} property for short) when the following condition holds:
\begin{itemize}
\item 
There exists a compact subset $K$ of $M$ such that for every $X$ in $\mathcal{A}$ and every irreducible component $\Gamma$ of $X$ we have $\Gamma \cap K \not=  \emptyset \hfill (NEI)$
\end{itemize}
\end{defn}

The fact that a subset $\mathcal{A}$ does not have the NEI property is then equivalent to the following, which explains our terminology. 
\begin{itemize}
\item 
There exists a sequence $(X_\nu)_{\nu \geq 0}$ in $\mathcal{A}$ and for each $\nu$ an irreducible component $\Gamma_\nu$ of $X_\nu$ such that  the sequence $(\Gamma_\nu)_{\nu \geq 0}$ {\em escapes to infinity}.
\end{itemize}

\parag{Remarks} 
\begin{enumerate}[(i)]
\item 
If a subset $\mathcal{A}$ of $\mathcal{C}_n^{loc}(M)$ has the NEI property, then it is contained in $\mathcal{C}_n^f(M)$.
\item 
Every subset of $\mathcal{A}$ has the NEI property if $\mathcal{A}$ has the NEI property.
\item 
A finite union of subsets having the NEI property has again the NEI property.
\item The singleton $\{\emptyset[n]\}$ has the NEI property. If $\mathcal{A}$ has the NEI property and contains $\{\emptyset[n]\}$, then $\mathcal{A} \setminus \{\emptyset[n]\}$ has again the NEI property.$\hfill \square$\\
\end{enumerate}

Let $X$ be a finite type $n$-cycle in a complex space $M$ in which there exists a sequence of irreducible $n$-cycles converging to the empty $n$-cycle $\emptyset[n]$. Then in any neighborhood of $X$ in $\mathcal{C}_n^{loc}(M)$, the property  $EI$ is true.\\
On the contrary, the lemma below shows that the topology defined on $\mathcal{C}_n^f(M)$ avoids this pathology.

\begin{lemma}\label{NEI 2} Any $X$ in $\mathcal{C}_n^f(M)$ has a neighborhood in $\mathcal{C}_n^f(M)$ which has the NEI property. So any compact subset $\mathcal{A}$ in $\mathcal{C}_n^f(M)$ has  the NEI property.
\end{lemma}

\parag{proof} For each $X$ in $\mathcal{C}_n^f(M)$ there exists a relatively compact open subset $W(X)$ of $M$ which intersects every irreducible component of $X$ and then any irreducible component of any $Y \in \Omega(W(X))$ meets $\overline{W(X)}$ and $\Omega(W(X))$ has the NEI property.\\
 When  $\mathcal{A}$ is compact the open subsets $\Omega(W(X)), X \in \mathcal{A}$ of $\mathcal{C}_n^f(M)$ form a covering of $\mathcal{A}$ and consequently there exist $X_1,\ldots,X_k$ in $\mathcal{A}$ such that $\mathcal{A}\subseteq\bigcup_{j=1}^k\Omega(W(X_j))$, since $\mathcal{A}$ is compact. Hence every irreducible component of any $Y$ in $\mathcal{A}$ intersects the compact subset $L := \bigcup_{j \in [1,k]} \overline W(X_j)$ of $M$.$\hfill \blacksquare$

\parag{Remark}It is important to notice that if a subset $\mathcal{A}$ in $\mathcal{C}_n^f(M)$ has the NEI property, its closure $\bar{\mathcal{A}}$ in $\mathcal{C}_n^f(M)$ may not have this property. This means that we might have a sequence $(X_\nu)_{\nu \geq 0}$ in $\partial\mathcal{A}$ and for each $\nu$ an irreducible component $\Gamma_\nu$ of $X_\nu$ such that the sequence $(\Gamma_{\nu})_{\nu\geq 0}$ escapes to infinity, as  is shown in the following example.

\parag{Example} Let $\mathcal{A}$ be  the subset of $\mathcal{C}_1^f(\mathbb{C}^2)$ which is defined by the family of  irreducible conics $X_{s,t} := \{(x,y) \in \mathbb{C}^2\ / \ x(sy+1) -t = 0 \}$ parameterized by $(s,t) \in D^*\times D^*$, where $D$ is the unit disc with center $0$ in $\mathbb{C}$. Then it is easy to see that $X_{s,t}$ contains the point $x= t, y= 0$ so that it meets the compact set $\bar D$. So the subset $\mathcal{A}$ has the NEI property. For any $s \in D^*$ the closure of $\mathcal{A}$ in $\mathcal{C}_1^f(\mathbb{C}^2)$ contains the $1$-cycle $X_{s, 0} := \{x= 0\} + \{y = -1/s\}$. The irreducible component $\{y = -1/s\}$ escapes at infinity when $s \to 0$ so $\bar{\mathcal{A}}$ does not have the NEI property.

\subsection{The Finite Type Limit  property}

Suppose $\mathcal{A}$ is a relatively compact subset of $\mathcal{C}_n^f(M)$. Then, as the natural injection $\mathcal{C}_n^f(M)\rightarrow  \mathcal{C}_n^{\rm loc}(M)$ is continuous and $\mathcal{C}_n^{\rm loc}(M)$ is a Hausdorff space, $\mathcal{A}$ satisfies the following condition:

\begin{itemize}
\item 
The closure of $\mathcal{A}$ in $\mathcal{C}_n^{\rm loc}(M)$ is contained in $\mathcal{C}_n^f(M)$. 
\hfill (FTL)
\end{itemize}

We shall say that a subset $\mathcal{A}$ of $\mathcal{C}_n^f(M)$ has the {\bf finite type limit} property (or the {\bf FTL} property for short) if it satisfies the condition above.

\medskip
The following example shows that a subset of $\mathcal{C}_n^f(M)$ (even a closed one) can have the properties BLV and NEI without without having property FTL.

\parag{Example} Let $\tau\colon M\rightarrow\C^2$ be the blow-up of  $\mathbb{Z}\times \{0\}$ in $\C^2$ and let $\pi\colon M \rightarrow \C$ be the composition of $\tau$ and the second projection 
$\C^2\rightarrow\C$. Denote $\mathcal{A}$ the the set of fibers of $\pi$ over $\bar D^* := \bar D \setminus \{0\}$. Then the closure of $\mathcal{A}$ in $\mathcal{C}_1^{loc}(M)$ contains the cycle $\pi^{-1}(0)$ which has infinitely many irreducible (compact) components. It is also easy to see that $\mathcal{A}$ is closed in $\mathcal{C}_1^f(M)$ (thanks to the uniqueness of  limits in 
$\mathcal{C}_n^{\rm loc}(M)$)  and satisfies the NEI property since every cycle in $\mathcal{A}$ is irreducible and intersects the compact subset $\tau^{-1}(\{1/2\}\times \bar D)$ of $M$. Now the sequence $(\pi^{-1}(1/n))_{ n \geq 1}$ in $\mathcal{A}$ has clearly no convergent subsequence in $\mathcal{C}_1^f(M)$ so $\mathcal{A}$ is not (relatively) compact in $\mathcal{C}_1^f(M)$.
\hfill{$\square$}

\begin{thm}\label{BLV+NEI + FTL}
A  subset of $\mathcal{C}_n^f(M)$ is relatively compact if and only if it has the BLV property, the NEI property and the FTL property. 
\end{thm}

\parag{proof}  We have already seen that every relatively compact subset of $\mathcal{C}_n^f(M)$ has these three properties.

Conversely, suppose that $\mathcal{A}$ is a subset of $\mathcal{C}_n^f(M)$ which has the three properties.
Take an arbitrary sequence $(X_\nu)_{\nu \geq 0}$ in $\mathcal{A}$ and let us show that it has a subsequence which converges in $\mathcal{C}_n^f(M)$.
 As $\mathcal{A}$ satisfies the BLV property, it  is relatively compact in $\mathcal{C}_n^{\rm loc}(M)$ and consequently $(X_\nu)_{\nu \geq 0}$ has a subsequence $(X_{\nu_j})_{j \geq 0}$ which converges in $\mathcal{C}_n^{loc}(M)$ to an $n$-cycle $X$. Moreover $X$ is in $\mathcal{C}_n^f(M)$ since $\mathcal{A}$ has the FTL property. To show that  $(X_{\nu_j})_{j \geq 0}$ converges to $X$ in $\mathcal{C}_n^f(M)$ it is enough to prove that, for every relatively compact open subset $V$ of $M$ which intersects every irreducible component of $X$, there exists $j_V\geq 0$ such that $X_{\nu_j}\in\Omega(V)$ for all  $j\geq j_V$.  To this end we fix such an open subset $V$ of $M$. As $\mathcal{A}$ has the  NEI property there exists a compact subset $K$ of $M$ which intersects every irreducible component of $X_{\nu_j}$ for all $j$, so if we take a relatively open subset $W$ of $M$ which contains $V$ and $K$, then $X\in\Omega(W)$ and $X_{\nu_j}\in\Omega(W)$  for all $j$.  Then, due to  Lemma \ref{1}, there exists an open neighborhood $\mathcal{U}$ of $X$ in $\mathcal{C}_n^{\rm loc}(M)$ having the property that
$$
\mathcal{U} \cap \Omega(V) = \mathcal{U}\cap \Omega(W).
$$
But $(X_{\nu_j})_{j \geq 0}$ converges to $X$ in $\mathcal{C}_n^{loc}(M)$ so there exists $j_V\geq 0$ such that $X_{\nu_j}\in\mathcal{U}$ and consequently $X_{\nu_j}\in \Omega(V)$ for all $j\geq j_V$. 
 $\hfill \blacksquare$\\
 
 We shall show now that if a subset of $\mathcal{C}_n^{f}(M)$ has the BLV property then the number of irreducible components of individual cycles belonging to this set cannot be unbounded without presenting some escape to infinity. This is precisely formulated in the proposition below.
 
\begin{prop}\label{FTI 1}
Let $\mathcal{A}$ be a subset of $\mathcal{C}_n^f(M)$ having properties BLV and NEI. Then there exists an integer $N$ such that every cycle $X$ in $\mathcal{A}$ has at most $N$ irreducible components (counted with multiplicities).
\end{prop}
 
The proof will use the following lemma.

\begin{lemma}\label{minimal volume}
Let $M$ be a complex space and $h$ be a continuous  hermitian metric on $M$. Let $K$ be a compact set in $M$ and  $W$ be a relatively compact open neighborhood of $K$ in $M$. Then there exists a positive number $\alpha(h,K,W)$ such that for any irreducible $n$-dimensional analytic set $\Gamma$ in $M$ which meets $K$ we have $\vol_h(\Gamma \cap W) \geq \alpha(h,K,W)$.
\end{lemma}

\parag{Proof} Assume the lemma is wrong. Then for each integer  $m \geq 1$ there exists an irreducible analytic subset $\Gamma_m$ of dimension $n$ which meets $K$ and satisfies 
$$
\vol_h(\Gamma_m \cap W) \leq 1/m
$$ 
In particular the sequence $(\Gamma_m\cap W)_{m\geq 1}$ has the BLV property so, by Theorem \ref{BLV 2}, it has a subsequence $(\Gamma_{m_j}\cap W)_{j\geq 1}$ which converges to an $n$-cycle $X$ in $\mathcal{C}_n^{\rm loc}(W)$. It follows that $\vol_h(X\cap W') = 0$ for every relatively compact open subset $W'$ of $W$ and consequently $X$ is the empty $n$-cycle in $W$. But, for every $j$, there exists a point $x_j$ in  $\Gamma_{m_j}\cap K$ and a subsequence of $(x_j)_{j\geq 1}$ which  converges to a point $x$ in $K$. This implies that $x\in |X|$ in contradiction to the fact that $|X|$ is empty. 
$\hfill \blacksquare$

\parag{Proof of Proposition \ref{FTI 1}} As $\mathcal{A}$ has the NEI property there exists a compact subset $K$ of $M$ such that every irreducible component of every cycle in $\mathcal{A}$ meets $K$.  Let $W$ be a relatively compact open neighborhood of $K$ in $M$ and $h$ be a continuous hermitian metric on $M$. 
Then  there exists a positive constant $C(h,W)$ such that, for every $X$ in $\mathcal{A}$, we have $\vol_h(X \cap W) \leq C(h,W)$ since $\mathcal{A}$ has the BLV property. Now, let $p(X)$ denote the number of  irreducible components of an $n$-cycle $X$ and let $\alpha(h,K,W)$ be the constant introduced in Lemma \ref{minimal volume}. Then, as every $X$ in $\mathcal{A}$ intersects $K$, we obtain the inequalities
$$  
p(X)\alpha(h,K,W) \leq \vol_h(X \cap W) \leq C(h,W)
$$
and consequently $p(X)\leq \frac{C(h,W)}{\alpha(h,K,W)}$ for all $X$ in $\mathcal{A}$.
\hfill{$\blacksquare$}\\

\begin{cor}\label{FTI 2}
Let $(X_\nu)_{\nu \geq 0}$ be a sequence in $\mathcal{C}_n^{f}(M)$ which converges in $\mathcal{C}_n^{\rm loc}(M)$ to a cycle $X$ which has infinitely many irreducible components. Assume moreover that $\{X_{\nu}\ /\ \nu\geq 0\}$ has the NEI property. Then there exists a subsequence $(X_{\nu_j})_{j\geq 0}$  and  for each $j$ an irreducible component $\Gamma_j$ of $X_{\nu_j}$ such that the  sequence $(\Gamma_j)_{j \geq 0}$ converges in $\mathcal{C}_n^{\rm loc}(M)$ to a cycle $Y \leq X$ which has also an infinite number of irreducible components.
\end{cor}

\parag{Proof} As $(X_\nu)_{\nu \geq 0}$ converges in $\mathcal{C}^{\rm loc}_n(M)$  the subset $\{X_{\nu}\ /\ \nu\geq 0\}$ has the BLV property. Then, thanks to Proposition \ref{FTI 1}, the number of  irreducible components of the cycles  $X_\nu$ is uniformly bounded. Hence, by taking a subsequence if necessary, we may assume that every  $X_\nu$  has exactly $N$ irreducible components (counted with multiplicities) so we can order the irreducible components and write $X_\nu = \Gamma_{\nu}^1 +\cdots +\Gamma_{\nu}^N$ for each $\nu$. Then, by taking a subsequence, we may assume that, for each $j$, the sequence $( \Gamma_{\nu}^j)_{\nu\geq 0}$ converges in $\mathcal{C}^{\rm loc}_n(M)$ to an $n-$cycle $Y_j$. It follows that $X = \sum_{j=1}^{N} Y_j$ and so there exists at least one  $j \in [1, N]$ such that $Y_j$ has infinitely many irreducible components. 
$\hfill \blacksquare$\\

Let us finish this paragraph by  giving two consequences of Theorem \ref{BLV+NEI + FTL} which are  easier tools to use for applications.

\begin{prop}\label{f-conv.}
Let $M$ be a reduced complex space and  $(X_{m})_{m \geq 0}$ be a sequence of finite type $n$-cycles  in $M$ having the following properties:
\begin{enumerate}[(i)]
\item 
There exists a compact subset $K$ in \ $M$ \ such that, for every $m\geq 0$,  every irreducible component of $X_{m}$  meets $K$.
\item 
The sequence $(X_{m})_{m \geq 0}$ converges in $\mathcal{C}_{n}^{\rm loc}(M)$ to an $n$-cycle \ $X$.
\item 
The cycle $X$ is in $\mathcal{C}_{n}^{f}(M)$.
\end{enumerate}
Then the sequence $(X_{m})_{m \geq 0}$ converges to $X$ in  $\mathcal{C}_{n}^{f}(M)$.
\end{prop}

\parag{Proof} Put $F := \{X_m\ /\  m \in \mathbb{N} \} $. Then $F$ is a relatively compact subset of  $\mathcal{C}_n^{\rm loc}(M)$ by (ii) and its closure in $\mathcal{C}_n^{\rm loc}(M)$  is  contained in $\mathcal{C}_n^f(M)$ by (iii). It follows that $F$ has the BLV and FTL properties. Moreover it satisfies the property NEI by  (i). Hence $F$ is a relatively compact subset of $\mathcal{C}_n^f(M)$ and, as any subsequence of $(X_{m})_{m \geq 0}$ which converges in $\mathcal{C}_n^f(M)$ must converges to $X$ by uniqueness of the limit in $\mathcal{C}_n^{\rm loc}(M)$, the conclusion follows.
$\hfill \blacksquare$

\begin{cor}\label{compact}
Let \ $M$ \ be a reduced complex space and let $\mathcal{A}$ be a subset of $\mathcal{C}_{n}^{f}(M) \setminus \{\emptyset[n]\}$. Then $\mathcal{A}$ is  compact in \ $\mathcal{C}_{n}^{f}(M)$ if and only if the following conditions are satisfied.
\begin{enumerate}[(i)]
\item
$\mathcal{A}$ is a compact subset of \ $\mathcal{C}_{n}^{\rm loc}(M)$.
\item 
There exists a compact subset $K$ of $M$ such that every irreducible component of every  $X$ in $\mathcal{A}$ intersects $K$.
\end{enumerate}
\end{cor}

\parag{Proof} Condition (ii) is equivalent to the NEI property so these conditions  are necessary.\\
Conversely, suppose that $\mathcal{A}$ satisfies conditions (i) and (ii). Condition (i) implies that every sequence in $\mathcal{A}$ has a subsequence which converges in $\mathcal{C}_{n}^{\rm loc}(M)$ to an $n$-cycle in $\mathcal{A}$, and thanks to Proposition \ref{f-conv.} condition (ii) implies that the subsequence converges in $\mathcal{C}_{n}^f(M)$ to this same $n$-cycle. Hence $\mathcal{A}$ is a sequentially compact subset of $\mathcal{C}_{n}^f(M)$ and consequently compact since $\mathcal{C}_{n}^f(M)$ is second-countable. 
\hfill$\blacksquare$

\chapter{$f$-Analytic Families of Cycles} 

\section{Introduction}

It is easy to define the notion of $f$-analytic family of $n$-cycles in a complex space $M$ parametrized by a Banach analytic set $S$ by adding to the classical definition a quasi-properness condition on the set theoretic graph $G$  of the family (when $S$ is not locally compact we ask that  the triple $(M, S, G)$ is quasi-proper). This is similar to the proper condition added in the case of compact cycles. But one crucial point in the compact case is the fact that all cycles nearby a given cycle $X_0$  may be described using a finite set of scales adapted to $X_0$.\\
In the finite type case, it is also possible to have a finite numbers of adapted scales such that each irreducible component of $X_0$ meets at least one center of these scales and then, there exists an open neighborhood of $X_0$, in the topology introduced on $\mathcal{C}_n^f(M)$, such that each irreducible component of any $X$ in this neighborhood meets at least one center of these scales.\\
The crucial point, in order that the weak analytic structure (defined mainly by the notion of $f$-analytic family introduced above) to  be close to the structure of  a Banach analytic set, is the Analytic Extension Theorem proved in section 3. It implies that, {\bf in the case where the parameter space  $S$ is a reduced complex space}, for a $f$-continuous family, the analyticity of the family in a finite set of adapted scales to the cycle $X_0$ (corresponding to $s_0 \in S$)  whose centers meet any irreducible component  of $X_0$, ensures the analyticity at the point $s_0$  of the family. Note that this demands the analyticity  at the point $s_0$ of the corresponding family of multigraph defined by this family {\bf in any $n$-scale on $M$ adapted to $X_0$}.\\
But it is important to be aware that this analytic extension theorem does not hold when the parameter space $S$ is not (locally)  finite dimensional, that is to say when the Banach analytic set $S$ is not a reduced complex space near $s_0$.\\
This is the reason for which the local classifying spaces we shall define  in analogy with these used  in the compact cycle case, are universal only for $f$-analytic families parametrized by reduced complex spaces.\\
Nevertheless it is interesting to define holomorphic maps from $\mathcal{C}_n^f(M)$ with values in a  Banach space and to consider analytic subsets defined (locally)  by such holomorphic  equations because by pull-back by a holomorphic map $\varphi : S \to \mathcal{C}_n^f(M)$ they will produce analytic subsets of $S$ even in the case where $S$ is a Banach analytic set.\\
But of course, for our main applications, $S$ will be assumed to be a  reduced  complex space.\\

This weak analytic structure, which is not locally isomorphic to a Banach analytic set in general,  will enable us to prove in section 4  a Direct Image Theorem for a semi-proper holomorphic map 
$$\varphi : S \longrightarrow \mathcal{C}_n^f(M)$$
 where $S$ is a reduced complex space, generalizing again Kuhlmann's Theorem. This is a major tool for applications.\\

We shall introduce, in section 5 below,  the notion of a meromorphic map from a reduced complex space $S$ to $\mathcal{C}_n^f(M)$ and this will be used to enlighten the notion of strongly quasi-proper map which will be introduced in the next chapter.

\section{Weak analytic structure on $\mathcal{C}_n^{f}(M)$} 

\subsection{Basic definitions}

\bigskip

In this section $M$ will be a reduced complex space and $i\colon\mathcal{C}_n^f(M)\hookrightarrow\mathcal{C}_n^{loc}(M)$ will denote the canonical injection.
 
 \bigskip
Although there is no hope to have, in general, even locally  a Banach analytic set  structure on the space $ \mathcal{C}_n^f(M)$, it is possible to define the notion of holomorphic map $S \to  \mathcal{C}_n^f(M)$ when $S$ is a Banach analytic set, and also the notion of a holomorphic map defined on an open set $\mathcal{U}$ of $ \mathcal{C}_n^f(M)$ with values in a Banach space $E$.\\

\begin{defn}\label{appl. hol. 1}
Let $S$ be a Banach analytic set and $(X_s)_{s \in S}$ be an analytic family of cycles in $M$. We say that the family is  $f${\bf\em-analytic} if $X_s$ is of finite type for every $s \in S$ and  the classifying map $\varphi : S \to \mathcal{C}_n^f(M)$ is continuous.\\
We say that a map $\psi : S \to \mathcal{C}_n^f(M)$ is {\bf \em holomorphic} when it classifies an f-analytic family.
\end{defn}

To decide whether or not an analytic family $(X_s)_{s \in S}$ of $n$-cycles is an $f$-analytic family of finite type $n$-cycles in $M$ it is sufficient to verify that the triple $(M, S, G)$  is quasi-proper, where $G \subset S \times M$ is the set-theoretic graph of the family.\\


We now show that the graph  of an $f$-analytic family of cycles in $M$ parameterized by an irreducible complex space is again of finite type.

 \begin{lemma}\label{finitude comp.}
 Let $M$ be a complex space and $N$ an irreducible complex space. Let $\varphi : N \to \mathcal{C}_{n}^{f}(M)$ be a holomorphic map. Let $G \subset N \times M$ be the set theoretic graph of the $f$-analytic family of $n$-cycles in $M$ classified by the map $\varphi$. Then $G$ has finitely many irreducible components.
 \end{lemma}
 
 \parag{Proof} The natural projection $\pi\colon G \rightarrow N$ is a holomorphic map which is both quasi-proper and open. As $\pi$ is equidimensional each irreducible component of $|G|$ is a union of irreducible components of the fibers of $\pi$ and consequently the restriction of $\pi$  to each irreducible component of  $G$ is a quasi-proper map.  Hence it is surjective since $N$ is irreducible.
Let us now fix one of the fibers of $\pi$. Then each irreducible component of $G$ contains at least one of the irreducible components of this fiber. But the map $\pi$ is quasi-proper so there exists a compact set in $G$ which intersects every irreducible component of the chosen fiber and intersects therefore every  irreducible component of  $G$.  It follows  that  $G$ has only  finitely many irreducible components.$\hfill \blacksquare$

\parag{{\bf IMPORTANT COMMENT}} Although we define above the notion of an holomorphic map $f : S \to \mathcal{C}_n^f(M)$ for $S$ a Banach analytic set, we shall now use only the case where $S$ is a reduced complex space in the definition of the weak analytic structure on $\mathcal{C}_n^f(M)$ (see Definition \ref{appl.hol. 2}). The reason is the following:\\
It may not be too difficult to prove the continuity of a map $f : S \to \mathcal{C}_n^f(M)$ even when $S$ is a Hausdorff topological space, but when $S$ is a Banach analytic set and when the cycles are not compact, to prove that  such a map $f$ composed with the inclusion $i : \mathcal{C}_n^f(M) \to \mathcal{C}_n^{\rm loc}(M)$ classifies an analytic family of cycles at a point $s_0 \in S$ is rather difficult. First because we cannot use the easy criterium when the parameter space is normal (see Theorem 4.3.12 in \cite{[e]}) because the Riemann extension theorem is valid only on non singular Banach analytic sets (that is to say for Banach manifolds) in the infinite dimensional case. So, in the  infinite dimensional singular case, to verify the analyticity of a family of non compact cycles at a point $s_0$ demands to check all scales in a covering of the cycle $X_{s_0}$. When $X_{s_0}$ is not compact, infinitely many scales will be necessary, and it is very difficult in practice to check the isotropy of the given family in each of these scales. Moreover, the intersection of the infinitely many neighborhood of $s_0$ in $S$ which appears in such a verification will intersect, in general, only in $\{s_0\}$, so that analyticity at this point will not implies the analyticity of the family at each point of an open neighborhood of $s_0$. \\
At the first glance, it seems that these difficulties have nothing to do with the fact that $S$ is a Banach analytic set or a reduced complex space. But in the case of a reduced complex space the situation will be much better thanks to the Analytic Extension Theorem which will be proved in the next section. It allows to show that it is enough to check the isotropy of the family in a finite set of scales. And so we can also deduced from this result that the analyticity condition is an {\bf open condition} on such an $S$.\\
But in the case of a Banach analytic set which is not a reduced complex space (and not a complex Banach manifold) this Analytic Extension Theorem is not true (a counter-example is described in \cite{[BM.2]} Chapter V, Proposition 2.4.2). This explains why we consider the "weak analytic structure" on $\mathcal{C}_n^f(M)$ as a functor defined only on the category of reduced complex spaces.$\hfill \square$\\

\begin{defn}\label{appl.hol. 2}
Let $\mathcal{U}$ be an open set in $\mathcal{C}_n^f(M)$ and $E$ be a Banach space. 
\begin{enumerate}[(i)]
\item
A mapping $g\colon \mathcal{U} \rightarrow E$ is called {\bf \em holomorphic} if it is continuous and if, for every reduced complex space $S$ and every  holomorphic map $\varphi\colon S \to \mathcal{U}$, the composed map
$$ 
g\circ\varphi\colon S \to E
$$
is holomorphic.
\item
We shall say that a subset $\mathcal{X}$  of \ $\mathcal{U}$ is {\bf\em analytic} if it is closed and if every point in $\mathcal{X}$ has an open neighbourhood $\mathcal{V}$ in $\mathcal{U}$ such that $\mathcal{X}\cap\mathcal{V}$ is the zero set of  a holomorphic map  with values in a Banach space. 
\item
Let $\mathcal{X}$ be an analytic subset of \ $\mathcal{U}$. A map $g\colon\mathcal{X}\rightarrow E$  is said to be {\bf\em holomorphic} if, for every point $x$ in  $\mathcal{X}$, there exists an open neighbourhood $\mathcal{V}$ of $x$ in $\mathcal{U}$ and a holomorphic map $G\colon\mathcal{V}\rightarrow E$ such that $g$ and $G$ coincide on $\mathcal{X}\cap\mathcal{V}$.\\
A {\bf\em holomorphic function} on $\mathcal{X}$ is a holomorphic map with values in $\C$.
\item
Let $\mathcal{X}$ be an analytic subset of \ $\mathcal{U}$ and  $T$ be a Banach analytic subset of an open subset of\ $E$. A map $g\colon\mathcal{X}\rightarrow T$  is called {\bf\em holomorphic} if it induces a holomorphic map $\mathcal{X}\rightarrow E$.
\item
Let $\mathcal{X}$ be an analytic subset of \ $\mathcal{U}$, $P$ be a reduced complex and $m$ be a natural number. We say that a map $g\colon\mathcal{X}\rightarrow \mathcal{C}_m^f(P)$ is {\bf\em holomorphic} if for every open subset $\mathcal{V}$ of $\mathcal{C}_m^f(P)$ and every holomorphic function $h$ on $\mathcal{V}$ the function $h\circ g_{|g^{-1}(\mathcal{V})}$ is holomorphic on $g^{-1}(\mathcal{V})$.
\item
Let $N$ be a complex space, $k$ a natural number and $\mathcal{V}$ an open subset of $\mathcal{C}_k^f(N)$. A map $g\colon\mathcal{V}\rightarrow\mathcal{C}_n^f(M)$ is {\bf\em holomorphic} if for every holomorphic map $\varphi\colon S\rightarrow\mathcal{V}$, where $S$ is a reduced complex space, the composition $g\circ\varphi$ is holomorphic.
\item Let $\mathcal{X}$ be an analytic subset of \ $\mathcal{U}$. We shall say that $\mathcal{X}$ is a  {\bf (reduced) complex subspace} of  \ $\mathcal{U}$ if, endowed with the sheaf of holomorphic functions which is defined above, it becomes a reduced complex space.
\end{enumerate}
\end{defn}

\parag{Remarks} \begin{enumerate} [(i)]
\item The sheaf of holomorphic functions  $\mathcal{O}_{\mathcal{C}_n^f(M)}$ on $\mathcal{C}_n^f(M)$ determines the weak holomorphic structure.
\item Let $\mathcal{X}$ be a closed subset of an open set $\mathcal{U}$ of $\mathcal{C}_n^f(M)$ such that $(\mathcal{X}, \mathcal{O}_{\mathcal{X}})$ is a reduced complex space, where $\mathcal{O}_{\mathcal{X}}$ is the restriction of $\mathcal{O}_{\mathcal{C}_n^f(M)}$ to $\mathcal{X}$, then $\mathcal{X}$ is an analytic subset of $\mathcal{U}$.
\end{enumerate}

Note a ajouter (JON).

\subsection{Some analytic subsets}

In this subsection we give the important  examples of  (closed) analytic subsets in $\mathcal{C}_{n}^{f}(M)$ which will be used in the sequel.

\begin{prop}\label{reduced}
Let $NR := \{C \in \mathcal{C}_{n}^{f}(M)\ | \  C \not= \vert C \vert \}$ be the subset of non reduced cycles. Then $NR$ is a (closed) analytic subset in $\mathcal{C}_{n}^{f}(M)$.
\end{prop}

\parag{Proof} As the empty $n$-cycle is an isolated point in $\mathcal{C}_{n}^{f}(M)$ it is enough to show that every non empty cycle in $\mathcal{C}_{n}^{f}(M)$ has an open neighborhood in which $NR$ is the zero set of a holomorphic function with values in a Banach space. Let $C_{0}$ be a non empty  cycle in $\mathcal{C}_{n}^{f}(M)$ and $C_0 = \sum_{i\in I}k_i\Gamma_i$ be the canonical expression of $C_0$. Choose, for each $i$ in $I$, an $n$-scale $E_{i} := (U_{i}, B_{i},j_{i})$ on $M$ adapted to $C_{0}$ such that the degree of $\vert C_{0}\vert$ and  $\Gamma_{i}$ in $E_{i}$ are equal to $1$. Remark that $C_{0}$ is reduced if and only if we have $k_{i} = 1$ for each $i \in I$. Let $W := \cup_{i \in I} \ j_{i}^{-1}(U_{i}\times B_{i})$ and $\mathcal{V} := \Omega(W) \cap \big(\cap_{i\in I} \ \Omega_{k_{i}}(E_{i})\big)$. Then a cycle $C \in \mathcal{V}$ is not reduced if and only if there exists at least one  $i \in I$ such that $C \cap  j_{i}^{-1}(U_{i}\times B_{i})$ is not reduced. As, for each $i$ in $I$, the natural map $\mathcal{V} \to H(\bar U_{i}, \Sym^{k_{i}}(B_{i}))$ is holomorphic, the proof is a consequence of the following claim :

\parag{Claim} 
{\em The subset of $H(\bar U, \Sym^{k}(B))$ corresponding to non reduced cycles in $U \times B$ is a closed analytic subset. It is empty for $k = 1$.}\\

For $k = 1$ the claim is clear. To prove the claim for $k \geq 2$ consider the discriminant map $\Delta_{0}\colon \Sym^{k}(\C^{p}) \rightarrow S^{k(k-1)}(\C^{p})$ defined by
 $$
 (x_{1}, \dots, x_{k}) \ \mapsto \prod_{1\leq i < j \leq k} (x_{i} - x_{j})^{2}.
 $$
 It is induced by a polynomial map $\oplus_{i=1}^{k} \ S^{i}(\C^{p}) \to S^{k(k-1)}(\C^{p})$ thanks to the standard Symmetric Function Theorem (see \cite{[BM.1]} Theorem I 4.2.7 or \cite{[e]}  Theorem 1.4.8), and so we have a holomorphic map 
 $$\Delta : H(\bar U, \Sym^{k}(B)) \to H(\bar U, S^{k(k-1)}(\C^{p}))$$
  given by $f \mapsto \Delta_{0}\circ f$. Of course, if $f \in H(\bar U, \Sym^{k}(B))$ defines a non reduced cycle in $U\times B$ we have $\Delta_{0}\circ f = 0$ in $H(\bar U, S^{k(k-1)}(\C^{p}))$.\\ 
Conversely, if $f \in H(\bar U, \Sym^{k}(B))$ defines a reduced cycle $X$ in $U\times B$, then there exists an open subset $V$ of $U$ above which $X$ is the union of mutually disjoint graphs of $k$ holomorphic functions $f_{1}, \dots, f_{k}\colon V \rightarrow B$. Thus, for all $t$ in $V$, we have
$$ 
(\Delta_{0}\circ f)(t) = \prod_{1\leq i < j \leq k} (f_{i}(t) - f_{j}(t))^{2} \neq 0
$$
and consequently  $\Delta_{0}\circ f$ is not identically zero on $U$. 
\hfill $\blacksquare$\\

Before we give more examples of analytic subsets of $\mathcal{C}_{n}^{f}(M)$ we have to prove the following lemma.

\begin{lemma}\label{rel.}
Let $U$ and $B$ be relatively compact open polydiscs in $\C^{n}$ and $\C^{p}$, and let $\pi : W \to F$ be a holomorphic map of an open neighborhood $W$ of $\bar U \times \bar B$ to a Banach space $F$. Then we have:
\begin{enumerate}[(i)]
\item
The multigraphs in $H(\bar U, \Sym^{k}(B))$\footnote{We identify each element in  $H(\bar U, \Sym^{k}(B))$ with the multigraph it determines in $U\times B$.}, which are contained in the analytic subset $\pi^{-1}(0)$  of $U\times B$, form a (closed) Banach analytic subset of $H(\bar U, \Sym^{k}(B))$.
\item
The subset $\mathcal{X}$ of $H(\bar U, \Sym^{k}(B))$ of multigraphs contained in a fiber of $\pi$ is a Banach analytic subset of $H(\bar U, \Sym^{k}(B))$. Moreover, for $k > 0$, the map $\varphi\colon \mathcal{X} \to F$ defined by sending  $X \in \mathcal{X}$  to the unique point in  $\pi(\vert X\vert)$,  is holomorphic.
\item
The multigraphs in $H(\bar U, \Sym^{k}(B))$, which have at least one irreducible component contained in the analytic subset $\pi^{-1}(0)$  of $U\times B$, form a Banach analytic subset $\mathcal{Y}$ of $H(\bar U, \Sym^{k}(B))$.

\end{enumerate}
\end{lemma}

\parag{proof} \ The only element in $H(\bar U, \Sym^{0}(B))$ is the empty multigraph which is contained in every fiber of $\pi$ and does not have any irreducible component so we will henceforth assume that $k > 0$.

\smallskip
For each $h \in [1,k]$ we have the holomorphic map
 $$
 N_{h}(\pi)\colon \Sym^{k}(W)  \longrightarrow S^{h}(F)
$$
given by the $h$-th Newton symmetric function $(z_{1}, \dots, z_{k}) \mapsto \sum_{j=1}^{k} \ \pi(z_{j})^{h}$, where $S^{h}(F)$ is the $h$-th symmetric power of $F$\footnote{That is to say the Banach space generated by the family $(x^{h})_{x\in F}$ in the Banach space of continuous homogeneous polynomials  of degree $h$ on the  dual  Banach space $F^*$ of $F$.}. Hence we get the holomorphic map
$$
\bigoplus_{h=1}^{k} N_{h}(\pi)\colon\Sym^{k}(W) \ \longrightarrow \ \bigoplus_{h=1}^{k} S^{h}(F)
$$

Each $f$ in $H(\bar U, \Sym^{k}(B))$ determines a map  $\tilde{f}$ in $H(\bar U, \Sym^{k}(W))$, defined by 
$$
\tilde{f}(t) := ((t,x_{1}), \dots, (t,x_{k}))\quad\text{if}\quad f(t) := (x_{1}, \dots, x_{k}),
$$
and it is easy to see that the map $f\mapsto \tilde{f}$ is holomorphic.

\medskip
(i) \ By compositon we obtain a holomorphic map
$$ 
\Phi\colon  H(\bar U, \Sym^{k}(B)) \ \longrightarrow \ H(\bar U, \oplus_{h=1}^{k} S^{h}(F)), \qquad f \mapsto \left(\oplus_{h=1}^{k} N_h(\pi)\right)\circ \tilde{f} 
$$ 
and the Banach analytic subset $\Phi^{-1}(0)$ consists of those multigraphs in $H(\bar U, \Sym^{k}(B))$, which are contained in the analytic subset $\pi^{-1}(0)$.

\medskip
(ii) \ 
Consider a fixed non empty open polydisc $U' \subset\subset U$ and let 
$$ 
\Psi\colon H\left(\bar U, \oplus_{h=1}^{k} S^{h}(F)\right) \ \longrightarrow \ H\left(\bar U', L(\C^{n}, \oplus_{h=1}^{k} S^{h}(F))\right) 
$$
be the holomorphic map which associates to each $g$ in $H(\bar U, \oplus_{h=1}^{k} S^{h}(F))$ the restriction to $\bar U'$ of the derivative of $g$. Then the Banach analytic subset $\Psi^{-1}(0)$ consists of all constant maps in  $H(\bar U,  \oplus_{h=1}^{k} S^{h}(F))$. 
Now fix a point $t_{0}$ in $U$ and let
$$
\chi\colon  H\left(\bar U, \oplus_{h=1}^{k} S^{h}(F)\right) \ \longrightarrow \ \bigoplus_{h=2}^{k} S^{h}(F)
$$
be the holomorphic map defined by 
$$
\chi(g) := \left(k^{h-1}N_{h}(g(t_{0})) - \left( N_{1}(g(t_{0}))\right)^{h}\right)_{h\in\{2,\ldots,k\}}
$$ 
Then the Banach analytic subset $\mathcal{Z} := \Psi^{-1}(0) \cap \chi^{-1}(0)$ of $ H(\bar U, \oplus_{h=1}^{k} S^{h}(F))$ consists of all constant maps $\bar U \rightarrow \oplus_{j=1}^{k} S^{h}(F))$ whose value is of the form $ k.a \oplus k.a^{2} \oplus \dots \oplus k.a^{k}$ for some $a \in F$. It follows that $\Phi^{-1}(\mathcal{Z})$ is exactly the subset $\mathcal{X}$ of $H(\bar U, \Sym^{k}(B))$.

\smallskip
To show that $\varphi\colon\mathcal{X}\rightarrow F$ is holomorphic, it is enough to notice that $\varphi$ is the restriction to $\mathcal{X}$ of the holomorphic map  $\frac{1}{k}.\left(ev_{1}\circ \Phi\right)$, where $\Phi$ is the map defined in (i) and  $ev_{1}\colon  H\left(\bar U, \oplus_{h=1}^{k} S^{h}(F)\right) \rightarrow F$ is defined by $ev_1(g_1,\ldots,g_k) := g_1(t_0)$.

\medskip
(iii) As above we let $\tilde{f}$ denote the element in $H(\bar U, \Sym^{k}(W))$ which is  determined by  $f$ in $H(\bar U, \Sym^{k}(B))$.  By composition with the  holomorphic map
$$
{\rm Nr}(\pi)\colon\Sym^k(W)\longrightarrow S^k(F),\qquad (z_1,\ldots,z_k)\mapsto \prod_{j=1}^k\pi(z_j)
$$
we then get a holomorphic map
$$
\Lambda\colon H(\bar U, \Sym^{k}(B))\longrightarrow H(\bar U,S^k(F)),\qquad f\mapsto {\rm Nr}(\pi)\circ\tilde{f}
$$
It follows that $\mathcal{Y} = \Lambda^{-1}(0)$ since ${\rm Nr}(\pi)\circ\tilde{f}$ is identically zero on $\bar{U}$ if and only if $\pi$ is identically zero on at least one of the irreducible components of $f$. 
\hfill $\blacksquare$\\

\begin{prop}\label{relatif}
Let $\pi\colon M \to N$ be a holomorphic map between two reduced complex spaces. Let $\mathcal{C}_{n}^{f}(\pi)$ be the subset of $\mathcal{C}_{n}^{f}(M)$ consisting of those $n$-cycles which are contained in a fiber of $\pi$\footnote{In other words  $\mathcal{C}_{n}^{f}(\pi)$ is the set of {\em $\pi$-relative $n$-cycles of finite type in $M$}.} and let $\mathcal{C}_{n}^{f}(\pi)^*$ be the (open) subset of all non empty cycles in  $\mathcal{C}_{n}^{f}(\pi)$, so
 $\mathcal{C}_{n}^{f}(\pi)^* := \mathcal{C}_{n}^{f}(\pi)\setminus\{\emptyset[n]\}$. 
\begin{enumerate}[(i)]
\item
The subset $\mathcal{C}_{n}^{f}(\pi)$ is a (closed) analytic subset of $\mathcal{C}_{n}^{f}(M)$. 
\item
The obvious map $p\colon \mathcal{C}_{n}^{f}(\pi)^*\rightarrow N$, which associates to each (non empty) $\pi$-relative cycle the unique point in $N$ whose fiber contains the cycle, is holomorphic.
\end{enumerate}
\end{prop}

\parag{Proof} First we show that the complement of $\mathcal{C}_{n}^{f}(\pi)$ is open. 
To do so we take a cycle $X_{0}$ in $\mathcal{C}_{n}^{f}(M)\setminus\mathcal{C}_{n}^{f}(\pi)$. Then $|X_0|$ contains two points $x$ and $y$ such that $\pi(x) \not= \pi(y)$ so there exist two $n$-scales $E$ and $E'$  adapted to $X_{0}$ with disjoint centers, one containing $x$ and the other $y$. It follows that the degrees $k := \deg_E(X_0)$ and $l := \deg_{E'}(X_0)$ are strictly positive and consequently  $\Omega^f_k(E)\cap\Omega^f_l(E')$ is an open neighborhood of $X_0$ in $\mathcal{C}_{n}^{f}(M)\setminus\mathcal{C}_{n}^{f}(\pi)$.
\smallskip
In order to obtain a local holomorphic equation for $\mathcal{C}_{n}^{f}(\pi)$ in $\mathcal{C}_{n}^{f}(M)$ we observe that, for every  $n$-scale  $E := (U,B,j)$ the natural map
$$ 
r_{E,k}\colon \Omega_{k}(E) \cap \mathcal{C}_n^f(M) \longrightarrow H(\bar U, \Sym^{k}(B)) 
$$
is holomorphic. Indeed, if $(X_s)_{s\in S}$ is an $f$-analytic family of $n$-cycles, parameterized by  a reduced complex space $S$, such that $E$ is adapted to $X_s$ for all $s$, then the natural map $S\rightarrow H(\bar U, \Sym^{k}(B))$ is holomorphic\footnote{This comes directly from the definition of an analytic family of cycles.}. Hence $r_{E,k}$ is holomorphic by (i) of Definition 4.2.2.  

\smallskip
Now let  $C_{0}$ be a cycle in $\mathcal{C}_{n}^{f}(\pi)$. As $\{\emptyset[n]\}$ is an open subset of $\mathcal{C}_{n}^{f}(\pi)$ we may assume that $C_0$ is not the empty $n$-cycle.
Then there exist finitely many arbitrary small $n$-scales,  $E_{1} = (U_1,B_1,j_1), \dots, E_{m} = (U_m,B_m,j_m)$,  on $M$ which satisfy the following conditions:
\begin{itemize}
\item
For each $i$, $E_i$ is adapted to $C_0$ and $k_{i} := \deg_{E_{i}}(C_{0}) > 0$.
\item
Every irreducible component of $C_0$ intersects the union of the centers of the scales, $W := \cup_{i=1}^mc(E_i)$.
\end{itemize}
It follows that  $\mathcal{U} := \left[\cap_{i \in [1,m]}\Omega_{k_{i}}(E_{i})\right] \cap\,\Omega(W)$ is an open neighborhood  of $C_0$ in $\mathcal{C}_{n}^{f}(M)$  and the holomorphic map 
$$
r\colon\mathcal{U}\ \longrightarrow\ \prod_{i \in [1,m]} H\left({\bar{U}_{i},\Sym^{k_i}(B_{i}})\right),
$$
induced by the product of the maps $\left(r_{E_{i},k_{i}}\right)_{i \in [1,m]}$, is injective.
As $C_0$ is in $\mathcal{C}_{n}^{f}(\pi)^*$ the set $\pi(|C_0|)$ is a singleton and  the $n$-scales can be chosen in such a way that $\pi$ maps their domains into a single chart on $S$. This means that there exists an open subset $S_0$ of $S$ and a holomorphic embedding $\rho\colon S_0\rightarrow\C^N$, for some $N$, such that $\pi^{-1}(S_0)$ contains the domains of the scales. By applying Lemma \ref{rel.} (ii) to the map $\rho\circ\pi$ with $F = \C^N$ we see that, for each $i$, the subset $\mathcal{X}_i$ of $H\left({\bar{U}_{i},\Sym^{k_i}(B_{i}})\right)$, consisting of all multigraphs which are contained in a single fiber of $\rho\circ\pi$, is analytic.
For each $j$ in $\{1,\ldots,m\}$ let \ $\phi_j\colon\mathcal{X}_j\rightarrow\C^N$  be the map which associates to $X$ in $\mathcal{X}_j$ the unique point in $(\rho\circ\pi)(|X|)$. By Lemma 4.2.4 (ii) this map is holomorphic and consequently the map
$$
\phi\colon\mathcal{X}_1\times\cdots\times\mathcal{X}_m\longrightarrow \left(\C^N\right)^m
$$
defined by $\phi := \phi_1\times\cdots\times\phi_m$ is also holomorphic.
Let $\Delta_m$ denote the small diagonal\footnote{ The subset of vectors $(x, x, \dots, x)$ for $x \in \mathbb{C}^N$.} of $\left(\C^N\right)^m$ and put $\mathcal{X} := \phi^{-1}(\Delta_m)$. We clearly have  $\mathcal{C}_n^f(\pi)\cap\mathcal{U} = r^{-1}(\mathcal{X})$ and $r^{-1}(\mathcal{X})$ is an analytic subset of $\mathcal{U}$. Hence we have proved that $\mathcal{C}_n^f(\pi)$ is an analytic subset of $\mathcal{C}_n^f(M)$.\\
The point  (ii) can be easily deduced from (ii) of Lemma \ref{rel.}.
$\hfill \blacksquare$\\

The next lemma is quite easy but it will be important in the sequel;

\begin{lemma}\label{Open}
Let $\pi\colon M \to N$ be a holomorphic map between two reduced complex spaces and let $V$  be an open subset in $N$. Note $\pi_V : \pi^{-1}(V) \to V$ the map induced by $\pi$. Then the obvious map
$$\xymatrix{ \mathcal{C}^f_n(\pi_V)^* \ar[r]^{i_V} \ar[d]^{\alpha_V} &  \mathcal{C}^f_n(\pi)^* \ar[d]^\alpha \\ V \ar[r] & N } $$
is an open embedding which induces an isomorphism between the corresponding weak analytic structures of  $\mathcal{C}^f_n(\pi_V)$ and of  the open set $\alpha^{-1}(V) \subset \mathcal{C}^f_n(\pi)^*$
\end{lemma} 

\parag{Proof} First we shall prove that $i_V$ is an homeomorphism onto its image which is the open set  $\alpha^{-1}(V)$ in $\mathcal{C}_n^f(\pi)^*$. This map is clearly open and bijective on $\alpha^{-1}(V)$. It is continuous because if $X \in \alpha^{-1}(V)$ and if $E$ is a scale adapted to $X$ we may always cover $\overline{c(E)}$ by finitely many  scales $(E_h)_{h \in H}$ on $\pi^{-1}(V)$ which are adapted to $X$ such that  $Y \in \cap_{h \in H} \Omega_{k_h}(E_h)\cap \mathcal{C}_n^f(\pi_V)$ implies $Y\in  \Omega_k(E)$. Moreover, if $W$ is a relatively compact open set in $M$ and  if $V_0$ is a relatively compact open set in $V$ we have the inclusion $\alpha_V^{-1}(V_0) \cap \Omega(W\cap \pi^{-1}(V_0)) \subset \Omega(W)$.\\
The fact that a  $f$-analytic family $(X_s)_{s \in S}$ of $\pi$-relative cycles in $M$ parametrized by a reduced complex space $S$  which are contained in $\alpha^{-1}(V)$ is a $f$-analytic family of 
 $\pi_V$-relative cycles in $\alpha^{-1}(V)$ is obvious because the quasi-properness of its graph of a family of  $\pi$-relative cycles  is a local property in $N$.$\hfill \blacksquare$\\

So when we shall consider a  holomorphic fiber map over an open set $V$  in $N$ (see for instance section V.2) we always may consider that it takes its values in $\alpha^{-1}(V) \subset \mathcal{C}_n^f(\pi)^*$ rather than in $\mathcal{C}_n^f(\pi_V)^*$.

\bigskip 

The next proposition gives an  analogous  result, but for  the inclusion of an analytic subset in $M$. 

\begin{prop}\label{une comp.}
Let $T $ be a closed analytic subset of the complex space $M$. 
\begin{enumerate}[(i)]
\item
The natural injection $\mathcal{C}_{n}^{f}(T)\hookrightarrow\mathcal{C}_{n}^{f}(M)$ is holomorphic and its image is an analytic subset of $\mathcal{C}_{n}^{f}(M)$.
\item
The subset $\mathcal{T}$ of $\mathcal{C}_{n}^{f}(M)$, consisting of all cycles having at least one irreducible component contained in $T$, is a (closed) analytic subset in $\mathcal{C}_{n}^{f}(M)$.
\end{enumerate}
\end{prop}

\parag{Proof} Since the singleton $\{\emptyset[n]\}$ is an open subset of $\mathcal{C}_{n}^{f}(M)$ it is enough to show  that every non empty cycle in $\mathcal{C}_{n}^{f}(M)$ has an open neighborhood in which  both subsets are analytic.
To this end  let $C_0$ be a non empty cycle in $\mathcal{C}_{n}^{f}(M)$ and consider a finite collection of  $n$-scales on $M$,  $E_{1} = (U_1,B_1,j_1), \dots, E_{m} = (U_m,B_m,j_m)$   with domains $V_1,\ldots,V_m$, satisfying the following conditions:
\begin{itemize}
\item
For each $i$, $E_i$ is adapted to $C_0$ and $k_{i} := \deg_{E_{i}}(C_{0}) > 0$.
\item
Every irreducible component of $C_0$ intersects the union of the centers of the scales, $W := \cup_{i=1}^mc(E_i)$.
\item
For each $i$, there exists a holomorphic map $g_i\colon V_i\rightarrow\C^{n_i}$ such that the set $T\cap V_i $ is equal to $ g_i^{-1}(0)$. 
\end{itemize}
Then  $\mathcal{U} := \left[\cap_{i \in [1,m]}\Omega_{k_{i}}(E_{i})\right] \cap\,\Omega(W)$ is an open neighborhood  of $C_0$ in $\mathcal{C}_{n}^{f}(M)$ and we let, for each $i$, 
$$
r_i\colon\mathcal{U}\ \longrightarrow H\left({\bar{U}_{i},\Sym^{k_i}(B_{i}})\right),
$$
denote the restriction of $r_{E_{i},k_{i}}$. Let $\Theta_i$ be the analytic subset in $H(\bar U, \Sym^k(B))$ of the multigraphs contained in $(g_i\circ j_i^{-1})^{-1}(0)$ (see Lemma \ref{rel.} point $(i)$). 
Then $r_i^{-1}(\Theta_i)$ is an analytic subset of $\mathcal{U}$, for each $i$, and,  to prove (i), we only have to observe that 
 $$\mathcal{U}\cap\mathcal{C}_n^f(T) = \bigcap\limits_{i = 1}^mr_i^{-1}(\Theta_i).$$
Let us now prove (ii).
By Lemma \ref{rel.}  point (iii), the subset $\mathcal{Y}_i$ of $H\left({\bar{U}_{i},\Sym^{k_i}(B_{i}})\right)$, consisting of those multigraphs which have at least one irreducible component contained in the zero set of $g_i\circ j_i^{-1}$, is  Banach analytic, so to finish the proof it is sufficient to show that 
$$
\mathcal{T}\cap\mathcal{U} = \bigcup_{i=1}^mr_i^{-1}(\mathcal{Y}_i).
$$ 
To this end suppose first that $C\in\mathcal{T}\cap\mathcal{U}$. Then $C$ has an irreducible component $\Gamma$ which is contained in $T$. Hence $\Gamma\cap W\neq \emptyset$ and it follows that  $\Gamma\cap c(E_i)\neq\emptyset$ for some $i$. Consequently the multigraph defined by $\Gamma$ in $U_i\times B_i$ is non-empty and contained in the zero set of $g_i\circ j_i^{-1}$. As every irreducible component of this multigraph is also an irreducible component of the multigraph defined by $C$ in $U_i\times B_i$ it follows that $C\in r_i^{-1}(\mathcal{Y}_i)$.

\smallskip
Conversely, if $C\in \cup_{i=1}^mr_i^{-1}(\mathcal{Y}_i)$ then $C\in r_i^{-1}(\mathcal{Y}_i)$ for some $i$ and at least one of the irreducible components of $C\cap V_i$ is contained in $T$. This component is contained in a unique irreducible component $\Gamma$ of $C$ and consequently $\Gamma\subseteq T$. Hence $C\in \mathcal{T}\cap\mathcal{U}$.
$\hfill \blacksquare$
\\

Note that Proposition \ref{une comp.} may not be true for an analytic family of cycles  which is not $f$-analytic as  the following example shows.

\parag{Example} Let $M := D = \{z \in \C \ / \ \vert z \vert < 1 \}$, $T := \{ 1 - \frac{1}{n}, n \in \mathbb{N}, n \geq 3\}$ and consider the family of $0$-cycles in $D$ parametrized by $D$ :
$$ X_{s} := \{ 1 - \frac{s+1}{s + m}, m \in \mathbb{N}, m \geq 3 \}\cap D  \quad {\rm for} \quad s \in D .$$
We have $X_{0} = T$ and a necessary and sufficient  condition on $s \in D$ in order that $X_{s}$ meets $T$ is that there exists $m,n \in \mathbb{N} \setminus \{0, 1,2\}$ with $ \frac{1}{n} = \frac{s+1}{s+m}$. This gives that $X_{s}$ meets $T$ if and only iff $s =\frac{p}{q}$ with $p \in \mathbb{Z}, \ q \in \mathbb{N}\setminus \{0, 1\}$ and $\vert \frac{p}{q} \vert < 1$. This is a dense set in $]-1,+1[$ !  
$\hfill \square$

\begin{prop}\label{universal.graph}
The graph of the tautological $f$-analytic family of $n$-cycles of finite type in $M$
$$
\mathcal{C}_n^f(M)\sharp M := \{(C,x)\in \mathcal{C}_n^f(M)\times M\ /\ x\in |C|\}
$$
is an analytic subset of $\mathcal{C}_n^f(M)\times M$\footnote{We leave  to the reader the definition of  a weak analytic structure on $\mathcal{C}_n^f(M)\times M$.}.
\end{prop}

\parag{Proof} Obviously $\mathcal{C}_n^f(M)\sharp M$ is a closed subset of $\mathcal{C}_n^f(M)\times M$. 

\smallskip
Suppose $(C,x)\in\mathcal{C}_n^f(M)\sharp M$ and let $E = (U,B,j)$ be an $n$-scale on $M$ adapted to $C$ such that $x\in c(E)$. Put $k := \deg_E(C)$. Then $\Omega^f_k(E)\times c(E)$ is an open neighborhood of $(C,x)$ in $\mathcal{C}_n^f(M)\times M$ and the map 
$$
r_{E,k}\times j\colon\Omega^f_k(E)\times c(E)\rightarrow H(\bar{U},\Sym^k(B))\times (U\times B)
$$ 
is holomorphic. By composing this map with the holomorphic map 
$$
H(\bar{U},\Sym^k(B))\times (U\times B)\rightarrow\Sym^k(B)\times B,\qquad (f,t,x)\mapsto( f(t),x)
$$
we obtain a holomorphic map $\Psi\colon\Omega_k(E)\times c(E)\rightarrow\Sym^k(B)\times B$.  As the subset
 $$\Sym^k(B)\sharp B := \{(\xi,x)\in\Sym^k(B)\times B\ /\ x\in\xi\}$$
  is  analytic in $\Sym^k(B)\times B$ it follows that 
  $$\Psi^{-1}(\Sym^k(B)\sharp B) = (\Omega_k(E)\times c(E))\cap(\mathcal{C}_n^f(M)\sharp M)$$
   is an analytic subset of $\mathcal{C}_n^f(M) \cap (\Omega_k(E)\times c(E))$.
\hfill{$\blacksquare$}\\

\subsection{Complements}

We begin this subsection by showing that the natural inclusion $\mathcal{C}_n(M)\hookrightarrow\mathcal{C}_n^f(M)$ is an open (holomorphic) embedding of the reduced complex space $\mathcal{C}_n(M)$ of compact analytic $n$-cycles in $M$ (see \cite{[BM.2]} Chapter V).

\begin{prop}\label{compact cycles}
Let $M$ be a complex space. Then $\mathcal{C}_n(M)$ is an open subset of $\mathcal{C}_n^f(M)$ and the induced topology coincides with the natural topology on $\mathcal{C}_n(M)$. Moreover its structure sheaf is the restriction of the sheaf of holomorphic functions on $\mathcal{C}_n^f(M)$ for the weak analytic structure defined in  \ref{appl.hol. 2} .
\end{prop}

\parag{Proof} Recall that the natural topology on $\mathcal{C}_n(M)$ is generated by all sets  of the form $\Omega_k(E)\cap\mathcal{C}_n(M)$, where $E$ is an $n$-scale on $M$ and $k$ is a natural number, and all sets of the form 
$$ 
\Omega_c(W) := \{ X \in \mathcal{C}_n(M) \ / \   |X| \subseteq W \} 
$$
where $W$ is a relatively compact open subset of $M$.\\
 Now fix a compact $n$-cycle $X_0$ in $M$, a relatively compact open subset $W$ of $M$ containing $X_0$ and a finite collection, $E_1, \dots, E_p$, of $n$-scales adapted to $X_0$. Put $k_j := \deg_{E_j}(X_0)$ and consider the open neighborhood $\mathcal{U} := \Omega_c(W) \cap \big(\cap_{j= 1}^p \Omega_{k_j}(E_j)\big) $ of $X_0$ in $\mathcal{C}_n(M)$. To complete the proof it is enough to 
construct an open neighborhood $\mathcal{V}$ of $X_0$ in $\mathcal{C}_n^f(M)$ such that $\mathcal{V} \subseteq \mathcal{U}$.\\
To do so we cover the compact set $\partial W$ by the centers of $n$-scales $F_1, \dots, F_q$ which are adapted to $X_0$ and satisfy the condition  $\deg_{F_h}(X_0) = 0$ for all $h$ in $\{1,\ldots, q\}$ and put 
$$
\mathcal{V} := \Omega(W)  \cap \left(\,\bigcap_{j= 1}^p \Omega_{k_j}(E_j)\right) \cap \left(\,\bigcap_{h=1}^q \Omega_0(F_h)\right).
$$
Now, if $X\in\mathcal{V}$, then $ \deg_{F_h}(X) = 0$ for all $h$ and it follows that $|X|\cap \partial W = \emptyset$. So every irreducible component of $X$ meets $W$ but not $\partial W$ and consequently $X$ is a compact cycle contained in $W$. Hence $X$ is in $\Omega_c(W)$.\\
The last assertion is obvious.\hfill{$\blacksquare$}\\

In general, it is difficult to check whether a given function on a given subset of $\mathcal{C}_n^f(M)$ is holomorphic with our definition (see \ref{appl.hol. 2} point $(i)$ ). So it is interesting to have a simple method to build, at least locally, holomorphic functions. This is the content of our next proposition.

\begin{prop}\label{hol. functions} Let $E := (U, B, j)$ be an $n$-scale on $M$, let $g : U \times B \to \mathbb{C}$ be a holomorphic function and  let $\varphi \in \mathscr{C}^\infty_c(U)^{(n, n)}$. Then the function defined by
\begin{equation*} G(X) := \int_X j^*(g.pr^*(\varphi))  = \int_U Trace_{X/U}(g).\varphi  \tag{@}
 \end{equation*}
  for $X \in \Omega^f_k(E)$ is holomorphic on the open set $ \Omega^f_k(E) $.
\end{prop}

In fact we shall prove that the map
$$\tilde{G} : H(\bar U, \Sym^k(B)) \times U \to \mathbb{C} , \quad (X, t) \mapsto Trace_{X/U}(g)(t) $$
is holomorphic. Then it is easy to conclude that for any holomorphic map 
$$h : S \longrightarrow \Omega^f_k(E) $$
  where $S$ is a Banach analytic set, the map $G$ is holomorphic using the holomorphy of an integral which depends holomorphically of a parameter\footnote{ see \cite{[BM.2]} Chapter V paragraph 2.4. Here the fact that the map $\tilde{G}$ is globally induced on 
$H(\bar U, \Sym^k(B))$ allows to apply the positive result in this case.}.
\parag{Claim}  The map $G$ is induced  on $H(\bar U, \Sym^k(B))\times U$ by a holomorphic map on the ambient Banach open set 
$ H(\bar U, \mathcal{U}_R)\times U \subset H(\bar U, \oplus_{h=1}^k S^h(\mathbb{C}^p)\times \mathbb{C}^n$
$$\mathcal{G} :  H(\bar U, \mathcal{U}_R)\times U \to \mathbb{C}.$$
 Here we assume that $B = B_R$ is the polydisc with center $0$ and radius $R$ in $\mathbb{C}^p$ and $\mathcal{U}_R$ is the open set in $\oplus_{h=1}^k S^h(\mathbb{C}^p)$ defined in Proposition V. 3.1.2 in \cite{[BM.2]}.\\
 To buid up the function $\mathcal{G}$ on $H(\bar U, \mathcal{U}_R)\times U$ let $N_{\alpha}(X)$ for $X \in H(\bar U, \oplus_{h=1}^k S^h(\mathbb{C}^p))$ be  the $\alpha$-component of the $\vert \alpha\vert$-th Newton function of $X$, for $\alpha \in \mathbb{N}^p$. Now write the Taylor expansion of $g$ at the point $(t, 0)$  for any $t$ fixed in $U$:
 $$ g(t, x) := \sum_{\alpha \in \mathbb{N}^p} g_{\alpha}(t).x^\alpha,\quad  {\rm for } \quad x \in B_R.$$
 Then define
 $$ \mathcal{G}(X, t) := \sum_{\alpha \in \mathbb{N}^p}  g_\alpha(t).N_{\alpha}(X), \quad {\rm for} \quad (X, t) \in H(\bar U, \mathcal{U}_R)\times U.$$
 Then the definition of the open set $\mathcal{U}_R$ (see Proposition 3.1.2 in \cite{[BM.2]}) gives the convergence of this series and the holomorphy of the function $\mathcal{G}$. It is easy to see that for $(X, t) \in H(\bar U, \Sym^k(B_R))$ we have $\mathcal{G}(X,t) = G(X, t)$ using the relation $(@)$.This completes the proof.$\hfill \blacksquare$\\
 
 \parag{Remark} It seems "a priori" that we may obtain locally  more holomorphic functions on $\mathcal{C}_n^f(M)$, using isotropy and the classifying spaces $\Sigma_{U, U'}(k)$  for isotropic morphism, by integration of $d''$-closed $\mathscr{C}^\infty(U' \times B)^{(n, n)}$-differential forms  with support in $K \times B$ for $K$ a compact set in $U'$. In fact, results in Chapter V of \cite{[BM.2]}  implies that this does not produce more local holomorphic functions than finite sums of functions obtained by using the previous proposition with enough adapted scales corresponding to several different linear projections of $U'\times B$ to $U$, sufficiently  near  the natural  (vertical) one.
 
 \parag{Example} Let $X_0$ be in $\mathcal{C}_n^f(M)$ and let $E := (U, B, j))$ be a scale on $M$ adapted to $X_0$ with $\deg_E(X_0) = k$. Let $t_0$ be a point in $U$. Then the  map
 $$ f : \Omega^f_E(k) \longrightarrow \Sym^k(B), \quad X \mapsto j_*(X) \cap (\{t_0\}\times B) $$
 is holomorphic.\\
  We leave to the reader the proof of the fact that this kind of holomorphic functions is obtained as uniform limits of holomorphic functions on $H(\bar U, \mathcal{U}_R)$ using an approximation of the Dirac mass at $t_0$ on $U$ by elements of $\mathscr{C}^\infty_c(U)^{(n, n)}$.$\hfill \square$
  
  \parag{Remark} Let $f$ be a  holomorphic map defined on a open set $\mathcal{U}$  in $\mathcal{C}_n^{\rm loc}(M)$. For any $f$-analytic family of $n$-cycles $(X_s)_{s \in S}$  with $X_s$ in $\mathcal{U}$ for each $s \in S$,  so a holomorphic map $\varphi : S \to \mathcal{U} \cap \mathcal{C}_n^f(M)$, $f$ induces, by composition with $\varphi$, a holomorphic map on $S$ and then $f$ is  a holomorphic map on $\mathcal{U} \cap \mathcal{C}_n^f(M)$.\\
  For instance, assuming that $M$ is a connected  $m$-dimensionnel complex manifold,  if $Y$ is a compact $q$-cycle in $M$, the subset  $\mathcal{U}(Y)$ in $\mathcal{C}_n^{\rm loc}(M)$ of $n$-cycles which cut properly $Y$ is open in $\mathcal{C}_n^{\rm loc}(M)$ and we have a holomorphic intersection map (see \cite{[BM.2]} chapter VII)  from $\mathcal{U}(Y)$ to $\mathcal{C}_d(M)$ where $d$ satisfies
  $$ m- d = m-n + m-q, \quad {\rm so} \quad d = n+q - m.$$
  So this defines a weakly holomorphic map from \, $\mathcal{U}(Y) \cap \mathcal{C}_n^f(M)$\,  to the reduced complex space $\mathcal{C}_d(M)$.

\newpage

\section{The Analytic Extension Theorem}

Let $S$ be a reduced complex space and $f$ be a function on $S$. Then, by definition, $f$ is analytic at a point $s_{0}$ when there exists an open neighborhood $S_{0}$ of $s_{0}$ in $S$ such that $f$ is holomorphic on $S_{0}$. But a family of $n$-cycles $(X_{s})_{s \in S}$ in a complex space $M$ which is analytic at a point $s_{0}$ is not necessarily analytic in a neighborhood of $s_0$ (even if the cycles are compact) as the following example shows. \\
On the other hand, if the family  $(X_{s})_{s \in S}$ is properly analytic at $s_0$, then it is analytic in an open neighborhood of $s_0$ (See Theoreme V.1.0.3 in \cite{[BM.2]}).  In the sequel we shall prove that this is also true if the family  $(X_{s})_{s \in S}$ is $f$-analytic at $s_0$. This result is part of the Analytic Extension Theorem below.

\parag{Example}
Let $D$ denote the open unit disk in $\C$ and put $X_s :=  (\{0\} + \{1- |s|\})\cap D$ for every $s$ in $\C$. Then the family $(X_s)_{s\in\C}$ of $0$-cycles in $D$ is analytic at $0$, but it is not analytic in any neighborhood of $0$. Note that in any scale on $D$ adapted to $X_0$ there exists an open neighborhood $V$ of $s = 0$ in $\mathbb{C}$ such that in the scale chosen the family coincides on $V$ with the constant family equal to $X_0$.$\hfill \square$

\begin{thm}\label{cycles}[{\bf Analytic Extension Theorem}]\\
Let $M$ be a complex space and  $n$ be a natural number. Consider an $f$-continuous family  $(X_{s})_{s\in S}$ of $n$-cycles in  $M$ parametrized by a reduced complex space  $S$. Fix a point $s_{0}$ in $S$ and assume that there exists an open subset $M'$ of $M$ meeting every irreducible component of $\vert X_{s_{0}}\vert$  and such that the family $(X_{s}\cap M')_{s \in S}$ is analytic at  $s_{0}$. Then there exists an open neighborhood $S_0$ of $s_0$ in $S$ such that  the family  $(X_{s})_{s\in S_0}$ is  $f$-analytic in $M$ at each point in $S_0$.
\end{thm}

Let us make explicit the situation of the previous theorem in terms of classifying maps : we have a continuous map $ \varphi : S \to \mathcal{C}_{n}^{f}(M) $ such that the composed map $r\circ \varphi$ is holomorphic at $s_{0}$, where $r : \mathcal{C}_{n}^{f}(M)  \to \mathcal{C}_{n}^{loc}(M') $ is obtained by restriction. Then the statement is that, assuming that the open set $M'$ meets each irreducible component of  $\vert X_{s_{0}}\vert$, the map  $\varphi$ is holomorphic on an open neighborhood of $s_0$.\\
Remark that the map $r$ is holomorphic\footnote{This means that for any holomorphic map $\psi : T \to \mathcal{C}_{n}^{f}(M) $ of a reduced complex space $T$ the composed map $r\circ\psi$ is holomorphic.} so that the holomorphy at $s_{0}$ of $r\circ \varphi$ is a necessary condition for the holomorphy of $\varphi$ at $s_{0}$. The theorem says that this condition is not only sufficient but  also open on $S$.\\

A key point in the proof of the previous theorem is given by the following analytic continuation result.

\begin{prop}\label{fonction}
Let $S$ be a reduced complex space and let  $U_{1} \subset U_{2}$ be two open  polydiscs in  $\C^{n}$ with $U_1 \not= \emptyset$. Let  $f: S \times U_{2} \to \C$ be a continuous function, holomorphic on  $\{s\}\times U_{2}$ for each $s \in S$ and assume also that the restriction of $f$ to   $S \times U_{1}$ is holomorphic. Then $f$ is holomorphic on $S \times U_{2}$.
\end{prop}

\parag{Proof} Consider first the case where $S$ is smooth. As the problem is local on $S$ it is enough to treat the case where $S$ is an open set in some $\C^{m}$. Fix then a relatively compact open polydisc $P$ in $S$. The function $f$ defines a map $ F : U_{2} \to \mathscr{C}^{0}(\bar P, \C)$, where $ \mathscr{C}^{0}(\bar P, \C)$ is the Banach space of continuous functions on $\bar P$, via the formula $F(t)[s] = f(s,t)$ for $t \in U_{2}$ et $s \in \bar P$. The map $F$ is holomorphic: this is an easy consequence of Cauchy's formula on a polydisc $U \subset\subset U_{2} $ with fixed $s \in S$ which computes the partial derivatives in $t := (t_{1}, \dots, t_{n})$:
$$
 \frac{\partial f}{\partial t_{i}}(s, t) = \frac{1}{(2i\pi)^{n}}\int_{\partial\partial U} f(s,\tau).\frac{d\tau_{1}\wedge \dots \wedge d\tau_{n}}{(\tau_{1} - t_{1})\dots (\tau_{i}-t_{i})^{2}\dots (\tau_{n}-t_{n})} \quad \forall t \in U \quad \forall i \in [1,n].
$$
This shows that $F$ is $\C-$differentiable and its differential at the point  $t \in U$ is given by  $h \mapsto \sum_{i=1}^{n} F_{i}(t).h_{i}, \ h \in \C^{n}$, where $F_{i} $ is the map associated to the function
$$ 
(s,t) \mapsto  \frac{\partial f}{\partial t_{i}}(s, t) \quad i \in [1,n]
$$
which is holomorphic for any fixed $s \in S$ thanks to the Cauchy formula above as $t \mapsto f(s,t)$ is holomorphic for each $s \in S$.\\
Let  $H(\bar P,\C)$ be the (closed) subspace of  $\mathscr{C}^{0}(\bar P, \C)$ of continuous functions which are holomorphic on $P$. Our assumption implies that the restriction of $F$ to $U_{1}$ takes its values in this subspace.\\
 Let us show that for each point $t \in U_{2}, F(t)$ is still in  $H(\bar P,\C)$.  Assume this is not true. Then there exists $t_{0} \in U_{2}$ with $F(t_{0})\not\in H(\bar P,\C)$, and so, by  the Hahn-Banach theorem, there exists a continuous linear form  $\lambda$  on $\mathscr{C}^{0}(\bar P, \C)$, vanishing on  $H(\bar P,\C)$ and such that $\lambda(F(t_{0})) \not= 0$. But the  function  $t \mapsto \lambda(F(t))$ is holomorphic on $U_{2}$ and vanishes on $U_{1}$. So it vanishes identically, contradicting the fact that $\lambda(F(t_{0})) \not= 0$. So $F$ is a holomorphic map with values in  $H(\bar P,\C)$ and  $f$ is holomorphic on $S\times U_{2}$ when $S$ is smooth.\\
The case where $S$ is a weakly normal complex space is then an immediate consequence of the smooth case, as the continuity of $f$ on $S \times U_{2}$ and the holomorphy of $f$ on $S_{reg}\times U_{2}$, obtained above, are enough to conclude.\\
When $S$ is a general reduced complex space the function $f$ is then a continuous meromorphic function on $S\times U_{2}$ which is holomorphic on $S \times U_{1}$. So the closed analytic subset $Y \subset S \times U_{2}$ of points at which $f$ is not holomorphic has empty interior in each $\{s\}\times U_{2}$. So the criterium  3.1.7 of analytic continuation of chapter IV  in 
\cite{[BM.1]} allows to conclude.
$\hfill \blacksquare$

\parag{Remarks} \begin{enumerate}
\item  It is an easy exercise to weaken the hypothesis of the previous proposition replacing the continuity of $f$ by the hypothesis ``$f$ is measurable and locally bounded  on $S \times U_{2}$''. In the first step the Banach space $\mathscr{C}^{0}(\bar P, \C)$ is replaced by the Banach space of bounded measurable functions on $\bar P$ and  in the second step  $S$ is assumed to be normal. Then the final conclusion is obtained following the same lines as above.
\item It is not difficult to extend the proposition above to the case where $S$ is a Banach open set. But this is not true, in general, for  singular  Banach analytic sets.\\
 The reader will find  in \cite{[BM.2]} p.33 a counterexample  where $S$ is a (infinite dimensional)  Banach analytic set  which has only one singular point (so it is a Banach manifold outside this point).$\hfill \square$
\end{enumerate}

\parag{Proof of  Theorem \ref{cycles}} Consider the graph $G \subset S \times M$ of the $f$-continuous family  $(X_{s})_{s\in S}$ and let $A$ be the open subset of $G$ consisting of all points  $(\sigma,\zeta)$ which satisfy the following condition:
\begin{itemize}
\item  
There exist an open neighborhood $S_{\sigma}$ of $\sigma$ in $S$ and an open neighborhood $M_{\zeta}$  of  $\zeta$ in $M$ such that the family $(X_{s} \cap M_{\zeta})_{s \in S_{\sigma}}$ is analytic at each point of $S_\sigma$.
\end{itemize} 

 Remark that our  assumption implies that the open set  $A$ meets every irreducible component of  $\{s_{0}\}\times \vert X_{s_{0}}\vert$. This point is consequence of the fact that isotropy at $s_0$ in an adapted scale implies analyticity of the family in the isotropy domain of the scale at any point in an open neighborhood of $s_0$ (see \cite{[BM.2]} Chapter V).\\
 We  prove now that $A$ contains $\vert X_{s_0}\vert$ which implies the analyticity of the family at the point $s_0$.\\
Assume now that there exists a smooth point of $\vert X_{s_{0}}\vert$  in the boundary of the set  $A \cap (\{s_{0}\}\times \vert X_{s_{0}}\vert)$.
Consider now such a point $(s_{0}, z_{0})$ and choose an $n$-scale $E := (U,B,j)$ which is adapted to $X_{s_{0}}$ and satisfies the following conditions:
 \begin{align*}
 & \deg_{E}(\vert X_{s_{0}}\vert) = 1,\quad    z_{0} \in j^{-1}(U\times B),\quad  j(z_{0}) := (t_{0}, 0),\\
 \\
 & j_{*}\left(X_{s_{0}}\cap j^{-1}(U\times B)\right) \ = \ k.(U \times \{0\}).
 \end{align*}
 Then we have a continuous classifying map $f : S_{1}\times U \to \Sym^{k}(B)$ where $S_{1}$ is an open neighborhood of $s_{0}$ in $S$. The map $f$ is holomorphic for each fixed $s \in S_{1}$. As the point  $(s_{0}, z_{0})$ is in the boundary of the open subset $A \cap (\{s_{0}\}\times \vert X_{s_{0}}\vert)$ of $\{s_{0}\}\times \vert X_{s_{0}}\vert$, there exists a (non empty) polydisc $U_{1} \subset U$ such that the restriction of $f$ to $S_{1}\times U_{1}$ is holomorphic near $s_{0}$. So, by shrinking $S_{1}$ if necessary, we can assume that $f$ is isotropic on $S_{1}\times U_{1}$. Applying  Proposition \ref{fonction} to each scalar component of $f$, we conclude that $f$ is isotropic on $S_{1}\times U$ (see \cite{[BM.2]} Chapiter V section 5). This contradicts the fact that $(s_{0}, z_{0})$ is in the boundary of $A \cap (\{s_{0}\}\times \vert X_{s_{0}}\vert)$.\\
 If the boundary of $A \cap (\{s_{0}\}\times \vert X_{s_{0}}\vert)$ is contained in the singular set of $\vert X_{s_{0}}\vert$, then we can apply the criterium (\cite{[BM.1]}  Chapter IV  Criterium  3.1.9)  to obtain directly that $A$ contains $\vert X_{s_{0}}\vert$ and the family is analytic at $s_0$.\\
 
 We shall prove now that there exists an open neighborhood $S_0$ of $s_0$ in $S$ such that the family is analytic at each point of $S_0$:\\
 Let $M''$ be a relatively compact open subset of $M'$ which intersect every irreducible component of $X_{s_0}$.  As $(X_s)_{s\in S}$ is $f$-continuous, $s_0$ admits an open neighborhood $S'$ in $S$ such that $M''$ meets every irreducible component of $X_s$ for all $s$ in $S'$. Then, by Theorem V.1.0.3 in \cite{[BM.2]}, there exists an open neighborhood $S_0 $ of $s_0$ in $S'$ such that $(X_s\cap M'')_{s \in S_0}$ is an analytic family of $n$-cycles in $M''$. Hence $(X_s)_{s\in S_0}$ is an $f$-analytic family of $n$-cycles in $M$ due to the first part of the proof applied to each point $s \in S_0$. $\hfill \blacksquare$\\
 
 To conclude this section, let us give an example of "weak analytic map" between $\mathcal{C}_n^f(M) $ and $\mathcal{C}_n^f(N)$.
 
 \begin{thm}\label{direct image}
 Let $q : M \to N$ be a proper holomorphic map between complex spaces. Then the map
 $$ q_* : \mathcal{C}_n^f(M) \to \mathcal{C}_n^f(N) $$
 given by the direct image of finite type cycles is holomorphic in the sense that for any holomorphic map $\varphi : S \to \mathcal{C}_n^f(M) $ where $S$ is a reduced complex space, the composed map $ q_*\circ \varphi $ is holomorphic.
 \end{thm}
 
 \parag{Proof} The only points which are not already contained in Direct Image Theorem  IV 3.5.3  in \cite{[BM.1]} is the fact that the direct image of a finite type cycle is a finite type cycle, and the fact that the direct image of a $f$-continuous family of cycles is f-continuous. The first point is trivial. To prove the second point, remark that if a compact set $K$ in $M$ meets every irreducible component of a cycle $X$ in $M$ then the compact set $q(K)$ meets each irreducible component of the cycle $q_*(X)$. $\hfill \blacksquare$\\

 \section{The semi-proper direct image theorem.}

The aim of this paragraph is to extend the Direct  Image Theorem  for semi-proper maps (see Theorem \ref{Kuhl-Banach})  to the case where the target space is $\mathcal{C}_n^f(M)$ for a given complex space $M$. This is not obvious because, as we have  already remarked above, the weak analytic structure on $\mathcal{C}_n^f(M)$  which is defined in the previous section is not, even locally, a structure of Banach analytic set. We give here an improvement of Theorem 5.0.5 in \cite{[B.13]}. It is also an opportunity to give a more elaborated  proof of this delicate result.\\

First we recall that  an analytic subset $\mathcal{X}$ of an open subset  $\mathcal{U}$ of $\mathcal{C}_{n}^{f}(M)$ is called a {\em reduced complex subspace} of $\mathcal{U}$ if, endowed with the sheaf of holomorphic functions of $\mathcal{C}_{n}^{f}(M)$, it becomes a reduced complex space (see definition \ref{appl.hol. 2} $(viii)$).

\begin{thm}\label{semi-proper direct image bis}
Let $M$  and $S$ be reduced complex spaces and $n$ be a natural number. Assume that we have a holomorphic map $\varphi\colon S \rightarrow \mathcal{U} $  which is semi-proper, 
 where $\mathcal{U}$ is an open set in $\mathcal{C}_{n}^{f}(M)$. Then  $\varphi(S)$ is a reduced complex subspace of $ \mathcal{U} \subset \mathcal{C}_{n}^{f}(M)$
\end{thm}

The main tools for the proof of Theorem \ref{semi-proper direct image bis} are Theorem \ref{Kuhl-Banach} and Theorem \ref{cycles} but  the following topological result is also needed.

\begin{lemma}\label{semi-proper direct image bis.0}
Let $\varphi\colon S\rightarrow T$ and $f\colon T\rightarrow Z$ be continuous maps between Hausdorff spaces. Suppose moreover that $S$ is locally compact, $\varphi$ is semi-proper and $f$ is injective. Then, for every $t_0$ in $\varphi(S)$, there exists an open neighborhood $T_0$ of $t_0$ in $T$ and an open subset $U$ of $Z$ such that $T_0 \subset f^{-1}(U)$ and such that the map
$$
\varphi^{-1}(T_0)\ \longrightarrow\ U,\qquad s\mapsto f(\varphi(s))
$$
is semi-proper and such that the map $\varphi(S)\cap T_0\rightarrow f(\varphi(S)\cap T_0)$, induced by $f$, is a homeomorphism.
\end{lemma}

\parag{Proof} Let $t_0$ be a point in $\varphi(S)$. As $\varphi$ is semi-proper and  $\varphi(S)$ is locally compact  $t_0$ admits a relatively compact open neighborhood $V_0$ in $\varphi(S)$. Hence the map  $V_0\rightarrow f(V_0)$, induced by $f$, is a homeomorphism and it follows that the map 
$$
\varphi^{-1}(V_0)\longrightarrow f(V_0),\qquad s\mapsto f(\varphi(s))
$$
is semi-proper and surjective. In particular $f(V_0)$ is locally compact and consequently locally closed in $Z$, so there exists an open subset $U$ of $Z$ such that $f(V_0)$ is a closed subset of $U$. Then for any open  neighborhood $T_0$ of $t_0$ in $T$, which satisfies $V_0 = \varphi(S)\cap T_0$, the map
$$
\varphi^{-1}(T_0)\ \longrightarrow\ U,\qquad s\mapsto f(\varphi(s))
$$
is semi-proper.  
\hfill{$\blacksquare$}

\parag{Proof of Theorem \ref{semi-proper direct image bis}} 
First recall that, for an $n$-scale $E = (U,B,j)$ on $M$, a relatively compact open polydisc $U'$ in $U$ and an integer $k\geq 0$, there exists a Banach analytic subset $\Sigma_{U,U'}(k)$ of a Banach open set and a holomorphic homeomorphism 
$$
\rho_{U,U'}\colon\Sigma_{U,U'}(k)\ \longrightarrow\ H(\bar{U},\Sym^k(B))
$$
having the following properties\footnote{These sets were introduced in \cite{[B.75]}   and are essential for the construction of the reduced complex space of compact analytic cycles in a given complex space. For a detailed discussion of these sets see Chapter V section 7.2 in \cite{[BM.2]}.}:
\begin{itemize}
\item[(a)]
For every holomorphic map from a reduced complex space into $\Sigma_{U,U'}(k)$, the corresponding family of multigraphs in $U\times B$ is an analytic family of cycles in $U'\times B$. 
\item[(b)]
If $(X_s)_{s\in S}$ is an analytic family of $n$-cycles in $M$ such that $X_s\in\Omega_k(E)$ for all $s$ in $S$, then the natural map $S\rightarrow\Sigma_{U,U'}(k)$, obtained by composing the classifying map $\Omega_k(E)\rightarrow H(\bar{U},\Sym^k(B))$ with the inverse of $\rho_{U,U'}$, is holomorphic.
\end{itemize}
Now fix $C_0$ in $\varphi(S)$ and let us show that $C_0$ admits an open neighborhood $\mathcal{V}$ in $ \mathcal{U}$ such that $\varphi(S)\cap\mathcal{V}$ is a reduced complex subspace of $\mathcal{V}$. The case $C_0 = \emptyset[n]$ being trivial we suppose that $C_0$ is not the empty $n$-cycle.

First we choose (as we did in the proof of Proposition \ref{une comp.}) finitely many $n$-scales, $E_1 = (U_1,B_1,j_1),\ldots,E_m = (U_m,B_m,j_m)$, which are all adapted to $C_0$ and let  $k_i  > 0$ be the degree of $C_0$ in $E_i$, for each $i$.  Assume also that every irreducible component of $C_0$ intersects the relatively compact open subset $\cup_{i=1}^mj_i^{-1}(U_i\times B_i)$ of $M$. Then we choose for each $i$ a relatively compact open polydisc $U_i'$ of $U_i$ in such a way that every irreducible component of $C_0$ meets $W := \bigcup\limits_{i=1}^mj_i^{-1}(U_i'\times B_i)$.

On the open subset
$$
\mathcal{W} := \Omega(W)\cap\left(\,\bigcap_{i=1}^m\Omega_{k_i}(E_i)\right)
$$
of $\mathcal{C}_n^f(M)$ the injective map $f\colon\mathcal{W}\ \rightarrow\  \prod_{i=1}^m\Sigma_{U_i,U_i'}(k_i)$, induced by the product of the classifying maps,  is holomorphic according to property (b) above. Then, by Lemma \ref{semi-proper direct image bis.0}, there exists an open neighborhood $\mathcal{V}$ of $C_0$ in $\mathcal{W}$ and an open subset $\mathcal{U}$ of $\prod_{i=1}^m\Sigma_{U_i,U_i'}(k_i)$ such that the map 
$$
\varphi^{-1}(\mathcal{V})\longrightarrow\mathcal{U},\qquad s\mapsto f(\varphi(s))
$$
is semi-proper and such that the map $\varphi(S)\cap\mathcal{V}\rightarrow f(\varphi(S)\cap\mathcal{V})$ induced by $f$ is a homeomorphism. As $\prod_{i=1}^m\Sigma_{U_i,U_i'}(k_i)$ is a Banach analytic subset of a Banach open set there exists an open subset $\tilde{\mathcal{U}}$ of the ambient Banach space such that $\mathcal{U} = \tilde{\mathcal{U}}\cap\prod_{i=1}^m\Sigma_{U_i,U_i'}(k_i)$. Then the restriction of $f$ to $\mathcal{V}$ composed with the natural inclusion $\mathcal{U}\hookrightarrow\tilde{\mathcal{U}}$ is an injective holomorphic map
$\tilde{f}\colon\mathcal{V}\ \longrightarrow\ \tilde{\mathcal{U}}$  and the map 
$$
\varphi^{-1}(\mathcal{V})\longrightarrow \tilde{\mathcal{U}},\qquad s\mapsto\tilde{f}(\varphi(s))
$$
is semi-proper. Hence $\tilde{f}(\varphi(S)\cap\mathcal{V})$ is a reduced complex subspace of $\tilde{\mathcal{U}}$, due to Theorem \ref{Kuhl-Banach}.  

To sum up the situation let us consider the commutative diagram
\begin{equation*}\tag{$@$}
\xymatrix{
\varphi(S)\cap\mathcal{V} \ar[d]_{\nu} \ar[r]^{g\ \ } & \tilde{f}(\varphi(S)\cap\mathcal{V}) \ar[d] \\
\mathcal{V} \ar[r]^{\tilde{f}} & \tilde{\mathcal{U}} } 
\end{equation*}
where $g\colon\varphi(S)\cap\mathcal{V}\rightarrow\tilde{f}(\varphi(S)\cap\mathcal{V})$ is the homeomorphism induced by $\tilde{f}$ and the vertical arrows are the natural inclusions. As $\varphi(S)\cap\mathcal{V}$ is a subset of $\mathcal{C}_n^f(M)$ the homeomorphism $g^{-1}$ defines an $f$-continuous family of $n$-cycles in $M$ and, by property (a) of the Banach analytic sets $\Sigma_{U_i,U_i'}(k_i)$, this family is analytic on $W$. But $W$ meets every irreducible component of every $C$ in $\varphi(S)\cap\mathcal{V}$ so we can apply Theorem \ref{cycles} and conclude that the family is analytic. This means that $\nu\circ g^{-1}\colon\tilde{f}(\varphi(S)\cap\mathcal{V})\rightarrow\mathcal{V}$ is a holomorphic map. We then deduce from ($@$) that $\varphi(S)\cap\mathcal{V}$, endowed with the sheaf  induced by the sheaf of holomorphic functions on $\mathcal{C}_{n}^{f}(M)$, is a reduced complex space and $g$ is an isomorphism between $\varphi(S) \cap \mathcal{V}$ and the reduced complex subspace $\tilde{f}(\varphi(S)\cap\mathcal{V})$ of $\tilde{\mathcal{U}}$.
\hfill{$\blacksquare$}
  
\newpage
 
 \section{Meromorphic maps to $\mathcal{C}_{n}^{f}(M)$.}
 
 In this section we introduce, for reduced complex spaces $M$ and $N$, the notion of a {\em meromorphic map} from $N$ to $\mathcal{C}_{n}^{f}(M)$. It is a generalization of   the notion of a meromorphic map between reduced complex spaces.  But to define the {\em graph} of a meromorphic map in this context we have to use the  Semi-Proper Direct Image Theorem proved in the previous section.
 
 \begin{defn}\label{mero}
 Fix a complex space $M$ and an integer $n$. Let $N$ be a reduced complex space and let  $\Sigma \subset N$ be  a nowhere dense analytic subset in $N$. We shall say that a holomorphic map $\varphi\colon N \setminus \Sigma \to \mathcal{C}_{n}^{f}(M)$ is {\bf\em meromorphic along $\Sigma$} (or more simply that $\varphi\colon N\dashrightarrow \mathcal{C}_{n}^{f}(M)$ is {\bf\em meromorphic}) if there exists a modification $\sigma\colon N_{1} \rightarrow N$ whose center is contained in $\Sigma$ and a holomorphic map $\varphi_{1}\colon N_{1} \to \mathcal{C}_{n}^{f}(M)$ extending the holomorphic map $\varphi\circ \sigma_{|\sigma^{-1}(N\setminus\Sigma)}$.
\end{defn}

To define the {\em graph} of such a meromorphic map we need the following corollary of Theorem \ref{semi-proper direct image bis}.

\begin{cor}\label{appli.}
 Fix a complex space $M$ and an integer $n\geq 0$. Let $N$ and $P$ be reduced complex spaces and $\varphi\colon N \rightarrow P \times \mathcal{C}_{n}^{f}(M)$ be a semi-proper holomorphic map. Then $\varphi(N)$ is a reduced complex subspace of $P \times \mathcal{C}_{n}^{f}(M)$.
 \end{cor}
 
 For the proof we use of the following lemma.
 
\begin{lemma}\label{appli.0}
Fix a complex space $M$ and a natural number $n$. Let $P$ be a reduced complex space. Denote $p\colon P\times M \rightarrow P$ and $q\colon P\times M \rightarrow M$ the canonical projections. Then the analytic subset $\mathcal{C}_{n}^{f}(p)^*$ of $\mathcal{C}_{n}^{f}(P \times M)^*$ is biholomorphic to the product $P \times \mathcal{C}_{n}^{f}(M)^*$.
\end{lemma}
  
  \parag{Proof}  Let $\alpha\colon \mathcal{C}_{n}^{f}(p)^* \rightarrow P$ be  the natural holomorphic map (see Proposition \ref{relatif}) and $\beta\colon \mathcal{C}_{n}^{f}(p) \rightarrow \mathcal{C}_{n}^{f}(M)$ be the map induced by the direct image of $n$-cycles by $q$ (see Theorem \ref{direct image}). Let us first show that $\beta$ is holomorphic.\\
   Indeed, every $n$-cycle in $\mathcal{C}_{n}^{f}(p)^*$ is of the form $\{x\}\times C$ where $C\in\mathcal{C}_{n}^{f}(M)^*$ and $q_*(\{x\}\times C) = C$. It follows  that, for every $f$-analytic family $(X_s)_{s\in S}$ of $p$-relative $n$-cycles in $P\times M$, the family $(q_*X_s)_{s\in S}$ of $n$-cycles in $M$  is well defined and  $f$-analytic. Hence  the map $(\alpha, \beta)\colon \mathcal{C}_{n}^{f}(p)^*  \rightarrow P \times \mathcal{C}_{n}^{f}(M)^*$ is  bijective and holomorphic. 
The inverse map $\gamma\colon P\times\mathcal{C}_{n}^{f}(M)^*\rightarrow \mathcal{C}_{n}^{f}(p)^*$,  given by $\gamma(p, C) := \{p\}\times C$, is also holomorphic  thanks to the product theorem for analytic families of cycles\footnote{This theorem, which was first proved in \cite{[B.75]}, is even not obvious in this simple case because an $n$-scale on $P \times M$ adapted to a cycle like  $\{p\}\times C$  is not necessarily given by the product of  an $n$-scale on $M$ with a local embedding for $P$.} 
(See Theorem 4.6.4 in \cite{[e]} or Theorem IV.6.2.3 in \cite{[BM.1]}).$\hfill \blacksquare$
  
\parag{Proof of Corollary \ref{appli.}} Note that $P \times \mathcal{C}_{n}^{f}(M)$ is the union of the disjoint open subsets $P \times \{\emptyset[n]\}$ and $P \times \mathcal{C}_{n}^{f}(M)^*$ and it follows that  the maps 
$$
\varphi^{-1}(P \times \{\emptyset[n]\})\rightarrow P \times \{\emptyset[n]\}\qquad\text{and}\qquad\varphi^{-1}(P \times \mathcal{C}_{n}^{f}(M)^*)\rightarrow P \times \mathcal{C}_{n}^{f}(M)^*
$$ 
induced by $\varphi$, are semi-proper and holomorphic. As Kuhlmann's theorem implies that $\varphi(N)\cap (P \times \{\emptyset[n]\})$ is a reduced complex subspace of $P \times \mathcal{C}_{n}^{f}(M)$  it is enough to consider the case where $\varphi$ takes its values in $P \times \mathcal{C}_{n}^{f}(M)^*$. But, with the notation of Lemma \ref{appli.0},  the spaces $P \times \mathcal{C}_{n}^{f}(M)^*$ and $\mathcal{C}_{n}^{f}(p)^*$ are biholomorphic so  Theorem \ref{semi-proper direct image bis} gives us the required result.
$\hfill \blacksquare$

\begin{prop}\label{graph mero 1}
Let $M$ and $N$ be reduced complex spaces, $\Sigma$ be a nowhere dense analytic subset of $N$ and  $n$ be a natural number. Let $\varphi\colon  N\setminus\Sigma\rightarrow \mathcal{C}_{n}^{f}(M)$ be a holomorphic map and let $\Gamma$ denote the closure of its graph in $N \times \mathcal{C}_{n}^{f}(M)$. Then $\varphi$ is meromorphic (along $\Sigma$) if and only if $\Gamma$ is an $N$-proper reduced complex subspace of $N \times \mathcal{C}_{n}^{f}(M)$.
\end{prop}

\parag{Proof} Suppose first that $\Gamma$ is a reduced complex subspace of $N \times \mathcal{C}_{n}^{f}(M)$ which is $N$-proper and let $\sigma\colon\Gamma\rightarrow N$ denote the natural projection. Then $\sigma$ is a holomorphic surjection and $\sigma^{-1}(\Sigma)$ is nowhere dense in $\Gamma$. It follows that $\sigma$ is a modification whose center is contained in $\Sigma$, and then the natural projection $\Gamma\rightarrow\mathcal{C}_{n}^{f}(M)$ is a holomorphic extension of $\varphi\circ\sigma_{|\sigma^{-1}(N\setminus\Sigma)}$. So $\varphi$ is meromorphic along $\Sigma$. 

\smallskip
Conversely, suppose that $\varphi$ is meromorphic along $\Sigma$. Then there exists a modification $\sigma\colon N_1\rightarrow N$, whose center is contained in $\Sigma$, and a holomorphic extension $$\varphi_1\colon N_1\rightarrow\mathcal{C}_{n}^{f}(M)$$
 of $\varphi\circ\sigma_{|\sigma^{-1}(N\setminus\Sigma)}$. The map $(\sigma,\varphi_1)\colon N_1\rightarrow N\times\mathcal{C}_{n}^{f}(M)$ is proper since it is the composition of the proper maps\footnote{Remember that a {\em proper map}  between Hausdorff spaces  is a closed map with compact fibers.} $(\id_{N_1},\varphi_1)\colon N_1\rightarrow N_1\times \mathcal{C}_{n}^{f}(M)$ and 
$$
\sigma\times\id_{\mathcal{C}_{n}^{f}(M)}\colon N_1\times \mathcal{C}_{n}^{f}(M)\rightarrow N\times\mathcal{C}_{n}^{f}(M).
$$ 
It follows that the image of $(\sigma,\varphi_1)$ is closed in $N\times\mathcal{C}_{n}^{f}(M)$ and consequently equal to $\Gamma$, as $N_1\setminus \sigma^{-1}(\Sigma)$ is dense in $N_1$. Then, by Corollary \ref{appli.}, the image of $(\sigma,\varphi_1)$ is a reduced complex subspace of $N\times\mathcal{C}_{n}^{f}(M)$. As the natural projection $\Gamma\rightarrow N$ is clearly proper the proof is completed.
$\hfill \blacksquare$

\begin{defn}\label{graph mero 2}
In the situation of Proposition \ref{graph mero 1} suppose that $\varphi$ is a meromorphic map. Then the $N$-proper reduced complex subspace $\Gamma$ of  $N \times  \mathcal{C}_{n}^{f}(M)$ with its natural projection on $N$ (which is a modification)  is called the {\bf\em graph} of the meromorphic map $\varphi$.
\end{defn}

\parag{Remark}
Under the hypotheses of Proposition \ref{graph mero 1} the mapping $\varphi$ is meromorphic along $\Sigma$ if and only if, for every open subset $V$ of $N$, the restriction $\varphi_{|V}$ is meromorphic along $\Sigma\cap V$.
\hfill{$\square$}

\section{Complements}

\subsection{Connected cycles}

\vspace{4mm}
In the cycle space $\mathcal{C}_n(M)$ of compact $n$-cycles of a complex space $M$ the connected cycles form a (closed) analytic subset (see Theorem 4.7.4 in \cite{[e]} or Theorem IV.7.2.1 in \cite{[BM.1]}), but this is no longer true for the space  $\mathcal{C}_n^f(M)$ of $n$-cycles of finite type in $M$. This is easily seen in the family of conics 
$$C_s := \{ x^2 + sy^2 = 1\}, s \in \mathbb{C}$$
 in $\mathcal{C}_1^f(\mathbb{C}^2)$ : for $s \not= 0$ the conic $C_s$ is smooth and connected and $C_0$ is the  disjoint union of two lines.\\
 In this paragraph we shall give a few results on connected cycles of finite type and give some  examples which show the difference between $\mathcal{C}_n(M)$ and 
$\mathcal{C}_n^f(M)$ in this regard.

\begin{defn}\label{weight 0}
For a finite type cycle $X$  in a complex space $M$ with canonical expression \ $X = \sum\limits_{i \in I} n_{i}.X_{i}$ \ the integer 
$$  
w(X) := \sum_{i\in I} n_{i} 
$$
will be called the {\bf\em  weight o}f $X$.
\end{defn}

The weight function on $\mathcal{C}_{n}^{f}(M)$  is characterized by the fact that it is additive and takes value $1$ on every irreducible cycle, with the convention $w(\emptyset[n]) = 0$.

\begin{prop}\label{weight 1}
Let $M$ be a complex space and let $n$ and $k$ be two non-negative integers. Then the subset
$$ 
F_{k} := \{ X \in \mathcal{C}_{n}^{f}(M) \ / \  w(X) \geq k \} 
$$
is a closed subset in $\mathcal{C}_{n}^{f}(M)$.
\end{prop}

\parag{Proof}It is enough to consider the case $k \geq 2$ because $F_0 = \mathcal{C}_n(M) = F_1 \cup \{\emptyset[n]\}$ and $\{\emptyset[n]\}$ is open and closed in $ \mathcal{C}_n(M)$.\\
  Let $X$ be a non empty  $n$-cycle of finite type in $M$ having the canonical expression $X = \sum_{i \in I} n_{i}.X_{i}$. We intend to show that there exists an open neighborhood of $X$ in $\mathcal{C}_{n}^{f}(M)$ in which every cycle is of weight at most $w(X)$. To this end choose, for each $i\in I$, an  $n$-scale $E_{i}$ on $M$ adapted to $X$ such that $\deg_{E_i}(X_i) = 1$ and $\deg_{E_i}(X_j) = 0$ for all $j\in I\setminus\{i\}$. Denote $c(E_i)$ the center of $E_i$ and put $U := \cup_{i \in I} c(E_{i})$. Then $U$ is a relatively compact open subset of $M$ and we claim that every cycle in the open neighborhood 
$$
\mathcal{U} \ := \ \Omega(U)\cap\left(\bigcap_{i\in I}\Omega_{n_i}(E_i)\right)
$$
of $X$ in $\mathcal{C}_{n}^{f}(M)$ is of weight at most $w(X)$. To see this, let $Y$ be a non empty  cycle in $\mathcal{U}$ with canonical expression $Y = \sum_{a \in A} n_{a}.Y_{a}$ and let $A_i$ denote the set of all $a$ in $A$ such that $Y_a$ intersects $c(E_i)$.  Then $\deg_{E_i}(Y) = \sum_{a\in A_i}n_a.\deg_{E_i}(Y_a) = n_i$ so that $\sum_{a\in A_i}n_a \leq n_i$. Hence $w(Y) = \sum_{a\in A}n_a \leq \sum_{i\in I} n_{i} = w(X)$.
\hfill{$\blacksquare$}\\

The set $F_{k}\cap\mathcal{C}_n(M) = \{ X \in \mathcal{C}_{n}(M) \ | \  w(X) \geq k \}$ is an analytic subset of the reduced complex space $ \mathcal{C}_n(M)$ (see Proposition 4.7.2 in \cite{[e]} or Proposition IV.7.1.2 in \cite{[BM.1]}). \\
{\bf This is no longer true for $F_k$} in general, as we can see from the following example.

\parag{Example} \ Let $N$ be a  complex connected manifold, $A$ be a closed subset of $N$ and consider the following open subset of $N\times\C$
$$
M \ := \ (N\times \C)\setminus (A\times \R) \ = \ (A\times(\C\setminus\R))\cup((N\setminus A)\times\C).
$$ 
Then the fibers of the natural projection $M\rightarrow N$ form an $f$-analytic family of (reduced) $1$-cycles $(X_s)_{s\in N}$. It is clear that $w(X_s) = 1$ for all $s\in N\setminus A$ and $w(X_s) = 2$ for all $s\in A$.\\

Remark that in this example the map $\pi : M \to N$ is a submersion between two complex connected manifolds if we assume that $N \setminus A$ is connected. \\
Note also that $A$ may be quite far from being an analytic subset in $N$.\\

In general, a limit of a convergent sequence of connected cycles in $\mathcal{C}_n^f(M)$ is not connected, as it may be seen on the example above. Nevertheless we have the following result.

\begin{prop}\label{connexe 1}
Let $(X_{\nu})_{\nu \geq 1}$ be a sequence in $ \mathcal{C}_{n}^{f}(M)$ converging to a cycle $X$ in $ \mathcal{C}_{n}^{f}(M)$. Let $M'$ be a relatively compact open set in $M$ such that every irreducible component of $X$ meets $M'$ and such that $\overline{X_{\nu}\cap M'}$ is connected for every $\nu$. Then $\overline{|X| \cap M'}$ is connected.
\end{prop}

\parag{Proof} Suppose, on the contrary, that $\overline{|X| \cap M'}$ is not connected. Then we can write $\overline{X \cap M'} = K_1\cup K_2$ where $K_1$ and $K_2$ are non-empty, compact and disjoint. It follows that there exist disjoint open neighbourhoods, $L_1$ of $K_1$ and $L_2$ of $K_2$, in $\overline{M'}$. Thus the set $K := \overline{M'}\setminus (L_1\cup L_2)$ is compact and does not meet  $\overline{|X| \cap M'}$. Since the set of cycles in $\mathcal{C}_{n}^{f}(M)$ which do not intersect the compact set $K$ form an open set, there exists $\nu_0$ such that $|X_{\nu}|\cap K = \emptyset$ for all \ $\nu\geq \nu_0$. This implies that there exists a subsequence $(X_{\nu_k})_k$ of $(X_{\nu})_{\nu \geq 1}$ having the property that $\overline{|X_{\nu_k}|\cap M'}$ is either contained in $L_1$ for all $k$ or in $L_2$ for all $k$. Let us show that this contradicts the fact that $(X_{\nu_k})_k$ converges to $X$:\\
So assume that $\overline{|X_{\nu_k}|\cap M'}$ is  contained in $L_1$ for all $k$. Then take a point $x_0$ in $\overline{X \cap M'}\cap K_2$. As there is a sequence of points in $X \cap M'$ converging to $x_0$ there exists a point $y_0$ in $X \cap M' \cap L'_2$ where $L'_2$ is the interior of $L_2$. Take a scale  $E$ on $M'\cap L'_2$ adapted to $X$ with $y_0 \in c(E)$. Then $\deg_E(X) \geq 1$ and for $k$ large enough the scale $E$ will be adapted to $X_{\nu_k}$ and we shall have $\deg_E(X_{\nu_k}) = \deg_E(X) \geq 1$. But this implies $X_{\nu_k} \cap M' \cap L'_2 \not= \emptyset$. Contradiction. So $\overline{|X| \cap M'}$ is connected. \hfill{$\blacksquare$}

\parag{Remarks}
\begin{enumerate}[(i)]
\item
In the above situation, $\overline{|X| \cap M'}$ is connected if  $|X| \cap M'$ is connected, but not vice versa. 
Observe also that the cycle $X$ is connected if $\overline{|X| \cap M'}$ is connected.
\item
Proposition \ref{connexe 1} is false if we replace $\overline{X_{\nu}\cap M'}$ by  $X_{\nu}\cap M'$ and $\overline{|X| \cap M'}$ by $|X| \cap M'$. This can be seen from the examples below.
\end{enumerate}

\parag{Example 1} Let $M$ be $\C^2$  and put 
$$
M' := (D\setminus\R)\times D,\qquad X_{\nu} := \C\times\{1/\nu\}\qquad\text{and}\qquad X := \C\times\{0\}
$$
where \ $\nu\in\N^*$ \ and $D$ is the unit disc in $\C$.

\parag{Example 2} Let $M := \mathbb{C}^2$ and $M' := \{\vert x \vert^2 + \vert y\vert^2 < 4 \}$ and consider the $1$-cycle
 $$C := \{(x, y)\in \mathbb{C}^2 \ / \ 4(x-1)^2 + y^2 = 4 \} .$$
 It is easy to see that $C \cap M'$ is smooth an connected so that $\overline{C \cap M'}$ is connected, but that $C \cap \bar M'$ has an isolated point $(2, 0)$, thanks to the following elementary computations, where we put $u := x-2$:
 \begin{align*}
 & y^2 = 4 - 4(1+u)^2 = -4u(2+u) \quad {\rm on } \ C \quad {\rm so}  \\
 & \vert 2+ u \vert^2 + 4\vert u\vert \vert 2+ u\vert -4 = \vert 2+u\vert \big[\vert 2+u\vert + 4 \vert u\vert\big] - 4 \geq \vert u\vert (4 - 3\vert u\vert) 
 \end{align*}
 and this is positive for $0 < \vert u\vert < 4/3$.\\
  It is easy to see that for $0 < \lambda \ll 1$ the conic $C_\lambda$ obtained from $C$ by the translation of vector $(\lambda, 0)$ does not meet $\bar M'$ in a neighborhood of the point $(2, 0)$. And of course when $\lambda \to 0$ the cycle $C_\lambda$ converges to $C$ in $\mathcal{C}_1^f(M)$.

\begin{cor}
Let $M$ be a reduced complex space and $M'$ be a relatively compact open subset of $M$. Then the set 
$$
\{X\in\Omega(M') \  / \ \overline{|X| \cap M'}\quad\text{\rm is connected}\}
$$
is a closed subset of $\Omega(M')$.
\end{cor}

\parag{Proof} This follows directly from Proposition \ref{connexe 1}.
\hfill{$\blacksquare$}

\parag{Example} Let $(z,w)$ denote the standard coordinates of $\C^2$. Then for any open neighborhood $\mathcal{U}$ of the connected cycle $X := 2.(\{0\}\times\C)$ in $\mathcal{C}_1^f(\C^2)$ there exists a complex number $a\neq 0$  such that, for every \ $n\in\N^*$,  the union of the complex line through the points $(a,0)$ and $(0,n)$ and the complex line through the points $(-a,0)$ and $(0,n)$, henceforth denoted by $Y_n$, is an element of $\mathcal{U}$. Moreover, the sequence $(Y_n)_{n\geq 1}$ of connected   $1$-cycles converges to the non-connected cycle 
$$
Y = \{z = a\} + \{z = -a\}
$$ 
in $\mathcal{C}_1^f(\C^2)$ and satisfy $w(Y_n) = w(Y) = 2$. This shows that the set of connected cycles in the closed subset $F_2$  is not locally closed. Hence the set of  connected cycles in $\mathcal{C}_1^f(\C^2)$ is not { \em locally closed}.\\

\subsection{Restriction}

We consider an irreducible complex space $M$ and a closed analytic subset $T$ which has no interior point in $M$. Then we define
$$\mathcal{C}^{f}_{n}(M, T) := \{ C \in \mathcal{C}^{f}_{n}(M) \ / \dim(\ \vert C\vert \cap T ) \leq n-1 \} .$$
Then for $C \in \mathcal{C}^{f}_{n}(M)$ the cycle $r^f(C) := C \setminus (\vert C\vert \cap T)$ is an element of $\mathcal{C}^{f}_{n}(M\setminus T)$.
Hence we have the following result.

\begin{lemma}\label{hol. inj.}
The subset $\mathcal{C}^{f}_{n}(M, T)$ is a  Zariski open subset in $\mathcal{C}^{f}_{n}(M)$\footnote{This means it is  the complement of a (closed) analytic subset in $ \mathcal{C}^{f}_{n}(M)$.} and the injective map
$$\rho^f:  \mathcal{C}^{f}_{n}(M, T) \to \mathcal{C}^{f}_{n}(M\setminus T), $$
induced by the restriction map $r^f$, is holomorphic.
\end{lemma}

\parag{Proof} The fact that the complement $\mathcal{T}$ of  $\mathcal{C}^{f}_{n}(M, T)$ is a closed analytic subset in $ \mathcal{C}^{f}_{n}(M)$ is proved in Proposition \ref{une comp.} point $(ii)$.\\
The only point to prove to obtain that $\rho^f$ is holomorphic is the continuity of $\rho^f$, because a scale on $M \setminus T$ is also a scale on $M$.\\
But a relatively compact open  subset in $M \setminus T$ is also open and relatively compact in $M$ and the fact that for $X \in \mathcal{C}^{f}_{n}(M, T)$ the irreducible components of $X$ are the closure of their intersection with $M \setminus T$ implies that for any relatively compact open set $W$ in $M \setminus T$ the inclusion $\Omega(W) \subset \mathcal{C}_n^f(M, T)$ in $\mathcal{C}_n^f(M)$ holds. $\hfill \blacksquare$\\

Note that the map $r^f$ is { \bf not continuous} at a point $C \in \mathcal{T}$ when $C$ is a limit of a sequence $(C_\nu)$ such that there exists, for each $\nu$, an irreducible component $\Gamma_\nu$ of $C_\nu$ with the property that the sequence $(\Gamma_\nu)$ converges in $\mathcal{C}_n^f(M)$  to a cycle contained in $T$.\\

\begin{thm}\label{res.2}
 Let $M$ be a compact complex space which is strongly K\"ahler (see \cite{[BM.2]}) and $T$ be a compact analytic subset with no interior point in $M$. Denote by ${\rm vol}_h$ the volume of $n$-cycles for the K\"ahler form $h$ of $M$\footnote{Recall that in this situation the volume of compact cycles is constant on a connected component of $\mathcal{C}_n(M)$.}. Then let $S_d$ be the union of the connected  components of $\mathcal{C}_n(M)$ such that the volume of their members is bounded by a constant $d > 0$ and note  $i : S_d  \to \mathcal{C}_n(M)$ the inclusion map. \\
Then  the image of  the map $\rho^f\circ i : S_d \cap \mathcal{C}_n(M, T) \to \mathcal{C}_n^f(M \setminus T)$ induces a closed analytic subset which is a reduced complex space in the open set  $$\mathcal{C}_n^f(M \setminus T) \setminus r^f(i(\mathcal{T} \cap S_d))$$ of $\mathcal{C}_n^f(M \setminus T)$.
\end{thm}

\parag{Proof} We begin by  proving that the subset $r^f(i(\mathcal{T} \cap S_d))$ of $\mathcal{C}_n^f(M \setminus T)$ is closed in $\mathcal{C}_n^f(M \setminus T)$.  \\
The subset $S_d \cap \mathcal{T}$ and then
 $i(S_d\cap \mathcal{T})$ are compact because $\mathcal{T}$ is closed and $i$ is continuous. But remember that $r^f$  is not continuous in general.\\
  So take  $C_0 \in \mathcal{C}_n^f(M \setminus T)$ and assume that $C_0$ is the limit in $\mathcal{C}_n^f(M \setminus T)$  of a sequence $r^f(i(X_\nu))$ where $i(X_\nu)$ are in $i(S_d\cap \mathcal{T})$. By taking a subsequence we may assume that the sequence $(X_\nu)$ converges to a cycle $X$ in $S_d \cap \mathcal{T}$ and all the cycles have the same volume. Fix an open relatively compact  subset $W$ in $M \setminus T$ such that any irreducible component of $C$ meets $W$. Let now 
$(Y_\nu)$ be the sequence obtained from $X_\nu$ by deleting all irreducible components of $X_\nu$ disjoint from $W$. In particular we delete at least the irreducible components of $X_\nu$ contained in $T$, and there exists at least one such component, so the h-volume of $Y_\nu$ is strictly less than the h-volume of $X_\nu$. Then $Y_\nu$ is in $S_d \setminus S_d \cap \mathcal{T}$ for each $\nu$. Again by taking a subsequence if necessary, we may assume that the sequence $(Y_\nu)$ has a fix h-volume and  converges to a cycle $Y$ in $S_d$.  Then it is clear that the sequence $r^f(i(Y_\nu))$ converges to $C$ in $\mathcal{C}_n^f(M \setminus T)$ because the limit of the sequence $r^f(i(Y_\nu))$ exists in $\mathcal{C}_n^f(M \setminus T)$ thanks to Corollary \ref{compact} and coincides with $C_0$ on $W$. \\
Remark that the volume of $Y$ is strictly less that the volume of $X$. So $Y \not= X$ and  $X$ has at least one irreducible component $\Gamma$ in $T$ which is not in $Y$. Then the cycle 
$Y+\Gamma$ is in $S_d$ and so are the cycles $Y_\nu + \Gamma$ which have the same volume than $Y + \Gamma$.\\
Then the sequence $(Y_\nu + \Gamma)$ is in $S_d \cap \mathcal{T}$ and $r^f(i(Y_\nu + \Gamma))$ converges to $C_0$ in $\mathcal{C}(M \setminus T)$  concluding the proof  that $r^f(i(S_d\cap \mathcal{T})$  is closed in $\mathcal{C}_n^f(M \setminus T) $.\\

So the target  set $ \mathcal{C}_n^f(M \setminus T) \setminus r^f(i((\mathcal{T} \cap S_d))$  is open in $\mathcal{C}_n^f(M \setminus T)$ and we may apply the Semi-Proper Direct Image Theorem \ref{semi-proper direct image bis} to the map $\rho^f\circ i$ if we are 	able to prove that it is a semi-proper map.\\

The main point is now to show that the map $\rho^f\circ i $ is a  semi-proper map. \\

We shall prove first  that the image of the map $\rho^f\circ i$ induces a closed  subset in  the target open set 
$$\mathcal{C}_n^f(M \setminus T) \setminus r^f(i(\mathcal{T} \cap S_d)).$$
 So let $(Y_\nu = \rho^f(i(X_\nu))) $ be a sequence of  cycles in this image converging to a cycle $Y$ in $\mathcal{C}_n^f(M \setminus T)$. By compactness of $S_d$ we may assume that the sequence $(X_\nu)$ converges to a cycle $X$ in $S_d$ and if $X$ is not in $S_d \cap \mathcal{T}$ the continuity of $\rho^f\circ i$ implies that $Y = \rho^f(i(X))$ and we are done. \\
 So assume that $X = Z + A$ where $A$ is a cycle in $S_d\cap \mathcal{T}$ and $Z$ is in $S_d \setminus S_d\cap \mathcal{T}$. Then the cycle $r^f(Z)$ is equal to $Y$ by uniqueness of the limit in $\mathcal{C}_n^{\rm loc}(M \setminus T)$ and this is impossible because we assumed that $Y$ is not in $r^f(S_d\cap\mathcal{T})$. So the image of $\rho^f\circ i $  is closed in $\mathcal{C}_n^f(M \setminus T) \setminus r^f(i(\mathcal{T} \cap S_d))$ and it is enough to check the semi-properness of the map at points in its image.\\

To complete the proof we shall use the following lemma.

\begin{lemma}\label{semi-proper}
Let $A$ be an analytic subset of $T$ of dimension at most $n-1$ and let $U$ a $(n-1)$-complete neighborhood of $A$\footnote{See \cite{[B.80]} where the existence of a basis of such open neihgborhood is proved. See also \cite{[BM.2]} Chapter XI Theorem  3.1.1.}. Let $W$ be an open neighborhood of $ T \setminus T \cap U$ and let $\mathcal{K}$ be the subset of $S_d$ of cycles disjoint from $W$. Then $\mathcal{K}$ is a compact set in $S_d \setminus (S_d\cap \mathcal{T})$.
\end{lemma}

\parag{Proof} First remark that $\mathcal{K}$ does not intersect $S_d \cap \mathcal{T}$ because if  a cycle $Z$ is disjoint from $W$ has an irreducible component in $T$ this component has to be inside $U$ which is $(n-1)$complete. This is not possible. Then,  as $S_d$ is compact and $\mathcal{K}$ is disjoint from $\mathcal{T}$ it is enough to prove that $\mathcal{K}$ is closed in $S_d$. So assume that we have a sequence $(Z_\nu)$ in $\mathcal{K}$ which converges to a cycle $Z$ in $S_d$. Then $Z$ is not in $\mathcal{K}$; it means that $Z$ meets $W$. But then we may find an $n$-scale $E$  in $W$ adapted to $Z$ with $\deg_E(Z) \geq 1$. And in the open neighborhood $\Omega_k(E)$  of $Z$ any cycle meets $W$;  so $Z_\nu$ is not in $\mathcal{K}$ for $\nu$ large enough. Contradiction.$\hfill \blacksquare$\\ 

\parag{End of proof of Theorem \ref{res.2}}To prove the semi-properness of the map $\rho^f\circ i $ let $Y$ be the restriction to $M \setminus T$ of a cycle $X \in S_d \setminus S_d\cap \mathcal{T}$ with $Y \not\in r^f(i(S_d\cap \mathcal{T}))$. Let $A := \vert X\vert \cap T$. This is an analytic subset of $T$ of dimension at most $(n-1)$. Then define $U$ and $W$ as in the previous lemma but small enough in order that any irreducible component of $Y$ does not meet $W$ but meets the relatively compact open set $M \setminus V$ in $M \setminus T$ where $V := W \cup U$. We want to prove that there is a neighborhood $\mathcal{V}$ of $Y$ in $\mathcal{C}_n^f(M \setminus T)$ such that any $Y' \in \mathcal{V} $ which is in the image of $\rho^f\circ i $ is in fact image of a cycle in $\mathcal{K}$.\\
So assume that $\Omega(V) \cap \rho^f(i(\mathcal{K}))$ is not a neighborhood of $Y$ in the image of $\rho^f\circ i$. Then there exists a sequence $(Y_\nu = \rho^f(i(X_\nu))$ in $\Omega(V)$ which converges to $Y$ in $\mathcal{C}_n^f(M \setminus T)$, where $X_\nu$ is in $S_d \setminus (S_d \cap \mathcal{T})$ for each $\nu$ and not in $\mathcal{K}$. Passing to a sub-sequence we may assume that  the sequence $X_\nu$ converges to a cycle $X$ in $S_d$ and we have $r^f(i(X)) = Y$ which implies that $X$ is not in $S_d\cap \mathcal{T}$ because $Y$ is not in $r^f(i(S_d\cap \mathcal{T}))$.\\
Then each irreducible component of $X$ is an irreducible component of $Y$ and then does not meet $W$. So for $\nu$ large enough any irreducible component of $X_\nu$ does not meet $W$ and then $X_\nu$ is in $\mathcal{K}$ for $\nu$ large enough. This contradicts our hypothesis.\\
Then the holomorphic map $\rho^f\circ i$ is semi-proper and the conclusion follows..$\hfill \blacksquare$\\

For instance, for $M := \mathbb{P}_m$  and let $h$ be the Fubini-Study metric on it; let  $T$ be  a hyperplane so that $M \setminus T$ is $\mathbb{C}^m$. We obtain  that algebraic cycles in $\mathbb{C}^m$  of dimension $n$ and  of degree at most  the integer  $d \geq 1 $ form a closed analytic subspace in $\mathcal{C}_n^f(\mathbb{C}^m)$ which is a reduced complex space. As the subset $S_d \cap\mathcal{T}$ in this case contains all $n$-cycles in $\mathbb{P}_m$ of dimension $n$ which are of degree at most $d-1$, we see that this reduced complex space in $\mathcal{C}_n^f(\mathbb{C}^m)$ is isomorphic \footnote{The fact that this holomorphic bijective map is an isomorphism may be obtained as in Theorem \ref{im.fiber map}  below.} via  the restriction  $\rho^f$ to the space of  $n$-cycles in $\mathbb{P}_n$ of degree $d$  with no irreducible component in $T$. \\

\parag{Remark} It is easy to extend such a result in the case where $M$ is a compact complex space of the class $\mathscr{C}$ of Fujiki-Varouchas (see \cite{[BM.2]} chapter XII). \marginpar{formulation precise ??}

\chapter{Geometrically $f$-flat maps and strongly quasi-proper maps}

 \section{Preliminaries}
 
 In this section we  give some technical results which will be used later on.\\
 
\begin{lemma}\label{very general}
Let $f : M \to N$ be a holomorphic map between reduced complex spaces and let $A$ be a closed analytic subset in $M$ with empty exterior. Then the set of points in $M$ such that $f^{-1}(f(x))$ has no irreducible component contained in $A$ is very general in $M$.
\end{lemma}

\parag{Proof} It is clearly enough to treat the case where $M$ is irreducible because a countable intersection of very general subsets is again very general (see Corollary 2.4.55 in \cite{[e]} for a proof). Then we may assume that $N$ is also irreducible\\
Let $q : M \times_N M \to M$ be the second projection. Then for each $x \in M$ we have
$$ q^{-1}(x) = f^{-1}(f(x))\times \{x\} \quad {\rm and} \quad  q^{-1}(x)\cap (A \times_N M) = (f^{-1}(f(x))\cap A)\times \{x\} .$$
So $A$ contains an irreducible component of $f^{-1}(f(x))$ if and only if $A\times_N M$ contains an irreducible component of $q^{-1}(x)$.\\
For each irreducible component $C$ of $M\times_N M$ let $T_C$ be the subset of points $x$  in $M$ such that $ (A \times_NM)\cap C$ contains an irreducible component of $q^{-1}(x)$. Again it is enough to prove that for each $C$ the subset $M \setminus T_C$ is very general in $M$.\\
So fix an irreducible component $C$ of $M \times_NM$. First remark that if $C$ does not dominate $M$ then $M \setminus q(C)$ is very general in $M$ thanks to Proposition 2.4.60 in \cite{[e]} and it follows that $M \setminus T_C $ is also very general  because  $T_C \subset q(C)$ .\\
So consider now the case where $C$ dominates $M$. Then, applying Corollary 2.4.61 of \cite{[e]} to the map $q_C : C \to M$, it is enough to prove that $(A \times_N M)\cap C$ has no interior point in $C$.\\
Assume on the contrary that there exists a point $(a, x) \in A \times_N M$ which  is an interior point of $C$.\\
 Define $k:= \min_{x \in M}  \dim_x f(f^{-1}(f(x)) $ and recall that $\Sigma_k(f) \setminus \Sigma_{k+1}(f)$ is a dense open  subset in $M$. Hence $\Sigma_k(q_C) \setminus \Sigma_{k+1}(q_C)$ is a dense  open  subset in $C$ as $C$ dominates $M$. Let $S(M)$ be the singular locus of  $M$. Then $M \setminus S(M)$ is a dense open subset in $M$ so we may assume that  $a$ and $x$ are smooth points of $M$, that $(a, x)$ is a smooth point in $C$  and that
$$\dim_a f^{-1}(f(a)) = \dim_x f^{-1}(f(x)) = k .$$
Then the first projection $ p_C : C \to M$ is an open map near $(a, x)$ because it is equidimensional and $a$ is a smooth point of $M$. Consequently $p_C((A\times_NM)\cap C)$ contains a neighborhood of $a$
and this contradicts the fact that $A$ has empty interior in $M$.$\hfill \blacksquare$\\

The following characterizations of a dominant map between irreducible complex spaces are useful.

\begin{lemma}\label{dominant}
Let $\pi\colon M \rightarrow N$ be a holomorphic map between two irreducible complex spaces. Then  the following conditions are equivalent:
\begin{enumerate}[(i)]
\item 
The map $\pi$ is dominant.
\item 
The pullback by $\pi$ of any dense subset in $N$ is dense in $M$.
\item  
There exists an open dense subset $M'$ of $M$ such that  the restriction of $\pi$ to $M'$ is an open map.
\item
The generic rank of $\pi$ is equal to $\dim N$.
\end{enumerate}
\end{lemma}

For the proof of this lemma uses the following  more or less standard result.

\begin{lemma}\label{rank max}
Let $\pi\colon M \rightarrow N$ be a holomorphic map between reduced complex spaces, where $N$ is irreducible. Assume that the generic rank of $\pi$ on every irreducible component of $M$ is strictly less than the dimension of $N$.  Then $\pi(M)$ has empty interior in $N$.
\end{lemma}

\parag{Proof} Remark first that, as the singular locus of $N$ is of empty interior in $N$, we may assume $N$ smooth. Hence it is enough to prove the result in the case where $N$ is an open subset of a numerical space. To do so we take an exhaustion by compact subsets, $M = \bigcup_{j \geq 1}K_j$. Then, thanks to Lemma \ref{Kuhl-banach.3}, each compact  subset $\pi(K_j)$ of $N$ is $b$-negligible for all $j$ and consequently $\pi(M) = \bigcup_{j \geq 1}f(K_j)$ has empty interior in $N$ since $N$ is a Baire space .\hfill{$\blacksquare$}

\parag{Proof of Lemma \ref{dominant}} Due to Lemma \ref{rank max}, (i) implies (iv), and by definition of  the {\em generic rank}, (iv) implies (iii).  To show that (ii) follows from (iii) let us take a dense subset $\Lambda$ of $N$. Then for any non-empty open subset $U$ in $M$ we get
$$
\pi(U\cap \pi^{-1}(\Lambda)) \supseteq \pi(M'\cap U\cap \pi^{-1}(\Lambda)) = \pi(M'\cap U)\cap \Lambda
$$ 
and $\pi(M'\cap U)\cap \Lambda$ is non empty since $\pi(M'\cap U)$ is a non-empty open subset of $N$. Finally, to show that  $(ii)$ implies $(i) $ we only have to  notice that the set  $N \setminus \pi(M)$  is dense in $N$ if $\pi$ is not dominant.$\hfill \blacksquare$\\
 
\begin{defn}\label{str.}
Let $\pi\colon M \rightarrow N$ be a holomorphic map between reduced complex spaces, where $N$ is irreducible, and consider a holomorphic map $f\colon Z \rightarrow N$ where $Z$ is an irreducible complex space. We define the {\bf \em strict transform of $\pi$ by $f$} as the holomorphic map $\tilde{\pi} : \tilde{Z} \to Z$ which is the natural projection onto $Z$ of the union $\tilde{Z}$ of those irreducible components of the fiber product $Z \times_{N}M$ which dominate $Z$. We shall denote $\tilde{Z}$ by $Z\times_{N,str}M$.
\end{defn}

In the situation above the space $Z\times_{N,str}M$ will be called the {\bf strict fiber product} (of $M$ over $N$  by $f$).  Note that $Z\times_{N,str}M$ is empty if and only if the image of the natural projection $Z\times_NM\rightarrow Z$ has empty interior, which is equivalent to $f^{-1}(\pi(M))$ being of empty interior in $Z$.\\
 Hence in general we have $Z\times_{N,str}M\subsetneq Z\times_NM$. 

\smallskip
One can easily find examples where $Z\times_NM\rightarrow Z$ is surjective and where the inclusion  $Z\times_{N,str}M\subset Z\times_NM$ is strict. For instance, if $\pi\colon M \to N$ is a modification which is not injective, the fiber product $M\times_{N} M$ has at least one  irreducible component which is not contained in the corresponding  strict fiber product. And the strict fiber product   is naturally isomorphic to $M$ in this case.

\smallskip
It should also be noted that there is a natural isomorphism between the fiber products $Z\times_NM$ and $M\times_NZ$ but the corresponding strict fiber products $Z\times_{N,str}M$ and $M\times_{N,str}Z$ need not be isomorphic.

\medskip 
It is an easy exercise left to the reader to prove that, in the case where $\pi\colon M \rightarrow N$ is a modification and $f\colon Z \rightarrow N$ a dominant map, the projection $Z \times_{N,str}M \rightarrow Z$ is a modification.

\medskip
Of course, each of the natural projections of $Z \times_{N,str}M$ onto $M$ and $Z$ factorizes the natural mapping $Z \times_{N,str}M \rightarrow N$, and when we write the commutative diagram 
$$
\xymatrix{Z \times_{N,str}M \ar[d]_{\tilde{\pi}} \ar[r] & M \ar[d]^{\pi} \\
Z \ar[r]^{f} & N } 
$$
the horizontal arrow above denotes the natural projection unless otherwise explicitly stated.

\medskip
If $M = \bigcup\limits_iM_i$ is the decomposition of $M$ into irreducible components, then $Z \times_{N,str}M $ is the union of the  $Z \times_{N,str}M_i$.

\begin{lemma}\label{trans.}
Consider the following diagram of holomorphic maps between irreducible complex spaces
$$
\xymatrix{\quad & \tilde{Z} \ar[d]_{\tilde{\pi}} \ar[r] & M \ar[d]^{\pi} \\
W \ar[r]^{g} & Z \ar[r]^{f} & N } 
$$
where $\tilde{Z} := Z\times_{N,str}M$ and  $g$ is a dominant map. Then we have a canonical identification between the strict transform of $\tilde{\pi}$ by $g$ and the strict transform of $\pi$ by $f\circ g$.
\end{lemma}

\parag{Proof} As there is a canonical identification between the fiber products : 
$$W\times_{N}M \simeq W \times_{Z}(Z\times_{N}M),$$
 it is enough to show that an irreducible component $A$ of \ $W\times_{N}M$ which dominates $W$ is an irreducible component of \ $W \times_{Z}\tilde{Z}$.\\
  Let $B$ be an irreducible component of $Z\times_{N}M$ containing the image of $A$ by the natural projection $W \times_{Z}(Z\times_{N}M)\rightarrow Z\times_{N}M$. Then $B$ dominates $Z$ because $A$ dominates $W$ and $g$ is assumed to be dominant. So $B$ is in $\tilde{Z}$ and $A$ is an irreducible component of $W\times_{Z}\tilde{Z}$ which dominates $W$.$\hfill \blacksquare$

\begin{lemma}\label{dominating.modification}
For finitely many modifications of a reduced complex space there exists a modification which factorizes through each of them. 
\end{lemma}

\parag{Proof} It is enough to prove the result for two modifications, so let $\tau_1\colon N_1\rightarrow N$ and $\tau_2\colon N_2\rightarrow N$ be modifications of a reduced complex space $N$. Then the strict transform $\tilde{\tau}_2\colon \tilde{N}_2\rightarrow N_1$ of $\tau_2$ by $\tau_1$ is a modification of $N_1$ and it follows that  $\tau_1\circ\tilde{\tau}_2\colon\tilde{N}_2\rightarrow N$ is a modification of $N$ which has the desired properties.
$\hfill \blacksquare$.

\section{ Holomorphic fiber maps and pull-back of cycles}

We recall here the notion of a geometrically flat map (see \cite{[BM.2]} chVI def. 4.6.7). \\

\quad Consider a holomorphic map $\pi: M \to N$ where $N$ is irreducible and define \\
 $n := \dim M - \dim N$. We say that $\pi$ is  {\bf geometrically flat} if it admits a holomorphic fiber map, that is to say a classifying map of an analytic family of cycles in $M$
$$ \varphi : N \to \mathcal{C}_n^{loc}(M)$$
which satisfies  $\vert \varphi(y)\vert = \pi^{-1}(y)$ for every $y \in N$ and, for $y$ very general\footnote{This implies that the cycle $\varphi(y)$ is reduced for all $y$ in a dense subset of $N$.} in $N$,  the cycle $\varphi(y)$ is reduced.\\
Then for an irreducible analytic subset $X \subset N$ of dimension $q$ the {\bf pull-back cycle $\pi^*(X)$} in $M$ is defined as the {\bf graph-cycle} (see \cite{[BM.1]}  or \cite{[e]}  Chapter IV section 3.3) of the analytic family of $n$-cycles in $M$ parametrized by $X$ and given by its classifying map $\varphi_{\vert X}$ which is the restriction of $\varphi$ to $X$. This is a $(n+q)$-cycle in $M$.\\

{ \em For a general $q$-cycle in $N$ the pull-back by $\pi$ is defined by additivity}  (note that $\pi^*(\emptyset[q]) = \emptyset[n+q]$).\\

Then the following result is proved in \cite{[BM.2]} ch.VI th. 4.8.1:

\begin{thm}\label{pull-back 1}
Let $\pi: M \to N$ be a geometrically flat  holomorphic map  between irreducible complex spaces and  let $(X_s)_{s\in S}$ be an analytic family of $q$-cycles in $N$ parametrized by a reduced complex space $S$. Then the family $(\pi^*(X_s))_{s \in S}$ of $(n+q)$-cycles in $M$ is analytic.$\hfill \blacksquare$
\end{thm}

We shall adapt this result to the case of quasi-proper maps  in Theorem \ref{pull-back 2} below, using the space of  finite type cycles.

\begin{defn} \label{ map}
Let $\pi : M \to N$ be a quasi-proper holomorphic map between reduced complex spaces with $M$ pure dimensionnel and  $N$ irreducible.
 Define  $n := \dim M - \dim N$. An {\bf\em $f$-fiber map} for $\pi$ is a holomorphic map $\varphi : N \to \mathcal{C}_n^f(M)$ which satisfies $\vert \varphi(y)\vert = \pi^{-1}(y)$ for all $y \in N$. We shall say that the $f$-fiber map $\varphi$  is {\bf\em reduced} when for $y$ {\em generic} in $N$ the cycle $\varphi(y)$ is reduced.
\end{defn} 

Note that under the hypothesis above the map $\pi$ is equidimensional and surjective when $M$ is not empty.\\
Remark that when a (holomorphic)  reduced $f$-fiber map exists it is unique.

\begin{lemma}\label{ft 1}
Let $\pi: M \to N$  a quasi-proper holomorphic map with $N$ irreducible and assume that there exists a reduced holomorphic $f$-fiber map for $\pi$,  $\varphi : N \to \mathcal{C}_n^f(M)$. Let $X$ be a finite type $q$-cycle in $N$. Then $\pi^*(X)$ is a finite type $(n+q)$-cycle in $M$.
\end{lemma}

\parag{Proof} It is enough to treat the case where $X$ is irreducible. Then $\pi^*(X) \to X$ is quasi-proper so $\pi^(X)$ has finitely many irreducible components.$\hfill \blacksquare$\\

The adapted version of Theorem \ref{pull-back 1} is now  easy.

\begin{thm}\label{pull-back 2}
Let $\pi: M \to N$  be a quasi-proper holomorphic map with $N$ irreducible and assume that there exists a  reduced holomorphic  $f$-fiber map for $\pi$,  $\varphi : N \to \mathcal{C}_n^f(M)$. Let 
$(X_s)_{s\in S}$ be an f-analytic family of $q$-cycles in $N$ parametrized by a reduced complex space $S$. Then the family $(\pi^*(X_s))_{s \in S}$ of $(n+q)$-cycles in $M$ is f-analytic.
\end{thm}

\parag{Proof}Thanks to Theorem \ref{pull-back 1},  the only point to prove is the continuity of the classifying map of the family $(\pi^*(X_s))_{s \in S}$ which takes its values in $\mathcal{C}_{n+q}^f(M)$ thanks to the previous lemma. Let $G$ be the set-theoretic graph of the family $(X_s)_{s\in S}$. Then $G \subset S\times N$ is quasi-proper over $S$. The set-theoretic graph $\Gamma \in S \times M$ of the family $(\pi^*(X_s))_{s \in S}$ is equal to $(\id_S\times \pi)^{-1}(G)$ and the point is to prove that $\Gamma$ is also quasi-proper over $S$. But this is an immediate consequence of Lemma \ref{open.quasipropre}.$\hfill \blacksquare$\\

\parag{Remark} We may reformulate the previous theorem as follows, using the definition of a $f$-GF map which will be given in the next section (see Definition \ref{GF-maps})  
\begin{itemize}
\item Let $\pi: M \to N$  be a $f$-GF map, $n$ the dimension of its fibers and $q \geq 0$ an integer. Then there exists a natural pull-bak map for finite type cycles
$$ \pi_q^* : \mathcal{C}_q^f(N) \to \mathcal{C}_{n+q}^f(M) $$
which is holomorphic for the weak analytic structures of these cycle spaces.\\
\end{itemize}

We  give now two example of $f$-fiber maps.

\parag{Example 1} We define
 $$N := \{ (x, y) \in \mathbb{C}^2 \ / \  x^2 = y^3\}\quad {\rm and} \quad  M := \{((x,y), z) \in N \times \mathbb{C} \ / \ z^2 = y \}.$$
 Then the natural projection
$\pi : M \to N$ gives a proper finite surjective map of degree $2$ which has a  reduced holomorphic $f$-fiber map $\varphi: N  \to \Sym^2(M) \subset \mathcal{C}^f_0(M)$ which associated to $(x, y) \in N$
the cycle of degree $2$
 $$\varphi(x, y) = \{(x,y,\sqrt{y})\} + \{(x,y, -\sqrt{y})\} \in N \times \Sym^2(\C) .$$
The holomorphy of this map is easily checked using the fact that the image of $\pi$ is contained in the closed analytic subspace $N \times \Sym^2(\mathbb{C}) \simeq N \times \mathbb{C}^2$
 of $\mathcal{C}^f_0(\pi)$ which is isomorphic to $N \times \mathbb{C}^2$ via the isomorphism $\Sym^2(\mathbb{C}) \simeq \mathbb{C}^2 $ given by the elementary symmetric functions $S_1$ and $S_2$ since
$$ S_1(\sqrt{y},- \sqrt{y}) = 0 \quad {\rm and} \quad  S_2(\sqrt{y},-\sqrt{y}) = -y .$$
Let $M_+ := \{(x,y), z) \in N\times \C \ / \ x = yz \}$ and $M_-  := \{(x,y), z) \in N\times \mathbb{C} \ / \ x =  -yz \}$. They are closed analytic subsets in $N\times \mathbb{C}$ and as we have $x^2 = y^2z^2$ in $M$ we see that $M_+ \cup M_- = M$ and they are the irreducible components of $M$. The projections $\pi_+ : M_+ \to N$ and $\pi_- = M_- \to N$ are holomorphic homeomorphisms and admits continuous inverses
given by $\psi_\pm(x,y) = (x,y,\pm x/y)$ which are not holomorphic at $(0,0)$, but are continuous meromorphic fiber maps for $\pi_\pm$ respectively.\\
Remark that $\pi_\pm$ does not admit a $f$-holomorphic fiber map, reduced or not.\\
This example show that even if  a (quasi-)proper  geometrically flat holomorphic fiber map  $\pi: M \to N$ with $N$ irreducible, has a reduced $f$-holomorphic  fiber map, it may happen that the restriction of $\pi$ to an irreducible component of $M$ does not admit a  holomorphic  map (reduced or not).

\parag{Example 2} We keep the same irreducible complex space $N$ as in Example 1 but we define now the non reduced complex space
$$M := \{((x,y), z) \in N \times \mathbb{C} \ / \  z^3 - 3yz + 2x = 0\}.$$
 Then again the natural projection
$\pi : M \to N$ gives a proper finite surjective map, which is now  of degree $3$ and  which has a non reduced  holomorphic $f$-fiber map 
$$\psi(x,y) = 2\{(x,y, x/y)\} +\{(x,y,-2x/y)\} \in N \times \Sym^3(\C).$$
The holomorphy of $\psi$ is easy to prove, as above, because the elementary symmetric functions of the cycle $ X(x,y) := 2\{x/y\} + \{-2x/y\} \in \Sym^3(\mathbb{C})$ are respectively equal to
$$ S_1(X(x,y)) = 0,\quad  S_2(X(x,y)) = -3y, \quad S_3(X(x,y) = -2x \quad {\rm for} \quad (x,y) \in N.$$
If we define the analytic subsets of $M$ by
 $$M_+ := \{(x,y), z) \in N\times \C \ / \ x = yz \}\quad {\rm and} \quad M_{-2} := \{(x,y), z) \in N\times \C \ / \ -2 x =  yz \}.$$
 We again find two irreducible components of $M$ which are, respectively the graph of the continuous meromorphic functions  $(x,y) \mapsto x/y $ and $(x,y) \mapsto -2x/y$ on $N$. This gives an example of a (quasi-)proper geometrically flat holomorphic map  $\pi: M \to N$ with $N$ irreducible but $M$ non reduced, which has a non reduced $f$-holomorphic fiber map, such that   $M$ does not admit a holomorphic reduced $f$-fiber map and such its irreducible components does not admit  any holomorphic $f$-fiber map, reduced or not reduced.\\
Note that the natural non reduced structure ( we have in $(M_+ \setminus \{(0,0,0\})$) the identity $$z^3 - 3yz + 2x = (z - x/y)^2(z + 2x/y) = (y - z^2)(z + 2x/y)$$  does not help !

 \section{Geometrically $f$-flat maps}

Among the holomorphic quasi-proper surjective maps, the equidimensional maps between two irreducible complex spaces are certainly the simplest. Suppose  that  we have such a map $\pi\colon M\rightarrow N$ and put   $n := \dim M - \dim N$.

\begin{defn}\label{GF-maps}
We shall say that a map $\pi\colon M \rightarrow N$  is {\bf\em geometrically f-flat} if the following conditions are satisfied:
\begin{enumerate}[(i)]
\item
$M$ is a reduced complex space, $N$ is an irreducible complex space and $\pi$ is holomorphic and surjective.
\item
There exists a reduced holomorphic  $f$-fiber map for $\pi$ (see Definition \ref{ map}) which will be called  {\bf\em the reduced $f$-fiber map} for $\pi$.
\end{enumerate}
\end{defn}
Geometrically   $f$-flat maps will often be called {\bf $f$-GF maps} for short.

The simplest examples of $f$-GF maps are given by the following lemma.

 \begin{lemma}\label{ map 2} 
 Let $\pi: M \to N$ be a quasi-proper holomorphic map between irreducible complex spaces. Let $n := \dim M - \dim N$ and assume  $N$ is normal and $\pi$ is equidimensional. Then there exists a {\em reduced} $f$-fiber map for $\pi$.
 \end{lemma}
 
 \parag{Proof}   This an immediate consequence of Theorem 3.4.1 of \cite{[BM.1]} ch.IV taking into account   the quasi-properness of $\pi$.$\hfill \blacksquare$

\parag{Remarks}  Let $\pi\colon M \rightarrow N$ be a holomorphic mapping from a reduced complex space to an irreducible complex space.
\begin{enumerate}[(i)]
\item
If the map $\pi$ is geometrically $f$-flat, then it is both quasi-proper and open. This is an immediate 
consequence of the continuity of a fiber map. Moreover,  if $M\neq\emptyset$, the map $\pi$ is surjective (remember that $\{\emptyset[n]\}$ is open and closed in $\mathcal{C}_n^f(M)$)  and the space $M$ has pure dimension.
\item 
If $M\neq\emptyset$ and $\pi$  is geometrically $f$-flat, then the reduced fiber map for $\pi$ takes its values in the analytic subset  $\mathcal{C}_n^f(\pi)^*$ and hence induces a holomorphic section of the natural map $\mathcal{C}_n^f(\pi)^*\rightarrow N$ which takes it generic values in the analytic subset of reduced cycles in $\mathcal{C}_n^f(\pi)^*$.

\item When $\pi$ is a quasi-proper equidimensional  map which does not have  a reduced $f$-fiber map we can always take the fiber product of $\pi$ with the normalization map $\nu\colon \tilde{N} \rightarrow N$ and obtain an $f$-GF map $\tilde{\pi}\colon \tilde{M}\rightarrow \tilde{N}$, where $\tilde{M} = \tilde{N}\times_NM$, and then the mapping $\tilde{\pi}$  admits a reduced  fiber map 
$\tilde{\varphi}\colon \tilde{N}\rightarrow \mathcal{C}_n^f(\tilde{M})$. Moreover, composing $\tilde{\varphi} $ with the direct image map $\tilde{\nu} _* : \mathcal{C}_n^f(\tilde{M}) \to \mathcal{C}_n^f(M)$ we obtain a {\em meromorphic}  reduced fiber map for $\pi$ parametrized by $\tilde{N}$. This will be explain later.
\item 
If $M$ is pure dimensional but not necessarily reduced and $\pi$ is flat (in the algebraic sense), then $\pi$ induces an analytic family $(X_y)_{y\in N}$ of $n-$cycles in $M$, where $n := \dim M - \dim N$, such that $|X_y| = \pi^{-1}(y)$ for all $y$ in $N$\footnote{This follows from Theorem X.3.3.5 in \cite{[BM.2]}}. Hence the map $\pi$ is geometrically $f$-flat if and only if it is quasi-proper and $M$ is generically reduced.\\
\end{enumerate}

The following lemma and its corollary give an easy way to recognize an $f$-GF map.

\begin{lemma}\label{hol. fiber map 0}
Let $\pi\colon M\rightarrow N$ be a holomorphic map between irreducible complex spaces. Put $n := \dim M-\dim N$ and assume that the canonical map $\mathcal{C}_n^f(\pi)^*\rightarrow N$ admits a holomorphic section $\varphi$ whose generic values are reduced cycles. Then the map $\pi$ is geometrically $f$-flat and $\varphi$ is the reduced  $f$-fiber map for $\pi$. 
\end{lemma}

\parag{Proof} Let $G$ denote the graph of $\pi$ in $M\times N$, let $\Gamma$ denote the graph cycle in $N\times M$ of the analytic family of $n$-cycles which $\varphi$ classifies and let $\rho\colon M\times N\rightarrow N\times M$ be the canonical isomorphism. Then $\Gamma = |\Gamma|$ and  $\rho(G)$ are (closed) analytic subsets of the same dimension in $N\times M$ and $\Gamma \subseteq \rho(G)$. Since $G$ is irreducible it follows that $\Gamma = \rho(G)$ and consequently  $\varphi$ is the reduced $f$-fiber map for $\pi$. 
$\hfill \blacksquare$

\begin{cor}\label{hol. fiber map}
Let $\pi\colon M\rightarrow N$ be a holomorphic map between irreducible complex spaces.  Put $n := \dim M-\dim N$ and assume that there exists a holomorphic map $\varphi\colon N \rightarrow \mathcal{C}^f_n(M)$ such that the restriction of $\varphi$ to a non-empty open set $N'$ in $N$ is the  reduced $f$-fiber map for the map $\pi^{-1}(N') \rightarrow N'$ induced by $\pi$. Then $\pi$ is an $f$-GF map and $\varphi$ is the reduced $f$-fiber map for $\pi$.
\end{cor}

\parag{Proof} The case where $M = \emptyset$ is trivial so we assume $M\neq\emptyset$. Then $\varphi(N')$ is a subset of $\mathcal{C}^f_n(\pi)^*$ and it follows that $\varphi(N)\subseteq\mathcal{C}^f_n(\pi)^*$ since $N$ is irreducible and $\mathcal{C}^f_n(\pi)^*$ is an analytic subset of $\mathcal{C}^f_n(M)$. Moreover, the composition of the canonical map $\mathcal{C}_n^f(\pi)^*\rightarrow N$ with $\varphi$ is a holomorphic map $N\rightarrow N$ which coincides with $\id_N$ on $N'$, so it is the identity map on $N$. Hence $\varphi$ induces a holomorphic section of the canonical map and Lemma \ref{hol. fiber map 0} allows us to conclude.
$\hfill \blacksquare$\\

The following easy consequence of the Direct Image Theorem \ref{semi-proper direct image bis}  shows that the reduced fiber map of an $f$-GF map $\pi : M \to N$  gives a realization  of $N$ as a reduced complex subspace of the  space of finite type cycles  in $M$.

\begin{thm}\label{im.fiber map}
Let $\pi: M \to N$ be a $f$-GF map and let $\varphi : N \to \mathcal{C}_n^f(M)$ be its reduced $f$-fiber map. Then $\varphi$ is a closed holomorphic embedding of the reduced complex space $N$\footnote{This means that the map $ N  \to \varphi(N)$ is an isomorphism of reduced complex spaces when $\varphi(N)$ is endowed with the sheaf of holomorphic functions induced from  $\mathcal{C}_n^f(M)$.}.
\end{thm}

\parag{Proof} As the map $\varphi$ induces a holomorphic section of the holomorphic map $\alpha: \mathcal{C}^f_n(\pi)^* \to N$ it is a closed topological embedding. It follows that $\varphi(N)$ is a complex subspace of $\mathcal{C}_n^f(M)$ due to the Direct Image Theorem\footnote{Here we consider only the case of  a proper map with finite fibers.}. Hence the map $N \to \varphi(N)$, induced by $\varphi$ is biholomorphic since its inverse is induced by $\alpha$, is holomorphic.$\hfill \blacksquare$\\

\section{Stability properties of $f$-GF maps}

Theorem \ref{pull-back 2} has the following corollary.

\begin{cor}\label{cor 1}
Let $\pi\colon M \rightarrow N$ and $\sigma\colon N \rightarrow P$ be two geometrically f-flat maps. Then the map $\sigma\circ\pi$ is geometrically f-flat.
\end{cor}

\parag{Proof} It is enough to apply the  theorem \ref{pull-back 2} to the map $\pi$ and to the reduced  $f$-fiber map for  $\sigma$ in order to show that $\sigma\circ \pi$ admits a holomorphic reduced  $f$-fiber map.$\hfill \blacksquare$

\begin{cor}\label{cor 3}
Let  $\pi\colon M\rightarrow N$ be a geometrically f-flat map and consider  a holomorphic map  $g\colon Z\rightarrow N$ where $Z$ is an irreducible complex space. Assume that $g(Z)$ is not contained in the subset of non-reduced fibers of $\pi$\footnote{Recall that this is an analytic subset with no interior points in $N$.}. Then the projection
 $\tilde{\pi}\colon Z \times_{N}M \rightarrow Z$ is geometrically f-flat. 
\end{cor}

\parag{Proof} Let $\varphi\colon N \rightarrow \mathcal{C}_n^f(M)$ be the reduced  $f$-fiber map of $\pi$. Then the map $\psi\colon Z \rightarrow \mathcal{C}_n^f(Z \times_{N}M) $ by $\psi(z) := \{z\}\times \varphi(g(z))$ is holomorphic , thanks to the Product Theorem (see Theorem 4.6.4 in \cite{[e]}). As we have $|\psi(z)| = \tilde{\pi}^{-1}(z)$ for all $z \in Z$,  $\psi$ is the reduced $f$-fiber map for $\tilde{\pi}$.$\hfill \blacksquare$\\

It should be noted that, without our hypothesis on $g(Z)$,  the corresponding $f$-fiber map $\psi$ is still an  $f$-fiber map for $\tilde{\pi}$ (but not  reduced in general). So, in the case where $Z$ is normal, it admits nevertheless a reduced fiber map using Theorem 4.2.12 in   \cite{[e]}.

\begin{defn}\label{geom.f-flattening}
Let $\pi\colon M\rightarrow N$ be a surjective holomorphic map from a pure dimensional reduced complex space to an irreducible complex space. A {\bf\em geometric $f$-flattening} (or simply {\bf\em $f$-flattening}) of $\pi$ is a modification $\tau\colon \tilde{N}\rightarrow N$ such that the {\em strict transform, $\tilde{\pi}\colon\tilde{M}\rightarrow\tilde{N}$, of $\pi$ by $\tau$}\footnote{See Definition \ref{str.} above} is an f-GF map.
\end{defn}

In the situation of Definition \ref{geom.f-flattening} it leads from (ii) of Lemma \ref{composition}, that the map $\pi$ is necessarily quasi-proper if it admits an $f$-flattening. On the other hand this condition is not sufficient as is shown in Example 1 below. In section 6 we study in detail the so-called, {\bf strongly quasi-proper maps}, which are exactly those  quasi-proper maps that have a geometric $f$-flattening. Moreover, we will show that a strongly quasi-proper map has a \lq\lq natural\rq\rq\  $f$-flattening. 

\parag{Example 1}  Let
 $$Y := \{((a,b),(x,y)) \in \C^{2}\times  \C^{2} \ / \  P(a, b, x, y) := a.x^{2} + b.x - a^{2}.y^{2}  = 0 \}$$
  and let $\pi : Y \to \C^{2}$ be the first projection onto $\C^{2}$. Then we have the following properties:
\begin{enumerate}[(i)]
\item 
The (algebraic) hypersurface $Y$ of $\C^{4}$ is irreducible (in fact normal and connected).
\item 
The map $\pi : Y \to \C^{2}$ is quasi-proper.
\item 
After blowing-up the origin in $\C^{2}$ the strict transform of $\pi$ is no longer quasi-proper.
\end{enumerate}

\parag{Proof of ${\rm (i)}$} The critical set of the polynomial $P(a,b,x,y)$ is given by the following equations
\begin{equation*}
\quad x^{2} - 2a.y^{2} = 0, \quad x  = 0, \quad 2a.x + b = 0, \quad 2a^{2}.y  = 0 . \tag{1}
\end{equation*}
So the subset $S := \{ a = b = x = 0\} \cup \{ x = y = b = 0\} $ which is one dimensional  is  the singular set of $Y$.  As it has codimension $2$ in $Y$, the hypersurface $Y$ is normal. We shall see below  that each fiber of $\pi$ is connected and then the existence of a holomorphic section \ $(\id, 0) : \C^{2} \to \C^{2}\times \{0\}$ of $\pi$ implies that $Y$ is connected. Hence $Y$ is irreducible.

\parag{Proof of ${\rm (ii)}$} First we shall describe the fibers of $\pi$ as subsets of $\C^{2}$. For $a.b \not= 0$   the fiber $\pi^{-1}(a,b)$ is a smooth conic containing the origin in $\C^{2}$. For  $a \not= 0$ and $b = 0$ the fiber $\pi^{-1}(a,b)$ is the union of two distinct lines through the origin. For $a = 0$ and $b \not= 0$ the fiber $\pi^{-1}(0,b)$ is the line $x  = 0$ which also contains the origin. Finally the fiber $\pi^{-1}(0,0)$ is $\C^{2}$. So each fiber is connected and contains the origin. Then the $\pi$-proper set $\C^{2}\times \{0\}$  meets every irreducible component of any fiber of $\pi$, so this map is quasi-proper.

\parag{Proof of ${\rm (iii)}$} Consider now the blow-up $\tau :  X \to \C^{2}$  of the (reduced) origin in $\C^{2}$. The complex manifold $X$ is the sub-manifold
 $$
 X := \{((a,b),(\alpha,\beta)) \in \C^{2}\times \mathbb{P}_{1}\ / \  a.\beta = b.\alpha\}.
 $$
 It will be enough to show that the strict transform of $\pi$ over the chart $\{ \beta \not = 0\}$ of $X$ is not quasi-proper to achieve our goal. So let $s := \alpha/\beta$. Then we have coordinates $(s,b) \in \C^{2}$ for this chart on $X$. The total transform of $Y$ is given by the equation
$$ s.b.x^{2} + b.x - s^{2}.b^{2}.y^{2}  = 0$$
and, as the function $b$ is not generically zero on the strict transform $\tilde{Y}$ of $Y$ by $\tau$. Then we have
$$\tilde{Y}_{\beta \not= 0} = \{\big((s,b),(x,y)\big) \in \C^{2}\times \C^{2} \ / \  x.(s.x + 1) - b.s^{2}.y^{2} = 0 \}.$$
So the fiber of the strict transform $\tilde{\pi}$  at the point $(s,0)$ is the union of the two lines $ \{x = 0\}$ and $\{x = - 1/s\}$ for $s \not= 0$. Then it is clear that this map is not quasi-proper as an irreducible component of the fiber at $(0,s), s\not= 0$ avoids any compact set in $\C^{2}$ when $s \not= 0$ goes to $0$.$\hfill \square$\\

\parag{Claim} The quasi-proper map $\pi$ in the previous example does not admit a $f$-flattening.\\

The proof is a consequence of the following useful criterium, which will be proved later (see Proposition \ref{characterization.SQP}).

\parag{Criterium} Let $\pi : M \to N$ be a quasi-proper surjective map between irreducible complex spaces and put $n := \dim M - \dim N$. Let $z$ a point in $N$ and assume that there exists a  sequence $(y_\nu)_{\nu \in \mathbb{N}}$ which  converges to $z$ and satisfies the following property:
\begin{itemize}
\item There does not exists a compact subset $K$ of $M$ which meets every irreducible component of  the set $\pi^{-1}(y_\nu)$ for all $\nu$.
\end{itemize}
Then the map $\pi$ does not admit an $f$-geometric flattening.$\hfill \square$

\parag{Proof of the Claim} Now,  consider the  double sequence $y_{\nu, q} := (1/q\nu, 1/q) \in \mathbb{C}^2 \setminus  \{(0, 0)\} $ for $(\nu, q) \in (\mathbb{N}^*)^2$. It converges to $(0, 0)$.\\
 The fiber of $\pi$ at $y_{\nu, q}$ is the smooth conic $ C_{\nu, q} := \{ x^2 + \nu x - y^2/q\nu = 0 \}$. Now the set 
$$ \overline{ \{ C_{\nu, q}\ / \  (\nu, q) \in  (\mathbb{N}^*)^2 \}} \subset \mathcal{C}_1^f(\mathbb{C}^2) $$
contains the sequence $(\{ x= 0\} + \{x = -1/\nu \})_{\nu \geq 1}$ of cycles  which has the sequence $(x = -1/\nu\})_{\nu \geq 1}$ of irreducible components  escaping at infinity when $\nu \to +\infty$. So this closure cannot be a compact subset in $\mathcal{C}_1^f(\mathbb{C}^2)$ (see III.3.1). \\
The following criterium whose proof of this fact is an immediate consequence of \ref{characterization.SQP}, implies that the quasi-proper map $\pi$ does not admit a  geometric $f$-flattening in any open neighborhood of the origin in $\mathbb{C}^2$.$\hfill \blacksquare$

\parag{Remark} Even though the blow-up of the origin in $\mathbb{C}^2$ gives a geometric flattening of $\pi$, it is not quasi-proper and so, it is not a geometric $f$-flattening.

The following example shows that a quasi-proper map can behave badly in another way.

\parag{Example 2} We shall give an example of a reduced hypersurface $M$ in  $\C^{3} \times \C^{2}$ which has two irreducible components $M_{1}$ and  $M_{2}$ such that the natural projection $p: M \to \C^{3}$ is quasi-proper but such that the restrictions $p_{1}$ and $p_{2}$ of $p$ to $M_{1}$ and $M_{2}$ are respectively quasi-proper and not quasi-proper.\\

Let 
 $$M_{1} := \{(x, y, z, u, v) \in \C^{3} \times \C^{2} \ / \  y.v = z.u \} $$
 and 
  $$M_{2}:= \{ (x, y, z, u, v) \in \C^{3} \times \C^{2} \ / \  x.u^{2}+ y.v^{2} + z.v - u = 0 \} $$
  and define $M := M_{1}\cup M_{2}$. Let us begin by the description of the fibers of $p_{1}$ and $p_{2}$ as subsets of $\C^{2}$.\\
  
  The fiber of $p_{1}$ at a point $(x, y, z)$ when $(y, z) \not= (0, 0)$ is a line passing through the origin. When $y = z = 0$ the fiber is equal to $\C^{2}$. So all fibers are irreducible and contain the origin in $\C^{2}$. Hence they meet the $0$-section of $p_{1}$ and $p_{1}$ is quasi-proper.\\
  
  The fiber of $p_{2}$ at a point  $(x, y, z)$ when $(y, z) \not= (0, 0)$ and  $x \not= 0$ is an irreducible conic passing through the origin. This is also the case if  $x = 0$ and  $y \not= 0$. For $x = y = 0$ the fiber is a line passing through the origin of $\C^{2}$. When $y = z = 0$ and  $x \not= 0$ the fiber is a couple of lines, one through the origin ($u = 0$) but the second  one ($u = 1/x$) is going to infinity\footnote{This means  that this line avoids any given compact set if $x$ is near enough to $0$ but not equal to $0$.} when $x$ goes to $0$. This shows that $p_{2}$ is not quasi-proper at the origin.\\
  
  For the map $p$ the fiber at a point $(x, y, z)$ when $(y, z) \not= (0, 0)$ and $x \not= 0$ is the union of an irreducible conic passing through the origin and of a line passing through the origin. For $x = y = 0$ and  $z \not= 0$ the fiber is the union of two lines passing through the origin  $\{u = 0 \}$ and  $\{ z.v = u \}$. For $y = z = 0$ the fiber is $\C^{2}$. So each irreducible component of  a fiber of $p$ contains the origin and so $p$ is quasi-proper.
$\hfill \blacksquare$
  
\parag{Remark} If we replace, in the above example, $M_{2}$ by
 $$
 M'_{2} := \{ (x, y, z, u, v) \in \C^{3} \times \C^{2} \ / \  x.u.v+ y.v^{2} + z.v - u = 0 \} 
 $$
 and denote $p'_2\colon M'_2\rightarrow\C^3$ the restriction of $p$ we obtain an similar example but where all fibers of the projection $p'_{2} $ are connected. 
  \hfill {$\square$}
  
\section{Strongly quasi-proper maps}
  
  \subsection{Definition and a characterization of strongly quasi-proper maps} 
  
  As is shown by Example 1 above, the notion of quasi-proper map is not stable by base change (even by a modification of the target space) in presence of "big fibers". We shall introduce in this section a stronger notion, called {\em strongly quasi-proper} maps (in short {\em SQP maps}) which has better functorial properties and is equivalent to quasi-properness when the map is equidimensional. It will be characterized by the fact that its maximal reduced fiber map (see definition below)  is a meromorphic  map from $N$  to $\mathcal{C}_n^f(M)$.\\
  
Our setting is now the following : we consider  quasi-proper surjective holomorphic maps. The lemma below shows that they are always  $f$-GF over a  dense Zariski  open set  $N'$ in the target space $N$ which is assumed to be irreducible.

\begin{lemma}\label{standart situation}
Let $\pi\colon M \rightarrow N$ be a  holomorphic quasi-proper and surjective  map between a pure dimensional complex space $M$ and an irreducible complex space $N$.  Then there exists a (closed) analytic subset $\Sigma$  with empty interior in $N$ such that the map induced  by $\pi$, 
$$
\pi'\colon M \setminus \pi^{-1}(\Sigma) \longrightarrow N\setminus \Sigma
$$
is an f-GF map, which admits a reduced $f$-fiber map.
\end{lemma}

\parag{Proof} Set  $n := \dim M - \dim N$ and let $S$ be the analytic subset  of points $x $ in $M$ such that the dimension at $x$ of the fiber $\pi^{-1}(\pi(x))$ is strictly bigger than $n$. As this analytic subset is a union of irreducible components of the fibers of $\pi$ which is assumed to be quasi-proper, 's theorem implies that the image $\Sigma_0$ of $S$ is a closed analytic subset in $N$. Moreover, $\Sigma_0$ is of empty interior since $\dim\Sigma_0 < \dim N$. Let $\Sigma_1$ be the set of non normal points in $N$ and put $\Sigma := \Sigma_0 \cup \Sigma_1$. Then, $\pi'\colon M\setminus\pi^{-1}(\Sigma)\rightarrow N\setminus\Sigma$ is a quasi-proper $n-$equidimensionnal map and $N\setminus\Sigma$ is a normal space, so  it is an f-GF map  due to Lemma \ref{ map 2}).$\hfill \blacksquare$\\

In the situation above we say that the reduced $f$-fiber map for $\pi'$ is a {\bf reduced $f$-fiber map for $\pi$ on $N' = N\setminus \Sigma$}.

\parag{Remark}
If in Lemma \ref{standart situation} we suppose that $\pi$ is quasi-proper and dominant, then $\pi$ is surjective and $\pi^{-1}(\Sigma)$ is a nowhere dense analytic subset of $M$. Moreover,  if we have two dense Zariski open subsets $N_1$ and $N_2$ of $N$ and, on each one of them, a reduced $f$-fiber map for $\pi$, then  these two $f$-fiber maps coincide on the intersection $N_1\cap N_2$. Consequently  there exists a {\em largest} dense Zariski open subset $N'$ of $N$ on which we have  a (unique) reduced $f$-fiber map $\varphi$ for $\pi$. Moreover, every reduced $f$-fiber map for $\pi$ on a dense Zariski open subset of $N$ is a restriction of $\varphi$. We call $\varphi$ {\bf the maximal reduced  $f$-fiber map for $\pi$}. Thanks to Lemma \ref{reduced} the cycle $\varphi(y)$ is equal to the reduced fiber $\pi^{-1}(y)$ for $y$ generic in $N'$.
\hfill{$\square$}\\ 

Even though the restriction of a quasi-proper holomorphic map to an irreducible component is in general not a quasi-proper map (see Example 2 following Corollary \ref{cor 3}), we still have the following result.

\begin{lemma}\label{quasi-proper.component}
Let $\pi\colon M\rightarrow N$ be a quasi-proper  holomorphic map between a reduced complex space $M$ and an irreducible complex space $N$,  and let $C$ be an irreducible component of $M$. Then we have:
\begin{enumerate}[(i)]
\item
The restriction $\pi_{|C}\colon C\rightarrow N$ is semi-proper.
\item
Suppose moreover that $\pi_{|C}$ is dominant and put $n := \dim C - \dim N$. Then $\pi_{\vert C}$ is surjective and  there exists a dense open subset $N'$ of $N$ and, on $N'$  a holomorphic $f$-fiber map $\varphi\colon N'\rightarrow\mathcal{C}_{n}^{f}(\pi)^*$ for $\pi_{\vert C}$ such that $\varphi(y)$ is reduced for all $y$ in $N'$.
\end{enumerate}
\end{lemma}

\parag{Proof} To prove (i) it is enough to show, thanks to Proposition \ref{semi-proper 1bis} $(a)$, that $\pi_{|C}$ is semi-proper at every point in $\overline{\pi(C)}$, so let us fix a point $y_0$ in $\overline{\pi(C)}$. As $\pi$ is quasi-proper,  there exists an open neighborhood $V$ of $y_0$ in $N$ and a compact subset $L$ of $M$ which intersects every irreducible component of $\pi^{-1}(y)$ for all $y$ in $V\cap \pi(C)$. So it is sufficient to show that  the compact subset $L\cap C$ intersects $\pi^{-1}(y)$ for all $y$ in $V\cap \pi(C)$.
Now, for every point $y$ in $V$, there exists a sequence $(y_{\nu})_{\nu\geq 1}$ which converges to $y$ and such that $C$ contains an irreducible component of $\pi^{-1}(y_{\nu})$ for all $\nu\geq 1$. Consequently there exists a sequence $(x_{\nu})_{\nu\geq 1}$ in $C\cap L$ such that, for all $\nu\geq 1$, $x_{\nu}\in\pi^{-1}(y_{\nu})\cap C\cap L$ and, by taking a subsequence, we may assume that $(x_{\nu})_{\nu\geq 1}$ converges to a point $x$ in $L\cap C$. Hence $\pi^{-1}(y)\cap L\cap C\neq\emptyset$, since $\pi(x) = y$ by continuity, and the proof  of $(i)$ is completed.

\medskip
For the proof of (ii) we observe first that the singular part of $N$ is $b$-negligible so we may assume that $N$ is smooth. Hence it is enough to prove the result in the case where $N$ is an open subset of a numerical space, since a subset of $N$ is $b$-negligible if and only if it is locally $b$-negligible in $N$.

Let $S(M)$ denote the singular part of $M$ and $T$ be the set of points $y$ in $N$ such that $S(M)$ contains at least one irreducible component of $\pi^{-1}(y)$. Now, for a fixed point $y_0$ in $N$, there exists an open neighborhood $V$ of $y_0$  in $N$ and a compact subset $K$ of $M$ which meets every irreducible component of $\pi^{-1}(y)$, for all $y$ in $V$, since $\pi$ is quasi-proper. Then $T\cap V$ is a $b$-negligible subset of $V$, thanks to Corollary \ref{Kuhl-banach.4}. Hence $T$ is a $b$-negligible subset of $N$ and consequently $N\setminus T$ is a dense open subset of $N$. As $\pi_{|C}$ is dominant, it follows that $C' := C\cap(M\setminus\pi^{-1}(T))$ is a dense open subset of $C$ and, for each $y$ in $N\setminus T$, the fiber $\pi_{|C'}^{-1}(y)$ is the union of those irreducible components of $\pi^{-1}(y)$ which intersect $C$. Thus $\pi_{|C'}\colon C'\rightarrow N\setminus T$ is a quasi-proper map and, and by Lemma \ref{standart situation},  there exists a dense open subset $N'$ of $N\setminus T$ and a holomorphic map $\varphi\colon N'\rightarrow\mathcal{C}_{n}^{f}(\pi_{\vert C})^*$ such that $\varphi(y)$ is reduced for all $y$ in $N'$.
\hfill{$\blacksquare$}

\begin{defn}\label{sqp}
We say that $\pi\colon M\rightarrow N$ is a {\bf\em strongly quasi-proper map} (an {\bf\em SQP map} for short) if  the following conditions are satisfied:
\begin{enumerate}[(i)]
\item
$M$ is a reduced complex space of pure dimension, $N$ is an irreducible complex space and $\pi$ is a holomorphic quasi-proper dominant map.
\item
The closure  in $N\times\mathcal{C}_n^f(M)$ of the graph of a reduced fibermap for $\pi$ over a dense Zariski  open subset $N'$ is proper over $N$. 
\end{enumerate}
\end{defn}

\parag{Remark} Suppose $\pi\colon M\rightarrow N$ is a map, which satisfies condition (i) of Definition \ref{sqp}, and let $\varphi\colon N'\rightarrow\mathcal{C}_n^f(M)$ be the reduced   fibermap for $\pi$ on $N'$. Denote  respectively $\Gamma$ the closure of the graph of $\varphi$ in $N\times\mathcal{C}_n^f(M)$ and $\overline{\varphi(N')}$ the closure of $\varphi(N')$ in $\mathcal{C}_n^f(M)$. As $\varphi$ induces a holomorphic section of the natural holomorphic map $\alpha\colon\mathcal{C}_n^f(\pi)^*\rightarrow N$, the canonical projection $N\times\mathcal{C}_n^f(M)\rightarrow\mathcal{C}_n^f(M)$ induces a homeomorphism from $\Gamma$ to $\overline{\varphi(N')}$. It follows that $\pi$ is an SQP map if and only if the map $\overline{\varphi(N')}\rightarrow N$ induced by $\alpha$ is proper.
\hfill{$\square$}\\

The following proposition gives a characterization of SQP maps. It is an improvement of the criterium given in \cite{[B.13]} for a holomorphic map  to be an SQP map.

\begin{prop}\label{SPQ dense}
Let $\pi\colon M \rightarrow N$ be a  dominant holomorphic map between a pure dimensional complex space $M$ and an irreducible complex space $N$. Define  $n := \dim M - \dim N$. Assume that there exists a  dense subset $\Lambda$ in $N$ such that for each $y \in \Lambda$ the fiber$\pi^{-1}(y)$ is non empty, and of pure dimension $n$ with finitely many irreducible components. Let $\gamma\colon\Lambda \rightarrow \mathcal{C}_{n}^{f}(M)$ be the map which associates to every $y$ in $\Lambda$ the reduced $n$-cycle $\pi^{-1}(y)$. Let $\Gamma$ be the graph of the map $\gamma$ and $\bar \Gamma$ be the closure of $\Gamma$ in $N\times \mathcal{C}_{n}^{f}(M)$. Our main assumption is now the following:
\begin{itemize}
\item 
The natural projection $\tau\colon\bar \Gamma \rightarrow N$ is proper.
\end{itemize}
Then the map $\pi$ is strongly quasi-proper.
\end{prop}

\parag{Proof} 
As $\tau(\bar \Gamma)$ is closed and contains $\Lambda$ we have $\tau(\bar \Gamma) = N$. Now recall that the singleton $\{\emptyset[n]\}$ is open (and closed) in $\mathcal{C}_{n}^{f}(M)$. Then by Proposition \ref{relatif} the set 
$$
\left\{(y,\xi)\in N\times\mathcal{C}_{n}^{f}(\pi)^* \ /\ \ |\xi| \subseteq \pi^{-1}(y) \right\}
$$ 
is closed in $N\times\mathcal{C}_{n}^{f}(M)$ and contains $\Gamma$. Hence it contains $\bar \Gamma$ also and it follows that $\pi$ is surjective.

Our second  step of the proof (which is infact the main step) we are going to show  that the map $\pi$ is quasi-proper. To do so let $p\colon N\times\mathcal{C}_{n}^{f}(\pi)\rightarrow \mathcal{C}_{n}^{f}(\pi)$ be the natural projection, let $y$ be an arbitrary point in $N$ and let $V$ be an open relatively compact neighborhood of $y$ in $N$. Fix $y' \in V$ and  choose an irreducible component $C$ of $\pi^{-1}(y')$\footnote{From the surjectivity of $\pi$ proved above, $\pi^{-1}(y')$ is not empty.}. Let $x'$ be a point in $C$ such that $x'$  does not belong to any other irreducible component of $\pi^{-1}(y')$. Then, as $\pi$ is dominant, $\pi^{-1}(\Lambda)$ is dense in $M$ (see  Lemma \ref{dominant}) and we can choose a sequence $(x_{\nu})_{\nu \geq 0}$ in $\pi^{-1}(\Lambda)$ converging to $x'$. For $\nu \gg 1$ we have $\pi(x_{\nu}) \in V$ so the cycles $\gamma(\pi(x_{\nu}))$ are in the compact subset $p(\tau^{-1}(\bar V))$ of $\mathcal{C}_{n}^{f}(\pi)$. 
By taking a subsequence, we can assume that the sequence $(\gamma(\pi(x_{\nu})))_{\nu\geq 0}$ converges to a cycle $\delta$ in 
$\mathcal{C}_{n}^{f}(\pi)$. As we have $x_{\nu}\in \gamma(\pi(x_{\nu}))$ for each $\nu$ we have $x' \in \vert \delta\vert$. Since $p(\tau^{-1}(\bar V))$  is a compact subset of \ $\mathcal{C}_{n}^{f}(M)$ there exists a compact subset $K$ in $M$ such that each irreducible component of every cycle in $p(\tau^{-1}(\bar V))$ meets $K$. So this is the case for each irreducible component of the cycle $\gamma(\pi(x_{\nu}))$, for every $\nu$, and for every irreducible component  of  $|\delta|$. Let $\delta_{0}$ be an irreducible component of $\delta$ containing $x'$. Then $\delta_{0}$ is contained in $C$ since $C$ is the only irreducible component of $\pi^{-1}(y')$ which contains $x'$. As $\delta_{0}$ meets $K$ so does $C$.

Thus we have proved that for all $y \in N$ there exists an open neighborhood  $V$ of $y$ in $N$ and a compact set $K$ in $M$ such that for every $y' \in V$ and any irreducible component $C$ of $\pi^{-1}(y')$ the intersection $C \cap K$ is not empty. This means that the map  $\pi$ is quasi-proper.

To prove that $\pi$ is strongly quasi-proper let 
$$ 
\varphi\colon N' \longrightarrow \mathcal{C}_{n}^{f}(M) 
$$
be the reduced fibermap for $\pi$ over $N'$. Then $\varphi(y)$ is reduced for generic $y$ in $N'$. It follows that $\varphi$ and $\gamma$ coincide on a dense subset of $N$. Hence the closure of the graph of $\varphi$ in $N \times \mathcal{C}_{n}^{f}(M) $ is contained in $\bar \Gamma$. Then, by definition,  the map $\pi$ is strongly quasi-proper. 
$\hfill \blacksquare$

\subsection{Basic properties of SQP-maps}

\begin{prop}\label{fat.fibers}
Let $\pi\colon M\rightarrow N$ be an SQP map, $\varphi\colon N'\rightarrow \mathcal{C}_{n}^{f}(\pi)$ be its reduced fibermap over $N'$ and $\Gamma$ be the closure of the graph of $\varphi$ in $N\times \mathcal{C}_{n}^{f}(M)$. Then  we have
$$
\bigcup_{\xi\in \Gamma_y}|\xi| = \pi^{-1}(y)
$$
for all $y$ in $N$\footnote{Recall that the equality $(\{y\}\times \mathcal{C}_{n}^{f}(M))\cap\Gamma = \{y\}\times\Gamma_y$ gives the definition of $\Gamma_y$.}. 
\end{prop}
\parag{Proof} For $y$ in $N'$ the assertion is obvious; so suppose $y\in N\setminus N'$ and $x\in\pi^{-1}(y)$. Then there exists a sequence $(x_\nu)$ in $\pi^{-1}(N')$ which converges to $x$ because $\pi$ is a dominant map. Since $\Gamma$ is proper over $N$ we may assume, by taking a subsequence if necessary, that the the sequence $(\varphi(\pi(x_\nu))$ of $n-$cycles in $M$ converges to an $n-$cycle $\xi\in \Gamma_y$. Then we have $x\in|\xi|$.  The other inclusion is obvious.\hfill{$\blacksquare$}

\begin{prop}\label{characterization.SQP}
Let $\pi\colon M\rightarrow N$ be a holomorphic map from a pure dimensional reduced complex space to an irreducible complex space which admits a fibermap $\varphi\colon N'\rightarrow\mathcal{C}_n^f(\pi)^*$, where $N'$ is a dense open set and  where $\varphi(y)$ is reduced for all $y$ in $N'$.
Let $\Gamma$ denote the closure of the graph of $\varphi$ in $N\times\mathcal{C}_n^{\rm loc}(M)$ and let 
$\tau\colon\Gamma\rightarrow N$ and $p\colon\Gamma\rightarrow\mathcal{C}_n^{\rm loc}(M)$ denote the natural projections. Then $\pi$ is an SQP map if and only if it satisfies the following condition:
\begin{itemize}
\item[\rm ($*$)]
For every compact subset $K$ of $N$ there exists a compact subset $L$ of $M$ which intersects every irreducible component of every $\xi$ in $\Gamma_y$ for all $y$ in $K$.
\end{itemize}
\end{prop}

\parag{Proof} Let $\Gamma^f$ denote the closure of the graph of $\varphi$ in  $N\times\mathcal{C}_n^{f}(M)$.\\
Suppose that $\pi$ is an SQP map.  Then $\Gamma^f$ is a closed subset of $\Gamma$ since $\Gamma^f$ is proper over $N$ and consequently $\Gamma_0 = \Gamma$. Now, let $K$ be a compact subset of $N$. As $\tau$ is a proper map  the subspace $p(\tau^{-1}(K))  = \cup_{y\in K}\Gamma_y $ of $\mathcal{C}_n^f(M)$ is compact and, due to Corollary \ref{compact}, there exists a compact subset $L$ of $M$ which intersects every irreducible component of every $\xi$ in 
$\Gamma_y$ for all $y$ in $K$. Hence ($*$) is satisfied. 

\smallskip
Conversely, suppose that condition ($*$) is satisfied and fix a compact subset $K$ of $N$. Let us show first that $p(\tau^{-1}(K))$ is a compact subset of $\mathcal{C}_n^{\rm loc}(M)$. We observe that $p(\tau^{-1}(K))$ is a closed subset of $\mathcal{C}_n^{\rm loc}(M)$, since the natural projection
 $$K\times\mathcal{C}_n^{\rm loc}(M)\rightarrow\mathcal{C}_n^{\rm loc}(M)$$
 is a closed map, so it is enough to show that $p(\tau^{-1}(K))$ is a relatively compact subset of $\mathcal{C}_n^{\rm loc}(M)$. By Theorem \ref{BLV 2} the set  $p(\tau^{-1}(K))$ is a relatively  compact subset of $\mathcal{C}_n^{\rm loc}(M)$ if for every continuous positive definite $(1,1)-$form $\omega$ on $M$ and every relatively compact open subset $M'$ of $M$  the function
\begin{equation*}\tag{$@$}
 \xi \mapsto \vol_{\omega}(\xi \cap M') := \int_{\xi\cap M'} \omega^{\wedge n} 
\end{equation*}
 is bounded on $p(\tau^{-1}(K))$. Now let us fix such  $\omega$ and $M'$ and show that the function $v$,  defined by $(@)$, is bounded on $p(\tau^{-1}(K))$. To do so we take a relatively compact open neighborhood $W$ of $K$ in $N$ and recall that, by Proposition 4.2.17 in \cite{[e]}, $v$ is continuous on  $\mathcal{C}_n^{\rm loc}(M)$  and moreover bounded on $p(\tau^{-1}(N'\cap W))$ by Theorem 3.6.6 in \cite{[e]}. As $\tau^{-1}(N'\cap W)$ is dense in $\tau^{-1}(\bar W)$ it follows that $v$ is bounded on $p(\tau^{-1}(\bar W))$ and consequently also on $p(\tau^{-1}(K))$.

It then follows  from ($*$) and Corollary \ref{compact},  that $p(\tau^{-1}(K))$ is a compact subset of $\mathcal{C}_n^{f}(M)$. This implies that the subset $\Gamma^f$ of $N\times\mathcal{C}_n^{f}(M)$ is equal to $\Gamma$ and is  proper over $N$. Hence $\pi$ is an SQP map.
\hfill{$\blacksquare$}

\parag{Example} In Example 1 following Corollary \ref{cor 3} the mapping $\pi\colon Y\rightarrow \C^2$ is quasi-proper and has a reduced fibermap $\varphi$ on $\C^2\setminus\{(0,0)\}$. Let $\Gamma$ denote the closure of its graph in $\C^2\times\mathcal{C}_1^f(Y)$ and let  $s\neq 0$ be a complex number. Then $\lim\limits_{t\to 0}\varphi(t,s.t)$ is the union of two lines given by $x = 0$ and $x = -\frac 1s$ \ in $\Gamma_{(0,0)}$. It follows that  no compact subset of $Y$ intersects every irreducible component of every $\xi$ in $\Gamma_{(0,0)}$, so the map $\pi$ is not strongly quasi-proper.
\hfill{$\square$}\\

The following result shows that strong quasi-properness of a map can be tested with any of its fiber maps, not only the reduced one.

\begin{cor}\label{indep.fiber.map}
Let $\pi\colon M\rightarrow N$ be a quasi-proper and dominant holomorphic map from a pure dimensional reduced complex space to an irreducible complex space and $\psi\colon N'\rightarrow\mathcal{C}_n^f(\pi)$ be a fiber map for $\pi$ (reduced or not) on an open dense subset $N'$ of $N$. Denote  $\Gamma_\psi$ the closure in $N\times\mathcal{C}_n^f(M)$ of the graph of $\psi$. If $\Gamma_\psi$ is proper over $N$, then $\pi$ is an SQP map.
\end{cor}
\parag{Proof} Let $\varphi$ denote the maximal reduced fiber map for $\pi$ and let $N''$ be a dense open subset of $N$ where both $\varphi$ and $\psi$ are defined and holomorphic. Let $\Gamma_\varphi$ denote the closure in $N\times\mathcal{C}_n^f(M)$ of the graph of $\varphi$. Let
$\tau_\psi\colon\Gamma_\psi\rightarrow N$ and $p\colon N\times\mathcal{C}_n^f(M)\rightarrow\mathcal{C}_n^f(M)$ denote the natural projections.  Notice first that we have:
\begin{enumerate}
\item[($*$)]
Every $y$ in $N''$ satisfies the inequality $\varphi(y)\leq\psi(y)$.
\end{enumerate}
Let $K$ be a compact subset of $N$. Then $\tau_\psi^{-1}(K)$ is compact and it follows from Corollary \ref{compact} that there exists a compact subset $L$ of $M$ which intersects every irreducible component of every cycle in $(\Gamma_\psi)_y$ for all $y$ in $K$. Now, for a fixed point $y$ in $K$ and a fixed cycle $\eta$ in $(\Gamma_\varphi)_y$, there exists a sequence $(y_\nu)$ in $N''$ which converges to $y$ and such that $\varphi(y_\nu)\to\eta$. As $\tau_\psi$ is proper we may assume, by taking a subsequence, that $(\psi(y_\nu))$ converges to a cycle $\xi$ in $(\Gamma_\psi)_y$. Thanks to ($*$) we then have $\eta\leq\xi$ and it follows that each irreducible component of $\eta$ is an irreducible component of $\xi$. Hence $L$ intersects every irreducible component of $\eta$ and consequently $\pi$ is an SQP map by Proposition \ref{characterization.SQP}.
\hfill{$\blacksquare$}

\begin{cor}
Let $\pi\colon M \rightarrow N$ be an SQP map and put $n := \dim M - \dim N$. Suppose $T$ is an irreducible subspace of $N$ which is not contained in $\pi(\Sigma_{n+1}(\pi))$\footnote{Recall that 
$\Sigma_{n+1}(\pi)$ is a union of fibers of $\pi$,  so its image by $\pi$, which is quasi-proper, is a closed analytic subset of $N$.}. Then the induced map $\pi^{-1}(T)\rightarrow T$ is an SQP map. 
\end{cor}
\parag{Proof} Let $\varphi\colon N'\rightarrow\mathcal{C}_n^f(M)$ be the reduced fibermap for $\pi$ and let $\Gamma$ denote the closure of its graph. The restriction of $\varphi$ to $T\setminus \pi(\Sigma_{n+1}(\pi))$ is then a (not necessarily reduced) fibermap for the induced map $\pi^{-1}(T)\rightarrow T$ and the closure of its graph is a closed subset of $\Gamma$. Thus $\pi^{-1}(T)\rightarrow T$ is an SQP map thanks to Corollary \ref{indep.fiber.map}.
\hfill{$\blacksquare$}\\

We would like to point out that Theorem \ref{general} below is a generalization of the above corollary to the case where $T$ is contained in the "big fibers locus" $\pi(\Sigma_{n+1}(\pi))$.

\begin{cor}\label{SQP.componentwise}
Let $\pi\colon M\rightarrow N$ be a dominant holomorphic map from a pure dimensional reduced complex space to an irreducible complex space. Suppose also that $M$ has only finitely many irreducible components. Then $\pi$ is an SQP map if and only if the restricion of $\pi$ to each irreducible component of $M$ is an SQP map.
\end{cor}
\parag{Proof} Put $n := \dim M -\dim N$ and let $M = M_1\cup\cdots\cup M_k$ be the decomposition of $M$ into irreducible components. 

Suppose first that $\pi$ is an SQP map. Then, thanks to (ii) of Lemma \ref{quasi-proper.component}, there exists a dense open subset $N'$ of $N$ and, for each $j \in [1, k]$, a holomorphic map  $\varphi_j\colon N'\rightarrow\mathcal{C}_n^f(M_j)$ such that $\varphi_j(y) = \pi_{|M_j}^{-1}(y)$ for all $y$ in $N'$. Hence 
$$
\varphi\colon N'\longrightarrow\mathcal{C}_n^f(M),\quad y\mapsto \varphi_1(y) + \cdots + \varphi_k(y)
$$
is a reduced fibermap for $\pi$.  Denote $\Gamma_1,\ldots,\Gamma_k$ and $\Gamma$  the closures in $N\times\mathcal{C}_n^f(M)$ of the graphs of $\varphi_1,\ldots,\varphi_k$ and $\varphi$. Let $K$ be any compact subset of $N$. If $\pi$ is an SQP map, then  there exists a compact subset $L$ of $M$ which intersects every irreducible component of every cycle in $\Gamma_y$ for all $y$ in $K$. Now, take a point $y$ in $K$ and a  cycle $\eta$ in $(\Gamma_j)_y$ for some $j$. Then there exists a sequence $(y_\nu)$ in $N'$ such that $(y_\nu,\varphi_j(y_\nu))$ converges to $(y,\eta)$. By taking a subsequence we may assume, since $\Gamma$ is proper over $N$, that $(\varphi(y_\nu))$ converges to a cycle $\xi$ in $\Gamma_y$. Then $\eta\leq\xi$ and it follows that $L$ cuts every irreducible component of $\eta$. Hence we have proved that the restrictions $\pi_{|M_1},\ldots,\pi_{|M_k}$ are all SQP maps. \\
The converse is proved in a similar way. \hfill{$\blacksquare$}\\

\section{Stability properties of SQP maps}

We begin by showing that SQP maps are stable by strict transform in any base change.

\begin{thm}\label{general}
Let $\pi\colon M \rightarrow N$ an SQP map and let $g\colon Z \rightarrow N$ be a holomorphic map where $Z$ is an irreducible complex space. Let $\tilde{\pi}\colon\tilde{Z} \rightarrow Z$ be the strict transform of $\pi$ by $g$ and consider the decomposition $\tilde{Z} = \cup_{j}\tilde{Z}^{(j)}$ where $\tilde{Z}^{(j)}$ is the union of all $j$-dimensional  irreducible components of $\tilde{Z}$.  Then, for each $j$ such that $\tilde{Z}^{(j)} \not= \emptyset$, the map $\tilde\pi_j\colon\tilde{Z}^{(j)}\rightarrow Z$, induced by $\pi$, is an SQP map.
\end{thm}

\parag{Proof}
Let us first prove the result in the case where $Z$ is a subspace of $N$ and $g$ is the natural injection. Then, by (i) of Lemma \ref{quasi-proper.component}, $\tilde{Z}$ is the union of those irreducible components of $\pi^{-1}(Z)$ which are mapped surjectively onto $Z$.

Consider a fixed  $j$ such that $\tilde{Z}^{(j)}\neq\emptyset$ and put $q := j - \dim Z$. As $\tilde{Z}^{(j)}$ has only finitely many irreducible components there exists, due to (ii) of Lemma \ref{quasi-proper.component}, a dense open subset $Z'$ of $Z$ and a holomorphic fibermap $\varphi_j\colon Z'\rightarrow\mathcal{C}_{q}^{f}(\tilde{\pi}_j)^*$ for $\tilde{\pi}_j$ over $Z'$ such that $\varphi_j(y)$ is reduced for all $y$ in $Z'$. Denote $\Gamma_j$ the closure of the graph of $\varphi_j$ in $Z\times\mathcal{C}_{q}^{f}(\tilde{Z}^{(j)})$. Then, thanks to Proposition \ref{characterization.SQP}, it is enough to prove that for any compact subset $K$ of $Z$ there exists a compact set $L_j$  of $\tilde{Z}^{(j)}$ having the following property:
\begin{itemize}
\item[]
If $(y_{\nu})_{\nu}$ is a sequence in $Z'$ which converges to a point $y_0$ in $K$ and such that $\varphi_j(y_{\nu})$ converges to a $q-$cycle $\xi$ in  $\mathcal{C}_{q}^{f}(\tilde{Z}^{(j)})$, then every irreducible component of $\xi$ meets $L_j$.
\end{itemize} 
Before proving this we put $n := \dim M - \dim N$ and denote $\varphi$ the reduced fibermap for $\pi$ over a dense open subset  $N'$ of $N$  and $\Gamma$ the closure of its graph in $N\times\mathcal{C}_{n}^{f}(M)$. Now let $K$ be a compact subset of $Z$. Then  there exists a compact subset $L$  of $M$ which,  for all $y\in K$,  intersects every irreducible component of every $n$-cycle $\beta\in \Gamma_y$, because $\pi$ is an SQP map.

Let $(y_{\nu})_{\nu}$ be a sequence in $Z'$, which converges to a point $y_0$ in $K$ and such that $\varphi_j(y_{\nu})$ converges to a $q-$cycle $\xi$ in  $\mathcal{C}_{q}^{f}(\tilde{Z}^{(j)})$, and $\xi_0$ be any irreducible component of $\xi$. Pick a point $x$ in $\xi_0$ which does not belong to any other irreducible component of $\xi$ and  choose, for each $\nu$, a point $x_\nu$  in $|\varphi_j(y_\nu)|$ such that $\lim_{\nu\to\infty}x_\nu = x$. 
Then there exists, for each $\nu$,   an $n$-cycle $\beta_\nu\in\Gamma_{y_\nu}$  which contains $x_{\nu}$, thanks to Proposition \ref{fat.fibers}. 
As $\Gamma$ is proper over $N$ we may assume, by taking a subsequence,  that the sequence $(\beta_{\nu})$ converges in  $\mathcal{C}_{n}^{f}(M)$ to an $n$-cycle $\beta$.  Moreover every irreducible component of $\beta$ intersects $L$. 

For each $\nu$ let $\gamma_\nu$ be an irreducible component of $\beta_\nu$ which contains $x_\nu$. Then we have  $\gamma_\nu\subseteq |\varphi_j(y_\nu)| = \tilde{\pi}_j^{-1}(y_\nu)$ for all $\nu$ and, due to Corollary \ref{compact},  we may assume, by taking a subsequence, that $(\gamma_\nu)_\nu$ converges in  $\mathcal{C}_{n}^{f}(M)$ to an $n$-cycle $\delta$ such that $\delta\leq\beta$. Now, for all $\nu$, we have $x_\nu\in\gamma_\nu\subseteq |\varphi(y_\nu)|$ and hence $x \in|\delta|\subseteq |\xi|$, because 
$(\varphi(y_\nu))_\nu$ converges to $\xi$ in $\mathcal{C}_{q}^{f}(\tilde{Z}^{(j)})$.\\
Let $\delta_0$ be an irreducible component of $|\delta|$ which contains $x$. Then $\delta_0\subseteq \xi_0$ and it follows that 
$\emptyset\neq\delta_0\cap L \subseteq \xi_0\cap L$, so we can set $L_j := L\cap\tilde{Z}^{(j)}$.

\medskip
The general case is now easily obtained from the case where $Z$ is a subspace of $N$. Indeed, the map $\id_{Z}\times \pi\colon Z\times M \rightarrow Z \times N$ is the strict transform of $\pi$ by the natural projection $Z \times N\rightarrow N$ and it is clearly an SQP map.  Then notice that  we can factorize $g$ by the canonical inclusion of its graph $G \hookrightarrow Z \times N$ for which the above case gives the result.
\hfill{$\blacksquare$}

\begin{lemma}\label{modif.quasiproper.1}
Let $\pi\colon M\rightarrow N$ be a quasi-proper and dominant holomorphic map between reduced complex spaces where $M$ is of pure dimension and $N$ is irreducible. Let $\tau\colon\tilde{M}\rightarrow M$ be a modification whose center is $\pi$-proper. Then there exists a nowhere dense analytic subset $\Sigma$ of $N$ and a reduced fibermap $\varphi$ for $\pi\circ\tau$ on $N\setminus\Sigma$ such that $\tau_*\circ\varphi$ is a reduced fibermap for $\pi$ on $N\setminus\Sigma$, where $\tau_*$ is the direct image map by $\tau$.
\end{lemma}

\parag{Proof} Set $n := \dim M - \dim N$ and let $C$ denote the center of $\tau$. As  $C$ has empty interior in $M$, $\dim C < n + \dim N$ and $\dim \tau^{-1}(C) < n + \dim N$ we have, thanks to Proposition 2.4.60 in \cite{[e]} and Kuhlmann's Theorem, that $\pi\left(\Sigma_n(\pi_{|C})\right)$ and $(\pi\circ\tau)\left(\Sigma_n(\pi\circ\tau_{|\tau^{-1}(C)})\right)$ are nowhere dense analytic subsets of $N$. Let $\Sigma'$ denote their union. Then, for every $y$ in $N\setminus\Sigma'$, no irreducible component of $\pi^{-1}(y)$ is contained in $C$ and no irreducible component of $\tau^{-1}(\pi^{-1}(y))$ is contained in $\tau^{-1}(C)$. Now, by Lemma \ref{modif.quasiproper.0}, $\pi\circ\tau$ is quasi-proper and consequently there exists a nowhere dense analytic subset $\Sigma''$ of $N$ and a reduced fibermap $\psi$ for $\pi\circ\tau$ on $N\setminus\Sigma''$. Put $\Sigma := \Sigma'\cup\Sigma''$ and let $\varphi$ denote the restriction of $\psi$ to $N\setminus\Sigma$. It follows that, for every $y$ in  $N\setminus\Sigma$, the map $\tau^{-1}(\pi^{-1}(y))\rightarrow \pi^{-1}(y)$ induced by $\tau$ is a modification and consequently $\tau_*(\varphi(y))$ is reduced if $\varphi(y)$ is reduced. Hence $\tau_*\circ\varphi\colon N\setminus\Sigma\rightarrow\mathcal{C}_n^f(\tilde M)$ is a reduced fibermap for $\pi$ since $\tau_*\colon\mathcal{C}_n^f(\tilde{M})\rightarrow\mathcal{C}_n^f(M)$ is a holomorphic map (see Theorem \ref{direct image}).
\hfill{$\blacksquare$}

\begin{thm}\label{modific.quasiproper}
Let $\pi\colon M\rightarrow N$ be a holomorphic map between reduced complex spaces where $M$ is of pure dimension and $N$ is irreducible. Let $\tau\colon\tilde{M}\rightarrow M$ be a modification whose center is $\pi$-proper. Then $\pi$ is an SQP map if and only if $\pi\circ\tau$ is an SQP map.
\end{thm}
\parag{Proof} Set $n := \dim M - \dim N$ and let $C$ denote the center of $\tau$. From Lemma \ref{modif.quasiproper.0} we have that $\pi$ is quasi-proper if and only if $\pi\circ\tau$ is quasi-proper so we may assume that $\pi$ is quasi-proper. Now let $\Sigma$ and $\varphi$ be as in Lemma \ref{modif.quasiproper.1} and let $\Gamma$, $\tilde\Gamma$ denote respectively  the closures of the graphs of $\tau_*\circ\varphi$ in $N\times\mathcal{C}_n^f(M)$ and $\varphi$ in $N\times\mathcal{C}_n^f(\tilde M)$. 

Suppose first that $\pi$ is an SQP map and let us show that $\pi\circ\tau$  is an SQP map. To do so we consider a compact subset $K$ of $N$. Then there exists a compact subset $L$ of $M$ which intersects every irreducible component of every $\gamma$ in $\Gamma_y$ for all $y$ in $K$. We shall show that the compact subset  $\tilde L := \tau^{-1}(L\cup(\pi^{-1}(K)\cap C))$ of $\tilde M$ intersects every irreducible component of every $\gamma$ in $\tilde\Gamma_y$ for all $y \in K$. So let's take a point $y$ in $K$ and an $n$-cycle $\gamma$ in $\tilde\Gamma_y$. Then $\tau_*(\gamma)$ is in $\Gamma_y$ since $\tau_*$ is continuous. Now let $\delta$ be an irreducible component of $\gamma$  which does not intersect $\tau^{-1}(\pi^{-1}(K)\cap C)$. Then $\delta$ does not intersect $\tau^{-1}(C)$ and consequently $\tau$ maps $\delta$ bi-holomorphically onto $\tau(\delta)$. Thus $\tau(\delta)$ is an irreducible component of $\tau_*(\gamma)$ and it follows that $\tau(\delta)\cap L\neq \emptyset$. Hence $\emptyset\neq\delta\cap\tau^{-1}(L) = \delta\cap\tilde L$.

\smallskip
Conversely, suppose that $\pi\circ\tau$  is an SQP map and consider a compact subset $K$ of $N$. Then there exists a compact subset $\tilde L$ of $\tilde M$ which intersects every irreducible component of every $n$-cycle in $\tilde\Gamma_y$ for all $y$ in $K$. Set $L := \tau(\tilde L)\cup(\pi^{-1}(K)\cap C)$. Take a point $y$ in $K$ and an $n$-cycle $\gamma$ in $\Gamma_y$. Then there exists a sequence $(y_\nu)_{\nu\geq 0}$ in $N\setminus\Sigma$ such that the sequence $(\tau_*(\varphi(y_\nu)))_{\nu\geq 0}$ converges to $\gamma$ in $\mathcal{C}_n^f(M)$. Since $\pi\circ\tau$ is strongly quasi-proper we may assume, by taking a subsequence, that the sequence $\varphi(y_\nu))_{\nu\geq 0}$ converges to an $n$-cycle $\tilde\gamma$ in $\tilde\Gamma_y$. Then by continuity we get $\tau_*(\tilde\gamma) = \gamma$ and with the same arguments as above it is clear that $L$ intersects every irreducible component of $\gamma$.
\hfill{$\blacksquare$}

\begin{prop}\label{SQP.compos.fGF}
Let $\pi_1 : M_1\rightarrow N$ be an SQP map and  $\pi_2 : M_2\rightarrow M_1$ \ be an $f$-GF map between reduced complex spaces. Then the composition $\pi_1\circ \pi_2$ is an SQP map.
\end{prop}
\parag{Proof} Denote respectively  $\Gamma_1$ and $\Gamma$ the closures of the graphs of the reduced fibermaps, $\varphi_1$ and $\varphi$, for $\pi_1$ and $\pi_1\circ\pi_2$. By Lemma \ref{open.quasipropre} the map $\pi_1\circ\pi_2$ is quasi-proper, so, thanks to Proposition \ref{characterization.SQP}, it is sufficient to show that, for every compact subset $K$ of $N$, there exists a compact subset $L_2$ of $M_2$ which intersects every irreducible component of every cycle in $\Gamma_y$ for all $y$ in $K$. Fix a compact subset $K$ of $N$. Then there exist a compact subset $L_1$ of $M_1$ such that ${\rm int}(L_1)$ intersects every irreducible component of every cycle in $(\Gamma_1)_y$ for all $y$ in $K$ and a compact subset $L_2$ of $M_2$ which intersects every irreducible component of $\pi_2^{-1}(x)$
for all $x$ in $L_1$. We are going to show that every irreducible component of every cycle in $\Gamma_y$ meets $L_2$ for all $y$ in $K$. Consider a cycle $\gamma$ in $\Gamma_y$. In a dense Zariski open subset of $N$ where both $\varphi_1$ and $\varphi$ are defined and holomorphic  we can find a sequence $(y_\nu)_{\nu\geq 0}$ which converges to $y$ and such that $\varphi_(y_\nu)$ converges  $\gamma$. By taking a subsequence we may also suppose that $\varphi_1(y_\nu)$ converges to a cycle $\beta$ in $(\Gamma_1)_y$. Then $|\gamma|\subseteq\pi_2^{-1}(|\beta|)$. Moreover  $|\gamma|$ and $\pi_2^{-1}(|\beta|)$ are of the same pure dimension so $|\gamma|$ is a union of irreducible components of $\pi_2^{-1}(|\beta|)$. As the induced map $\pi_2^{-1}(|\beta|)\rightarrow|\beta|$ is open and quasi-proper it maps each irreducible component of $\gamma$  onto an irreducible component of $|\beta|$. Hence for an irreducible component $\delta$ of $\gamma$ we get $\pi_2(\delta)\cap {\rm int}(L_1)\neq\emptyset$, so there exists a point $x$ in $\pi_1^{-1}(y)$ such that $\delta$ contains an irreducible component of $\pi_2^{-1}(x)$. It follows that $\delta$ \ intersects $L_2$.\hfill{$\blacksquare$}\\

Note that the composition of any two SQP maps is not an SQP map in general. Even in the simple case where $\pi : M \to N$ is an $f$-GF map and $\tau : \tilde{M} \to M$ is a modification, the composition $\pi\circ \tau$ is not necessarily an SQP map, without  the assumption that the center of $\tau$ is $\pi$-proper. This is the content of the following example.

 In fact a proper modification of $M$ for a f-GF map $\pi$  with center $C$ not $\pi$-proper is enough to give a counter-example, as follows:\\

\parag{Example}  Put $M := \mathbb{C}^2$,  $N := \mathbb{C}$  and  let $\pi: M \to N$ be the projection $\pi(x, y) = x$. Let $\tau : \tilde{M} \to M$ be the blow-up of $M$ at each point of $\{0\}\times \mathbb{Z}$. Then $\tau$ is a (proper) modification and $\pi$ is clearly a $f$-GF map. But $\pi\circ\tau : \tilde{M} \to N$ is not an SQP map  because its fiber at the origin is not a finite type cycle.
 So $\pi\circ\tau$ is even  not quasi-proper (but  it is equidimensional ! )\\

\section{Further characterizations of SQP maps}

\parag{Notation}
In this section we use the following notation for a continuous map  $\pi\colon M\rightarrow N$. If $V$ is an open subset of $N$ we write $M_V$ instead of $\pi^{-1}(V)$ and we denote $\pi_V\colon M_V\rightarrow V$ the map induced by $\pi$.\hfill{$\square$}\\

 The main result of this section is the following  characterization of SQP maps in terms of $f$-GF maps  which is a variant of  Theorem 2.4.4 of \cite{[B.13]} (see also   \cite{[Mt.00]}).

\begin{thm}\label{SQP et GF}
Let $\pi\colon M \rightarrow N$ be a quasi-proper holomorphic dominant  map between reduced  complex spaces where $M$ is of pure dimension and $N$ is irreducible. Then the following three conditions are equivalent.
\begin{enumerate}[(i)]
\item
The map $\pi$ is strongly quasi-proper.
\item
The maximal reduced fiber map of $\pi$ is meromorphic map from $N$ to  $\mathcal{C}_n^f(M)$\footnote{See Definition \ref{mero}.}).
\item
The map $\pi$ admits an f-flattening.
\end{enumerate}
\end{thm}

The proof is given in Subsection 7.2 below.

\subsection{The local $f$-flattening theorem}

\begin{lemma}\label{induction step}
Let $\pi\colon M\rightarrow N$ be a surjective holomorphic map between reduced complex spaces of pure dimensions and put $n := \dim M - \dim N$. Let $y_0$ be a point in $N$ such that the fiber 
$\pi^{-1}(y_0)$ has only finitely many irreducible components and $\dim\pi^{-1}(y_0) > n$. Then there exists an open neighborhood $V$ of $y$ and a modification $\tau \colon\tilde{V} \to V$ such that the strict transform $\tilde{\pi}_V\colon\tilde{M}_{\tilde{V}}\rightarrow\tilde{V}$ of $\pi_V$ by $\tau$ has the property that $\dim\tilde{\pi}_V^{-1}(\tilde{y}) < \dim\pi^{-1}(y_0)$ for all $\tilde{y}\in\tau^{-1}(y_0)$.
\end{lemma}

\parag{Proof} Write $\dim\pi^{-1}(y_0) = n+k$ with $k > 0$ and let $\Gamma_{1}, \dots, \Gamma_{N}$ be the irreducible components of $\pi^{-1}(y_0)$ which has dimension $n + k$. For each $i \in [1, N]$ choose a  point $x_{i}$ in $\Gamma_{i}$, which is a smooth point of $\dim\pi^{-1}(y_0)$, and an $(n+k)-$scale $E_{i} = (U_{i},B_{i},j_{i})$ adapted to $\pi^{-1}(y_{0})$  in such a way that the following conditions are satisfied:
\begin{enumerate}[i)]
\item  
$x_{i}\in j_{i}^{-1}(U_{i}\times B_{i})$ and $j_{i}(x_{i}) = (0,0)$.
\item 
$ j_{i} (\Gamma_{i}\cap j_{i}^{-1}(U_{i}\times B_{i})) = U_{i}\times \{0\}$ 
\item
$\deg_{E_{i}}(\Gamma_{i}) = 1$ and $\deg_{E_i}(\Gamma_j) = 0$  if $j \not=  i$.
\end{enumerate}
As the compact set $\cup_{i=1}^{N} \  j_{i}^{-1}(\bar U_{i}\times \partial B_{i})$ does not meet $\pi^{-1}(y_{0})$ there exists an open neighborhood $V_0$ of $y_{0}$ in $N$ such that $\pi^{-1}(V_0) \cap j_{i}^{-1}(\bar U_{i}\times \partial B_{i}) = \emptyset$ for each $i$.
Put $W_{i} := \pi^{-1}(V_0) \cap j_{i}^{-1}(U_{i}\times B_{i})$. Let $pr_i\colon U_i\times B_i\rightarrow B_i$ be the natural projection and let $\theta_{i}\colon W_{i} \rightarrow V_0 \times U_{i}$ be the map induced by 
$(\pi, pr_{i}\circ j_{i})$. Clearly $\theta_{i}$ is a proper map with finite fibers and consequently  $\theta_i(W_i)$ is a nowhere dense analytic subset of $V_0\times U_i$ since $\dim W_i < \dim V_0\times U_i$. Then, due to Proposition 3.6.5 in \cite{[e]} (or Proposition III.6.1.5 in \cite{[BM.1]}), there exists an open neighborhood $V_i$ of $y_0$ in $V_0$ and a modification $\tau_i\colon\tilde{V}_i\rightarrow V_i$ such that the fibers of the strict tranform $\tilde{W}_i\rightarrow\tilde{V}_i$ are at most of dimension $n+k-1$. Moreover, we may assume that $\tau_1,\dots,\tau_N$ are modifications of the same neighborhood $V$ of $y_0$ in $V_0$. Hence, for each $i$, we then have the commutative diagrams

$$
\xymatrix{
W_i \ar[d]_{\rm can. incl.} \ar[r]^{\theta_i \ } & V\times U_i\ar[ldd]^{\rm can. proj.} \\
M_V \ar[d]_{\pi_V} & \\
V&
} 
\qquad\text{and}\qquad
\xymatrix{
\tilde{W}_i \ar[d]_{\rm can. incl.} \ar[r]^{\tilde{\theta}_i \ } & \tilde{V}_i\times U_i\ar[ldd]^{\rm can. proj.} \\
\tilde{M}_i \ar[d]_{\tilde{\pi}_i} & \\
\tilde{V}_i&
} 
$$
where the latter is obtained by taking the strict transform of the former by $\tau_i$. Let $C$ be an irreducible component of $\tilde{\pi}_i^{-1}(\tilde{y})$ where $\tilde{y}$ is a point  in $\tau_i^{-1}(y_0)$. Then from the above we see that $\dim C < n+k$ if $C\cap \tilde{W}_i\neq\emptyset$. By Lemma \ref{dominating.modification} there exists a modification $\tau\colon\tilde{V}\rightarrow V$ such that, for each $i$, we have a factorization 
\ $\tau\colon\xymatrix{\tilde{V}\ar[r]^{\sigma_i}&\tilde{V}_i\ar[r]^{\tau_i}& V}$. Hence the commutative diagram
$$
\xymatrix{
\tilde{M}_V\ar[r]\ar[d]_{\tilde{\pi}_V}&\tilde{M}_{i} \ar[d]_{\tilde{\pi}_{i}} \ar[r] & M_V\ar[d]^{\pi_V} \\
\tilde{V}\ar[r]^{\sigma_i}&\tilde{V}_i \ar[r]^{\tau_i} & V 
} 
$$
where $\tilde{\pi}_V\colon\tilde{M}_V\rightarrow\tilde{V}$ is the strict transform of $\pi_V\colon M_V\rightarrow V$ by $\tau$. Now, let $\tilde{y}$ be a point in $\tau^{-1}(y)$ and let us show that every irreducible component of $\tilde{\pi}^{-1}(\tilde{y})$ is at most of dimension $n+k-1$. We argue by contradiction and assume that there is an irreducible component $C$ of $\tilde{\pi}^{-1}(\tilde{y})$ such that $\dim C = n+k$. Then the proper map $\tilde{M}_V\rightarrow M_V$ maps $C$ biholomorphically onto an irreducible component of $\pi^{-1}(y_0)$, say $\Gamma_i$. Consequently the image of $C$ in $\tilde{M}_i$ is an irreducible component of $\tilde{\pi}_i^{-1}(\sigma_i(\tilde{y}))$ and intersects $\tilde{W}_i$. This is a contradiction since every irreducible component of $\tilde{\pi}_i^{-1}(\sigma_i(\tilde{y}))$ which intersects $\tilde{W}_i$ is at most of dimension $n+k-1$
\hfill{$\blacksquare$}

\parag{Remark}
In the situation above $\tilde{\pi}_V^{-1}(\tilde{y})$ can have infinitely many irreducible components, even when $\pi$ is quasi-proper.
\hfill{$\square$}

\parag{Example} Consider $\mathbb{C}^3$ with coordinates $(x, y, z)$ and define the two smooth hyper-surfaces  $M_1 := \{ y = 0\}$ and $M_2 := \{ y = x^2 \}$ in $\mathbb{C}^3$. Note that the set-theoretic intersection $D := M_1 \cap M_2$ is the line $\{x = y = 0\}$ and that $M_1$ and $M_2$ are tangent at each point of this line.\\
Let $\tau : Z \to \mathbb{C}^3$ be the blow-up of $\mathbb{C}^3$ at each point $(0, 0, n)$ with $n \in \mathbb{Z}$. Denote respectively  $\tilde{M}_1$ and $\tilde{M}_2$ the strict transforms of $M_1$ and $M_2$. Then $\tau_i : \tilde{M}_i \to M_i$ is the the blow-up of $M_i$ at each point $(0, 0, n)$ with $n \in \mathbb{Z}$ for $i = 1, 2$. Moreover the intersection $\tilde{M}_1 \cap \tilde{M}_2$ in $Z$ is the union of the strict transform of $D$ with  the exceptional $\mathbb{P}_1$ in $M_1$ (or in $M_2$) over the points $(0, 0, n)$ with $n \in \mathbb{Z}$.\\
Let $\tilde{M} := \tilde{M}_1 \cup \tilde{M}_2$ and let $\pi : \tilde{M} \to \mathbb{C}$ be the holomorphic function defined by $p_2\circ \tau_{\vert \tilde{M}}$ where $p_2(x, y, z) := y$.\\
The fiber of $\pi$ at the point $y_0 \not= 0$ is the pull-back by $\tau$ of the couple of lines given by the  equations  $y = x^2, y = y_0$ in $\mathbb{C}^3$. The fiber at $0$ of $\pi$ is equal to $\tilde{M_1}$ which is irreducible (smooth and connected) of dimension $2$.\\
Then $\pi$ is quasi-proper because the analytic subset $X := \tau^{-1}(\{z= 1/2\})\cap \tilde{M}_2$ in $\tilde{M}$  is proper (and finite of degree 2) on $\mathbb{C}$ via $\pi$ and it meets every irreducible component of each fiber of $\pi$ : for $\pi^{-1}(y_0)$ with $y_0 \not= 0$ this is clear as it contains the points $(\pm \sqrt{y_0}, y_0, 1/2)$ and for $y_0 = 0$ it contains the point $(0, 0, 1/2) \in \tilde{M}_1$.\\
But the strict transform of $\pi$ by the blow-up of the origin in $\mathbb{C}$ (which is the identity map) is the (equidimensional) holomorphic map 
$$ p_2\circ \tau_{\vert \tilde{M}_2} : \tilde{M}_2 \to \mathbb{C} $$
whose fiber at the origin has infinitely many irreducible components: the strict transform of the line $D$ and each exceptional $\mathbb{P}_1$ over the points $(0, 0, n), n \in \mathbb{Z}$.$\hfill \square$\\

The following theorem is originally due to D. Mathieu (see \cite{[Mt.00]}).

\begin{thm}\label{local.f-flattening}
Let $\pi\colon M\rightarrow N$ be an SQP map and  $y$ be a point in $N$. Then there exists an open neighborhood of $y$ in $N$ such that $\pi_V\colon M_V\rightarrow V$ admits an f-flattening. 
\end{thm}
 
\parag{Proof} Define  $n := \dim M - \dim N$. By iterated use of Lemma \ref{induction step} we get an open neighborhood $V_1$ of $y$ and a modification $\tau_1\colon\tilde{V_1}\rightarrow V_1$ such that, for all $\tilde{y}$ in $\tau_1^{-1}(y)$, the fiberover $\tilde{y}$ of the strict transform, $\tilde{\pi}_{V_1}\colon\tilde{M}_{V_1}\rightarrow \tilde{V}_1$, is of (pure) dimension $n$. Let $T$ denote the image of $\Sigma_{n+1}(\tilde{\pi}_{V_1})$ by $\tilde{\pi}_{V_1}$. Then $T\cap\tau_1^{-1}(y) = \emptyset$ and, due to Theorem \ref{general}, $\tilde{\pi}_{V_1}$ is an SQP map. Hence $T$ is an analytic subset of $\tilde{V}_1$ and consequently $\tilde{V}_1\setminus T$ is an open neighborhood of $\tau_1^{-1}(y)$. As $\tau_1$ is a proper map, there exists an open neighborhood of $V$ of $y$ such that $\tau_1^{-1}(V)\subseteq \tilde{V}_1\setminus T$ and the induced map $\tau_1^{-1}(V)\rightarrow V$ is a modification. By composing this modification with the normalization map $\tilde{V}\rightarrow \tau_1^{-1}(V)$ we get a modification $\tau\colon\tilde{V}\rightarrow V$ which has the property that the strict transform, $\tilde{\pi}\colon\tilde{M}_V\rightarrow\tilde{V}$, of $\pi_V$ by $\tau$ is an equidimensional SQP map. It follows that $\tilde{\pi}_V$ is an $f$-GF map since $\tilde{V}$ is normal.
\hfill{$\blacksquare$}

\subsection{Proof of the main theorem}

 For the  proof of Theorem \ref{SQP et GF} we need some technical results.

\begin{lemma}\label{graph.image}
Let $\pi\colon M \rightarrow N$ be a quasi-proper holomorphic dominant  map between reduced  complex spaces where $M$ is of pure dimension and $N$ is irreducible. Let $n$ be the relative dimension of $\pi$, let $\alpha\colon\mathcal{C}_n^f(\pi)^*\rightarrow N$ be  the natural map and let $\varphi$ be the maximal  reduced $f$-fiber map for $\pi$ over $N' := N \setminus \Sigma$. Then $\varphi$ is meromorphic along \ $\Sigma$ \  if and only if the closure of the image of $\varphi$ in $\mathcal{C}_n^f(\pi)^*$  is an $\alpha$-proper reduced complex subspace of $\mathcal{C}_n^f(\pi)^*$.
\end{lemma}

\parag{Proof}
Let $\Gamma_{\varphi}$ be the closure of the graph of $\varphi$ in $N\times\mathcal{C}_n^f(M)$ and $\Gamma$ be the closure of the image of $\varphi$ in $\mathcal{C}_n^f(M)$ (which is the also its closure in $\mathcal{C}_n^f(\pi)^*$). Then, by Proposition \ref{graph mero 1}, the map $\varphi$ is meromorphic if and only if $\Gamma_{\varphi}$ is an $N$-proper reduced complex subspace of $N\times\mathcal{C}_n^f(M)$ and, due to the remark following Definition \ref{sqp}, $\Gamma_{\varphi}$ is proper over $N$ if and only if $\Gamma$ is $\alpha$-proper. Now, the canonical projection $N\times\mathcal{C}_n^f(M)\rightarrow \mathcal{C}_n^f(M)$ is a holomorphic map, which induces a homeomorphism $\Gamma_{\varphi}\rightarrow\Gamma$, and  its inverse $\Gamma\rightarrow\Gamma_{\varphi}$ is induced by 
the holomorphic map $\left(\alpha,\,\id_{\mathcal{C}_n^f(\pi)^*}\right)\colon\mathcal{C}_n^f(\pi)^*\rightarrow N\times\mathcal{C}_n^f(\pi)^*$. From Theorem \ref{semi-proper direct image bis} it then follows that $\Gamma$ is an $\alpha$-proper reduced complex subspace of $\mathcal{C}_n^f(\pi)^*$ if and only if $\Gamma_{\varphi}$ is an $N$-proper reduced complex subspace of $N\times\mathcal{C}_n^f(M)$ (and in that case $\Gamma_{\varphi}$ and $\Gamma$ are biholomorphic).\hfill{$\blacksquare$}

\parag{Proof of Theorem \ref{SQP et GF}} Define $n := \dim M - \dim N$. Let $\alpha\colon\mathcal{C}_n^f(\pi)^*\rightarrow N$ be the natural map, $\varphi$ be the maximal  reduced $f$-fiber map for 
$\pi$ and $\Gamma$ denote the closure of the image of $\varphi$  in $\mathcal{C}_n^f(\pi)^*$. 

We have a natural identification $\mathcal{C}_n^f(\pi_V) = \alpha^{-1}(V)$ so we denote $\alpha_V\colon\mathcal{C}_n^f(\pi_V)\rightarrow V$ the map induced by $\alpha$ and put $\Gamma_V := \Gamma\cap\mathcal{C}_n^f(\pi_V) = \left(\alpha_{\vert\Gamma}\right)^{-1}(V)$. 

\smallskip
To prove that (i) implies (ii) we fix a point $y_0$ in $N$. By Proposition \ref{local.f-flattening} we get an open neighborhood $V$ of $y_0$ and an $f$-flattening $\tau_1\colon\tilde{V}\rightarrow V$ of $\pi_V$.  Hence the commutative diagram
$$
\xymatrix{
\tilde{M}_V \ar[d]_{\tilde{\pi}_V} \ar[r]^{q} & M_V \ar[d]^{\pi_V}\\ 
\tilde{V}\ar[r]^{\tau} & V \\
}
$$
where $q$ is the natural projection. Let $\Sigma_V$ denote the center of $\tau$ and $\psi\colon\tilde{V}\rightarrow\mathcal{C}_n^f(\tilde{\pi}_V)$ be the reduced fibermap for $\tilde{\pi}_V$. As the map $\tilde{V}\setminus\tau^{-1}(\Sigma_V)\rightarrow V\setminus\Sigma_V$ induced by $\tau$ is biholomorphic, the holomorphic map $q_*\circ\psi$ induces the reduced fibermap for $\pi_V$ on  $V\setminus\Sigma_V$ and consequently we have $(q_*\circ\psi)(\tilde{V}) = \Gamma_V$ by continuity.

Consider the commutative diagram
$$\xymatrix{
& \mathcal{C}_n^f(\pi_V)^* \ar[d]^{\alpha_V}\\ 
\tilde{V}\ar[ur]^{q_*\circ\psi} \ar[r]^{\tau} & V \\}$$
Observe that $q_*\circ\psi$ is proper, since $\alpha_V\circ q_*\circ\psi = \tau$ is proper, so $\Gamma_V$ is an $\alpha_V$-proper reduced complex subspace of $\mathcal{C}_n^f(\pi_V)^*$ due to Theorem \ref{semi-proper direct image bis}. It follows that $\Gamma$ is an $\alpha$-proper reduced complex subspace of $\mathcal{C}_n^f(\pi)^*$ and $\varphi$ is meromorphic.

\smallskip
In order to  show that (ii) implies (iii) we observe that, by Lemma \ref{graph.image}, $\Gamma$ is an $\alpha$-proper reduced complex subspace of $\mathcal{C}_n^f(\pi)^*$. Let $\tilde{M}$ denote the set-theoretic graph in $\Gamma\times M$ of the $f$-analytic family of $n$-cycles defined by $\Gamma\hookrightarrow\mathcal{C}_n^f(M)$ and consider the commutative diagram 
$$
\xymatrix{
\tilde{M}\ar[d]_{\tilde{\pi}}\ar[r]&M \ar[d]^{\pi}\\ 
\Gamma\ar[r]^{\alpha_{|\Gamma}} & N \\
}
$$
where $\tilde{\pi}$ and $\tilde{M}\rightarrow M$ are induced by the natural projections $\Gamma\times M\rightarrow\Gamma$ and $\Gamma\times M\rightarrow M$. Then $\alpha_{|\Gamma}\colon\Gamma\rightarrow N$ is a modification and $\tilde{\pi}\colon\tilde{M}\rightarrow\Gamma$ is the strict transform of $\pi$ by $\alpha_{|\Gamma}$. As the map $\tilde{\pi}$ is geometrically $f$-flat, the proof is completed.

\smallskip
Finally we prove that (iii) implies (i). To this end let $\sigma\colon\tilde{N}\rightarrow N$ be a modification such that the strict transform $\tilde{\pi}\colon\tilde{M}\rightarrow\tilde{N}$ of $\pi\colon M\rightarrow N$ by $\sigma$ is an $f$-GF map and let $\psi\colon\tilde{N}\rightarrow\mathcal{C}_n^f(\tilde{M})$ be the reduced $f$-fibermap for $\tilde{\pi}$. Then we have a commutative diagram of holomorphic maps (see (i) implies (ii))
$$
\xymatrix{
& \Gamma \ar[d]^{\alpha_{\vert\Gamma}}\\ 
\tilde{N}\ar[ur] \ar[r]^{\tau} & N \\
}
$$
where $\tilde{N}\rightarrow\Gamma$ is surjective. It follows that $\alpha_{\vert\Gamma}\colon\Gamma\rightarrow N$ is proper since $\tau$ is proper. Hence $\pi$ is an SQP map.
$\hfill \blacksquare$

\parag{Remark} Under the hypotheses of Theorem \ref{SQP et GF} suppose that $\pi$ is an SQP map. Then (using the same notation as above) the $f$-flattening $\alpha_{\vert\Gamma}\colon\Gamma\rightarrow N$ is \lq\lq optimal\rq\rq\ in the sense that every $f$-flattening of $\pi$ factorizes through $\alpha_{\vert\Gamma}$. More precisely, it has the following property (see (iii) implies (i)):

\begin{itemize}
\item
If $\tau\colon\tilde{N}\rightarrow N$ is a modification such that the strict transform of $\pi$ by $\tau$ is geometrically $f$-flat, then $\tau = \alpha_{\vert\Gamma}\circ q_*\circ\psi$.
\end{itemize}

\begin{defn}\label{appli. fiber}
Let $\pi\colon M \rightarrow N$ be an SQP map, put $n:= \dim M - \dim N$ and let $\varphi$ be the maximal reduced $f$-fiber map of $\pi$.  The closure of the image of $\varphi$ in  $\mathcal{C}_{n}^{f}(M)$ will henceforth be denoted by $N_{\pi}$ and the $f$-analytic family of $n$-cycles in $M$ classified by the canonical injection $N_{\pi}\hookrightarrow\mathcal{C}_{n}^{f}(M)$
will be called the  {\bf\em (meromorphic) family of fibers of $\pi$}.
\end{defn}

\begin{lemma}\label{propre et fini}
Let $\pi : M \to N$ be a dominant and quasi-proper holomorphic map from a pure dimensional  to an irreducible complex space. Let $\theta: M_1 \to M$ a proper finite and surjective map. Then $\pi\circ \theta : M_1 \to N$ is SQP if and only if $\pi$ is an SQP map.
\end{lemma}

\parag{Proof} First assume that $\pi$ is a SQP map. Let $K$ be a compact set in $N$. Then there exists a compact set $L$ in $M$ such that any irreducible component of any $n$-cycle in $\mathcal{C}_n^f(\pi)$ which is in the meromorphic family of fibers of $\pi$ and over $K$ has to meet $L$. Let $\Gamma$ be an irreducible component of a cycle in $\mathcal{C}_n^f(\pi\circ\theta)$ and which is over a point in $K$. The image by $\theta$ of $\Gamma$ is an irreducible component of a $n$-cycle in $\mathcal{C}_n^f(\pi)$ which is in the meromorphic family of fibers of $\pi$ because the direct image of cycles by $\theta$ is continuous and the generic fibers of $\pi\circ \theta$ are the pull-back by $\theta$ of generic fibers of $\pi$. Then $\Gamma$ has to meet the compact set $\theta^{-1}(L)$ in $M_1$. This gives the compactness of the subset of $\mathcal{C}_n^f(\pi\circ\theta)$ which is the closure of generic fibers over $K$  of the map
 $\pi\circ \theta$. This implies that $\pi\circ \theta$ is SQP.\\
The  converse is proved in a similar way.  $\hfill \blacksquare$

\begin{cor}\label{stab.}
Let $\pi\colon M \rightarrow N$ and $g\colon N \rightarrow P$ be two SQP maps. Let $\Sigma \subset N$ denote the locus of big fibers of $\pi$ and assume that it is $g$-proper. Then the map $g\circ\pi\colon M \rightarrow P$ is also an SQP map.
\end{cor}

\parag{Proof} Consider first the case where $N$ is normal.  By Theorem \ref{SQP et GF} there exists a modification $\tau\colon \tilde{N} \rightarrow N$ whose center is $\Sigma$ such that the strict transform $\tilde{\pi}\colon \tilde{M} \rightarrow \tilde{N}$ of $\pi$ by $\tau$ is an $f$-GF map. Then we have the following commutative diagram
$$
\xymatrix{
\tilde{M} \ar[d]^{\tilde{\pi}} \ar[r]^{\rho} & M \ar[d]^{\pi}\\ 
\tilde{N}\ar[r]^{\tau}\ar[dr]_{g\circ\tau} & N \ar[d]^{g} \\
 & P
}
$$
where $\rho$ is a modification since $\pi$ is a dominant map. Now thanks to Proposition \ref{SQP.compos.fGF} the map $g\circ\tau$ is strongly quasi-proper so $g\circ\tau\circ\tilde{\pi} = g\circ\pi\circ\rho$ is equally an SQP map by Theorem \ref{modific.quasiproper}.\\
 Hence $g\circ\pi$ is an SQP map due to Theorem \ref{modific.quasiproper}.\\
 When $N$ is not normal,  let $\nu : N_1 \to N$ the normalization of $N$ and $\pi_1 : M_1 \to N_1$ the strict transform of $\pi$ by $\nu$. Then $\pi_1$ is a SQP map Thanks to  Lemma \ref{propre et fini} the map $g_1 := g\circ \nu$ is an SQP map. Now $\nu^{-1}(\Sigma)$ is the locus of big fibers for $\pi_1$. We conclude from the previous case that $g_1\circ \pi_1$ is an SQP map. Let $q : M_1 \to M$ be the natural projection. It is a finite modification. Then $g\circ\pi$ is an  SQP map again by Lemma \ref{propre et fini}.  \hfill{$\blacksquare$}

\parag{Remark} In the case where the locus $\Sigma \subseteq N$ of big fibers of $\pi$ is compact it is of course proper over $P$. In particular this is the case when $\pi$ is equidimensional.

 \begin{cor}\label{cor.1}
 Let $\tau : \tilde{M} \to M$ be a modification with center $C$ of an irreducible complex space $M$ and let $(X_{s})_{s \in S}$ be an $f$-analytic family of $n$-cycles in $M$ parametrized by a reduced complex space $S$. 
 Assume moreover that the following conditions are satisfied:
\begin{enumerate}[(i)]
\item 
For $s$ generic in $S$  the $n$-cycle $X_{s}$ is reduced and has no irreducible component contained in $C$.
\item 
The projection $G \cap (S \times C) \to S$ is proper, where $G \subset S \times M$ is the graph of the family $(X_{s})_{s \in S}$
\end{enumerate}
Then there exists a   modification $\theta\colon\tilde{S} \rightarrow S$ and an $f$-analytic family of $n$-cycles $(Y_{\tilde{s}})_{\tilde{s} \in \tilde{S}}$ in $\tilde{M}$ parametrized by $\tilde{S}$ such that 
\begin{enumerate}
\item 
for $\tilde{s}$ generic in $\tilde{S}$ the cycle $Y_{\tilde{s}}$ is the strict transform by $\tau$ of the cycle $X_{\theta(\tilde{s})}$,
 \item 
 for each $\tilde{s} \in \tilde{S}$ we have $\tau_{*}(Y_{\tilde{s}}) = X_{\theta(\tilde{s})}$.
 \end{enumerate}
 \end{cor}
  
 \parag{Proof} Remark first that normalizing $S$ we may replace our initial family by a finite sum (may be with multiplicities)  of $f$-analytic families having irreducible generic cycles. We may also assume then that $S$ is irreducible.
 So it is enough to consider the case where the graph $G$ is irreducible. Let $\tilde{G} \subset S \times \tilde{M}$ be the strict transform of the natural projection $G\rightarrow M$ by the  modification 
 $\tau$. Then $\tilde{G}\rightarrow G$ is a modification whose center  is $G \cap (S \times C)$. The restriction of the natural projection  $p: G \to S$ to $G \cap (S \times C)$ is proper by assumption (ii) and for $s$ generic in $S$ the fiber of $p$ is not contained in $S \times C$ by assumption (i). Then by Theorem \ref{modific.quasiproper} the natural projection $\tilde{p}\colon \tilde{G} \rightarrow S$ is an SQP map. 
Now, using the notation introduced in Definition \ref{appli. fiber}, we put $\tilde{S} := S_{\tilde p}$ and let $\theta\colon\tilde{S} \rightarrow S$ and $\psi\colon\tilde{S} \rightarrow\mathcal{C}_n^f(\tilde{G})$ be the natural projections. Then we get $\psi(\tilde s) = \{\theta(\tilde s)\}\times Y_{\tilde s}$, where $(Y_{\tilde{s}})_{\tilde\in \tilde S}$ is an $f$-analytic family in $\tilde{M}$ which has the required properties. $\hfill \blacksquare$\\
 
 \subsection{Extendable cycles}

We shall discuss now the converse of the restriction problem considered in chapter V subsection 6.2.\\

Recall the classical important theorem of E. Bishop (see \cite{[Bish.64]} Theorem 3 p.299).

\begin{thm}\label{extend.2}
Let $M$ be a complex space, $n$ be a non negative integer  and fix a closed analytic subset $T \subset M$ with no interior point in $M$. Fix a continuous hermitian metric $h$ on $M$. Let $X$ be a pure $n$-dimensional analytic subset in $M \setminus T$. Assume that for each point  $t_0$ in $T$ there exists a relatively compact open subset $V(t_0)$ of $t_0$ in $M$ such that the integral 
$$\int_{V(t_0)\cap X} h^n < +\infty .$$
Then $\bar X$, the closure of $X$ in $M$,  is a complex analytic subset of pure dimension $n$  in $M$.
\end{thm}

 In the situation of the theorem, since $\bar X$ has pure dimension $n$, it has no irreducible component contained in $T$.

\begin{defn}\label{extend.1}
Let $M$ be a complex space, $n$ be a non negative integer  and fix a closed analytic subset   $T \subset M$ with no interior point in $M$. Fix a continuous hermitian metric $h$ on $M$. Let $S$ be a   pure dimensional complex space with dimension $\sigma$ and let $(X_s)_{s \in S}$ be a f-analytic family of $n$-cycles in $M \setminus T$. We say tha the family $(X_s)_{s \in S}$ is {\bf\em pre-extendable along $T$} if for each $s_0 \in S$ and each $t_0 \in T$ there exist relatively compact open neighborhoods $U(s_0)$ and $V(t_0)$ respectively of $s_0$ in $S$ and of $t_0$ in  $M$ and a constant $C > 0$ such that
$$ \int_{V(t_0)\cap X_s}  h^n \leq C \quad \forall s \in U(s_0).$$
\end{defn}

\parag{Remarks}\begin{enumerate}
 \item  The condition above is automatic for $t_0 \not\in T$ because in a  $n$-scale $E := (U, B, j)$ on $M \setminus T$ adapted to $X_{s_0}$ the fact that for $s$ near $s_0$ the scale $E$ is still adapted to $X_s$ with $\deg_E(X_s) = \deg_E(X_{s_0})$ implies that the volume of $X_s$ in the relatively compact open set $j^{-1}(U\times B)$ is uniformly bounded in a neighborhood of $s_0$ in $S$ (see  Lemma 4.2.3 in Chapter IV of  \cite{[e]}) .
 \item As an obvious consequence of Bishop's result, in the situation of the definition above, for each $s \in S$ the closure  of  $\overline{\vert X_s\vert}$ is a pure $n$-dimensional  analytic set in $M$. Note $Y_s$ the $n$-cycle in $M$ defined by $\vert Y_s\vert := \overline{\vert X_s\vert}$ and such that, for each irreducible component $\Gamma$ of $X_s$,  the multiplicity of $\bar \Gamma$ in $Y_s$ is equal to the multiplicity of $\Gamma$ in $X$.
 Then, in general, the family $(Y_s)_{s \in S}$ is not $f$-continuous. See the simple example below. \\
 But if  the family $(Y_s)_{s \in S}$  is $f$-continuous on an open set $S'$ in $S$, then it is $f$-analytic on $S'$ thanks to Analytic Extension Theorem \ref{cycles}.
 \end{enumerate}

\parag{Example} Let $M$ be the blow-up of the origin in $\mathbb{C}^2$ and  let $T$ be the exceptional divisor. Let $(X_s)_{s \in \mathbb{C}}$ be the family of lines obtained by translating a line $X_{s_0}$ through the origin in $\mathbb{C}^2$. Then the limit of $Y_s$ when $s$ goes to $s_0$, $s \not= s_0$, is the total transform of $ X_{s_0}$ which contains the exceptional divisor. But $Y_{s_0}$ is only the strict transform of the line through the origin.\\

\begin{thm}\label{extend.3}
Let $M$ be a complex space, $n$ be an integer  and fix a compact analytic subset $T \subset M$ with no interior point in $M$. Let  $(X_s)_{s \in S}$ be an $f$-analytic family of $n$-cycles in $M \setminus T$ which is pre-extendable. Then there exists a modification $\tau : \tilde{S} \to S$ and an f-analytic family $(Z_{\tilde{s}})_{\tilde{s}} \in \tilde{S}$ of $n$-cycles in $M$  with the following properties:
\begin{enumerate}[(i)]
\item For each $\tilde{s} \in \tilde{S}$ we have $Z_{\tilde{s}} \cap (M\setminus T) = X_{\tau(\tilde{s})}$.
\item For $\tilde{s}$ generic in $\tilde{S}$ we have $Z_{\tilde{s}}  = \overline{X_{\tau(\tilde{s})}}$.
\end{enumerate}
\end{thm}

\parag{Proof} Let $G$ be the  graph-cycle of the family $(X_s)_{s \in S}$ in $S \times (M \setminus T)$. Then our hypothesis implies, thanks to Bishop's theorem recalled above, that the closure of $\vert G\vert$ in $S \times M$ is a complex analytic set of pure dimension $\sigma + n$, and the compactness of $T$ implies that its projection to $S$ is strongly quasi-proper, thanks also to the theorem of \cite{[B.78]}(see Theorem  3.6.6 \cite{[e]})  to bound the volume of its generic fibers. Then the existence of a geometric $f$-flattning (see paragraph V.7.2) for the projection onto $S$ of the closure of $\vert G\vert$ in $S \times M$ allows to conclude.$\hfill \blacksquare$\\

 \marginpar{Il faut encore estimer le volume près de l'infini !! c'est faisable via le degré voir probleme Ya. Ajout du 8/3/23.}

\parag{Example} Let $(X_s)_{s \in S}$ be an $f$-analytic family of $n$-cycles in $\mathbb{C}^q$ parametrized by an irreducible complex space $S$. Then this family is pre-extendable to $\mathbb{P}_q$ if and only if each cycle is algebraic. Then the previous result explains that, up to a modification of $S$, we obtain a proper  family of compact cycles in $\mathbb{P}_q$ which is given, on a dense Zariski open set of $S$, by the closure of these cycles.\\

We conclude this subsection by giving a simple case where we have a nice restriction map without assuming that the cycles are compact. It is an obvious corollary of subsection 6.3.

\begin{cor}\label{res.4}
Let $\pi : M \to N$ be an  $f$-geometrically flat map between reduced complex spaces and let $T$be  a closed analytic subset in $M$. Let $\Theta$ be the analytic subset of points $y$  in $N$  such that $T$ contains at least one irreducible component of the fiber $\pi^{-1}(y)$. Then the induced map $\pi_T : M \setminus (\pi^{-1}(\Theta)\cup T) \to N \setminus \Theta$ is a geometrically f-flat map.$\hfill \blacksquare$\\
\end{cor}

Note that when $\Theta = \emptyset$ the map $\pi_T : M \setminus T \to N$ is geometrically f-flat.

 \chapter{ Applications}

\section{Application to meromorphic quotients}

Let $M$ be an irreducible complex space. Classically, an analytic equivalence relation on $M$ is defined by its graph which is an analytic subset $\mathcal{R} \subset M \times M$.  In his fundamental paper \cite{[C.60]}  Henri Cartan studies the case of a proper analytic equivalence relation, which is the case where  the first projection $p_1 : \mathcal{R} \to M$ is a proper holomorphic map. In his article he gives a necessary and sufficient condition for the existence of a holomorphic quotient. Such an existence  means that the quotient space endowed with the sheaf of invariant holomorphic functions is a complex space. But  this condition is not always true even assuming that $M$ is compact.\\
Nevertheless, under this compactness condition, assuming for instance that $\mathcal{R}$ is irreducible, there is always an irreducible complex space which is an "almost" quotient for such an equivalence relation using the reduced complex space of compact analytic cycles in $M$ as follows:\\
Let $n := \dim \mathcal{R} - \dim M$. Then we have, thanks to the fact that $p_1: \mathcal{R} \to M$ is proper and surjective, an analytic subset $\Sigma$ in $M$ and a holomorphic fiber map
$$ \varphi : M \setminus \Sigma \to \mathcal{C}_n(M) $$
classifying the compact $n$-cycles in $M$ which are the (generically reduced)$n$-dimensional fibers of $p_1$. Moreover this map is meromorphic along $\Sigma$ which means that there exists a modification 
$\tau: \tilde{M} \to M$ with center in $\Sigma$  and a holomorphic map $\tilde{\varphi}: \tilde{M} \to \mathcal{C}_n(M)$ which coincides with $\varphi$ on $M \setminus \Sigma \simeq \tilde{M} \setminus \tau^{-1}(\Sigma)$. Assuming now that $M$ is compact (or more generally that $\tilde{\varphi}$ is proper) Remmert's  Direct Image Theorem ensures that $\tilde{\varphi}(\tilde{M}) = \overline{\varphi(M \setminus \Sigma)}$ is an analytic subset $Q$ in the reduced complex space $\mathcal{C}_n(M)$. Then, it is clear that there exists an open and dense subset $Q'$ in $Q$ which  is in bijection with the generic  equivalence classes for the given analytic equivalence relation.\\
These considerations motivate the introduction of  the notion of meromorphic quotient.\\
As the reader may see, the tools introduced in the previous chapters are precisely those which allow us to generalize to strongly quasi-proper analytic equivalence relations this point of view, in order to obtain an existence theorem for meromorphic quotients in a rather large context.\\

A simple way to produce an analytic equivalence relation on an irreducible complex space $M$ is to look at a holomorphic action of a complex Lie group $G$ on $M$. In such a case we often encounter the following situation: \\
 there exists an open dense set $\Omega$ which is $G$-stable and in which the $G$-orbits are closed (in $\Omega$), but these orbits are not closed in $M$ in general. This is already the case for the obvious action of $G := \mathbb{C}^*$ on $\mathbb{C}^n$. This kind of situation and the existence of meromorphic quotients for some of these actions are studied in \cite{[G.86]} and  in \cite{[B.18]}. They motivate the definitions of meromorphic equivalence relations and of meromorphic quotients which are given below, although they may seem a little more complicated than necessary after  the comments  we give above in the  case of a proper equivalence relation.  We do  not present the results of \cite{[B.18]} in this book,  but nevertheless it seems interesting to treat the general situation in which the tools introduced in the previous chapters may be used with success.\\

\subsection{Holomorphic quotient}

In this paragraph  we collect some basic facts on  holomorphic quotients with respect to analytic  equivalence relations defined by  holomorphic maps.\\
  First let us  recall the basic definitions concerning holomorphic quotients. 

\medskip
In the sequel $M$ will always be a reduced complex space. 

\begin{defn}
We say that an equivalence relation,  $R\subseteq M\times M$,  on $M$ is  {\bf\em analytic}  if  $R$ is an analytic subset of $M\times M$.
\end{defn}

Let $R$ be an analytic equivalence relation on $M$ and $Q$ be the topological quotient of $M$ by $R$.
For every open subset  $U$ of $Q$  we let $\mathcal{A}(U)$ denote the $\C-$algebra of all functions $g$ on $U$ such that $g\circ q$ is holomorphic on $q^{-1}(U)$. Then $\mathcal{A}$ is a sheaf of $\C$-algebras   on $Q$.
If the ringed space $(Q,\mathcal{A})$ is a reduced complex space, then we call it the {\bf holomorphic quotient of $M$ by $R$}. In this case we say that  {\bf $M$ admits a holomorphic quotient with respect to $R$}.

\medskip
Now suppose we have a reduced complex space $T$ and a holomorphic map $f\colon M\rightarrow T$. Let $R_f$ be the equivalence relation determined by $f$ and let $(Q,\mathcal{A})$ be the corresponding ringed space (as defined above). Then we say that $f$ is a {\bf holomorphic quotient map} if the canonical morphism of ringed spaces from $(Q,\mathcal{A})$ to $(T,\mathcal{O}_T)$ is an  isomorphism. In other words $f$ is a holomorphic quotient map if $(Q,\mathcal{A})$ is a reduced complex space and moreover isomorphic to $(T,\mathcal{O}_T)$.

\medskip
The proposition below is proved in \cite{[G.86]} without being explicitly stated there. It gives a necessary topological condition for $M$ to admit a holomorphic quotient with respect to an analytic equivalence relation.

\begin{prop}\label{semi-proper.quotient}
Let   $R$ is an analytic equivalence relation on $M$ and denote respectively $Q := M/R$ and $q\colon M\rightarrow Q$ the corresponding toplogical quotient and quotient map.
Then the map $q$ is semi-proper if and only if $Q$ is a first countable Hausdorff space.
\end{prop}  
\parag{Proof} We first notice that the saturation with respect to $R$ of a compact subset of $M$ is closed. Indeed, if $p_1$ and $p_2$ denote the first and second projections of $M\times M$ onto $M$ and $K$ is a compact subset of $M$, then $p_1((M\times K)\cap R)$ is the saturation of $K$.  But, as $R$ is closed, the map $(M\times K)\cap R\rightarrow M$ induced by $p_1$ is proper and consequently $p_1((M\times K)\cap R)$ is closed. 

\smallskip
It follows that $q(K)$ is closed in $Q$ for  every compact subset of $M$. 

\smallskip
Now, suppose that $q$ is semi-proper and let $x$ and $y$ be two distinct points in $Q$. Then there exists a compact subset $L$ of $M$ such that $q(L)$ is a neighborhood of $\{x,y\}$. Let $V_x$ and $V_y$ be disjoint open neighborhoods of the compact subsets $q^{-1}(x)\cap L$ and $q^{-1}(y)\cap L$. Then $K_x := L\setminus V_x$ and $L\setminus V_y$ are compact subsets of $M$ such that $q^{-1}(x)\cap K_x = \emptyset$ and $q^{-1}(y)\cap K_y = \emptyset$. Hence $q(K_x)$ and $q(K_y)$ are closed subsets of $Q$ and it follows that $q(L)\setminus q(K_x)$ and $q(L)\setminus q(K_y)$ are disjoint neighborhoods of $x$ and $y$ in $Q$.

In fact we have shown that, for every neighborhood $V$ of $q^{-1}(x)\cap L$ in $L$, there exists an open neighborhood $W$ of $x$ in $Q$ such that $q^{-1}(W)\cap L\subseteq V$. This implies that $x$ has a countable basis of neighborhoods in $Q$ since $q^{-1}(x)\cap L$ has a countable basis of neighborhoods in $L$. Hence $Q$ is first countable. 

\smallskip
Conversely, suppose that $Q$ is a first countable Hausdorff space and let us prove by contradiction that $q$ is semi-proper. So assume that $q$ is not semi-proper. Then there exists a point $x$ in $Q$ which is not an interior point of $q(L)$ for any compact subset $L$ of $M$.  Let $(L_n)_{n\in \N}$ be an exhaustion of $M$ by compact subsets such that $x\in q(L_n)$ for all $n$, and let  $(W_n)_{n\in \N}$ be a decreasing neighborhood basis of $x$ in $Q$. Then, for each $n\in\N$, there exists a point $x_n$ in $W_n$ such that $q^{-1}(x_n)\cap L_n = \emptyset$. It follows that $M\setminus\bigcup\limits_{n\in\N}q^{-1}(x_n)$ is an open saturated subset of $M$. This contradicts the hypothesis that $q$ is a quotient map because $q\left(M\setminus\bigcup\limits_{n\in\N}q^{-1}(x_n)\right) = Q\setminus\{x_n\ |\ n\in\N\}$ is not an open subset of $Q$.
\hfill{$\blacksquare$}\\

The following result is a direct consequence of Proposition \ref{semi-proper.quotient}.

\begin{cor}
In the situation of Proposition \ref{semi-proper.quotient} assume that $q$ is a semi-proper map. Then $Q$ is a locally compact Hausdorff space.
\hfill{$\blacksquare$}
\end{cor}

\bigskip
Now suppose that we have a surjective holomorphic map $f\colon M\rightarrow T$. Then, by Proposition \ref{semi-proper.quotient}, a necessary condition for the map $f$ to be a holomorphic quotient map is that it  is semi-proper. Moreover we have the following result.

\begin{lemma}\label{semi-proper.top.quot}
Let  $f\colon M\rightarrow T$ be a semi-proper surjective holomorphic map. Let $R_f$ be the analytic equivalence relation determined by $f$ and let $(Q,\mathcal{A})$ be the corresponding ringed space. Then the canonical morphism of ringed spaces $\tilde{f}\colon Q\rightarrow T$ induces a homeomorphism.
\end{lemma}

\parag{Proof} We have a commutative diagram of ringed spaces
\begin{equation*}
\xymatrix{M \ar[d]_{q}\ar[r]^{f}  &T\\
Q\ar[ur]_{\tilde{f}} } 
\end{equation*}
where $\tilde{f}$ is bijective and continuous so  it is enough to show $\tilde{f}$ is a closed map. To this end let $K$ be a compact subset of $T$. As $f$ is semi-proper there exists a compact subset $L$ of $M$ such that $f(L) = K$ and consequently $\tilde{f}^{-1}(K) = q(L)$. It follows that $\tilde{f}$ is proper and hence a homeomorphism.
\hfill{$\blacksquare$}

\begin{prop}\label{quotient.weak.noramlization}
Let  $f\colon M\rightarrow T$ be a dominant semi-proper  holomorphic map between reduced complex spaces and suppose moreover that $M$ is weakly normal. Let $R_f$ be the analytic equivalence relation determined by $f$ and let $(Q,\mathcal{A})$ be the corresponding ringed space. Then $(Q,\mathcal{A})$ is the weak normalization of $(T,\mathcal{O}_T)$. 
\end{prop}
\parag{Proof} Since the the normalization map  $\nu\colon \tilde{M}\rightarrow M$ is obviously a holomorphic quotient map it is not restrictive to assume $M$ normal and, due to Lemma \ref{semi-proper.top.quot}, we may identify the topological spaces $Q$ and $T$. Thus the proof consists of showing that $\mathcal{A}$ is the sheaf of continuous meromorphic  functions on $T$.

Let us first prove that, in the case where $T$ is a connected manifold, we have $\mathcal{A} = \mathcal{O}_T$.  Put $n := \dim T$ and let $S(M)$ denote the singular locus of $M$. Let $A$ be the set of those $t$ in $T$ which satisfy   $f^{-1}(t)\subseteq S(M)$ and  let $B$ be the set of all $t$ in $T\setminus A$ such that $f$ is of rank strictly less than $n$ at every point in $f^{-1}(t)$. Due to Lemma \ref{Kuhl-banach.3},  the set $A$ is $b$-negligible in $T$ and the set $B$ is $b$-negligible in $T\setminus A$, since  $f$ is semi-proper. It follows that $A\cup B$ is a $b$-negligible subset of $T$. Now let $U$ be an open subset of $T$ and $g$ be a holomorphic function on $f^{-1}(U)$ which is constant on every fiberof $f$ over $U$. As $f$ is a topological quotient map there exists a (unique) continuous function $\tilde{g}\colon U\rightarrow\C$ satisfying $\tilde{g}\circ f = g$. But, for every $y$ in $U\setminus A\cup B$, the map $f$ admits a holomorphic section $\sigma_y$ in an open neighborhood $V_y$ in $U$ and consequently $\tilde{g} = g\circ\sigma_y$ on $V_y$. It follows that  $\tilde{g}$ is holomorphic on $U\setminus A\cup B$ and hence on $U$ as $A\cup B$ is a $b$-negligible subset of $T$. This shows that $\mathcal{A} = \mathcal{O}_T$. 

So, in the general case, the sheaves $\mathcal{A}$ and  $\mathcal{O}_T$ are identical on $T\setminus S(T)$. 

Now let $\tilde{T}$ denote the weak normalization of $T$ and $U$ be an open subset of $T$. If $g\in\mathcal{A(U)}$, then $g$ is continuous on $U$ and holomorphic on $U\setminus S(T)$. It follows that $\mathcal{A}(U)\subseteq\mathcal{O}_{\tilde T}(U)$.  
Conversely, suppose that $g\in \mathcal{O}_{\tilde T}(U)$. Then $g\circ f$ is holomorphic on $f^{-1}(U)\setminus f^{-1}(S(T))$ and continuous on  $f^{-1}(U)$. It follows that $g\circ f$ is holomorphic on $f^{-1}(U)$ since $M$ is normal and $f$ is dominant. Hence $g\in\mathcal{A(U)}$. This shows that  $\mathcal{O}_{\tilde T}(U)\subseteq\mathcal{A}(U)$.
\hfill{$\blacksquare$}

\begin{cor}\label{holom.quot.weak.norm}
In the situation of Proposition \ref{quotient.weak.noramlization} suppose moreover that $T$ is weakly normal. Then $f\colon M\rightarrow T$ is a holomorphic quotient map. 
\end{cor}
\parag{Proof} Due to Proposition \ref{quotient.weak.noramlization}, the induced map $Q\rightarrow T$ is the weak normalization of $T$ and hence an isomorphism since $T$ is weakly normal.
\hfill{$\blacksquare$}

\begin{cor}\label{holom.quot.}
Let  $f\colon M\rightarrow T$ be a dominant semi-proper  holomorphic map between reduced complex spaces. Then $M$ admits a holomorphic quotient with respect to $R_f$.
\end{cor}

\parag{Proof} Define the sheaf $\mathcal{A}$ as above and let $\tilde{\mathcal{O}}_T$ denote the sheaf of continuous meromorphic functions on $T$. Then we have $\mathcal{O}_T\subseteq\mathcal{A}\subseteq\tilde{\mathcal{O}}_T$ and $\mathcal{A}$ is an $\mathcal{O}_T-$algebra of finite type. Consequently $\mathcal{A}$ is $\mathcal{O}_T-$coherent and consequently the ringed space $(T,\mathcal{A})$ is a reduced complex space.
\hfill{$\blacksquare$}

\parag{Remark}
Under the hypotheses of Corollary \ref{holom.quot.} the map $f$ is in general not a holomorphic quotient map.

 \subsection{Meromorphic equivalence relations and meromorphic quotients}
 
We first give the main definitions.
 
 \begin{defn}\label{mero-equiv.} Let $M$ be an irreducible  complex space and let $\mathcal{R} \subset M \times M$ be an analytic subset and  $\Omega$ be a dense open  set  in $M$.
We shall say that $(\mathcal{R}, \Omega)$ is a {\bf \em meromorphic equivalence relation} if it  satisfies the following conditions:
\begin{enumerate}[(i)]
\item The subset $\mathcal{R} \cap (\Omega\times \Omega)$ is an equivalence relation on $\Omega$.
\item  There exists a dense subset $\Omega'$ in $\Omega$ such that for each $x \in \Omega'$ we have $\overline{\Omega_x} = \mathcal{R}_x$ where $\Omega_x$ is the equivalence class of $x$ in $\Omega$ and where $\mathcal{R}_x$ is defined by the relation 
$\{x\}\times \mathcal{R}_x = \mathcal{R} \cap (\{x\}\times M) $. 
\end{enumerate}
\end{defn}

As the choice of the dense set $\Omega'$ in $\Omega$ is not so important (but its existence is important) we often omit it in the definition of a meromorphic equivalence relation.

\parag{Remark}
Let $(\mathcal{R}, \Omega)$ be a meromorphic equivalence relation on an irreducible complex space $M$. By density of $\Omega'$ in $M$ the set  $\mathcal{R}$ is both reflexive and symmetric. In particular we have for all $x$ in $M$ the equality $\mathcal{R} \cap (M\times\{x\}) = R_x\times\{x\}$.

\parag{Example} Let $g: M \dasharrow N$ be a meromorphic map where $M$ is an irreducible complex space and $N$ is a reduced complex space, or a Banach analytic set or  $\mathcal{C}_n^f(P)$   where $P$ is a reduced  complex space (see subsection IV.2). Then let $\tau : \tilde{M} \to M$ be the modification given by the graph of $g$ and let $\tilde{g} :\tilde{M} \to N$ be the natural projection. Let $\tilde{\mathcal{R}} \subset \tilde{M} \times \tilde{M}$ be the graph of the analytic equivalence relation given by $\tilde{g}$ and $\mathcal{R}$ be the image of
 $\tilde{\mathcal{R}}$ on $M\times M$ by $\tau\times \tau$. Then $\mathcal{R}$ is a closed analytic subset of $M \times M$ thanks to Remmert's Direct Image Theorem and if $\Sigma$ is the center of $\tau$ the open dense set $\Omega := M \setminus \Sigma$ satisfies the condition $(i)$ of  Definition \ref{mero-equiv.}.. \\
  Moreover the set of point $x$  in $\Omega$ where the condition $(ii)$ is satisfied contains a dense subset because  $\tau^{-1}(\Sigma)$ is a closed analytic subset with empty interior in $\tilde{M}$ 
  and due to Lemma \ref{very general}, there exists a dense set of points $y$  in $\tilde{M}$ such that $\tilde{g}^{-1}(\tilde{g}(y)) \cap \Omega$ is dense in  $\tilde{g}^{-1}(\tilde{g}(y))$, so  condition $(ii)$ is 
  also satisfied.$\hfill \square$

\begin{defn}\label{mero.quot.}
Let $(\mathcal{R}, \Omega)$ be a meromorphic equivalence relation on an irreducible complex space $M$. We say that  $(\mathcal{R}, \Omega)$ admits {\bf\em a meromorphic quotient} if there exists a modification $\tau\colon \tilde{M} \rightarrow M$ with center $\Sigma\subset M \setminus \Omega $ and a holomorphic quotient map $q\colon\tilde{M}\rightarrow Q$ which satisfy the following condition:
\begin{itemize}
\item[{\rm $(iii)$}]
There exists a dense open subset $Q'$ of $Q$ such that for every $y \in \Omega' \cap q^{-1}(Q')$ the set $\tau^{-1}(R_{\tau(y)})\cap \Omega$ is dense in $q^{-1}(q(y))$
\end{itemize}
In this case we  say that $q\colon\tilde{M}\rightarrow Q$ (or simply $q : M \dasharrow Q$) is {\bf\em a meromorphic quotient of $M$ by $(\mathcal{R}, \Omega)$}.\\
\end{defn}

\begin{defn} \label{SQP-equiv.} We say that a meromorphic equivalence relation $(\mathcal{R}, \Omega)$ on an irreducible complex space $M$ is {\bf\em strongly quasi-proper} (resp. {\bf\em geometrically $f-$flat}) if the natural projection $p_1 : \mathcal{R}\rightarrow M$ is an SQP map (resp. an f-GF map). \\
\end{defn} 

Our main existence result for meromorphic quotients is the next theorem.

\begin{thm}\label{quotient}
Let $(\mathcal{R}, \Omega)$ be an SQP meromorphic equivalence relation on an irreducible complex space $M$. Then $(\mathcal{R}, \Omega)$ admits a meromorphic quotient $q\colon \tilde{M}\rightarrow Q$, where $q$ is an $f$-GF map. Moreover, in the case where $(\mathcal{R}, \Omega)$ is geometrically $f$-flat the modification $\tilde{M}\rightarrow M$ is finite.
\end{thm}

\parag{Proof} Put $n := \dim \mathcal{R} - \dim M$. Let $\pi_1\colon \mathcal{R} \rightarrow M$ and $\pi_2\colon \mathcal{R} \rightarrow M$  be  the natural projections onto the first and second components of $M\times M$,  $\varphi\colon M'\rightarrow\mathcal{C}_n^f(\pi_1)$ be the maximal  reduced fiber map for $\pi_1$ and $\Gamma$ be  the closure of the image of $\varphi$ in $\mathcal{C}_n^f(\pi_1)$. Let $\gamma\colon\Gamma\rightarrow M$ be the restriction of the natural holomorphic  map $\alpha\colon\mathcal{C}_n^f(\pi_1)^*\rightarrow M$ and $\pi\colon\Gamma\rightarrow\mathcal{C}_n^f(M)$ be the restriction to $\Gamma$ of the direct image morphism
 $$(\pi_2)_*\colon\mathcal{C}_n^f(\mathcal{R})\rightarrow\mathcal{C}_n^f(M).$$
 Then, by assumption, $\gamma$ is a modification and we shall now show that $\pi$ is semi-proper. To do so we first observe that, by identifying $\mathcal{C}_n^f(\pi_1)^*$ with the analytic subset 
$$
\{(x,C)\in M\times\mathcal{C}_n^f(M)^*\ /\ |C|\subseteq \mathcal{R}_x \},
$$
$\pi\colon\Gamma\rightarrow\mathcal{C}_n^f(\pi_1)$ is induced by the natural projection $M\times\mathcal{C}_n^f(M)^*\rightarrow\mathcal{C}_n^f(M)^*$.\\
 Next we notice that $\Gamma$ is contained in the analytic subset 
$$
M\,\sharp\,\mathcal{C}_n^f(M) := \{(x, C) \in M \times \mathcal{C}_n^f(M) \ / \  x \in |C| \}
$$
since $x\in|\varphi(x)|$ for all $x$ in $M'$. For $A\subseteq M$ and $\mathcal{B}\subseteq\mathcal{C}_n^f(M)$ we put
$$
A\,\sharp\,\mathcal{B} := (A\times\mathcal{B})\cap(M\,\sharp\,\mathcal{C}_n^f(M))
$$
and let $p_1\colon M\,\sharp\,\mathcal{C}_n^f(M)\rightarrow M$ and $p_2\colon M\,\sharp\,\mathcal{C}_n^f(M)\rightarrow \mathcal{C}_n^f(M)$ be  the natural projections. Observe that, for $\emptyset\neq W\subset\subset M$, we have $p_2(\bar{W})\,\sharp\,\Omega(W) = \Omega(W)$ and that, for every $C\in \mathcal{C}_n^f(M)$, we have $p_2^{-1}(C) = |C|\times\{C\}$. As for any $(x,C)$ in $\Gamma$ and any $y$ in $|C|$ there exists a sequence $(x_\nu)$ in $M'$ such that the sequence $(x_\nu,\varphi(x_\nu))$ tends to $(y,C)$ in $M \times \mathcal{C}_n^f(M)$ the subset  $\Gamma$ of $M\,\sharp\,\mathcal{C}_n^f(M)$ is $p_2$-saturated, i.e. $p_2^{-1}(p_2(\Gamma))= \Gamma$. Now, let $C\in \mathcal{C}_n^f(M)$ and let $W$ be a relatively compact open subset of $M$ which intersects every irreducible component of $C$. Then, keeping in mind that $\pi = {p_2}_{|\Gamma}$, we get
\begin{eqnarray*}
\pi(\gamma^{-1}(\bar{W}))\cap\Omega(W) &=& p_2(\Gamma\cap(\bar{W}\,\sharp\,\mathcal{C}_n^f(M)))\cap\Omega(W) \\
&=& p_2(\Gamma\cap(\bar{W}\,\sharp\,\Omega(W))) =  \pi(\Gamma)\cap\Omega(W)
\end{eqnarray*} 
The last equality being valid because $\Gamma$ is $p_2$-saturated and $p_2(\bar{W}\,\sharp\,\Omega(W)) = \Omega(W)$. Hence $\pi$ is semi-proper.

Due to Theorem \ref{semi-proper direct image bis} it then follows that  $T := \pi(\Gamma)$ is a reduced complex subspace of $\mathcal{C}_n^f(M)$ and $\Gamma$ is the graph of the $f$-analytic family of $n$-cycles classified by the natural inclusion $T\hookrightarrow\mathcal{C}_n^f(M)$. Now, let $\tilde{M}$ and $Q$ denote the weak normalizations of $\Gamma$ and $T$, and let $\tau\colon\tilde{M}\rightarrow M$ and   $q\colon\tilde{M}\rightarrow Q$ be the holomorphic maps which are determined by $\gamma$ and $\pi$. Then $\tau$ is a modification and $q$ is geometrically $f$-flat. Moreover, $q$ is a quotient map by Proposition \ref{quotient.weak.noramlization}. As $\tau$ and $q$ clearly satisfy condition ($iii$) of Definition \ref{mero.quot.} it follows that $q\colon\tilde{M}\rightarrow Q$ is a meromorphic quotient of $M$ by $\mathcal{R}$.\\
In the case where $\mathcal{R}$ is geometrically f-flat the modification $\tilde{M}\rightarrow M$  constructed above is clearly finite.\hfill{$\blacksquare$}\\

The next result shows that in the situation of the previous theorem, and in particular when an $f$-GF meromorphic quotient exists, the meromorphic quotient factorizes any holomorphic map on $M$ which is constant on the equivalence classes in $\Omega$.

\begin{prop}\label{Mero.quot.fGF}
In the situation of Definition \ref{mero.quot.} and with the same notation, suppose that $(\mathcal{R}, \Omega)$ is strongly quasi-proper and that $q\colon\tilde{M}\rightarrow Q$ is an $f$-GF meromorphic quotient of $M$ by $\mathcal{R}$. Then, for every reduced complex space $N$ and every holomorphic map $\pi\colon M\rightarrow N$ which is constant on $\mathcal{R}_{\tau(x)}\cap\Omega$ for all $x$ in $\Omega' \cap q^{-1}(Q')$, there exists a unique holomorphic map $g\colon Q\rightarrow N$ such that $g\circ q = \pi\circ\tau$.
$$
\xymatrix
{\tilde{M}\ar[r]^q\ar[d]_\tau& Q\ar[d]^g\\
M\ar[r]^{\pi}& N}
$$ 
\end{prop} 

\parag{Proof} Let $N$ be a reduced complex space and $\pi\colon M\rightarrow N$ be a holomorphic map which is constant on $R_{\tau(x)}\cap\Omega$ for all $x$ in $\Omega' \cap q^{-1}(Q')$. As $q$ is a holomorphic quotient map it is enough to show that $\pi\circ\tau$ is constant on every fiber of $q$. 
By assumption the set $\tau^{-1}(R_{\tau(x)}\cap\Omega)$ is dense in $ q^{-1}(q(x))$ for all $x$ in  $\Omega' \cap q^{-1}(Q')$ so $\pi\circ\tau$ \ is constant on $q^{-1}(q(x))$ for all $x$ in  $\Omega' \cap q^{-1}(Q')$. Hence by continuity the conclusion follows.$\hfill \blacksquare$\\


 \section{Reparametrization of an $f$-analytic family}
 
 An (analytic) equivalence relation on a reduced complex space $M$ may be seen as a collection of (closed analytic) subsets\footnote{We give here a translation in english of a part of the introduction of \cite{[B.08]} which gives some light on the relations between "reparametrization" and existence of some quotients in complex geometry.} 
 $(X_s)_{s \in M}$ parametrized by the set $M$ itself. The corresponding quotient is then obtained by identifying two points $s$ and $s'$ when the subset $X_s$ and $X_s'$ are the same.
 With this point of view the fact that the subsets $X_s$ define a partition of $M$ looks useless. Moreover, in the case of a meromorphic equivalence relation we no longer have the condition that the subsets $\vert X_s \vert$ are mutually disjoint.\\
 This point of view highlights the fact that the space $M$ plays two very different roles:
 \begin{itemize}
 \item $M$ is the ambient space in which live the closed analytic sets $\vert X_s\vert$.
 \item $M$ is the parameter space for the family $(X_s)_{s \in M}$.
 \end{itemize}
 In what follows, we shall keep $M$ in its first role (as the ambiant complex space) and we shall introduce a reduced complex space $S$, in general without any relation with $M$,  to parametrize the
 analytic subsets in $M$ (in fact finite type $n$-cycles in $M$). The graph $G \subset S \times M$ of the family will play the role of the graph of the equivalence relation, requiring that $G$ satisfies  the following  condition:\\
 The projection $\pi : G \to S$ is quasi-proper equidimensional map whose fibers define an $f$-analytic family of $n$-cycles in $M$ parametrized by $S$. Recall that this last condition is automatic when $S$ is normal. Then we want to find a reduced  complex space which is an analytic quotient of $S$ by the equivalence relation associated to the holomorphic map $\varphi: S \to \mathcal{C}_n^f(M)$ classifying the $f$-analytic family $(X_s)_{s \in S}$ of $n$-cycles in $M$.\\
 So as a set, such a quotient is the image $\varphi(S) \subset \mathcal{C}_n^f(M)$ and we look for a structure of reduced complex space on $\varphi(S)$, compatible with the weak analytic structure defined above on $\mathcal{C}_n^f(M)$. The key tool for such a result is of course the generalization of Khulmann's Direct Image Theorem \ref{semi-proper direct image bis} proved in chapter IV.\\

 \begin{thm}\label{reparam.} Let $M$  and $S$ be  reduced complex spaces and let $(X_s)_{s \in S}$ be a f-analytic family of $n$-cycles in $M$ parametrized by $S$. Assume that the classifying map
 $\varphi : S \to \mathcal{C}_n^f(M)$ of this family is semi-proper. Then the image $T := \varphi(S) $ is a reduced complex space (endowed with the structure sheaf induced from the weak analytic structure of $\mathcal{C}_n^f(M)$) and the restriction to $T$ of the tautological family of $\mathcal{C}_n^f(M)$   has the following universal property:
 \begin{itemize} 
 \item For any f-analytic family $(Y_v)_{v \in V}$ of $n$-cycles in $M$ parametrized by a reduced complex space $V$ such that for each $v \in V$ there exists a $s \in S$ with $Y_v = X_s$, there exists a holomorphic map $h : V \to T$ such that we have $Y_v = X_{h(v)}$ for each $v \in V$.
 \end{itemize}
 \end{thm} 
 
 \parag{Proof} Theorem \ref{semi-proper direct image bis} gives the fact that $T= \varphi(S)$ is a closed analytic subset in $\mathcal{C}_n^f(M)$ which is  reduced complex space. Then the classifying map of the f-analytic family $(Y_v)_{v \in V}$ takes its values in $T$ so defines the holomorphic map $h$.$\hfill \blacksquare$

 \parag{Remark} In fact, the previous theorem gives the existence of a  weak quotient\footnote{It may be necessary to normalize weakly $S$ and $T$  to get a holomorphic quotient in the sense defined in section VI 1.} of $S$ in the category of reduced complex spaces for the equivalence relation associated to the holomorphic classifying map $\varphi : S \to \mathcal{C}_n^f(M)$.\\
  If $G \subset S \times M$ is the set theoretic graph of this family, we may also look at this result as a existence of a weak quotient of $G$ in the category of complex spaces for the analytic equivalence relation defined by the holomorphic map $p_S\circ\varphi : G \to \mathcal{C}_n^f(M)$.\\
  
 If we begin with a f-meromorphic family of $n$-cycles in $M$  parametrized by the reduced complex space $S$, assuming that the map $\tilde{\varphi} : \tilde{S} \to \mathcal{C}_n^f(M)$ is semi-proper (where the modification $\tau : \tilde{S} \to S$ is the projection on $S$ of the graph of the meromorphic classifying map $\varphi : S \dasharrow \mathcal{C}_n^f(M)$ of our family) we obtain a f-GF meromorphic quotient of $S$ by the meromorphic analytic equivalence relation defined by the set theoretic graph of the family, the open dense set $\Omega$ in $S$ corresponding to the complement of the  polar set  of  $\tilde{\varphi}$.

\newpage
\section{Generalized Stein factorization }

In the first paragraph of this section we give an extension to the case of an SQP map of a weak version of the Stein factorization of a proper  holomorphic map.
As in this context the factorization obtained does not give, in general, a map with irreducible generic  fibers, we give in the second paragraph a necessary and sufficient condition for the existence of a factorization with a map having irreducible generic fibers.

\subsection{The general case}

Let us begin by considering the following weak version of the classical {\em Stein factorization} for a proper holomorphic map.

\begin{thm}\label{Stein.classic}
Let $\pi\colon M\rightarrow N$ be a proper surjective holomorphic map between connected normal spaces. Then there exists, a surjective holomorphic map $g\colon M \rightarrow T$ with connected fibers, where $T$ is a normal complex space, and  a proper holomorphic map   $h\colon T \rightarrow N$ with finite fibers such that  $\pi = h \circ g$. Moreover the generic fibers of $g$ are irreducible.
\hfill{$\blacksquare$}
\end{thm}

Notice that the surjectivity condition on $\pi$ can be skipped since $\pi(M)$ is a reduced complex subspace of $N$ by Remmert's Direct Image Theorem. Observe that both $g$ and $h$ are proper maps since their composition is proper and $g$ is surjective. In addition, it is easy to see that the triple $(g,T,h)$ is unique in the  sense that it is determined by a universal property (see Theorem \ref{Stein.SQP} below).\\

The following result is an extension of the previous theorem to the case of a SQP map.

\begin{thm}\label{Stein.SQP}
Let $\pi\colon M\rightarrow N$ be an SQP map between connected normal complex spaces. Then there exists, an SQP map $g\colon M \rightarrow T$, where $T$ is a (connected) normal space, and  a proper holomorphic map $h\colon T \rightarrow N$ with finite fibers such that  $ \pi = h \circ g$. Moreover, the triple $(g,T,h)$ has the following universal property:
\begin{itemize}
\item
If $\xymatrix{M \ar[r]^{g_1} & T_1 \ar[r]^{h_1} & N}$ is a factorization of $\pi$, where $T_1$ is normal and $h_1$ is a proper surjective map with finite fibers, then there exists a unique morphism $\theta\colon T\rightarrow T_1$ such that  $g_1 = \theta\circ g$ and $h = h_1\circ \theta$.
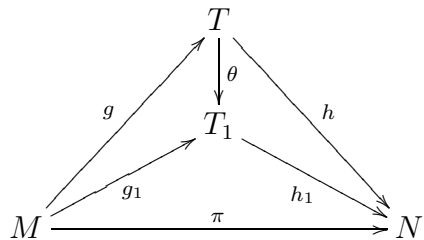
$$
\xymatrix{
&&T\ar[d]^{\theta}\ar[rrdd]^h&\\
&& T_1 \ar[rrd]_{h_1} & \\
M \ar[rruu]^{g}\ar[rru]_{g_1}\ar[rrrr]^{\pi}&&&&N
}
$$
\end{itemize}
\end{thm}

\medskip
Before proving the theorem it is interesting to point out that for a map $\pi$ which is either proper or SQP  we get the \lq\lq optimal\rq\rq\ factorization $\pi = h\circ g$, but the main difference between these two cases is, that in the SQP case the fibers of $g$ are not necessarily connected (See the remark following the proof).

\medskip
The factorization $\pi = h\circ g$ shall henceforth be called the {\bf (generalized) Stein factorization} of $\pi$.

\bigskip
For the proof of the theorem we need the following lemma.

\begin{lemma}\label{Compos.irr.diagon}
Let $\pi\colon M\rightarrow N$ be a holomorphic map between irreducible complex spaces and let $R_{\pi}$ be the equivalence relation which $\pi$ defines on $M$. Then exactly one of the irreducible components of  $R_{\pi}$ contains the diagonal. 
\end{lemma}

\parag{Proof} Denote respectively $\Delta_M$ and $\Delta_N$ the diagonals of $M\times M$ and $N\times N$. Then $R_{\pi} = (\pi\times\pi)^{-1}(\Delta_N)$. 

\smallskip
Since $M$ is irreducible $\Delta_M$ is also irreducible and consequently contained in at least one irreducible component of $R_{\pi}$. 

\smallskip
To prove that $\Delta_M$ is contained in exactly one irreducible component of $R_{\pi}$ it is not restrictive to assume that $\pi(M)$ is not contained in $S(N)$. Otherwise we may replace $N$ by $S(N)$ and so on. Next we observe that, for every smooth point $x$ in $M$, the point $(x,x)$ is a smooth point of $M\times M$ so the points $(x,x)$ in $\Delta_M$ such that $(x,x)$ is a smooth point of $M\times M$ and $(\pi(x),\pi(x))$ is a smooth point of $N\times N$ form a dense open subset $U$ of $\Delta_M$. Now, for every $(x,x)$ in $U$ such that $\pi$ is of maximal rank at $x$, the map $\pi\times\pi$ is of maximal rank at $(x,x)$. Hence there exits a point in $\Delta_M$ at which $\pi\times\pi$ is of maximal rank and such point is a smooth point of $R_{\pi}$. It follows that this point is contained in exactly one irreducible component of $R_{\pi}$ and the proof is completed.
\hfill{$\blacksquare$}

\parag{Proof of Theorem \ref{Stein.SQP}} Put $n := \dim M - \dim N$ and let $R\subseteq M\times M$ denote the analytic equivalence relation defined by $\pi$, i.e. $R := M\times_NM$. Let $R_1$ denote the union of those irreducible components of $R$ which dominate $M$ by the projection  $M\times M\rightarrow M$ onto the first factor. The map $R_1\rightarrow M$, induced by the projection onto the first factor, is the strict transform of $\pi$ by $\pi$ so it is an SQP map, due to Theorem \ref{modific.quasiproper}. Exactly one irreducible component $R_0$ of  $R$ contains the diagonal of $M\times M$, due to Lemma \ref{Compos.irr.diagon}, and we denote $\pi_1\colon R_0\rightarrow M$ the map induced by the projection onto the first factor. Then, thanks to Corollary \ref{SQP.componentwise}, $\pi_1$ is an SQP map.  

\smallskip
We shall now prove the theorem in two steps.

\medskip
{\sc First step.} Suppose that $\pi$ is an $f$-GF map.  As $\pi_1$ is an equidimensional map and $M$ is normal it follows from Corollary 4.3.13 in \cite{[e]}, that $\pi_1$ is an $f$-GF map. Let $\varphi_0\colon M\rightarrow\mathcal{C}_n^f(\pi_1)^*$ denote the reduced fibermap for $\pi_1$ and let $\Gamma$ denote its image. Now, following the same line as in the proof of Theorem \ref{quotient}, we identify $\mathcal{C}_n^f(\pi_1)^*$ with the analytic subset
$$
\left\{(x,C)\in M\times\mathcal{C}_n^f(M)^*\ /\ |C|\subseteq (R_0)_x\right\}.
$$
Then, as it is shown in the above mentioned proof, the direct image morphism $\mathcal{C}_n^f(\pi_1)^*\rightarrow\mathcal{C}_n^f(M)^*$ is induced by the natural projection $M\times\mathcal{C}_n^f(M)^*\rightarrow\mathcal{C}_n^f(M)^*$ and its restriction $\gamma\colon\Gamma\rightarrow\mathcal{C}_n^f(M)^*$ is semi-proper. Hence $\gamma(\Gamma)$ is a reduced complex subspace of $\mathcal{C}_n^f(M)^*$ thanks to Theorem \ref{semi-proper direct image bis}. Let $T$ denote the normalization of $\gamma(\Gamma)$. Then, as $M$ is normal, there is a unique holomorphic map $g\colon M\rightarrow T$ which makes the diagram 
$$
\xymatrix{& T\ar[d]\\M\ar[ru]^{g}\ar[r]_{\gamma\circ\varphi_0}& \gamma(\Gamma)}
$$ 
commutative and we observe that $g$ is surjective. Hence we obtain the commutative diagram 
$$
\xymatrix{M\ar[rr]^{g}\ar[rd]_{\pi}& &T\ar[ld]^{h}\\& N&}
$$ 
where $h$ is the composition of the normalization $T\rightarrow \gamma(\Gamma)$ and the restriction to $\gamma(\Gamma)$ of the natural map $\mathcal{C}_n^f(\pi)^*\rightarrow N$.

Now, let $y$ be a point in $N$ and let $\pi^{-1}(y) = \bigcup_{i\in I}C_i$ be its decomposition into irreducible components. For each $i$ we have $C_i\times C_i\subseteq R_0$  so $g$ is constant on $C_i$ and consequently $g(\pi^{-1}(y))$ is a finite subset of $T$.  Hence the map $h$ has finite fibers and since it is also quasi-proper, by Lemma \ref{composition}
, $h$ is proper.  It is then easily seen that $g$ is a quasi-proper map and consequently an $f$-GF map since $g$ is equidimensional and $T$ is normal.

\smallskip
Let us now show that this factorization has the universal property. To do so suppose that $\xymatrix{M \ar[r]^{g_1} & T_1 \ar[r]^{h_1} & N}$ is another such  factorization of $\pi$. Then the equivalence relation defined by $g_1$ on $M$ contains $R_0$ and it follows that $g$ is constant on the fibers of $g_1$. Thus there exist a unique holomorphic map $\theta\colon T\rightarrow T_1$ such that $\theta\circ g = g_1$ since $g_1$ is a holomorphic quotient map by Corollary \ref{holom.quot.weak.norm}. Obviously this implies $h_1\circ\theta = h$.

\medskip
{\sc Second step} In the general case there exists, thanks to Theorem \ref{SQP et GF}, a modification $\tau\colon\tilde{N}\rightarrow N$ such that the strict transform $\tilde{\pi}\colon\tilde{M}\rightarrow \tilde{N}$ of $\pi$ by $\tau$ is an $f$-GF map. Hence the commutative diagram
$$
\xymatrix{\tilde{M}\ar[d]_{\sigma}\ar[r]^{g_1} & T_1\ar[r]^{h_1} & \tilde{N}\ar[d]^{\tau}\\ 
M\ar[rr]^{\pi}&& N}
$$ 
where $\tilde{\pi} = h_1\circ g_1$, is the Stein factorization, whose existence is assured by the first step of the proof. As  $\tau\circ h_1$ is a proper map it admits a Stein factorization $\xymatrix{T_1\ar[r] & T\ar[r]^{h} & N}$, due to Theorem \ref{Stein.classic}, and consequently we get the commutative diagram 
$$
\xymatrix{\tilde{M}\ar[d]_{\sigma}\ar[r]^{\tilde{g}} & T\ar[d]^{h}\\ 
M\ar[r]^{\pi}& N}
$$ 
where $\tilde{g}$ is quasi-proper and $h$ is proper with finite fibers. Now, let $\tilde{\Gamma}$ be the graph of $\tilde{g}$ in $\tilde{M}\times T$ and put $\Gamma := (\sigma\times\id_T)(\tilde{\Gamma})$. Then $\Gamma$ is a reduced complex subspace of $M\times T$. Moreover, the natural projection $\Gamma\rightarrow M$ is a modification with finite fibers and consequently a biholomorphic map since $M$ is normal. Hence there exists a unique holomorphic map $g\colon M\rightarrow T$ whose graph is $\Gamma$. It follows that $h\circ g = \pi$  is a (generalised) Stein factorization of $\pi$. 
\hfill{$\blacksquare$}

\parag{Remark}
Assume that in the example preceding Proposition \ref{connexe 1}
\ the space  $N$ is irreducible and $\emptyset\neq A\neq N$. Then the natural projection $\pi\colon M\rightarrow N$  is clearly an SQP map and the equivalence relation defined by $\pi$ has only one irreducible component. Consequently the (generalized) Stein factorization of $\pi$ is trivial, i.e. $\pi = \id_N\circ\,\pi$, even though all the fibers of $\pi$ over $A$ have two connected components.  \\

Based on the same idea one can even construct SQP maps with trivial Stein factorization but such that the number of connected components in the  fibers is not limited.

\subsection{The useful case}

A necessary condition in order that the theorem proved in the previous paragraph gives a map  $g : M \to T$ which has irreducible  generic fibers is that there exists an integer $k \geq 1$ such that the initial map $\pi$ has  generic fibers with exactly $k$ irreducible components. The following theorem shows that this condition is sufficient for an f-GF map between connected normal complex spaces.

\begin{thm}\label{Stein.fact.}
Let $\pi\colon M\rightarrow N$ be an f-GF map between connected normal spaces and let 
$$
\xymatrix{M \ar[r]^{g} & T \ar[r]^{h} & N}
$$ 
be the generalized Stein factorization of $\pi$.  Assume, moreover, that there exists a dense  Zariski open  subset $\Omega$ of $N$ and an integer $k \geq 1$ such that, for every fiber$y$ in $\Omega$, the fiber
$\pi^{-1}(y)$  has exactly $k$ irreducible components. Then  there exists a dense open subset $T'$ of $T$ such that $g^{-1}(t)$ is irreducible for all $t$ in $T'$.
\end{thm}

\parag{Remark} We leave to the reader as an exercise on the f-flattening theorem,  the proof of the analog result for SQP map.\\

For the proof of the theorem we need the following lemmas.

\begin{lemma}\label{Stein.fact.1}
Let $\pi\colon M\rightarrow N$ be an $f$-GF map between normal connected spaces. Then there exists a nowhere dense analytic subset $A$ of $N$, a nowhere dense analytic subset $B$ of $N\setminus A$ and an open subset $M_0$ of $M$ which satisfy the following conditions:
\begin{itemize}
\item
If $x\in M_0$, then $\pi(x)$ is a smooth in $N$ and $\pi$ is a submersion at $x$.
\item
For all $y$ in $N\setminus (A\cup B)$, every irreducible component of $\pi^{-1}(y)$ intersects $M_0$.
\end{itemize}
\end{lemma}

\parag{Proof} Let $A_1$ be the set of all $y$ in $N$ such that $S(M)$ contains an irreducible component of $\pi^{-1}(y)$. Then $A_1$ is a nowhere dense analytic subset of $N$, due to Proposition \ref{une comp.}, since $\pi$ is an $f$-GF map. Put $A := A_1\cup S(N)$. Now let $\Sigma$ be the analytic subset of $M\setminus (S(M)\cup\pi^{-1}(A))$ consisting of all of points where $\pi$ is not a submersion 

Put $n := \dim M - \dim N$ and let $\varphi\colon N\rightarrow\mathcal{C}_n^f(M)$ be a fiber map for $\pi$. Then $\varphi_{|_{N\setminus A}}$ induces (by restriction of cycles) a holomorphic map 
$$
\psi\colon N\setminus A\rightarrow\mathcal{C}_n^f(M\setminus (S(M)\cup\pi^{-1}(A)))
$$
thanks to Corollary \ref{res.4}. The cycles in $\mathcal{C}_n^f(M\setminus (S(M)\cup\pi^{-1}(A)))$ which have at least one irreducible component contained in $\Sigma$ form an analytic subset $\mathcal{B}$ of $\mathcal{C}_n^f(M\setminus (S(M)\cup\pi^{-1}(A)))$ so if we put
$B := \psi^{-1}(\mathcal{B})$ and 
$$
M_0 := M\setminus (S(M)\cup\pi^{-1}(A\cup B)\cup\Sigma)
$$
the sets $A$, $B$ and $M_0$ have the required properties.
\hfill{$\blacksquare$}

\begin{lemma}\label{somme}
Let $M$ be a reduced complex space and $(Y_{\nu})$ be a sequence in $\mathcal{C}_n^f(M)$ which converges to $X$. Assume, moreover, that $X$ is reduced and has exactly $k$ irreducible components and that $Y_{\nu}$ is also reduced and has exactly $k$ irreducible components for every $\nu$. Write, for each $\nu$,
$$
Y_{\nu} = \sum_{j=1}^k\Gamma_j^{\nu}
$$
If for each $j\in\{1,\ldots,k\}$ the sequence $(\Gamma_j^{\nu})$converges in $\mathcal{C}_n^f(M)$ to a cycle $\Gamma_j$, then each $\Gamma_j$ is an irreducible component of $X$ and $X =  \sum_{j=1}^k\Gamma_j$.
\end{lemma}

\parag{Proof} As the addition of cycles is continuous we obtain $X =  \sum_{j=1}^k\Gamma_j$. It follows that each $\Gamma_j$ is irreducible since $X$ is reduced and has exactly $k$ irreducible components.
\hfill{$\blacksquare$}

\parag{Proof of Theorem \ref{Stein.fact.}} 
Then $\Omega_1 = N\setminus Z$ where $Z$ is a closed analytic   subset of empty interior in $N$. Let $A$, $B$ and $M_0$ be as in Lemma \ref{Stein.fact.} and put $\Omega' := N\setminus (A\cup Z\cup B)$ and $M' := \pi^{-1}(\Omega')\cap M_0$. Then $M'$ and $\Omega'$ are complex connected manifolds and $\pi$ induces a surjective submersion $\pi'\colon M'\rightarrow\Omega'$.\\
Now, as in the proof of  the Generalized Stein factorization Theorem \ref{Stein.SQP}, we let $R$ denote the analytic equivalence relation defined by $\pi$ and $R_0$ be the unique irreducible component of $R$ which contains the diagonal $\Delta_M$.\\
We observe that $M'$ is constructed from $M$  by first taking away an analytic subset with empty interior and then from the dense Zariski open subset which is left as the complement of  a nowhere dense closed  analytic subset. Hence $(M'\times M)\cap R_0$ is  dense in $R_0$. \\
For every $x$ in $M'$, the set $M'\cap\pi^{-1}(\pi(x)) = M'\cap R_x$ is the union of $k$ mutually disjoint connected manifolds so $x$ is contained in exactly one of them. It follows that $x$ is contained in exactly one irreducible component of $R_x$ which we shall denote by $C_x$.

\parag{Claim}  For every $x$ in $M'$ we have $(R_0)_x = C_x$.
\parag{Proof of the claim} It is clear that $C_x\subseteq (R_0)_x$ for all $x\in M'$. \\
To prove that $C_x = (R_0)_x$, for all $x\in M'$, it is enough to show that $(C_x)_{x\in M'}$ is an analytic family of $n$-cycles in $M$, because then the graph $G$ of this family is an analytic subset of $M'\times M$, which is contained in $R_0$ and of the same dimension as $R_0$. As $(M'\times M)\cap R_0$ is irreducible  it follows  that $G = (M'\times M)\cap R_0$ and consequently $(R_0)_x = C_x$ for all $x\in M'$.\\
Let us now show that $(C_x)_{x\in M'}$ is an analytic family of $n$-cycles. \\
The main point here is to prove that the family is continuous\footnote{The error in \cite{[B.13]} is at this point.}. To do so it is enough to prove that every sequence $(x_n)$  in $M'$ converging to a point $z$ in $M'$ has a subsequence  $(x_{n_k})_k$ such that the sequence of $n$-cycles $(C_{x_{n_k}})_k$ converges to $C_z$.

Let $(x_\nu)$ be a sequence in $M'$ converging to a point $z$ in $M'$ and write, for every $\nu$, 
$$
\pi^{-1}(\pi(x_\nu)) = \Gamma^\nu_1 + \cdots + \Gamma^\nu_k
$$ 
where $C_{x_\nu} = \Gamma_\nu^1$. For each $j$, every subsequence of $(x_\nu)$ has a subsequence $(x_{n_k})_k$ such that the sequence $(\Gamma_{\nu_k}^j)_k$ converges in $\mathcal{C}_n^f(M)$, because  $\pi$  is a f-GF map.  It follows that $(x_\nu)$ has a subsequence $(x_{n_k})_k$ such that, for every $j$,  the sequence $(\Gamma_{\nu_k}^j)_k$ converges in $\mathcal{C}_n^f(M)$ to an $n$-cycle $\Gamma_j$.  Due to lemma \ref{somme} it then follows that $C_{x_{\nu_k}} = \Gamma_{\nu_k}^1$ converges to an irreducible component of $\pi^{-1}(\pi(z))$ when $k$ goes to infinity. But $x_{n_k}$ is in $C_{x_{\nu_k}}$ for all $k$ and converges to $z$ as $k$ goes to infinity so $C_{x_{\nu_k}}$ converges to $C_z$ when $k$ goes to infinity because $C_z$ is the unique irreducible component of $\pi^{-1}(\pi(z))$  which contains $z$.  Hence we have proved that $(C_x)_{x\in M'}$ is an $f$-continuous family of $n-$cycles in $M$.

Now the analyticity of this family is an easy consequence of the fact that the family of fibers of $\pi$ is an f-analytic family of cycles in $M$: \\
It is enough to consider the analyticity of this family near a point $(z_0, x_0) \in R_0$ using the criterium IV 3.1.9 in \cite{[BM.1]} . And in this case this family locally coincides with the family of fibers of $\pi$.\hfill{$\blacksquare$}

\parag{Proof of Theorem \ref{Stein.fact.} continued}  As $R_0$ is an analytic subset of $M\times M$ such that $R_0\cap(M'\times M')$ is an analytic equivalence relation on $M'$ and such that $R_0\cap(\{x\}\times M')$ is dense in $R_0\cap(\{x\}\times M)$  for all $x\in M'$ it follows that $R_0$ is a meromorphic equivalence relation on $M$. Moreover, the canonical projection $R_0\rightarrow M$ is an $f$-GF map so, due to Theorem \ref{quotient}, it admits a meromorphic quotient $q\colon \tilde{M}\rightarrow Q$ where $\tilde{M}\rightarrow M$ is a modification and $q$ is an $f$-GF map. But $\pi\colon M\rightarrow N$ is an $f$-GF map so this modification is finite and hence an isomorphism since $M$ is a normal space. It then follows from Lemma \ref{Mero.quot.fGF} that there exists a unique holomorphic map $h\colon Q\rightarrow N$ such that $h\circ q = \pi$
$$
\xymatrix
{
\tilde{M}\ar[d]_q\ar[r]^\pi& N\\
Q\ar[ur]_h&
}
$$ 
and this completes the proof.
\hfill{$\blacksquare$}\\

The following example shows that Theorem \ref{Stein.fact.} is not correct if the normality condition on $M$ is skipped.

\parag{Example}
Consider the surface $M := \{(u,v,w)\in\C^3\ /\ uw^2 = v^2\}$ \ and let $\pi\colon M\rightarrow\C$ be the projection $(u,v,w)\mapsto u$. Then the following is easily seen.
\begin{itemize}
\item
The surface $M$ is irreducible and its normalization is 
$$
\nu\colon \C \times\C\longrightarrow M,\qquad (x,y)\mapsto (x^2,xy,y)
$$
\item 
The map $\pi$ is geometrically f-flat and $\pi^{-1}(u)$ consists of two reduced lines which intersect at the origin if $u\neq 0$ and $\pi^{-1}(0)$ is a double line.
\end{itemize}

\parag{Remark}
In the classical case (proper case) of the Stein factorization the conditions of Theorem \ref{Stein.fact.} are always fulfilled so, if $M$ and $N$ are normal, the irreducible components of the fibers are generically disjoint.

\end{document}